\documentclass[psamsfonts,reqno]{memo-l}
\usepackage{amscd,amssymb,eucal,latexsym}

\newcommand{\ii}[1]{\emph{#1}}

\theoremstyle{plain}

\newtheorem{lemma}{Lemma}[section]
\newtheorem{theorem}[lemma]{Theorem}
\newtheorem{proposition}[lemma]{Proposition}
\newtheorem{corollary}[lemma]{Corollary}
\newtheorem{claim}{Claim}
\newtheorem*{sclaim}{Claim}

\newtheorem*{stat}{\name}
\newcommand{\name}{testing}

\theoremstyle{definition}
\newtheorem{definition}[lemma]{Definition}
\newtheorem{example}[lemma]{Example}
\newtheorem*{problem}{Problem}

\theoremstyle{remark}
\newtheorem{remark}[lemma]{Remark}
\newtheorem{notation}[lemma]{Notation}
\newtheorem*{note}{Note}

\newenvironment{all}[1]{\renewcommand{\name}{#1}\begin{stat}}
                          {\end{stat}}

\newcommand{\qedc}{{\qed}~{\rm Claim~{\theclaim}.}}
\newcommand{\qedsc}{{\qed}~{\rm Claim.}}

\newenvironment{cproof}
{\begin{proof}[Proof of Claim.]}
{\qedc\renewcommand{\qed}{}\end{proof}}

\newenvironment{scproof}
{\begin{proof}[Proof of Claim.]}
{\qedsc\renewcommand{\qed}{}\end{proof}}

\numberwithin{equation}{section}

\newcommand{\Mh}{$(\mathrm{M_{ht}})$}
\newcommand{\Ml}{$(\mathrm{M_{lh}})$}

\DeclareMathOperator{\wdt}{wd}
\DeclareMathOperator{\supp}{supp}
\newcommand{\two}{\mathbf{2}}
\newcommand{\On}{\mathbf{On}}
\newcommand{\Cn}{\mathbf{2}_{\infty}}
\newcommand{\Zn}{\mathbb{Z}_{\infty}}
\newcommand{\Rn}{\mathbb{R}_{\infty}}
\newcommand{\Nh}{\mathcal{N}}
\newcommand{\eps}{\varepsilon}
\renewcommand{\iff}{if and only if}
\newcommand{\ol}[1]{\,\overline{\!#1}}
\newcommand{\oll}[1]{\overline{#1}}
\newcommand{\trim}{\leq_{\mathrm{trim}}}
\newcommand{\id}{\mathrm{id}}
\newcommand{\es}{\varnothing}
\newcommand{\rem}{\ll_{\mathrm{rem}}}
\newcommand{\srem}{\lesssim_{\mathrm{rem}}}
\newcommand{\sd}{\smallsetminus}
\newcommand{\conc}{\mathbin{{}^{\frown}}}
\newcommand{\scal}[2]{\left\langle\nobreak#1\nobreak\mid\nobreak#2
\nobreak\right\rangle}

\newcommand{\di}[1]{\tfrac{\mathstrut\textstyle{#1}}{\scriptstyle{\infty}}}
\newcommand{\DI}[1]{#1|_{\infty}}

\newcommand{\annr}{\mathrm{ann_r}}
\newcommand{\pdi}[1]{\left(\di{#1}\right)}
\newcommand{\CC}{\mathbf{C}}
\newcommand{\lsc}{lower semicontinuous}
\newcommand{\usc}{upper semicontinuous}
\newcommand{\Sbul}{S^{\bullet}}
\newcommand{\blp}{+^{\bullet}}
\newcommand{\Ref}[1]{\widetilde{#1}}
\newcommand{\mm}{\mathsf{m}}
\newcommand{\nn}{\mathsf{n}}
\newcommand{\Bo}{B_{\omega}}
\newcommand{\Co}{C_{\omega}}
\newcommand{\lA}{\boldsymbol{A}}
\newcommand{\lL}{\boldsymbol{L}}
\newcommand{\tlA}{\widetilde{\lA}}

\newcommand{\CCx}{\mathbb{C}}
\newcommand{\tCCx}{\widetilde{\CCx}}
\newcommand{\lBo}{\boldsymbol{B_{\omega}}}
\newcommand{\lCo}{\boldsymbol{C_{\omega}}}
\DeclareMathOperator{\rB}{B}
\DeclareMathOperator{\rE}{E}

\newcommand{\tvi}{\vrule height 12pt depth 6pt width 0pt}

\newcommand{\fa}{\mathfrak{a}}
\newcommand{\fb}{\mathfrak{b}}

\newcommand{\fm}{\mathfrak{m}}

\newcommand{\la}{\boldsymbol{a}}
\newcommand{\lb}{\boldsymbol{b}}
\newcommand{\lc}{\boldsymbol{c}}

\newcommand{\lx}{\boldsymbol{x}}
\newcommand{\lp}{\boldsymbol{p}}
\newcommand{\ly}{\boldsymbol{y}}
\newcommand{\lz}{\boldsymbol{z}}

\newcommand{\lu}{\boldsymbol{u}}

\newcommand{\lga}{\boldsymbol{\alpha}}
\newcommand{\lgb}{\boldsymbol{\beta}}

\newcommand{\lgc}{\boldsymbol{\gamma}}

\newcommand{\dx}{\dot{x}}
\newcommand{\dy}{\dot{y}}

\newcommand{\fin}{\mathrm{fin}}
\newcommand{\mf}{\mathrm{mf}}

\renewcommand{\SS}{\mathfrak{S}}
\newcommand{\bv}[1]{\left\|#1\right\|}
\newcommand{\vbv}[1]{{[\![}#1{]\!]}}

\newcommand{\I}{\mathrm{I}}
\newcommand{\II}{\mathrm{II}}
\newcommand{\III}{\mathrm{III}}

\newcommand{\pup}[1]{\textup{(}{#1}\textup{)}}

\newcommand{\NN}{\mathbb{N}}
\newcommand{\ZZ}{\mathbb{Z}}

\newcommand{\RR}{\mathbb{R}}

\DeclareMathOperator{\cc}{cc}

\DeclareMathOperator{\BB}{Proj}
\newcommand{\BBp}[1]{\BB^*{#1}}

\newcommand{\cm}{commutative mon\-oid}
\newcommand{\pcm}{partial commutative mon\-oid}
\newcommand{\PCM}{Partial Commutative Mon\-oid}

\newcommand{\prm}{partial refinement mon\-oid}

\newcommand{\DD}{\Delta}
\DeclareMathOperator{\Dim}{Dim}
\DeclareMathOperator{\Drng}{Drng}

\newcommand{\set}[1]{\{#1\}}
\newcommand{\Set}[1]{\left\{#1\right\}}

\newcommand{\setm}[2]{\set{#1\mid#2}}
\newcommand{\Setm}[2]{\Set{#1\mid#2}}

\newcommand{\famm}[2]{(#1)_{#2}}

\font\ssfont=cmss8 scaled\magstep1

\newcommand{\Lat}{\mathcal{L}}
\newcommand{\ModR}{\mbox{\ssfont Mod-}R}
\newcommand{\addA}{\mbox{\ssfont add-}A}
\newcommand{\addR}{\mbox{\ssfont add-}R}
\newcommand{\Hom}{\mathrm{Hom}}
\newcommand{\J}{\mathrm{J}}
\newcommand{\NSIR}{\mbox{\ssfont NSI-}R}
\newcommand{\End}{\mathrm{End}}
\newcommand{\IIf}{${\II}_{\mathrm{f}}$}
\newcommand{\calB}{\mathcal{B}}
\DeclareMathOperator{\barotimes}{\overline\otimes}

\author{K.\,R. Goodearl}
\address{Department of Mathematics\\
University of California\\
Santa Barbara, CA 93106\\
U.S.A.}
\email{goodearl@math.ucsb.edu}
\urladdr{http://www.math.ucsb.edu/\~{}goodearl/}

\author{F. Wehrung}
\address{CNRS, UMR 6139\\
D\'epartement de Math\'ematiques\\
Universit\'e de Caen\\
14032 Caen Cedex\\
France}
\email{wehrung@math.unicaen.fr}
\urladdr{http://www.math.unicaen.fr/\~{}wehrung}

\keywords{Partial monoid, commutative monoid, continuous dimension scale,
espalier, projection, directly finite, purely infinite, Boolean-valued, lower
embedding, lattice, dimension theory, complete, orthogonality, extremally
disconnected}

\subjclass[2000]{Primary 06B15, 06C15, 06B23, 06B30, 06E15, 06C10, 08A55.
Secondary 06B35, 06C20, 06A05, 16D50, 16D70, 16D90, 16E20, 16E50, 18E10,
19A49, 19K14, 28B10, 28A60, 46L51}

\makeindex

\begin{document}

\title[Complete dimension theory]
{The complete dimension theory of\\
partially ordered systems with\\
equivalence and orthogonality\\
}

\begin{abstract}
We develop dimension theory for a large class of structures of the form
\linebreak
$(L\nobreak,\leq\nobreak,\perp\nobreak,\sim\nobreak)$,
where $(L\leq)$ is a
partially ordered set, $\perp$ is a binary relation on~$L$, and $\sim$ is
an equivalence relation on
$L$, subject to certain axioms. We call these structures
\emph{espaliers}.\index{espalier} For $x$, $y$, $z\in L$, we say that
$z=x\oplus y$ holds, if $x\perp y$ and $z$ is the supremum of $\set{x,y}$.
The \emph{dimension theory} of $L$ is the universal $\sim$-invariant
homomorphism from $(L,\oplus,0)$ to a
\pcm\index{partial commutative monoid}\
$S$. We say that $S$ is the \emph{dimension range}\index{dimension range} of
$L$. Particular examples of espaliers are the following:
\begin{enumerate}
\item Let $B$ be a \index{Boolean algebra!complete ---} complete Boolean
algebra. For $x$, $y\in B$, we say that
$x\perp y$ if $x\wedge y=0$, and we take $\sim$ to be any zero-separating,
unrestrictedly additive and refining equivalence relation on $B$ (for
instance, equality).

\item Let $R$ be a right self-injective von Neumann regular ring.
\index{ring!Von Neumann regular ---}\index{ring!right self-injective ---}%
We denote by~$L$ the lattice\index{lattice}
of all direct summands of a given nonsingular
injective right $R$-module, for instance, the lattice of finitely generated
right ideals of $R$. For $A$, $B\in L$, we say that $A\perp B$ if
$A\cap B=\set{0}$, and $A\sim B$ if $A\cong B$.

\item More generally,
let $L$ be a complete, meet-continuous, complemented, modular
\index{lattice!complemented modular ---} lattice.
For $x$, $y\in L$, we say that $x\perp y$ if $x\wedge y=0$, and $x\sim y$ if
$x$ and $y$ are projective
\index{projectivity by decomposition@projectivity by decomposition
($\approxeq$)} by (finite) decomposition.

\item Let $A$ be an \index{AW*-algebra}AW*-algebra. We denote by $L$ the
lattice of projections of~$A$, and take the standard orthogonality and
equivalence relations on~$L$. For $p$, $q\in L$, then, $p\perp q$ if $pq=0$,
and
$p\sim q$ if $p$ and $q$ are \emph{Murray-von~Neumann equivalent}, that is,
there exists $x\in A$ such that $p=x^*x$ and $q=xx^*$.

\end{enumerate}
We prove that the dimension range\index{dimension range} of any
espalier\index{espalier}
$(L,\leq,\perp,\sim)$ is a \emph{lower interval} of a \cm\ of the form
   \begin{equation}
   \CC(\Omega_{\I},\ZZ_\gamma)\times\CC(\Omega_{\II},\RR_\gamma)
   \times\CC(\Omega_{\III},\two_\gamma),\tag{*}
   \end{equation}
where $\Omega_{\I}$, $\Omega_{\II}$, and $\Omega_{\III}$ are complete
Boolean spaces,\index{Boolean space}
and where we put, for every ordinal $\gamma$,
   \begin{align*}
   \ZZ_\gamma&=\ZZ^+\cup\setm{\aleph_\xi}{0\leq\xi\leq\gamma},\\
   \RR_\gamma&=\RR^+\cup\setm{\aleph_\xi}{0\leq\xi\leq\gamma},\\
   \two_\gamma&=\set{0}\cup\setm{\aleph_\xi}{0\leq\xi\leq\gamma},
   \end{align*}
\index{Zzzgamma@$\ZZ_\gamma$}\index{Rzzgamma@$\RR_\gamma$}%
\index{Tzzgamma@$\two_\gamma$}%
endowed with their interval topology\index{interval topology}
and natural addition operations.
Conversely, we prove that every lower interval of a monoid of the form
(*) can be represented as the dimension range\index{dimension range} of an
espalier\index{espalier} arising from each of the contexts (i)--(iv) above.
The context of W*-algebras requires the spaces
$\Omega_{\I}$, $\Omega_{\II}$, and $\Omega_{\III}$ to be
\emph{hyperstonian}, and no further restriction is needed.

This subsumes many earlier dimension-theoretic results,
and, in applications, completes theories developed for examples such as
(i)--(iv) above.
\end{abstract}

\maketitle

\tableofcontents

\chapter{Introduction}\label{Ch:Intro}

\section{Background}\label{S:Backgr}

The central theme of this paper is, as indicated by the title,
\emph{dimension theory}. Basically, an equivalence relation $\sim$ is given
on a structure $L$, and our goal is to elucidate the quotient structure
$L/{\sim}$. We shall be interested in cases where $L$ is a \emph{complete
lattice}\index{lattice!complete ---} endowed with a notion of orthogonality
subject to a number of axioms. To make it clear that the motivations for
this are widespread and by no means confined to lattice theory, we start by
discussing what is known about some fundamental examples. In
Chapter~\ref{Ch:ClEsp}, we will apply our general theory to these examples,
thus showing the improvements that it brings to them.

\subsection{Abstract measure theory}\label{Su:AMeasTh}
One of the most basic examples of what could be called a ``dimension theory''
arises from measure theory. To pick a favorite, we first consider
the Lebesgue measure\index{Lebesgue measure}
$\mm$ on the real line $\RR$. It is defined on
the Boolean algebra\index{Boolean algebra} $\mathcal{B}$ of all
Lebesgue-measurable subsets of
$\RR$. However, it fails total additivity of measure, for every subset of
$\RR$ is the union of singletons, which have Lebesgue measure zero. To bring
back total additivity, the standard way is to say that $\mm$ is defined not on
$\mathcal{B}$, but on the quotient algebra $B=\mathcal{B}/\mathcal{N}$,
where $\mathcal{N}$ is the ideal of null sets.
The Boolean algebra\index{Boolean algebra} $B$ and the resulting map from
$B$ to $[0,+\infty]$, which we still denote by $\mm$, have the following
properties:
\begin{itemize}
\item[(a)] $B$ is a \emph{complete}
\index{Boolean algebra!complete ---}Boolean algebra.

\item[(b)] The map $\mm$ is \emph{unrestrictedly additive},
\index{unrestrictedly additive!--- map|ii} that is,
the following equality holds:
   \[
   \mm\left(\bigvee_{i\in I}x_i\right)=\sum_{i\in I}\mm(x_i),
   \]
for any \emph{disjoint} family $\famm{x_i}{i\in I}$ of elements of $B$.
The notation $\bigvee_{i\in I}x_i$ stands for the \emph{join} (i.e.,
supremum) of the set $\setm{x_i}{i\in I}$ in $B$.

\item[(c)] For all $x$, $y\in B$, if $y$ is a translate of $x$ (that is,
$y=\alpha+x$ for some real number $\alpha$), then
$\mm(x)=\mm(y)$.
\end{itemize}

Rule (b) above seems somehow puzzling at first glance, because of the
apparent possibility of an uncountable index set $I$. However, since the
Boolean algebra\index{Boolean algebra} $B$ is countably saturated, all
infinite joins in $B$ are, really, \emph{countable} joins, so that, in~(b),
all the $x_i$-s are majorized by the join of countably many of them.

For $x$, $y\in B$, we define the relation $x\sim y$ to hold, if there
are disjoint families $\famm{x_i}{i\in I}$ and $\famm{y_i}{i\in I}$ of
elements of $B$ such that $x=\bigvee_{i\in I}x_i$ and
$y=\bigvee_{i\in I}y_i$, and $y_i$ is a translate of $x_i$, for all
$i\in I$. It is not difficult to verify that $\sim$ is an \emph{equivalence
relation} on $B$. Furthermore, by (b) and (c) above, $x\sim y$ implies
that $\mm(x)=\mm(y)$, for all $x$, $y\in B$.

It is harder to verify that the converse of the above fact
also holds, namely: $\mm(x)=\mm(y)$ implies that $x\sim y$, for
all $x$, $y\in B$. This fact is due to S. Banach and A. Tarski,
see \cite{BaTa24},\index{Banach, S.}\index{Tarski, A.}
or \cite[Theorem~9.17]{Wago85}.\index{Wagon, S.} Hence the quotient 
set $B/{\sim}$ is
isomorphic, \emph{via} the measure~$\mm$, to the interval $[0,+\infty]$. A
moment's reflection shows that $B/{\sim}$ can be endowed with a \emph{partial
addition}, defined by the rule
   \[
   [x]+[y]=[x\vee y],\quad\text{for all disjoint }x,\,y\in B,
   \]
that endows it with a structure of \emph{\pcm}\index{partial commutative
monoid} (see Definition~\ref{D:PartCM}), and that the measure $\mm$ factors
through an isomorphism of partial monoids between $B/{\sim}$ and
$[0,+\infty]$. We see in this particular case that $B/{\sim}$ is a
\emph{total} monoid, that is, the addition of $B/{\sim}$ is defined
everywhere.

\smallskip

Now let us consider the \emph{converse} of the above paragraph. That is, we
are given a Boolean algebra\index{Boolean algebra} $B$, endowed with an
equivalence relation $\sim$, and we wish to find the structure of
$B/{\sim}$. While this problem in full generality can lead to almost any
structure, we focus the study by making the following assumptions on $B$ and
$\sim$, that are satisfied for the example above:
\begin{itemize}
\item[(1)] $B$ is a complete\index{Boolean algebra!complete ---} Boolean
algebra.

\item[(2)] $x\sim 0$ implies that $x=0$, for all $x\in B$.

\item[(3)] (see Axiom~(L6) of Definition~\ref{D:MeasChLatt}) The relation
$\sim$ is \emph{unrestrictedly refining},
\index{unrestrictedly refining relation} that is, for every $a\in L$ and
every disjoint family $\famm{b_i}{i\in I}$ of elements of~$L$, if
$a\sim\bigvee_{i\in I}b_i$, then there exists a decomposition
$a=\bigvee_{i\in I}a_i$, with $\famm{a_i}{i\in I}$ disjoint, such that
$a_i\sim b_i$ for all $i\in I$.

\item[(4)] (see Axiom~(L7) of Definition~\ref{D:MeasChLatt}) The relation
$\sim$ is \emph{unrestrictedly additive},
\index{unrestrictedly additive!--- relation} that is, for all disjoint
families $\famm{a_i}{i\in I}$ and $\famm{b_i}{i\in I}$ of elements
of~$L$, if $a_i\sim b_i$ for all $i\in I$, then
$\bigvee_{i\in I}a_i\sim\bigvee_{i\in I}b_i$.
\end{itemize}

The most basic example of this situation is for
$B=\mathfrak{P}(\Omega)$, the powerset algebra of an infinite set $\Omega$,
where $\sim$ is the relation of \index{equipotency|ii}\emph{equipotency}
on subsets of $\Omega$, that is, $X\sim Y$ \iff\ there exists a bijection
from $X$ onto $Y$. If $\gamma$ is the unique ordinal such that
$|\Omega|=\aleph_\gamma$, then $L/{\sim}$ is isomorphic to the monoid
\index{Zzzgamma@$\ZZ_\gamma$}
   \[
   \ZZ_\gamma=\ZZ^+\cup\setm{\aleph_\xi}{0\leq\xi\leq\gamma},
   \]
endowed with the addition that extends the natural addition of the set
$\ZZ^+$ of nonnegative integers and such that
$n+\aleph_\beta=\aleph_\alpha+\aleph_\beta=\aleph_\beta$, for all $n\in\ZZ^+$
and all ordinals $\alpha$, $\beta$ such that $\alpha\leq\beta\leq\gamma$, see
page~\pageref{Pg:ZR2gam}. So, if $\mu\colon B\to\ZZ_\gamma$ is the map
defined by the rule $\mu(X)=|X|$, for all $X\in B$, then $\mu$ factors
through $\sim$, thus defining an isomorphism from $B/{\sim}$ onto
$\ZZ_\gamma$.

As we shall see in this paper, it is still possible, in the general case, to
obtain a ``measure'' $\mu$ on $B$ such that $B/{\sim}$ is isomorphic to the
range of $\mu$. The range of the measure $\mu$ is not necessarily
$[0,+\infty]$ and not even some $\ZZ_\gamma$ (as in the example above), but
rather a certain set of continuous functions from a complete Boolean space
\index{Boolean space|ii}
(i.e., extremally disconnected compact Hausdorff topological space)
$\Omega$ to a monoid of the form
$\RR^+\cup\setm{\aleph_\xi}{0\leq\xi\leq\gamma}$ (or a submonoid of this
monoid). A similar result is achieved by D.~Maharam
in~\cite{Maha49},\index{Maharam, D.} in a slightly different context---for
instance, all sums and joins are countable joins, while~$B$ satisfies the
countable chain condition. This is not the only restriction imposed in
Maharam's work, as, for example, Axiom~III, page~281 in
\index{Maharam, D.}\cite{Maha49}, that rules out what we will call later the
``Type~III'' case.

\subsection{Nonsingular injective modules over self-injective regular
rings}\label{Su:RSIReg}
Let $R$ be a (von~Neumann) regular, right
self-injective ring,\index{ring!Von Neumann regular ---}
\index{ring!right self-injective ---} let $M$ be a nonsingular injective
right
$R$-module. We order the set $L$ of all direct summands of $M$ by inclusion
and we endow it with the relation of isomorphism, $\cong$. The
\emph{dimension theory} of
$M$ is the study of the structure of $L/{\cong}$. We say that a family
$\famm{X_i}{i\in I}$ of elements of $L$ is \emph{orthogonal}, if the sum of
the submodules $X_i$ is a direct sum. We recall some fundamental properties
of $L$ and $\cong$ (references will be given in Section~\ref{S:RSIReg}):
\begin{itemize}
\item[(1)] $L$ is a complete lattice,\index{lattice!complete ---|ii} that is,
every subset of $L$ has a supremum.

It is known that the infimum of a family $\famm{X_i}{i\in I}$ of
elements of $L$ is their \emph{intersection}, $\bigcap_{i\in I}X_i$.

\item[(2)] $L$ is \emph{complemented},
\index{lattice!complemented ---|ii}
that is, every element $X$ of $L$
has a complement (that is, an element $Y$ of $L$ such that $X\oplus Y=M$).

\item[(3)] $L$ is \emph{meet-continuous},
\index{lattice!meet-continuous ---|ii}
that is, for every $X\in L$ and
every upward directed family $\famm{Y_i}{i\in I}$ of elements of $L$,
the following equality holds:
   \[
   X\cap\bigvee_{i\in I}Y_i=\bigvee_{i\in I}(X\cap Y_i).
   \]

\item[(4)] $L$ is \emph{modular}, \index{lattice!modular ---|ii}
that is, the equality
   \[
   X\cap(Y\vee Z)=(X\cap Y)\vee Z
   \]
holds, for all $X$, $Y$, $Z\in L$ such that $X\supseteq Z$.

\item[(5)] (see Axiom~(L6) of Definition~\ref{D:MeasChLatt}) The relation
$\cong$ is \emph{unrestrictedly refining},
\index{unrestrictedly refining relation} that is, for every $X\in L$ and
every orthogonal family $\famm{Y_i}{i\in I}$ of elements of~$L$, if
$X\cong\bigvee_{i\in I}Y_i$, then there exists an orthogonal decomposition
$X=\bigvee_{i\in I}X_i$ such that $X_i\cong Y_i$ for all $i\in I$.

\item[(6)] (see Axiom~(L7) of Definition~\ref{D:MeasChLatt})  The relation
$\cong$ is \emph{unrestrictedly additive},
\index{unrestrictedly additive!--- relation} that is, for all orthogonal
families $\famm{X_i}{i\in I}$ and $\famm{Y_i}{i\in I}$ of elements of~$L$, if
$X_i\cong Y_i$ for all $i\in I$, then
$\bigvee_{i\in I}X_i\cong\bigvee_{i\in I}Y_i$.

\end{itemize}

We observe that the supremum in $L$ of a family $\famm{X_i}{i\in I}$ of
elements of $L$ is not given by the sum of submodules $\sum_{i\in I}X_i$, but
by its\index{injective hull} \emph{injective hull},\index{ErzzE@$\rE(M)$}
$\rE\left(\sum_{i\in I}X_i\right)$ (which can be identified with a unique
submodule of $M$ because $M$ is injective and nonsingular).

As in Subsection~\ref{Su:AMeasTh}, the quotient set $L/{\cong}$ can be endowed
with a structure of \emph{\pcm},\index{partial commutative monoid} under
the addition given by the rule
   \[
   [X]+[Y]=[X\oplus Y]\quad\text{if }X\cap Y=\set{0},
   \]
for all $X$, $Y\in L$.

Essentially by using Axioms (1)--(6) above, the structure of $L/{\cong}$
has been completely elucidated in several particular cases. For example, in
case $M$ is \emph{directly finite}\index{directly finite}
(i.e., $M$ is not isomorphic to any
proper direct summand of itself), $L/{\cong}$ is isomorphic to a lower
subinterval (with respect to the componentwise ordering) of a monoid of the
form
   \begin{equation}\label{Eq:MonI,II}
   M=\CC(\Omega_{\I},\ZZ^+\cup\set{\infty})\times
   \CC(\Omega_{\II},\RR^+\cup\set{\infty}),
   \end{equation}
where $\Omega_{\I}$ and $\Omega_{\II}$ are complete Boolean spaces;
\index{Boolean space} see \index{Goodearl, K.\,R.}\index{Boyle, A.\,K.}%
Chapter~11 in K.\,R.~Goodearl and A.\,K.~Boyle \cite{GoBo}. In the general
case, there are a monoid~$M$ of the form given by \eqref{Eq:MonI,II} and a
direct power~$N$ of a monoid of the form
$\set{0}\cup\setm{\aleph_\xi}{\xi<\gamma}$ (for a certain ordinal $\gamma$)
such that $L/{\cong}$ embeds into $M\times N$, see Chapters 12 and~13
\index{Goodearl, K.\,R.}\index{Boyle, A.\,K.}%
in~\cite{GoBo}, and the variations in \index{Goodearl, K.\,R.}
\cite[Chapter 12]{GvnRR}. Further results along these lines were obtained by
C.~Busqu\'e \index{Busqu{\'e}, C.}
\cite{Busq90}, who showed, in particular, that the second factor of the
embedding above, namely the map $L/{\cong} \rightarrow N$, actually sends
$L/{\cong}$ to $\CC(\Omega, \set{0}\cup\setm{\aleph_\xi}{\xi<\gamma})$
for a suitable complete Boolean space\index{Boolean space} $\Omega$
(containing $\Omega_{\I}\sqcup\Omega_{\II}$)
\cite[Proposition~4.7]{Busq90}.\index{Busqu{\'e}, C.} However, these
embeddings do not provide an isomorphism of
$L/{\cong}$ onto a \emph{lower subset} of a monoid of continuous functions.
The difficulties are already visible in case $R$ is a complete Boolean algebra
\index{Boolean algebra!complete ---}
(viewed as a ring), due to an example of K.~Eda \index{Eda, K.}\cite{Eda84}:
There exists a complete Boolean algebra\index{Boolean algebra!complete ---}
$R$ such that the injective hull\index{injective hull} of the free
$R$-module of rank $\aleph_0$ contains a direct sum of $\aleph_1$ copies of
itself (see the discussion of Problem 18 in \cite[p.~374]{GvnRR}).
\index{Goodearl, K.\,R.} Here the
image of the embedding obtained from \cite{GoBo},
\index{Goodearl, K.\,R.}\index{Boyle, A.\,K.}%
\cite{GvnRR}, and
\cite{Busq90}\index{Busqu{\'e}, C.} contains a function with all values at
least $\aleph_2$, but not the constant function with value $\aleph_1$.

\subsection{Conditionally complete, meet-continuous, sectionally complemented,
modular lattices}\label{Su:CMSMLatt}
For elements $a$, $b$, and $c$ in a lattice $L$ with zero, we say that
$c=a\oplus b$,
\index{cabozzplus@$c=a\oplus b$ (in a lattice)|ii}
if $c=a\vee b$ and $a\wedge b=0$. We say that $L$
\index{lattice!sectionally complemented ---|ii}
is \emph{sectionally complemented}, if for all $a$, $b\in L$ such that
$a\leq b$, there exists $x\in L$ such that $a\oplus x=b$. If $L$ is
\emph{modular}, that is, the implication
   \[
   x\geq z\Longrightarrow x\wedge(y\vee z)=(x\wedge y)\vee z
   \]
holds, for all $x$, $y$, $z\in L$, then the partial
operation $\oplus$ gives~$L$ a structure of
\emph{\pcm}.\index{partial commutative monoid} Completeness and
meet-continuity of~$L$ are defined as in (1) and (3) of
Subsection~\ref{Su:RSIReg}. So, in particular, if $M$ is a nonsingular
injective right module over a right self-injective regular ring
\index{ring!Von Neumann regular ---}\index{ring!right self-injective ---}%
$R$, then the lattice of all direct summands of $M$ is complete,
\index{lattice!sectionally complemented modular ---}
meet-continuous, sectionally complemented, and modular. The classical
von~Neumann\index{continous geometry}
\emph{continuous geometries}, see J. von~Neumann
\index{von Neumann, J.}\cite{JvNbook} or F. Maeda
\cite{FMae55},\index{Maeda, F.} are obtained by adding the conditions that
$L$ has a unit (that is, a largest element) and is join-continuous.

At this point, we seem to be stymied because of the following problem. We
cannot claim outright that our lattice-theoretical context could lead to
generalizations of Subsection~\ref{Su:RSIReg}, for there is no such
thing \emph{a priori} as ``isomorphism of submodules'' between the elements of
$L$. In the case of\index{continous geometry} continuous geometries, it is
easy to remedy this by replacing isomorphism by
\index{perspectivity@perspectivity ($\sim$)}\emph{perspectivity}, which
turns out to be
\emph{transitive} (this is a difficult result, due to J.
von~Neumann~\cite{JvNbook}).\index{von Neumann, J.} Elements $a$ and $b$ of
a lattice~$L$ are
\emph{perspective}, in notation \index{perspectivity@perspectivity ($\sim$)}
$a\sim b$, if there exists
$x\in L$ such that $a\wedge x=b\wedge x$ and $a\vee x=b\vee x$.
For\index{continous geometry} continuous geometries, the structure of
$L/{\sim}$ is completely understood, see \cite{JvNbook}
\index{von Neumann, J.} and, for the general, reducible case,
\index{Iwamura, T.}T.~Iwamura \cite{Iwam44}---namely,
$L/{\sim}$ is isomorphic to a lower segment of the positive cone of a
Dedekind complete lattice-ordered group. The paper J.~Harding and
M.\,F.~Janowitz~\cite{HaJa97}
\index{Harding, J.}\index{Janowitz, M.\,F.}%
shows how a reducible\index{continous geometry} continuous
geometry can be represented as the space of continuous sections of a bundle
of irreducible continuous geometries, thus shedding more light on the
transition from irreducible continuous geometries\index{continous geometry}
to reducible ones.

However, for a general complete, meet-continuous, sectionally complemented,
modular
\index{lattice!sectionally complemented modular ---}
lattice $L$, the relation of
perspectivity\index{perspectivity@perspectivity ($\sim$)} $\sim$ on $L$ is
not transitive as a rule---see, for example, the obvious case where $L$ is
the subspace lattice of an infinite-dimensional vector space over a field.
Hence, we have to find a better candidate than
$\sim$ to replace isomorphism of submodules. A natural guess is of course
the \emph{transitive closure} $\approx$ of $\sim$
\index{projectivity@projectivity ($\approx$)}
(usually called \emph{projectivity}), but this relation
fails to be additive, as defined in Axiom~(L7) of
Definition~\ref{D:MeasChLatt}, and as isomorphism of submodules should be.
The final answer is, in fact, nontrivial, and it follows from the theory of
\emph{normal equivalences}\index{normal equivalence}
introduced by the second author in
Chapters 10--13 of\index{Wehrung, F.} \cite{WDim}. Namely, there is
(fortunately!)
\emph{exactly one} ``reasonable'' candidate for isomorphism, and it is the
binary relation $\equiv$ on $L$ defined by the rule
   \begin{multline*}
   a\equiv b\text{ if there are decompositions }a=x_0\oplus x_1\text{ and }
   b=y_0\oplus y_1,\\
   \text{ with }x_0\sim y_0\text{ and }x_1\sim y_1,
   \end{multline*}
for all $a$, $b\in L$.
The transitivity of $\equiv$ is proved in Theorem~13.2 of\index{Wehrung, F.}
\cite{WDim}, while the complete additivity of $\equiv$ (Axiom~(L7) of
Definition~\ref{D:MeasChLatt}, see also item~(6) of
Subsection~\ref{Su:RSIReg}) is proved as in Proposition~13.9
of\index{Wehrung, F.} \cite{WDim} by replacing countable families by
arbitrary families. The quotient
$L/{\equiv}$ is then a lower subset of the so-called \emph{dimension monoid}
\index{dimension monoid} $\Dim L$\index{dzzimL@$\Dim L$}
of $L$, which, as its name indicates, is a (commutative) monoid.
The dimension theory of $L$ is elucidated here in
Theorem~\ref{T:DimMonRCMLatt}. In the context of
Subsection~\ref{Su:RSIReg}, that is, $L$ is the lattice of all direct
summands of a given nonsingular injective right module over a right
self-injective regular ring,
\index{ring!Von Neumann regular ---}\index{ring!right self-injective ---}%
it is then the case that $\equiv$ is identical to submodule isomorphism on
$L$, see Lemma~10.2 and Theorem~13.2 of\index{Wehrung, F.} \cite{WDim}.

Of crucial importance for all the proofs of these results is
a result of I. Halperin and J.~von~Neumann \cite{HaNe40}
\index{Halperin, I.}\index{von Neumann, J.}%
that states that
$x\approx y$ and $x\wedge y=0$ implies that $x\sim y$, for all $x$,
$y\in L$. This result is extended in\index{Wehrung, F.} \cite{WDim} to
countably meet-continuous lattices, where it is used to prove that the
quotient $L/{\equiv}$ is then a so-called \emph{generalized cardinal
algebra}, see Proposition~13.10 of\index{Wehrung, F.}
\cite{WDim}. However, even in case the lattice-theoretical version of
direct finiteness (see Subsection~\ref{Su:RSIReg}) holds in $L$, no
analogue of an embedding into monoids of the form \eqref{Eq:MonI,II} had been
found before the present work.

\subsection{Lattices of projections of W*- and AW*-algebras}
\label{Su:ProjAW*}\index{W*-algebra}
We recall that an \emph{AW*-algebra}\index{AW*-algebra|ii}
is a C*-algebra $A$ such that the right annihilator of any subset $X$ of $A$
has the form $pA$, for a projection $p$ of $A$ (a \emph{projection} of $A$
is an element $p$ of~$A$ such that $p=p^2=p^*$). We denote by $L$ the set of
projections of $A$. Let $a\leq b$ hold, if $ab=a$ (equivalently, $ba=a$),
for all $a$, $b\in L$. Thus $\leq$ is a partial ordering on $L$.
Orthogonality of any projections $a$ and $b$, in notation
$a\perp b$, is defined by $ab=0$, and (Murray-von~Neumann) equivalence is
defined by the rule
   \[
   a\sim b\text{ if there exists }x\in A\text{ such that }a=x^*x\text{ and }
   b=xx^*.
   \]
Much of the structure of $L$ was developed axiomatically by I. Kaplansky
in his monograph \cite{Kapl68}.\index{Kaplansky, I.} Some of the axioms and
methods we use were inspired by Kaplansky's work, as was the structure
theory for nonsingular injective modules constructed by Goodearl and Boyle
\cite{GoBo}.
\index{Goodearl, K.\,R.}\index{Boyle, A.\,K.}%
Readers familiar with either of those works will recognize the parallels
below.

Again, the quotient $L/{\sim}$ can be endowed with a structure of
\pcm,\index{partial commutative monoid} where addition is given by the rule
   \[
   [a]+[b]=[a+b],\quad\text{if }ab=0,\text{ for all }a,\,b\in L,
   \]
where $[p]$ denotes the $\sim$-equivalence class of a projection $p$ of $A$.
The amount of known general information on the structure of $L/{\sim}$ is more
fragmentary than for the examples considered in previous sections, due to
fewer axioms satisfied. For example, the analogues of properties (3)
(meet-continuity) and (4) (modularity) considered in Subsection~\ref{Su:RSIReg}
fail for projections of AW*-algebras as a rule. In F.\,J.
Murray and J. von~Neumann \cite{MuNe36},
\index{Murray, F.\,J.}\index{von Neumann, J.}%
a $[0,+\infty]$-valued ``dimension function'' is
constructed on the projections of any \emph{W*-factor} (i.e.,
indecomposable von Neumann algebra);\index{factor (W*-)|ii}
Kaplansky showed that the same
construction could be carried out for AW*-factors. Still in the indecomposable
case, it is known that the closed two-sided ideals are well-ordered, see
F.\,B. Wright \index{Wright, F.\,B.}\cite{Wrig58}. Most of what was known
about $L/{\sim}$ in the general case could be obtained from more general,
often lattice-theoretical works that we shall discuss now.

\subsection{Lattice-theoretical generalizations}\label{Su:LattGen}
A common feature of the structures considered in Subsections
\ref{Su:AMeasTh}--\ref{Su:ProjAW*} is that they all involve a complete,
sectionally complemented lattice $L$,
\index{lattice!complete ---}\index{lattice!sectionally complemented ---}%
a binary relation $\perp$ on $L$, and
an equivalence relation $\sim$ on $L$. It has been observed early that even
apart from the classical study of\index{continous geometry} continuous
geometries, the dimension theory of a given structure could be done by just
studying the associated structure
$(L,\perp,\sim)$. Furthermore, these structures will be \emph{ordered}
structures, so that we shall write $(L,\leq,\perp,\sim)$ instead of
$(L,\perp,\sim)$:
\begin{itemize}
\item[---] It is in S. Maeda \cite{SMae55}\index{Maeda, S.} that the most
general axiomatization of the structures $(L,\leq,\perp,\sim)$ is given.
It holds for all the examples considered in Subsections
\ref{Su:AMeasTh}--\ref{Su:ProjAW*}, and this allows to construct ``dimension
functions''---corresponding to the measures of
Subsection~\ref{Su:AMeasTh}---on $L$ that, in the ``finite'' case, separate
the elements of $L$.

\item[---] In L.\,H. Loomis \cite{Loom55},\index{Loomis, L.\,H.} another
axiomatization is used, that involves an \emph{orthocomplementation} on $L$,
thus it does not apply to the examples considered in
Subsections~\ref{Su:RSIReg} and \ref{Su:CMSMLatt}.

\item[---] In P.\,A. Fillmore \index{Fillmore, P.\,A.}\cite{Fill65}, a
further axiomatization of the structures $(L,\perp\nobreak,\sim\nobreak)$ is
introduced, that does not assume completeness of $L$ but rather
\emph{countable} completeness, and that assumes an orthocomplementation
(thus, again, it does not encompass Subsections~\ref{Su:RSIReg} and
\ref{Su:CMSMLatt}). One of the main results is that the structure $L/{\sim}$
is a generalized cardinal algebra (as in Subsection~\ref{Su:CMSMLatt}).
Furthermore, under some countability assumptions, $L$ is complete and
$L/{\sim}$ is isomorphic to a lower subset of a monoid of the form
\eqref{Eq:MonI,II}, see \cite[Theorem~3.12]{Fill65}.\index{Fillmore, P.\,A.}
\end{itemize}

Nevertheless, in each class of examples considered in Subsections
\ref{Su:AMeasTh}--\ref{Su:ProjAW*}, some of the dimension-theoretical
properties that one could have expected to hold were still missing from the
known results. For example, there has been no general treatment of the
reducible Type~III case; it was seemingly not even clear whether or not it had
to be treated as a pathology.

\section{Results and methods}\label{S:ResMeth}

In view of the various examples presented in Section~\ref{S:Backgr} and of
what is known about them, the main goals of this paper are the
following:
\begin{itemize}
\item[(1)] To capture in a convenient set of axioms the various properties of
the structures $(L\nobreak,\leq\nobreak,\perp\nobreak,\sim\nobreak)$
encountered in these examples. This set of axioms should be sufficient to
develop a \emph{complete} dimension theory of these structures, that is, a
complete description of the structures $L/{\sim}$, without additional
assumptions such as finiteness or chain conditions.

\item[(2)] To develop a set of monoid-theoretical axioms that should be
satisfied by the structures $L/{\sim}$.

\item[(3)] Although the set of axioms obtained in (2) is quite complicated,
our third goal will be to give a simple description of the structures
satisfying the axioms of (2) in terms of continuous functions on complete
Boolean spaces.\index{Boolean space}
\end{itemize}

We shall now give some details about our road to these goals.

\subsection*{Espaliers}\index{espalier} The relevant structures
$(L,\leq,\perp,\sim)$ will be called \emph{espaliers}, see
Definition~\ref{D:MeasChLatt}. The axiom system defining espaliers is
stronger than the axiom system $(1,\alpha)$,
$(1,\beta)$, \dots, $(1,\zeta)$,
$(2,\alpha)$, \dots, $(2,\zeta)$ considered by S.~Maeda in\index{Maeda, S.}
\cite{SMae55}. Nevertheless, these axioms are sufficient for our
purposes---for instance, all the examples considered in
Section~\ref{S:Backgr} are espaliers.\index{espalier} The only drastic
generalization that we will introduce is to state that the underlying
partial ordering of an espalier
$(L,\leq,\perp,\sim)$ defines a
\emph{partial}, as opposed to total,
\index{lattice!partial ---} lattice, so that for elements $a$ and
$b$ of $L$, the meet (i.e., infimum)
$a\wedge b$ of $\set{a,b}$ always exists, but the join (i.e., supremum) of
$\set{a,b}$ exists only in case $\set{a,b}$ is majorized. This small
generalization affects neither the proofs nor even the results---the
structures $L/{\sim}$ are partial structures anyway---and it paves the way
for further algebraic constructions on espaliers, such as \emph{amalgamation}.

\subsection*{Continuous dimension scales}
\index{continuous dimension scale}
In parallel to this, we shall develop a system
of monoid-the\-o\-ret\-i\-cal axioms, (M1)--(M6) (see
Definition~\ref{D:DimInt}), that captures the structures (\pcm s)
$L/{\sim}$, for an espalier\index{espalier} $L$. This axiom system is
rather complicated, but it completely
isolates what monoid theory we need to understand the structures
$L/{\sim}$. Among these axioms is a variant of \emph{conditional completeness}
(see Axiom~(M2)), that is, every nonempty subset admits an infimum for the
algebraic (pre)ordering (see Definition~\ref{D:AlgPr}), but there are other,
less natural-looking axioms, such as (M6).

The \pcm s satisfying Axioms (M1)--(M6) will be called \emph{continuous
dimension scales}. They are unrelated to H. Lin's ``continuous scales''
\index{Lin, H.} introduced in \cite{Lin88,Lin91}.

The relation between espaliers\index{espalier} and
continuous dimension scales\index{continuous dimension scale} is then
given by the following (see Theorem~\ref{T:DimEsp}).

\begin{all}{Theorem A}
Let $(L,\leq,\perp,\sim)$ be an espalier.\index{espalier} Then the
\pcm\index{partial commutative monoid}\ $L/{\sim}$ of all
$\sim$-equivalence classes of elements of $L$ is a continuous dimension
scale.\index{continuous dimension scale}
\end{all}

\subsection*{Descriptions of continuous dimension
scales}\index{continuous dimension scale}
At first glance, Theorem~A may
appear as the ultimate goal of this paper. However, it provides only a
list of properties of the partial monoids
$L/{\sim}$, without giving any representation in terms of known structures.
Moreover, although the axioms describing the structure of espalier seem to
be almost the weakest possible to obtain a complete dimension theory, and
thus, in some sense, unavoidable, this might not seem to be
the case \emph{a priori} for the axioms describing continuous dimension
scales.\index{continuous dimension scale} We counter this by proving that
there are no ``missing'' axioms for continuous dimension
scales\index{continuous dimension scale} relative to espaliers.

\begin{all}{Theorem B}
A partial commutative monoid $S$ is a continuous dimension scale if and
only if $S\cong L/{\sim}$ for some espalier\index{espalier}
$(L,\leq,\perp,\sim)$.
\end{all}

\noindent Theorem~B follows from the fact that most
of our classes of examples of espaliers are
\emph{universal} in the sense that arbitrary continuous dimension
scales\index{continuous dimension scale} can be represented
(isomorphically) as lower subsets of the continuous dimension scales
$L/{\sim}$ arising from these examples---see, for example, Theorems
\ref{T:MeasDUniv}, \ref{T:meetcontinDuniv}, \ref{T:L(R)Duniv},
\ref{T:L(A)Duniv}.

As for a concrete representation of continuous dimension
scales\index{continuous dimension scale}, we exhibit them as lower
subsets (for the algebraic preordering) of product spaces of the form
   \[
   \CC(\Omega_{\I},\ZZ_\gamma)\times\CC(\Omega_{\II},\RR_\gamma)
   \times\CC(\Omega_{\III},\two_\gamma),
   \]
where $\Omega_{\I}$, $\Omega_{\II}$, and $\Omega_{\III}$ are complete Boolean
spaces\index{Boolean space} and, for any ordinal $\gamma$, the monoids
$\ZZ_\gamma$,\index{Zzzgamma@$\ZZ_\gamma$|ii}
$\RR_\gamma$,\index{Rzzgamma@$\RR_\gamma$|ii} and
$\two_\gamma$\index{Tzzgamma@$\two_\gamma$|ii} are defined as
   \begin{align*}
   \ZZ_\gamma&=\ZZ^+\cup\setm{\aleph_\xi}{0\leq\xi\leq\gamma},\\
   \RR_\gamma&=\RR^+\cup\setm{\aleph_\xi}{0\leq\xi\leq\gamma},\\
   \two_\gamma&=\set{0}\cup\setm{\aleph_\xi}{0\leq\xi\leq\gamma},
   \end{align*}
endowed with the natural addition and ordering, together with
the\index{interval topology} interval topology (see
Section~\ref{S:NotTerm}). See page~\pageref{Pg:ZR2gam} for more details.

\begin{all}{Theorem C}
Let $S$ be a \pcm.\index{partial commutative monoid}
Then $S$ is a continuous dimension scale\index{continuous dimension scale} \iff\ it can be
embedded as a lower subset into a product monoid of the form
   \[
   \CC(\Omega_{\I},\ZZ_\gamma)\times\CC(\Omega_{\II},\RR_\gamma)
   \times\CC(\Omega_{\III},\two_\gamma),
   \]
where $\Omega_{\I}$, $\Omega_{\II}$, and $\Omega_{\III}$ are complete Boolean
spaces.\index{Boolean space}
\end{all}

A more precise version of Theorem~C is formulated in
Theorem~\ref{T:EmbDimInt}. The concrete version of Theorem~B is that any lower
subset of a monoid of the form
   \[
   \CC(\Omega_{\I},\ZZ_\gamma)\times\CC(\Omega_{\II},\RR_\gamma)
   \times\CC(\Omega_{\III},\two_\gamma),
   \]
can be represented as $L/{\sim}$, for a suitable espalier\index{espalier}
$L$. More precisely, we show that~$L$ may arise from each of the above
contexts---abstract measure theory, nonsingular injective modules over
self-injective regular rings,
\index{ring!Von Neumann regular ---}\index{ring!right self-injective ---}%
meet-continuous complemented modular lattices,
\index{lattice!complemented modular ---}
and projection lattices of\index{AW*-algebra}
AW*-algebras, see Sections~\ref{S:AMeasTh}--\ref{S:ProjAW*}. For projection
lattices of\index{W*-algebra} W*-algebras, there is an additional
restriction on the spaces $\Omega_{\I}$, $\Omega_{\II}$,
$\Omega_{\III}$---namely, they are
\emph{hyperstonian},\index{hyperstonian} see
Corollary~\ref{C:W*hyperDuniv}. In addition, it is worth noticing that
although the embedding in Theorem~C is not unique as a rule, it is
determined by the condition that it ``commutes with projections'' and its
value at the elements of a \emph{finitary unit} of $S$
(Definition~\ref{D:FinUn}), see Theorem~\ref{T:Uneps}. Finally, all this
extends to ``continuous dimension scales''\index{continuous dimension scale} that are no
longer \emph{sets}, but rather \emph{proper classes}. The corresponding
common extensions of the abovementioned ``existence'' and ``uniqueness''
statements hold, and they are presented in Theorem~\ref{T:GenEmbDI}.

In order to make the results and methods of this paper accessible to the
widest audience, we have avoided the use of forcing and Boolean-valued
models for most proofs. Exceptions to this rule are the proofs of
D-universality\index{D-universal} for the classes of Boolean
espaliers\index{espalier!Boolean ---} (Theorem~\ref{T:MeasDUniv}) and
espaliers of projections of AW*-algebras\index{AW*-algebra}
(Theorem~\ref{T:bigAW*Drng}), as reasonable ``forcing-free'' proofs do not
seem to be available.

\section{Notation and terminology}\label{S:NotTerm}

Disjoint unions of sets will be denoted by $\sqcup$, $\bigsqcup$,
\index{unczzqsup@$\sqcup$, $\bigsqcup$|ii} so that, for
example, $X=\bigsqcup_{i\in I}X_i$ means that $X=\bigcup_{i\in I}X_i$ and
that $X_i\cap X_j=\es$, for all distinct $i$, $j\in I$.

Following the usual set-the\-o\-ret\-i\-cal terminology, we denote by
\index{ozzmega@$\omega$|ii}
$\omega$ the set of all natural numbers. We identify any natural number
$n$ with $\set{0,1,\dots,n-1}$.
Any ordinal $\alpha$ is identified with the set of all ordinals less than
$\alpha$. A cardinal is an initial ordinal. Following well-established
set-the\-o\-ret\-i\-cal practice, for an ordinal~$\alpha$, the notations
$\omega_\alpha$ and $\aleph_\alpha$ both denote the $\alpha$-th infinite
cardinal, except that the first one is viewed as an ordinal while the second
one is viewed as a cardinal.

If $P$ is a partially preordered set, a subset $X$ of $P$ is a \emph{lower
subset} (resp., \emph{upper subset}) of $P$ if $x\leq y$ and $y\in X$ (resp.,
$x\in X$) implies that $x\in X$ (resp., $y\in X$), for all $x$, $y\in P$.
For an element $a$ of $P$, we denote by $(a]$ (resp., $[a)$) the lower subset
(resp., upper subset) of $P$ generated by $a$. A subset $X$ of $P$ is
\emph{coinitial},\index{coinitial (subset)|ii}
if $[X)=\bigcup_{a\in X}[a)$ is equal to $P$. If $P$ has a
least element $0$, a subset $X$ of $P$ is \index{dense|ii}\emph{dense} in
$P$, if $X\setminus\set{0}$ is coinitial in $P\setminus\set{0}$. We say that
$X$ is an \emph{antichain}\index{antichain|ii} of $P$, if $0\notin X$ and
$(a]\cap(b]=\set{0}$ for any distinct $a$, $b\in X$.
If $X$ and $Y$ are subsets of $P$, we abbreviate the statement
   \[
   \forall (x,y)\in X\times Y,\ x\leq y
   \]
by $X\leq Y$. If $X=\set{a}$ (resp., $Y=\set{a}$), we write $a\leq Y$
(resp., $X\leq a$).

The \emph{interval topology}\index{interval topology|ii} on $P$ is the least
topology of $P$ for which all the intervals of the form $(a]$ or $[a)$, for
$a\in P$, are closed.

We shall consider the interval topology only in the totally ordered, complete
case. The relevant result is then the following, see, for example,
\index{Birkhoff, G.}\cite[\S X.12]{Birk}.

\begin{proposition}[O. Frink]\label{P:Frink}
\index{Frink, O.}
Let $(E,\leq)$ be a totally ordered set. We suppose that $E$ is
\emph{complete}, that is, every subset of $E$ has an infimum in $E$.
Then the interval topology of $E$ is compact Hausdorff.
\end{proposition}

If $G$ is a partially ordered group, $G^+$
\index{gzzplus@$G^+$|ii}
denotes the positive cone of $G$. We put\index{nzzn@$\NN$|ii}
$\NN=\ZZ^+\setminus\set{0}$.
We say that $G$ is
\index{partially ordered abelian group!directed ---|ii} \emph{directed},
if is upward directed as a partially
ordered set; we say that $G$ satisfies the
\index{interpolation property|ii} \emph{interpolation property}, if for all
$a_0$, $a_1$, $b_0$, $b_1\in G$ such that $a_0,a_1\leq b_0,b_1$, there exists
$x\in G$ such that $a_0,a_1\leq x\leq b_0,b_1$. We say that $G$ is
\index{partially ordered abelian group!Dedekind complete ---|ii}
\emph{Dedekind complete}, if it is directed and every nonempty majorized
subset of $G$ has a supremum. It is well-known that every Dedekind complete
partially ordered group is abelian, see, for example,\index{Birkhoff, G.}
\cite[Theorem~28]{Birk}. We shall write such groups using additive
notation.

For any point $x$ in a topological space $\Omega$, we denote by
$\Nh^\Omega(x)$\index{Nzzh@$\Nh^\Omega(x)$|ii}
(or $\Nh(x)$ if $\Omega$ is understood) the set of all open
neighborhoods of $x$ in $\Omega$. For a subset $X$ of $\Omega$, we denote by
$\overset{\,\circ}{X}$ the \emph{interior}\index{interior|ii}
\index{Xzzinterior@$\overset{\,\circ}{X}$|ii} of $X$ and by
$\ol{X}$ the \emph{closure}\index{closure|ii}
\index{Xzzclosure@$\ol{X}$|ii} of $X$ in $\Omega$. If $K$ is a
totally ordered set, endowed with its interval topology, a map
$f\colon\Omega\to K$ is \index{semicontinuous!lower ---|ii}
\emph{\lsc} (resp., \index{semicontinuous!upper ---|ii}\emph{\usc}), if the
set $\setm{x\in\Omega}{f(x)\leq\alpha}$ (resp.,
$\setm{x\in\Omega}{\alpha\leq f(x)}$) is closed, for every $\alpha\in K$.

A topological space $\Omega$ is \emph{extremally disconnected}, if the
closure of every open subset of $\Omega$ is open.
\index{extremally disconnected|ii}
We use the terminology \emph{complete Boolean space}\index{Boolean space} as
a synonym for \emph{extremally disconnected compact Hausdorff topological
space}. See Section~\ref{S:cmpBooleanspc} for more detail on these concepts.
Complete Boolean spaces\index{Boolean space} are also called
\emph{Stone spaces}\index{Stone space|ii} (or \emph{stonian spaces}) in the
literature.

\chapter{Partial commutative monoids}\label{Ch:PartMon}

\section[Basic results]
{Basic results about partial commutative monoids}\label{S:BasicResults}

\subsection{Partial commutative monoids}
Many monoid-theoretical objects we shall deal with through this paper are not
monoids, but just \emph{partial} monoids. The following fundamental example
provides us with a large supply of partial monoids.

\begin{example}\label{Ex:FundPartMon}
Let $(M,+,0)$ be a \cm. For a subset $S$ of~$M$ satisfying the two
following properties
\begin{enumerate}
\item $0\in S$;

\item $x+y\in S$ implies that $x$, $y\in S$, for all $x$, $y\in M$,
\end{enumerate}
we endow $S$ with the partial addition $+_S$ defined by
   \[
   a+_Sb=c,\text{ only in case }c\in S,
   \]
for any $a$, $b\in S$. We call $S$ a \emph{partial submonoid} of $M$.
\end{example}

Observe that we do not merely consider \emph{all} subsets of $M$, but only
those that satisfy the conditions (i) and (ii) above---they are exactly the
nonempty lower subsets of $M$ for the \emph{algebraic preordering} of $M$,
see Definition~\ref{D:AlgPr}.

It turns out that the properties of partial submonoids of \cm s are captured
by the following definition.

\begin{definition}\label{D:PartCM}
A \emph{\pcm}\index{partial commutative monoid|ii}
is a structure $(S,+,0)$, where $+$ is a partial binary
operation on $S$ which satisfies the following properties:
\begin{itemize}

\item[(a)] $+$ is \emph{associative}, that is, for all $a$, $b$, $c\in S$,
$(a+b)+c$ is defined \iff\ $a+(b+c)$ is defined, and then, both have the
same value.

\item[(b)] $+$ is \emph{commutative}, that is, for all $a$, $b\in S$,
$a+b$ is defined \iff\ $b+a$ is defined, and then, both have the same
value.

\item[(c)] There exists an element, denoted by $0$ (\emph{necessarily
unique}), of~$S$ such that $a+0=a$, for all $a\in S$.
\end{itemize}
\end{definition}

We generalize to this context the classical definition of the algebraic
preordering on a \cm.\index{algebraic preordering|ii}

\begin{definition}\label{D:AlgPr}
Let $(S,+,0)$ be a \pcm. The
\emph{algebraic preordering} on $S$ is the (reflexive, transitive) binary
relation $\leq$ defined on $S$ by the rule
   \[
   a\leq b\text{ \iff\ }a+x=b,\text{ for some }x\in S.
   \]
An element $u\in S$ is an \emph{order-unit}, if every element of $S$ lies
below $nu$ (defined), for some $n\in\ZZ^+$.\index{order-unit|ii}
\end{definition}

The following definition is of course a direct generalization of
Example~\ref{Ex:FundPartMon}.

\begin{definition}\label{D:PartSubm}
A \emph{partial submonoid} of a \pcm\index{partial commutative monoid}~$S$
is a lower subset $T$ of~$S$ (for the algebraic preordering of~$S$)
containing $0$ as an element, endowed with the partial addition defined by
   \[
   a+b=c\text{ \iff\ }a+b=c\text{ in }S\text{ and }c\in T,
   \text{ for all }a,\,b\in T.
   \]
\end{definition}

We omit the trivial proof of the following result.

\begin{proposition}
Every partial submonoid \pup{as in Example~\textup{\ref{Ex:FundPartMon}}} of
a \pcm\index{partial commutative monoid}\ is a \pcm.
\end{proposition}

The following class of embeddings will be of special interest.

\begin{definition}\label{D:LowEmb}
Let $A$ and $B$ be \pcm s,\index{partial commutative monoid} and
let\linebreak $\varphi\colon A\to B$. We say that
$\varphi$ is a \emph{lower embedding},
\index{lower embedding|ii} if the following conditions hold:
\begin{enumerate}
\item $\varphi$ is a homomorphism of partial monoids.

\item $\varphi$ is one-to-one, and $\varphi(x)\leq\varphi(y)$ implies that
$x\leq y$, for all $x$, $y\in A$.

\item The range of $\varphi$ is a lower subset of $B$, with respect to the
algebraic preordering of $B$.
\end{enumerate}
\end{definition}

Hence, a lower embedding from $A$ into $B$ identifies $A$ with a lower subset
(with respect to the algebraic preordering) of $B$, endowed with the
structure of partial submonoid as in Definition~\ref{D:PartSubm}.

The following result shows that all \pcm s\index{partial commutative
monoid} can be obtained from Example~\ref{Ex:FundPartMon}.

\begin{proposition}\label{P:PcmCm}
Every \pcm\index{partial commutative monoid}\ admits a lower embedding into
a \cm.
\end{proposition}

\begin{proof}
Let $(S,+,0)$ be a \pcm.\index{partial commutative monoid} Let $\infty$ be
any object such that
$\infty\notin S$, and put\index{Szzbul@$\Sbul$|ii}
$\Sbul=S\cup\set{\infty}$. We define on
$\Sbul$ the binary operation $\blp$\index{abzzzlp@$a\blp b$|ii} defined by
the rule
   \[
   a\blp b=\begin{cases}
   a+b,&\text{if }a,\,b\in S\text{ and }a+b\text{ is defined in }S,\\
   \infty,&\text{otherwise},
   \end{cases}
   \quad\text{for all }a,\,b\in\Sbul.
   \]
It is easy to verify that $(\Sbul,\blp,0)$ is a \cm\ and that the inclusion
map from $S$ into $\Sbul$ is a lower embedding.
\end{proof}

\begin{remark}
A noticeable effect of Proposition~\ref{P:PcmCm} is to make computations in
\pcm s\index{partial commutative monoid} much more convenient. For example,
suppose that we have to prove that an equality of the form $A=B$ holds in a
given \pcm\index{partial commutative monoid}\ $S$, \emph{via} a sequence of
equalities
$A=C_0=C_1=\dots=C_n=B$, where $A$, $B$, and the $C_i$ are finite sums of
elements of~$S$. We assume in addition that the sum defining $A$ is defined
in $S$. Instead of having to verify that all the terms $C_i$ are defined in
$S$ and pairwise equal, it is sufficient to argue in $\Sbul$ that
$A=C_0=C_1=\dots=C_n=B$, without having to worry about undefined terms.
\end{remark}

This applies, in particular, to the following Lemmas
\ref{L:leqStillDef}, \ref{L:InvPermPlus}, and \ref{L:SetAssoc}.

\begin{lemma}\label{L:leqStillDef}
Let $(S,+,0)$ be a partial monoid, with algebraic preordering $\leq$. Let
$a$, $b$, $a'$, $b'\in S$. If $a+b$ is defined and $a'\leq a$ and
$b'\leq b$, then $a'+b'$ is defined, and $a'+b'\leq a+b$.
\end{lemma}

In any given \pcm\index{partial commutative monoid}\ $S$, we define
inductively the statement
$a=\sum_{i<n}a_i$ to hold, for $n<\omega$, $a$, $a_0$, \dots,
$a_{n-1}\in S$, as follows:
\begin{enumerate}
\item $a=\sum_{i<0}a_i$ \iff\ $a=0$.

\item $a=\sum_{i<n+1}a_i$ \iff\ $a=\left(\sum_{i<n}a_i\right)+a_n$.
\end{enumerate}

If the operation of~$S$ is denoted by $\oplus$, then we shall write
$\oplus_{i<n}a_i$ instead of~$\sum_{i<n}a_i$.

\begin{lemma}\label{L:InvPermPlus}
Let $(S,+,0)$ be a \pcm.\index{partial commutative monoid}
For all $n<\omega$, all $a$, $a_0$,\dots, $a_{n-1}\in S$, and every
permutation $\sigma$ of $n$,
   \[
   a=\sum_{i<n}a_i\qquad\text{\iff}\qquad a=\sum_{i<n}a_{\sigma(i)}.
   \]
\end{lemma}

By Lemma~\ref{L:InvPermPlus}, for a finite set $I$ and elements $a$, $a_i$
(for $i\in I$) of~$S$, we can define unambiguously the statement
$a=\sum_{i\in I}a_i$ to hold, if $a=\sum_{j<n}a_{\sigma(j)}$, where $n$
is the cardinality of $I$ and $\sigma$ is any bijection from $n$ onto $I$.

\begin{lemma}\label{L:SetAssoc}
Let $(S,+,0)$ be a \pcm.\index{partial commutative monoid} Let $I$ and $J$
be finite sets, let
$\pi\colon I\twoheadrightarrow J$ be a surjective map, let
$\famm{a_i}{i\in I}$ be a family of elements of~$S$, and let $a\in S$.
Then the following are equivalent:

\begin{enumerate}
\item $a=\sum_{i\in I}a_i$.

\item For all $j\in J$, the term $\sum_{i\in\pi^{-1}\set{j}}a_i$ is defined,
and, if we denote its value by~$b_j$, then $a=\sum_{j\in J}b_j$.
\end{enumerate}

\end{lemma}

\subsection{Partial refinement monoids}
\index{refinement monoid!partial ---}
\begin{definition}\label{D:RefPpty}
We say that a \pcm\index{partial commutative monoid}\ $(S,+,0)$ has the
\emph{refinement property}, or is a
\emph{\prm}, if for all $a_0$, $a_1$, $b_0$, $b_1\in S$ such that
$a_0+a_1=b_0+b_1$, there are elements $c_{i,j}$ in $S$, for $i$, $j<2$, such
that the equalities $a_i=c_{i,0}+c_{i,1}$ and $b_i=c_{0,i}+c_{1,i}$ hold,
for all $i<2$.
\end{definition}

The information contained in the equalities $a_i=c_{i,0}+c_{i,1}$ and
$b_i=c_{0,i}+c_{1,i}$ for all $i<2$ will often be condensed in the format of
a \emph{refinement matrix} as follows:
   \[
   \begin{tabular}{|c|c|c|}
   \cline{2-3}
   \multicolumn{1}{l|}{} & $b_0$ & $b_1$\tvi\\
   \hline
   $a_0$ & $c_{0,0}$ & $c_{0,1}$\tvi\\
   \hline
   $a_1$ & $c_{1,0}$ & $c_{1,1}$\tvi\\
   \hline
   \end{tabular}
   \]
These notations can be easily generalized to refinement matrices of
arbitrary, finite or even infinite, dimensions.
These notations are also very widely used in\index{Wehrung, F.} \cite{WDim}.

Define a \emph{refinement monoid}\index{refinement monoid|ii} as a
\cm\ satisfying the refinement property. In a spirit similar to
Proposition~\ref{P:PcmCm}, we shall now prove (see
Proposition~\ref{P:PrmRm}) that every
\prm\ can be obtained as a lower subset of a refinement monoid. The proof of
Proposition~\ref{P:PcmCm} does not apply for this result, because
$\Sbul$ fails in general to satisfy refinement even if $S$ has refinement. We
shall use instead a procedure adapted to refinement monoids.

\begin{proposition}\label{P:PrmRm}
Every \prm\ $S$ admits a lower embedding\index{lower embedding} into a
refinement monoid\index{SzzRef@$\Ref S$|ii} $\Ref S$. In addition, one can
take $\Ref S$ to be generated by $S$ as a monoid, and such that the canonical
embedding from $S$ into $\Ref S$ is universal among the homomorphisms of
partial monoids from $S$ to \cm s.
\end{proposition}

\begin{proof}
The following construction is a particular case of the construction presented
in Chapter~4 of\index{Wehrung, F.} \cite{WDim}---with the notations used
there, $\Ref S=\Dim(S,+,=)$. However, in this context, a direct verification
is easy, so we give an outline here.

Let $S$ be a \prm. We endow the set $\mathbb{S}$ of all finite,
nonempty sequences of elements of~$S$ with the binary relation $\equiv$
defined by the rule
   \begin{multline*}
   \famm{a_i}{i<m}\equiv\famm{b_j}{j<n}\text{ if there are }c_{i,j}\in S,
   \text{ for }i<m\text{ and }j<n,\text{ such that }\\
   a_i=\sum_{j<n}c_{i,j},\text{ for all }i<m,\text{ and }
   b_j=\sum_{i<n}c_{i,j},\text{ for all }j<n.
   \end{multline*}
By using the refinement property in $S$, it is not difficult to verify that
$\equiv$ is an equivalence relation on $\mathbb{S}$. For any
$s\in\mathbb{S}$, we denote by $[s]$ the equivalence class of $s$ modulo
$\equiv$. We endow the quotient $\Ref S=\mathbb{S}/{\equiv}$ with the binary
addition $+$ defined by\index{stczzonc@$s\conc t$|ii}
   \[
   [s]+[t]=[s\conc t],\text{ for all }s,\,t\in\mathbb{S},
   \]
where $s\conc t$ denotes the \emph{concatenation} of $s$ and $t$. It is
straightforward to verify that $\Ref S$, endowed with $+$, is
a\index{refinement monoid} refinement monoid. For any $a\in S$,
we denote by $j(a)$ the equivalence class modulo
$\equiv$ of the finite sequence $(a)$ of length one. Then $j(0)$ is the zero
element of $\Ref S$, and $j$ is a lower embedding\index{lower embedding}
from $S$ into $\Ref S$.

For the remainder of the proof, we shall identify $S$ with its image
in $\Ref S$ under the embedding $j$. Thus the elements of $\Ref S$ are exactly
the finite sums $\sum_{i<m}a_i$, where $m\in\NN$ and
$a_0$, \dots, $a_{m-1}\in S$, and the equality $\sum_{i<m}a_i=\sum_{j<n}b_j$
holds \iff\ there exists a refinement matrix of the form
   \[
   \begin{tabular}{|c|c|}
   \cline{2-2}
   \multicolumn{1}{l|}{} & $b_j\ (j<n)$\tvi\\
   \hline
   $a_i\ (i<n)$ & $c_{i,j}$\tvi\\
   \hline
   \end{tabular}
   \]
for some elements $c_{i,j}$ (for $i<m$ and $j<n$) of~$S$.
Obviously $\Ref S$ is generated by $S$ as a monoid.

Now we verify the second assertion of Proposition~\ref{P:PrmRm}. Let $M$ be a
\cm\ and let $f\colon S\to M$ be a homomorphism of partial monoids. Let
$\oll{g}\colon\mathbb{S}\to M$ be the map defined by the rule
   \[
   \oll{g}(\famm{a_i}{i<n})=\sum_{i<n}f(a_i),\text{ for all }
   \famm{a_i}{i<n}\in\mathbb{S}.
   \]
Then $\oll{g}(j(0))=0_M$,
$\oll{g}(s\conc t)=\oll{g}(s)+\oll{g}(t)$, and
$s\equiv t$ implies that $\oll{g}(s)=\oll{g}(t)$, for all $s$,
$t\in\mathbb{S}$. Hence $\oll{g}$ can be factored through
$\equiv$, thus yielding a homomorphism $g\colon\Ref S\to M$ that extends $f$.
Since $S$ generates $\Ref S$ as a monoid, $g$ is the only homomorphism with
this property.
\end{proof}

For any \prm\ $S$, we take $\Ref S$ to be the refinement
monoid\index{refinement monoid} having all the properties described in
Proposition~\ref{P:PrmRm}. By the given universal property, $\Ref S$ is
unique up to isomorphism. We will also identify $S$ with its canonical image
in $\Ref S$.

Our next lemma collects some basic information about $\Ref S$. For
$n\in\NN$, we put
   \[
   nS=\Setm{\sum_{i<n}x_i}{x_0,\ldots,x_{n-1}\in S}\subseteq\Ref S.
   \]

\begin{definition}\label{D:Conical}
A \prm\ $S$ is \emph{conical},\index{conical (\prm)|ii}
if $x+y=0$ implies that $x=y=0$,
for all $x$, $y\in S$. In other words, $x\leq 0$ implies that $x=0$, for
all $x\in S$.
\end{definition}

\begin{lemma}\label{L:BasicRefS}
Let $S$ be a \prm. Then the following assertions hold:
\begin{enumerate}
\item $nS$ is a lower subset of $\Ref S$, for all
$n\in\NN$.

\item If $S$ is\index{cancellative|ii} cancellative
\pup{i.e., $a+c=b+c$ in $S$ implies that $a=b$}, then $\Ref S$ is
cancellative.

\item If $S$ is conical, then $\Ref S$ is conical.
\end{enumerate}
\end{lemma}

\begin{proof}
(i) is an easy consequence of refinement in $\Ref S$.

(ii) Folklore. A proof can be found in Lemma~3.6 in\index{Wehrung, F.}
\cite{WDim}.

(iii) Let $x$, $y\in\Ref S$ such that $x+y=0$. Write $x=\sum_{i<m}x_i$ and
$y=\sum_{j<n}y_j$ for some $m$, $n\in\NN$ and $x_0$, \dots, $x_{m-1}$,
$y_0$, \dots, $y_{n-1}\in S$. Then, for all $i<m$, $x_i\leq x\leq x+y=0$ in
$\Ref S$, thus, since $S$ is a lower subset of $\Ref S$, $x_i\leq 0$ in
$S$. Hence, as $S$ is conical, $x_i=0$, so $x=0$. Hence $y=0$.
\end{proof}

\section[Direct decompositions]{Direct decompositions of partial refinement
monoids}\label{S:DirectDecompRefMon}
\index{refinement monoid!partial ---}
\begin{quote}
\em In this section, we shall fix a conical \prm~$S$. We shall denote by
$\leq$ the algebraic preordering of~$S$.
\end{quote}

\begin{definition}\label{D:Ideal}
An \emph{ideal} of~$S$\index{ideal (of a \prm)|ii}
is a nonempty subset $I$ of~$S$ such that $a+b\in I$
\iff\ $a\in I$ and $b\in I$, for all $a$, $b\in S$ such that $a+b$ is
defined.
\end{definition}

We define elements $a$ and $b$ of~$S$ to be \emph{orthogonal}, in
notation, $a\perp b$,
\index{orthogonal (elements)|ii}\index{abozzrthab@$a\perp b$|ii}%
if $x\leq a,b$ implies that $x=0$, for
all $x\in S$. If $X$ and $Y$ are subsets of~$S$, then we define $X\perp Y$
to hold if $x\perp y$ for all $(x,y)\in X\times Y$.
We shall put\index{Xzzbot@$X^\bot$|ii}
   \begin{equation}\label{Eq:DefXbot}
   X^\bot=\setm{s\in S}{s\perp x,\text{ for all }x\in X},
   \qquad\text{for all }X\subseteq S.
   \end{equation}
If $X=\set{x}$, a singleton, then we write $x^\bot$ instead of
$\set{x}^\bot$. In particular, $X\perp Y$ \iff\ $Y\subseteq X^\bot$, \iff\
$X\subseteq Y^\bot$.

\begin{lemma}\label{L:XbotId}\hfill
\begin{enumerate}
\item $a\perp c$ and $b\perp c$ implies that $a+b\perp c$,
for all $a$, $b$, $c\in S$ such that $a+b$ is defined.

\item The set $X^\bot$ is an ideal of~$S$, for all $X\subseteq S$.

\item $a\perp b$ and $a$, $b\leq c$ implies that $a+b$ is defined and
$a+b\leq c$.
\end{enumerate}

\end{lemma}

\begin{proof}
(i) Let $x\leq c,a+b$. By refinement, there are $a'$,
$b'\in S$ such that $a'\leq a$, $b'\leq b$, and $x=a'+b'$. So $a'\leq a,c$,
whence
$a'=0$. Similarly,
$b'=0$, so $x=0$, thus proving $a+b\perp c$.

(ii) is an obvious consequence of (i).

(iii) By the definition of $\leq$, there are $a'$, $b'\in S$ such that
$a+a'=b+b'=c$. By applying refinement to the equality $a+a'=b+b'$ and by
using the assumption that $a\perp b$, we obtain $t\in S$ such that $a'=b+t$
and $b'=a+t$. Since $a+a'$ is defined, $a+b$ is defined, and
$c=a+b+t\geq a+b$.
\end{proof}

\begin{notation}\label{Not:SumSubsets}
For $n\in\NN$ and $X_0$, \dots, $X_{n-1}\subseteq S$, we put
   \[
   X_0+\cdots+X_{n-1}=\Setm{x\in S}
   {\exists(x_0,\ldots,x_{n-1})\in X_0\times\cdots\times X_{n-1},\quad
   x=\sum_{i<n}x_i}.
   \]
\end{notation}
\index{XzzplusY@$X_0+\cdots+X_{n-1}$, $nX$|ii}

We shall also write $\sum_{i<n}X_i$ instead of $X_0+\cdots+X_{n-1}$.
If $X_i=X$ for all $i$, then we shall abbreviate this further by $nX$.
If $X_i\perp X_j$ for all $i\neq j$, then we shall write
$X_0\oplus\cdots\oplus X_{n-1}$, or $\bigoplus_{i<n}X_i$, instead of
$\sum_{i<n}X_i$, and we shall say that the sum of the $X_i$ is
\emph{orthogonal}.

\begin{lemma}\label{L:OplusId}
Let $n\in\NN$ and let $S_i$, for $i<n$, be nonempty subsets of~$S$ such that
$S=\bigoplus_{i<n}S_i$. Then the following hold:

\begin{enumerate}

\item $S_i=\left(\bigoplus_{j\neq i}S_j\right)^\bot$, for all $i<n$.
In particular, $S_i$ is an ideal of~$S$.

\item For all $x\in S$, there exists a unique decomposition
$x=\sum_{i<n}x_i$ such that $x_i\in S_i$ for all $i<n$.

\end{enumerate}

\end{lemma}

\begin{proof}
(i) The sum of all the $S_j$ is orthogonal, thus so is the sum of all $S_j$,
for $j\neq i$. Furthermore, $S_i\perp S_j$, for all $j\neq i$, so
$S_j\subseteq S_i^\bot$. Hence, by using Lemma~\ref{L:XbotId}(ii),
$\bigoplus_{j\neq i}S_j\subseteq S_i^\bot$. Conversely, let $x\in S_i^\bot$.
By assumption, there exists a decomposition $x=\sum_{j<n}x_j$, where
$x_j\in S_j$, for all $j<n$. But $x_i\in S_i$, thus $x\perp x_i$; whence
$x_i=0$, so $x\in\bigoplus_{j\neq i}S_j$. Hence
$S_i=\left(\bigoplus_{j\neq i}S_j\right)^\bot$. By Lemma~\ref{L:XbotId}(ii),
it follows that $S_i$ is an ideal of $S$.

(ii) Suppose $x=\sum_{i<n}x_i=\sum_{i<n}y_i$, with elements
$x_i$, $y_i\in S_i$, for all $i<n$. Since $S$ satisfies refinement, there
exists a refinement matrix of the form
   \[
   \begin{tabular}{|c|c|}
   \cline{2-2}
   \multicolumn{1}{l|}{} & $y_j\ (j<n)$\tvi\\
   \hline
   $x_i\ (i<n)$ & $z_{i,j}$\tvi\\
   \hline
   \end{tabular}
   \]
with elements $z_{i,j}\in S$, for all $i$, $j<n$. But if $i\neq j$, then
$S_i\perp S_j$, whence $z_{i,j}=0$. Hence, $x_i=z_{i,i}=y_i$, for all $i<n$.
\end{proof}

\begin{remark}\label{R:EmbdSSi}
The direct product $\prod_{i<n}S_i$ can be naturally endowed with a
structure of partial monoid, by defining the addition componentwise. In the
context of Lemma~\ref{L:OplusId}, we obtain a map
   \[
   \varphi\colon S\to\prod_{i<n}S_i,\quad x\mapsto\famm{x_i}{i<n}
   \text{ with }x=\sum_{i<n}x_i,\quad x_i\in S_i\text{ for all }i<n.
   \]
This map is a one-to-one homomorphism of partial monoids. However, it is
\emph{not}, in general, surjective: for arbitrary $x_i\in S_i$, for $i<n$,
the sum $\sum_{i<n}x_i$ may not be defined. But of course, if $S$ is a
(total) monoid, then $\varphi$ is an isomorphism.
\end{remark}

\begin{proposition}\label{P:EmbdSSi}
In the context of Remark~\textup{\ref{R:EmbdSSi}}, $\varphi$ is a lower
embedding\index{lower embedding} from $S$ into $\prod_{i<n}S_i$.
\end{proposition}

\begin{proof}
Only part (iii) of the definition of a lower embedding\index{lower
embedding} is not completely trivial.

Let $x\in S$ and $\famm{y_i}{i<n}\in\prod_{i<n}S_i$ such that
$\famm{y_i}{i<n}\leq\varphi(x)$.
Put $\varphi(x)=\famm{x_i}{i<n}$, so $y_i\leq x_i$, for all $i<n$. By the
definition of $\varphi$, $\sum_{i<n}x_i$ is defined, and equal to $x$. By
Lemma~\ref{L:leqStillDef}, $\sum_{i<n}y_i$ is also defined. Denote its value
by $y$. By the definition of $\varphi$, $\famm{y_i}{i<n}=\varphi(y)$.
\end{proof}

\section{Projections of partial refinement monoids}\label{S:ProjPCM}
\index{refinement monoid!partial ---}
\begin{quote}
\em Standing hypothesis: $S$ is a conical \prm. We denote again by $\leq$ the
algebraic preordering of~$S$.
\end{quote}

\begin{definition}\label{D:Proj}
A \emph{projection}\index{projection (of a \prm)|ii}
of~$S$ is an endomorphism $p$ of $(S,+,0)$ such that
   \[
   x\in p(x)+(pS)^\bot,\qquad\text{for all }x\in S.
   \]
\end{definition}

In particular, if $p$ is a projection of~$S$, then $p(x)\leq x$, for all
$x\in S$. Thus $p(0)=0$. Furthermore, $p$ preserves the algebraic preordering
of~$S$, see Definition~\ref{D:AlgPr}.

\begin{proposition}\label{P:CharProj}
Let $p$ be an endomorphism of~$S$. Then the following are equivalent:
\begin{enumerate}
\item $p$ is a projection of~$S$.

\item There are ideals $S_0$ and $S_1$ of~$S$ such that
\begin{enumerate}
\item $S=S_0\oplus S_1$.

\item $p(x_0+x_1)=x_0$, for all $(x_0,x_1)\in S_0\times S_1$ such that
$x_0+x_1$ is defined.
\end{enumerate}
\end{enumerate}

\end{proposition}

\begin{proof}
(ii)$\Rightarrow$(i) is easy.

(i)$\Rightarrow$(ii) Assume (i).
We put $S_0=pS$ and $S_1=(pS)^\bot$. By the definition of a projection,
$S=S_0+S_1$. Since $S_0\perp S_1$, it follows that $S=S_0\oplus S_1$. In
particular, $S_0$ and $S_1$ are ideals of~$S$ (see Lemma~\ref{L:OplusId}(i)).
For $x\in S$, let $y\in S_1$ such that $x=p(x)+y$. If $x=x_0+x_1$ in $S$ such
that $x_i\in S_i$ for all $i<2$, then, by Lemma~\ref{L:OplusId}(ii),
$p(x)=x_0$ and $y=x_1$.
\end{proof}

\begin{corollary}
Every projection of~$S$ is idempotent.
\end{corollary}

In the context of Proposition~\ref{P:CharProj}(ii), we observe that $S_0=pS$
while $S_1=p^{-1}\set{0}=S_0^\bot=(pS)^\bot$. In particular, $p$ is determined
by $S_0$ alone, so we shall call $p$ the \emph{projection of~$S$ onto~$S_0$}.

\begin{definition}
A \emph{direct summand}\index{direct summand (of a \prm)|ii}
of~$S$ is a subset $X$ of~$S$ such that
$S=X\oplus Y$, for some $Y\subseteq S$.
\end{definition}

Of course, by Lemma~\ref{L:OplusId}, $X$ is then an
ideal of~$S$, and $Y=X^\bot$, so $S=X\oplus X^\bot$.

It follows that the direct summands of~$S$ are exactly the ranges of
the projections of~$S$.

Furthermore, by exchanging the roles of $S_0$ and $S_1$, we obtain another
projection, which we shall denote by
\index{pzzbot@$p^\bot$ ($p$ projection)|ii}
$p^\bot$. Formally, $p^\bot$ is the
unique projection of~$S$ such that $p^\bot S=(pS)^\bot$ and
$(p^\bot)^{-1}\set{0}=pS$. We observe that $p^{\bot\bot}=p$.

\begin{notation}
Let $\BB{S}$\index{pzzroj@$\BB{S}$|ii}
denote the set of projections of~$S$. We shall also often use
the notation
\index{pzzrojst@$\BBp{S}$|ii} $\BBp{S}=\BB{S}\setminus\set{0}$.
\end{notation}

We shall now study the structure of $\BB{S}$, towards
Proposition~\ref{P:ProjBa}.

\begin{lemma}\label{L:MeetCompPr}
Let $p$, $q\in\BB{S}$. The the following holds:

\begin{enumerate}
\item $S=(pS\cap qS)\oplus(pS\cap q^\bot S)
\oplus(p^\bot S\cap qS)\oplus(p^\bot S\cap q^\bot S)$.

\item Let $r$ denote the projection from $S$ onto $pS\cap qS$. Then
$r=pq=qp$.
\end{enumerate}

\end{lemma}

\begin{proof}
(i) It is obvious that all ideals $pS\cap qS$, $pS\cap q^\bot S$,
$p^\bot S\cap qS$, and $p^\bot S\cap q^\bot S$ are pairwise orthogonal.
Let $x\in S$. Since $S=pS+(pS)^\bot$, there exists a decomposition
$x=x_0+x_1$, where $x_0\in pS$ and $x_1\in(pS)^\bot$. For $i<2$,
$x_i\in qS+(qS)^\bot$, thus $x_i=x_{i,0}+x_{i,1}$, for some $x_{i,0}\in qS$
and $x_{i,1}\in(qS)^\bot$. Since $pS$ and $(pS)^\bot$ are ideals of~$S$,
$x_{0,0}\in pS\cap qS$, $x_{0,1}\in pS\cap q^\bot S$,
$x_{1,0}\in p^\bot S\cap qS$, and $x_{1,1}\in p^\bot S\cap q^\bot S$.
Observe that $x=x_{0,0}+x_{0,1}+x_{1,0}+x_{1,1}$.

(ii) Since both $p$ and $q$ act as the identity on $pS\cap qS$, so does
$pq$. Furthermore, $q(x)=0$ for all $x\in q^\bot S$ and $pq(x)=p(x)=0$
for all $x\in p^\bot S\cap qS$, so, $pq=r$. By symmetry, $r=qp$.
\end{proof}

We shall put $p\wedge q=pq=qp$, for all
\index{pzzroj@$\BB{S}$} $p$, $q\in\BB{S}$.

\begin{corollary}
The structure $(\BB{S},\wedge)$ is a semilattice
\pup{i.e., an idempotent \cm}.
\index{semilattice|ii}\index{pzzroj@$\BB{S}$}%
\end{corollary}

We endow $\BB{S}$\index{pzzroj@$\BB{S}$} with the partial ordering $\leq$
defined by
   \[
   p\leq q\text{ \iff\ }p\wedge q=p,\qquad\text{for all }p,\,q\in\BB{S}.
   \]
For this partial ordering, $p\wedge q$ is, of course, the infimum of
$\set{p,q}$. The least element of $\BB{S}$\index{pzzroj@$\BB{S}$} is $0$
(the zero map), while the greatest element of
$\BB{S}$\index{pzzroj@$\BB{S}$} is $\id_S$ (the identity on~$S$).

\begin{lemma}\label{L:pLeqOrtq}
Let $p$, $q\in\BB{S}$. Then\index{pzzroj@$\BB{S}$} the following holds:
\begin{enumerate}
\item $p\leq q$ \iff\ $pS\subseteq qS$ \iff\ $p(x)\leq q(x)$ holds, for all
$x\in S$.

\item $p\wedge q=0$ \iff\ $q\leq p^\bot$.

\end{enumerate}
\end{lemma}

\begin{proof}
(i) If $p\leq q$, then $pS=qpS\subseteq qS$.

Suppose now that $pS\subseteq qS$. Let $x\in S$.
The inequality $p(x)\leq q(x)$ holds for all $x\in pS$
(because then $p(x)=x=q(x)$) and for all $x\in p^\bot S$
(because then $p(x)=0\leq q(x)$), so it holds for all $x\in S$ since
$S=pS+p^\bot S$.

If $p(x)\leq q(x)$ for all $x\in S$, then $p(x)\in qS$ since $qS$ is an
ideal of~$S$, so $qp(x)=p(x)$. Hence $p\leq q$.

(ii) By Lemma~\ref{L:MeetCompPr}, $p\wedge q=0$ \iff\ $pS\cap qS=\set{0}$.
Hence, $p\wedge q=0$ \iff\ $qS\subseteq(pS)^\bot=p^\bot S$, \iff\
$q\leq p^\bot$ by (i) above.
\end{proof}

\begin{corollary}\label{C:BotAA}
The map $p\mapsto p^\bot$ is an involutive
anti-automorphism\index{pzzroj@$\BB{S}$} of
$(\BB{S},\leq)$.
\end{corollary}

\begin{proof}
We already know that $p^{\bot\bot}=p$, for all\index{pzzroj@$\BB{S}$}
$p\in\BB{S}$. Furthermore, by Lemma~\ref{L:pLeqOrtq}, $p\leq q$ implies that
$q^\bot\leq p^\bot$, for all $p$, $q\in\BB{S}$.
\end{proof}

Since $(\BB{S},\leq)$\index{pzzroj@$\BB{S}$} is a meet-semilattice, we thus
obtain the following.

\begin{corollary}\label{C:BSLatt}
$(\BB{S},\leq)$\index{pzzroj@$\BB{S}$} is a lattice.
\end{corollary}

So we denote by $p\vee q$ the supremum of $\set{p,q}$, for
all\index{pzzroj@$\BB{S}$} $p$,
$q\in\BB{S}$. We can strengthen Corollary~\ref{C:BSLatt} right away.

\begin{proposition}\label{P:ProjBa}
$(\BB{S},\leq)$\index{pzzroj@$\BB{S}$} is a\index{Boolean algebra} Boolean
algebra.
\end{proposition}

\begin{proof}
By Corollary~\ref{C:BSLatt}, $(\BB{S},\leq)$\index{pzzroj@$\BB{S}$} is a
lattice. Furthermore,
$p^\bot$ is a complement of $p$, for all\index{pzzroj@$\BB{S}$}
$p\in\BB{S}$. Hence, to conclude the proof, it suffices to prove
distributivity. The argument below is classical, and it can be traced back
to Glivenko's work, see, for example,
\cite[\S V.11]{Birk}.\index{Birkhoff, G.}

So, let\index{pzzroj@$\BB{S}$} $p$, $q$, $r\in\BB{S}$. We put
   \[
   s=p\wedge(q\vee r)\qquad\text{and}\qquad t=(p\wedge q)\vee(p\wedge r).
   \]
Then $t^\bot\wedge p\wedge q=t^\bot\wedge p\wedge r=0$, which implies, by
Lemma~\ref{L:pLeqOrtq}(ii), that $t^\bot\wedge p\leq q^\bot$ and
$t^\bot\wedge p\leq r^\bot$, thus, meeting both inequalities,
$t^\bot\wedge p\leq q^\bot\wedge r^\bot$. By Corollary~\ref{C:BotAA},
$q^\bot\wedge r^\bot=(q\vee r)^\bot$, so it follows that
$t^\bot\wedge p\wedge(q\vee r)=0$, that is, by Lemma~\ref{L:pLeqOrtq}(ii),
$t^\bot\leq s^\bot$; whence $s\leq t$. But the converse inequality $s\geq t$
is obvious, thus $s=t$.
\end{proof}

\begin{notation}
For\index{pzzroj@$\BB{S}$} $p$, $q$, $r\in\BB{S}$, let $r=p\oplus q$ hold
just in case $r=p\vee q$ and $p\wedge q=0$.
\end{notation}

\begin{lemma}\label{L:JoinDisj}
Let $p$, $q\in\BB{S}$ such that $p\wedge q=0$. Then
   \[
   (p\vee q)(x)=p(x)+q(x),\qquad\text{for all }x\in S.
   \]
\end{lemma}

\begin{proof}
Let $x\in S$, and put $r=p\vee q$.
We apply the definition of a projection to $p$ and to $q$. So there are
$u\in(pS)^\bot$ and $v\in(qS)^\bot$ such that
$r(x)=pr(x)+u=qr(x)+v$. We observe that $pr(x)=p(x)$ and $qr(x)=q(x)$.
By applying the refinement property to the equality $p(x)+u=q(x)+v$ and by
observing that $p(x)\perp q(x)$, we obtain $t\in S$ such that $u=q(x)+t$
and $v=p(x)+t$. On the one hand, $t\leq u,v$, thus
$t\in p^\bot S\cap q^\bot S=r^\bot S$, see Lemma~\ref{L:MeetCompPr}. On
the other hand, $t\leq r(x)$. Hence, $t=0$, so $r(x)=p(x)+u=p(x)+q(x)$.
\end{proof}

\begin{notation}\label{Not:WeVe}
For $x$, $y$, $z\in S$, $z=x\wedge y$ is the statement
   \[
   z\leq x,y\qquad\text{and}\qquad
   \forall t,\quad t\leq x,y\Rightarrow t\leq z.
   \]
We define, dually, the statement $z=x\vee y$. Note that $z$ is uniquely
defined by either statement only in case $\leq$ is antisymmetric.

Similarly, one can define the notations $a=\bigwedge_{i\in I}a_i$ and
$a=\bigvee_{i\in I}a_i$.
\end{notation}

\begin{proposition}\label{P:pvwq(x)}
Let\index{pzzroj@$\BB{S}$} $p$, $q\in\BB{S}$, let $x\in S$. Then the
following statements are satisfied:
   \[
   (p\wedge q)(x)=p(x)\wedge q(x),\qquad
   (p\vee q)(x)=p(x)\vee q(x).
   \]
\end{proposition}

\begin{proof}
We put $r=p\wedge q$ and $s=p\vee q$.

By Lemma~\ref{L:pLeqOrtq}, $r(x)\leq p(x),q(x)$.
Let $y\in S$ such that $y\leq p(x),q(x)$. Since $pS$ and $qS$ are ideals of
$S$ (see Lemma~\ref{L:OplusId}(i)), $y\in pS\cap qS$, so $p(y)=q(y)=y$. Thus
$y=pq(y)=r(y)\leq r(x)$. Hence $r(x)=p(x)\wedge q(x)$.

By Lemma~\ref{L:pLeqOrtq}, $p(x),q(x)\leq s(x)$.
Let $y\in S$ such that $p(x),q(x)\leq y$.
Thus, \emph{a fortiori}, $p(x),(p^\bot\wedge q)(x)\leq y$.
Since $\BB{S}$\index{pzzroj@$\BB{S}$} is a Boolean algebra,
\index{Boolean algebra!complete ---}
$s=p\oplus(p^\bot\wedge q)$.
It follows, by Lemma~\ref{L:JoinDisj}, that
$s(x)=p(x)+(p^\bot\wedge q)(x)$, thus $s(x)\leq y$ by
Lemma~\ref{L:XbotId}(iii). Hence $s(x)=p(x)\vee q(x)$.
\end{proof}

If $S$ is a \emph{total} (as opposed to partial) monoid, then the
projections of~$S$ correspond to direct decompositions of~$S$, thus, they
preserve arbitrary suprema and infima. For our partial structures, the
corresponding result still holds.

\begin{lemma}\label{L:ProjCont}
Let $p$ be a projection of~$S$. For every family $\famm{a_i}{i\in I}$
of elements of~$S$ and every $a\in S$,
\begin{enumerate}
\item If $I\ne\es$, then $a=\bigwedge_{i\in I}a_i$ implies that
$p(a)=\bigwedge_{i\in I}p(a_i)$.

\item Suppose that any two elements of $S$ have a meet. Then
$a=\bigvee_{i\in I}a_i$ implies that $p(a)=\bigvee_{i\in I}p(a_i)$.
\end{enumerate}
\end{lemma}

\begin{note}
The natural settings of Lemma~\ref{L:ProjCont} are in situations where $S$ is
antisymmetric as well. However, that condition is not, strictly speaking,
necessary, if we use the interpretation of the symbols $\bigwedge$ and
$\bigvee$ given in Notation~\ref{Not:WeVe}.

\end{note}

\begin{proof}
(i) Of course, $p(a)\leq p(a_i)$, for all $i$. Let $b$ be a minorant of
$\setm{p(a_i)}{i\in I}$. Since $I$ is nonempty, $b$ belongs to $pS$, so
$p(b)=b$. Since $b$ is also a minorant of $\setm{a_i}{i\in I}$, $b\leq a$.
Hence, $b=p(b)\leq p(a)$.

(ii) Of course, $p(a_i)\leq p(a)$, for all $i$. Let $b$ be a majorant of
$\setm{p(a_i)}{i\in I}$. Since $p(a)$ is also a majorant of
$\setm{p(a_i)}{i\in I}$ and by assumption, $c=p(a)\wedge b$ exists and it is
a majorant of $\setm{p(a_i)}{i\in I}$. {}From $c\leq p(a)$ it follows that
$c+p^\bot(a)$ is defined. Furthermore,
$a_i=p(a_i)+p^\bot(a_i)\leq c+p^\bot(a)$, for all $i\in I$, thus
$a\leq c+p^\bot(a)$. Therefore, $p(a)\leq c\leq b$.
\end{proof}

\begin{definition}\label{D:LsDiff}
Suppose that $S$ is antisymmetric. For $a$, $b$, $c\in S$, let $c=b\sd a$
\index{cabszzmsms@$c=b\sd a$|ii}
mean that $c$ is the least $x\in S$ such that $b\leq a+x$. We say that $c$ is
the \emph{least difference} of $b$ and $a$.
\end{definition}

\begin{lemma}\label{L:LstDiff}
Suppose that $S$ is antisymmetric, let $a\leq b$ in $S$. If $b\sd a$
exists, then $b=a+(b\sd a)$.
\end{lemma}

\begin{proof}
Put $c=b\sd a$.
Since $a\leq b$, there exists $d\in S$ such that $b=a+d$, thus, by the
definition of the least difference, $c\leq d$, whence $b\leq a+c\leq a+d=b$.
Therefore, since $S$ is antisymmetric, $b=a+c$.
\end{proof}

\begin{lemma}\label{L:SdProj}
Suppose that $S$ is antisymmetric and that $b\sd a$ exists for all $a$,
$b\in S$ such that $a\leq b$. Then $p(b\sd a)=p(b)\sd p(a)$, for all $a$,
$b\in S$ such that $a\leq b$ and all\index{pzzroj@$\BB{S}$} $p\in\BB{S}$.
\end{lemma}

\begin{proof}
{}From $b\leq a+(b\sd a)$ it follows that $p(b)\leq p(a)+p(b\sd a)$, thus
$p(b)\sd p(a)\leq p(b\sd a)$. Conversely, put $d=p(b)\sd p(a)$. Then
$p(b)\leq p(a)+d$ by definition, while we also have
$p^\bot(b)\leq p^\bot(a)+p^\bot(b\sd a)$, thus, adding the two inequalities
together (and observing that, since $d\leq p(b\sd a)$, $d+p^\bot(b\sd a)$ is
defined), we obtain the inequality $b\leq a+d+p^\bot(b\sd a)$. It follows
that $b\sd a\leq d+p^\bot(b\sd a)$, whence, by applying $p$, we obtain that
$p(b\sd a)\leq d$. Therefore, $d=p(b\sd a)$.
\end{proof}

\section{General comparability}\label{S:GenComp}

\begin{quote}
\em Standing hypothesis: $S$ is a conical \prm. We denote again by $\leq$ the
algebraic preordering of~$S$.
\end{quote}

\begin{definition}\label{D:GenCompMon}
We say that $S$ has \emph{general comparability},
\index{general comparability|ii} if for all $x$, $y\in S$,
there exists $p\in\BB{S}$\index{pzzroj@$\BB{S}$} such that $p(x)\leq p(y)$
and $p^\bot(x)\geq p^\bot(y)$.
\end{definition}

We give a sufficient condition that implies general comparability.

\begin{lemma}\label{L:GenCompAx}
Suppose that $S$ satisfies the following axioms:
\begin{enumerate}
\item
$\forall a,b$, $\exists c,x,y$ such that $a=c+x$, $b=c+y$, and $x\perp y$.

\item $S=a^\bot+a^{\bot\bot}$, for all $a\in S$.
\end{enumerate}

Then $S$ satisfies general comparability.
\end{lemma}

\begin{proof}
Let $a$, $b\in S$. Consider $c$, $x$, $y$ as in (i).
By (ii), there exists $p\in\BB{S}$\index{pzzroj@$\BB{S}$} such that
$pS=x^\bot$ and
$p^\bot S=x^{\bot\bot}$. So $p(x)=0$ and $p^\bot(y)=0$, whence
$p(a)\leq p(b)$ and $p^\bot(b)\leq p^\bot(a)$.
\end{proof}

\begin{lemma}\label{L:meetjoinS}
Suppose that $S$ has general comparability and that the algebraic
preordering of~$S$ is antisymmetric. Let $a$, $b\in S$.
The following assertions hold:
\begin{enumerate}
\item The pair $\set{a,b}$ has an infimum.

\item If the pair $\set{a,b}$ is majorized, then it has a supremum.
\end{enumerate}
\end{lemma}

A partially ordered set satisfying (i) and (ii) above is sometimes called a
\emph{chopped lattice}.\index{lattice!chopped ---}

\begin{proof}
By general comparability, there exists $p\in\BB{S}$\index{pzzroj@$\BB{S}$}
such that $p(a)\leq p(b)$ and $p^\bot(a)\geq p^\bot(b)$.
So $c=p(a)+p^\bot(b)$ is defined, and $c\leq a$, $b$. Furthermore, it is
easy to verify that $c=a\wedge b$.

Similarly, if the pair $\set{a,b}$ is majorized by an element $e$, then
$d=p(b)+p^\bot(a)$ is defined, and $d\leq e$. It is easy to verify that
$d=a\vee b$.
\end{proof}

\begin{notation}
For $a$, $b\in S$, let $a\ll b$ hold, if $a+b=b$.
We also say that $b$ \emph{absorbs} $a$.
\end{notation}

We recall the following axiom, see \index{Wehrung, F.}\cite{Wehr92a,Wehr92b}:
   \begin{align*}
   &\text{\textbf{The pseudo-cancellation property:}}\\
   &\forall a,b,c,\quad
   a+c=b+c\Rightarrow\exists d,\ \exists u,v\ll c\text{ such that }
   a=d+u\text{ and }b=d+v.
   \end{align*}
\index{pseudo-cancellation property|ii}
If $a+c=b+c$, $a=d+u$, $b=d+v$, and $u$, $v\ll c$, then $b+u$ is defined
(because $u\leq c$ and $b+c$ is defined) and $b+u=d+u+v\geq a$. Hence we
obtain the following weaker version of the pseudo-cancellation property:
   \[
\forall a,b,c,\quad
   a+c\leq b+c\Rightarrow\exists x\ll c\text{ such that }\ a\leq b+x.
   \]

\begin{lemma}\label{L:GCimpPC}
Suppose that $S$ has general comparability. Then $S$ satisfies the
pseudo-cancellation property.
\end{lemma}

\begin{proof}
Suppose $a+c=b+c$ in $S$. By using refinement, we find a refinement matrix
as follows:
   \[
   \begin{tabular}{|c|c|c|}
   \cline{2-3}
   \multicolumn{1}{l|}{} & $b$ & $c$\tvi\\
   \hline
   $a$ & $t$ & $a'$\tvi\\
   \hline
   $c$ & $b'$ & $c'$\tvi\\
   \hline
   \end{tabular}
   \]
By general comparability, there exists $p\in\BB{S}$\index{pzzroj@$\BB{S}$}
such that
$p(a')\leq p(b')$ and $p^\bot(b')\leq p^\bot(a')$. Let $u$,
$v\in S$ such that $p(a')+v=p(b')$ and $p^\bot(b')+u=p^\bot(a')$. Then
   \[
   \begin{aligned}
   p(c)&=p(b')+p(c')&&=p(a')+p(c')+v&&=p(c)+v,\\
   p^\bot(c)&=p^\bot(a')+p^\bot(c')&&=p^\bot(b')+p^\bot(c')+u&&=p^\bot(c)+u,
   \end{aligned}
   \]
It follows that $u+c=v+c=c$. Put $d=t+p(a')+p^\bot(b')$. By using
Lemma~\ref{L:JoinDisj}, we obtain the equalities
   \begin{align*}
   a&=t+a'=t+p(a')+p^\bot(a')=d+u\\
   b&=t+b'=t+p(b')+p^\bot(b')=d+v.\tag*{\qed}
   \end{align*}
\renewcommand{\qed}{}
\end{proof}

\begin{corollary}\label{C:GCSep}
Suppose that $S$ has general comparability. Then $S$ is
\emph{separative},\index{separative (\pcm)|ii}
that is, it satisfies the statement
   \[
   \forall a,b,c,\quad(a+c=b+c\text{ and }c\leq a,b)\Rightarrow a=b.
   \]
\end{corollary}

\begin{proof}
Let $a$, $b$, $c\in S$ such that $a+c=b+c$ and $c\leq a$, $b$.
By Lemma~\ref{L:GCimpPC}, there are $d$, $a'$, $b'\in S$ such that $a=d+a'$,
$b=d+b'$, and $a'$, $b'\ll c$. In particular, $a'$, $b'\leq c$, thus, since
$a+c$ and $b+c$ are defined, $a+b'$ and $b+a'$ are defined,
see Lemma~\ref{L:leqStillDef}. Note that
$a+b'=b+a'$. However, $c\leq a$, $b$, thus, since $a'\ll c$, we obtain
that $a'\ll b$, so $b+a'=b$. Similarly, $a+b'=a$. Therefore,
$a=b$.
\end{proof}

\begin{definition}\label{D:DirFinMon}
An element $c$ of~$S$ is
\begin{itemize}
\item[---] \emph{directly finite|ii},\index{directly finite} if $x+c=c$
implies that $x=0$, for all $x\in S$,

\item[---] \emph{cancellable}, if\index{cancellable (element)|ii}
$x+c=y+c\in S$ implies that $x=y$, for
all $x$, $y\in S$.
\end{itemize}
We say that $S$ is \emph{stably finite}, if\index{stably finite|ii} every
element of $S$ is\index{directly finite} directly finite. We denote
by~$S_\fin$\index{Szzfin@$S_\fin$|ii} the subset of $S$ consisting of all
directly finite elements.
\end{definition}

It is obvious that every cancellable element is\index{directly finite}
directly finite. By Lemma~\ref{L:GCimpPC}, we obtain immediately the
following converse.

\begin{lemma}\label{L:DfimCanc}
Suppose that $S$ has general comparability. Then every directly finite
\index{directly finite}element of~$S$ is cancellable.
\end{lemma}

\section{Boolean-valued partial refinement monoids}\label{S:BoolValPRM}
\begin{quote}
\em Standing hypothesis: $S$ is a conical \prm. We denote again by $\leq$ the
algebraic preordering of~$S$.
\end{quote}

For elements $a$ and $b$ of~$S$, it follows from Proposition~\ref{P:pvwq(x)}
that the set of all projections $p$ of~$S$ such that $p(a)\leq p(b)$ is
closed under finite join. We shall now consider a stronger statement.

\begin{definition}\label{D:BoolVal}
For $a$, $b\in S$, we shall denote by $\bv{a\leq b}$
\index{abzzvaleb@$\bv{a\leq b}$|ii}
the largest projection
$p$ of~$S$ such that $p(a)\leq p(b)$ if it exists.
Hence,\index{pzzroj@$\BB{S}$}
$\bv{a\leq b}\in\BB{S}$.

We say that $S$ is\index{refinement monoid!Boolean-valued partial ---|ii}
\emph{Boolean-valued}, if the Boolean value $\bv{a\leq b}$
is defined, for all $a$, $b\in S$.
\end{definition}

\begin{notation}\label{Not:bva=b}
For $a$, $b\in S$, if both $\bv{a\leq b}$ and $\bv{b\leq a}$ are defined, we
put\index{abzzvaeqb@$\bv{a=b}$|ii}
$\bv{a=b}=\bv{a\leq b}\wedge\bv{b\leq a}$.
\end{notation}

\begin{lemma}\label{L:cc(a)}
Assume that $S$ has general comparability.
Let $a\in S$, and suppose that $\bv{a=0}$ is defined. Then the following
assertions hold:
\begin{enumerate}
\item $a^\bot=\bv{a=0}S$.

\item $a^{\bot\bot}=\bv{a=0}^\bot S$.

\item $S=a^\bot\oplus a^{\bot\bot}$.
\end{enumerate}
\end{lemma}

\begin{proof}
(i) Put $p=\bv{a=0}$. For $x\in pS$ such that $x\leq a$, we have
$x=p(x)\leq p(a)=0$. Hence $pS\subseteq a^\bot$.

Conversely, let $x\in a^\bot$. By general comparability, there exists
$q\in\BB{S}$\index{pzzroj@$\BB{S}$} such that $q(a)\leq q(x)$ and
$q^\bot(x)\leq q^\bot(a)$. Since
$a\perp x$, the equalities $q(a)=q^\bot(x)=0$ hold. It follows from the
definition of $p$ that $q\leq p$. Therefore, $p^\bot(x)\leq q^\bot(x)=0$,
so $x\in pS$.

(ii) follows immediately from (i), while (iii) follows immediately from (i),
(ii), and the fact that $\bv{a=0}$ is a projection of $S$.
\end{proof}

\begin{lemma}\label{L:bva=0}
Assume that $S$ has general comparability.
Let $a\in S$. Then $\bv{a=0}$ exists \iff\ $S=a^\bot\oplus a^{\bot\bot}$.
\end{lemma}

\begin{proof}
If $\bv{a=0}$ exists, then $S=a^\bot\oplus a^{\bot\bot}$ by
Lemma~\ref{L:cc(a)}(iii). Conversely, suppose that
$S=a^\bot\oplus a^{\bot\bot}$. So
there exists a unique projection $p$ of $S$ such that $pS=a^\bot$. {}From
$p(a)\leq a$ and $p(a)\in a^\bot$ it follows that $p(a)=0$. Let $q\in\BB{S}$
\index{pzzroj@$\BB{S}$} such that $q(a)=0$. We claim that $q(x)\perp a$, for
any $x\in S$. Indeed, let $y\in S$ such that $y\leq q(x),a$. {}From $y\leq
q(x)$ it follows that
$q(y)=y$, thus $y=q(y)\leq q(a)=0$, so $y=0$, thus establishing our claim. So
$qS\subseteq a^\bot=pS$, whence $q\leq p$ by Lemma~\ref{L:pLeqOrtq}(i).
Therefore, $p=\bv{a=0}$.
\end{proof}

\begin{definition}\label{D:Remov}
Let $a$, $b\in S$. We say that $a$ is \emph{removable} from $b$,
\index{removable (in \prm s)|ii} and we write $a\rem b$,
\index{abrzzem@$a\rem b$|ii} if the following conditions hold:
\begin{enumerate}
\item $a\leq b$.

\item $b\leq a+x$ implies that $b\leq x$, for all $x\in S$.
\end{enumerate}
\end{definition}

In particular, we observe that $a\rem b$ implies that $a+b\leq b$ (the
converse does not hold as a rule). In particular, in case $S$ is
antisymmetric, $a\rem b$ implies that $a\ll b$.

\begin{lemma}\label{L:TrLeqTr}
Let $a$, $b$, $c\in S$ such that either $a\rem b\leq c$ or
$a\leq b\rem c$. Then $a\rem c$.
\end{lemma}

\begin{proof}
In both cases, it is trivial that $a\leq c$.

Suppose that $a\rem b\leq c$. Let $x\in S$ such that $c\leq a+x$.
So $b\leq a+x$, thus, since $a\rem b$, $b\leq x$, that is, $x=b+y$ for some
$y$. Hence $c\leq a+x=a+b+y$. But $a\rem b$, thus $a+b\leq b$, so
$c\leq b+y=x$. So $a\rem c$.

Suppose now that $a\leq b\rem c$. Let $x\in S$ such that $c\leq a+x$.
So $c\leq b+x$, thus (since $b\rem c$) $c\leq x$. So, again, $a\rem c$.
\end{proof}

\begin{lemma}\label{L:GCIVPtr}
Suppose that $S$ is antisymmetric \pup{that is, the algebraic preordering of
$S$ is antisymmetric} and that $S$ has pseudo-cancellation. For all $a$,
$b$, $c\in S$, the following assertions hold:
\begin{enumerate}
\item $a\leq b\leq a+c$ implies that there exists $x\leq c$ such that
$b=a+x$.

\item If $a\leq b$ in $S$, then $a\rem b$ \iff\ $b=a+x$ implies that $b=x$,
for all $x\in S$.
\end{enumerate}

\end{lemma}

\begin{proof}
(i) Since $a\leq b$, there exists $y\in S$ such that $b=a+y$. Hence
$a+y\leq a+c$, thus, by pseudo-cancellation, there exists $u\ll a$ such that
$y\leq u+c$. By refinement, there are $v\leq u$ and $x\leq c$ such that
$y=v+x$. Since $S$ is antisymmetric, $v\ll a$. Hence, $b=a+y=a+x$, with
$x\leq c$.

(ii) We prove the nontrivial direction. So, suppose that $b=a+x$ implies
$b=x$, for all $x\in S$. Now let $x\in S$ such that $b\leq a+x$. By (i)
above, there exists $y\leq x$ such that $b=a+y$. By assumption, $b=y$;
whence $b\leq x$.
\end{proof}

\begin{lemma}\label{L:ProjTr}
Suppose that $S$ is antisymmetric and satisfies general comparability. Let
$a$, $b\in S$.
\begin{enumerate}
\item If $a\rem b$, then $p(a)\rem p(b)$, for all\index{pzzroj@$\BB{S}$}
$p\in\BB{S}$.

\item Let $\famm{p_i}{i\in I}$ be a family of projections
of~$S$. We assume that both $\oll{a}=\bigvee_{i\in I}p_i(a)$ and
$\oll{b}=\bigvee_{i\in I}p_i(b)$ are defined. If $p_i(a)\rem p_i(b)$
for all $i\in I$, then $\oll{a}\rem\oll{b}$.

\item Let $n<\omega$, let $\famm{p_i}{i<n}$ be a finite
sequence of projections of~$S$, and let $p=\bigvee_{i<n}p_i$. If
$p_i(a)\rem p_i(b)$ for all $i<n$, then $p(a)\nobreak\rem\nobreak p(b)$.
\end{enumerate}
\end{lemma}

\begin{proof}
(i) It is clear that $p(a)\leq p(b)$.
Now let $x\in S$ such that $p(a)+x=p(b)$. So,
   \begin{align*}
   b&=p(b)+p^\bot(b)\\
   &=p(a)+p^\bot(b)+x\\
   &=p(a)+p^\bot(a)+p^\bot(b)+x&&(\text{because }p^\bot(a)\ll p^\bot(b))\\
   &=a+p^\bot(b)+x.
   \end{align*}
Since $a\rem b$, it follows that $p^\bot(b)+x=b$; whence $p(b)=p(x)=x$.
By Lemma~\ref{L:GCIVPtr}(ii), $p(a)\rem p(b)$.

(ii) Observe first that $\oll{a}\leq\oll{b}$. Let $x\in S$ such that
$\oll{a}+x=\oll{b}$. Observe that $p_i(a)\leq\oll{a}\leq a$
for all $i\in I$; hence $p_i(\oll{a})=p_i(a)$. Similarly,
$p_i(\oll{b})=p_i(b)$. Therefore, $p_i(a)+p_i(x)=p_i(b)$, for all $i$,
hence, since $p_i(a)\rem p_i(b)$, $p_i(b)=p_i(x)\leq x$. This holds for all
$i$, whence $\oll{b}\leq x$, so $x=\oll{b}$.
The conclusion follows from Lemma~\ref{L:GCIVPtr}(ii).

(iii) By Proposition~\ref{P:pvwq(x)}, $p(a)=\bigvee_{i<n}p_i(a)$ and
$p(b)=\bigvee_{i<n}p_i(b)$. The conclusion follows then from (ii).
\end{proof}

\begin{corollary}\label{C:2-5.?}
Suppose that $S$ is antisymmetric and satisfies general comparability.
For all $a$, $b$, $c\in S$, if $a\rem b,c$, then $a\rem b\wedge c$.
\end{corollary}

\begin{proof}
By general comparability, there exists $p\in\BB{S}$\index{pzzroj@$\BB{S}$}
such that
$p(b)\leq p(c)$ and $p^\bot(c)\leq p^\bot(b)$. Then, as in the proof of
Lemma~\ref{L:meetjoinS}, $b\wedge c=p(b)+p^\bot(c)$. By
Lemma~\ref{L:ProjTr}(i), $p(a)\rem p(b)=p(b\wedge c)$ and
$p^\bot(a)\rem p^\bot(c)=p^\bot(b\wedge c)$. Hence,
Lemma~\ref{L:ProjTr}(iii) implies that $a\rem b\wedge c$.
\end{proof}

\begin{definition}\label{D:PurInf}
An element $a$ of~$S$ is \emph{purely infinite}, if $2a=a$.

We denote by $\DI{S}$\index{SzzDI@$\DI{S}$|ii}
the set of all purely infinite elements of~$S$.
\end{definition}

We observe that the only element of~$S$ which is both
\index{directly finite} directly finite and purely infinite is $0$.

\begin{lemma}\label{L:2-5.?}
Suppose that $S$ is antisymmetric and satisfies general comparability. Then
$\DI{S}$ is closed under finite infima and suprema.
\end{lemma}

\begin{proof}
This is clear from the descriptions of pairwise infima and suprema given in
the proof of Lemma~\ref{L:meetjoinS}.
\end{proof}

\begin{lemma}\label{L:a+b=b(pi)}
Suppose that $S$ is antisymmetric.
Let $a$, $b\in S$ such that $a\leq b$. If either $a$ or $b$ is purely
infinite, then $a\ll b$.
\end{lemma}

\begin{lemma}\label{L:BVtr}
Suppose that $S$\index{refinement monoid!Boolean-valued partial ---}
is antisymmetric, Boolean-valued, and that it has general
comparability. Let $a\in\DI{S}$ and $b\in S$ such that $a\leq b$. Put
$p=\bv{b\leq a}$. Then $p^\bot(a)\rem p^\bot(b)$.
\end{lemma}

\begin{proof}
By the definition of $p$, $p(b)\leq p(a)$. Since $S$ is antisymmetric,
$p(a)=p(b)$. Furthermore, $p^\bot(a)\leq p^\bot(b)$ (because $a\leq b$).

Let $x\in S$ such that
   \begin{equation}\label{Eq:pbotba}
   p^\bot(b)=p^\bot(a)+x.
   \end{equation}
By general comparability, there
exists $q\in\BB{S}$\index{pzzroj@$\BB{S}$} such that
   \begin{align}
   q(x)&\leq q(a)\label{Eq:qxqa}\\
   q^\bot(a)&\leq q^\bot(x)\label{Eq:qbotaqbotx}
   \end{align}
By applying $q$ to \eqref{Eq:pbotba}, we obtain that
   \begin{equation}\label{Eq:qpbotba}
   qp^\bot(b)=qp^\bot(a)+q(x).
   \end{equation}
However, $x\leq p^\bot(b)$, thus $p^\bot(x)=x$, so
$q(x)=qp^\bot(x)=p^\bot q(x)\leq p^\bot q(a)=qp^\bot(a)$.
Hence, by Lemma~\ref{L:a+b=b(pi)},
$qp^\bot(a)+q(x)=qp^\bot(a)$, so, by \eqref{Eq:qpbotba},
$qp^\bot(b)=qp^\bot(a)$. By the definition of $p$, $qp^\bot\leq p$, thus,
since $qp^\bot\leq p^\bot$, $qp^\bot=0$, that is, $q\leq p$. Hence
$p^\bot\leq q^\bot$, thus, by \eqref{Eq:qbotaqbotx},
$p^\bot(a)\leq p^\bot(x)=x$. Hence, by \eqref{Eq:pbotba} and by
Lemma~\ref{L:a+b=b(pi)}, $x=p^\bot(b)$. We conclude the proof by
Lemma~\ref{L:GCIVPtr}(ii).
\end{proof}

We now introduce a useful definition.

\begin{definition}\label{D:cc(a)}
For $a\in S$, the \emph{central cover}\index{central cover|ii} of $a$,
denoted by $\cc(a)$, \index{czzca@$\cc(a)$|ii} is
defined as $\bv{a=0}^\bot$.
\end{definition}

Note that $a^\bot=\cc(a)^\bot S$, by Lemma~\ref{L:cc(a)}(i).

\begin{corollary}\label{C:remequiv}
Suppose that $S$ is\index{refinement monoid!Boolean-valued partial ---}
antisymmetric, Boolean-valued, and that it has general
comparability. Let $a\in\DI{S}$ and $b\in S$ such that $a\leq b$. Then
$a\rem b$ if and only if $q(b) \nleq q(a)$ for all nonzero projections
$q\leq \cc(b)$.
\end{corollary}

\begin{proof}
Assume first that $a\rem b$, and let $q\leq\cc(b)$ be a
projection such that $q(b)\leq q(a)$. Then $b\leq q(a)+q^\perp(b)\leq
a+q^\perp(b)$, and it follows from the assumption $a\rem b$ that
$b\leq q^\perp(b)$. Thus $q(b)=0$, so $q\wedge\cc(b)=0$, and hence $q=0$.

Conversely, assume that $q(b) \nleq q(a)$ for all nonzero projections
$q\leq \cc(b)$, set $p=\bv{b\leq a}$, and observe that $p\perp \cc(b)$. By
Lemma \ref{L:BVtr}, $p^\perp(a)\rem p^\perp(b)$, and so $\cc(b)(a)\rem
\cc(b)(b)$. Therefore $a\rem b$.
\end{proof}

\begin{lemma}\label{L:Basiccc}
Assume that $S$ is\index{refinement monoid!Boolean-valued partial ---}
antisymmetric, Boolean-valued, and satisfies general
comparability. Let $a\in S$ and let\index{pzzroj@$\BB{S}$} $p\in\BB{S}$.
\begin{enumerate}
\item $\cc(a)\leq p$ \iff\ $a\in pS$.

\item $\cc(p(a))=p\wedge\cc(a)$.

\item Suppose that $a=\bigvee_{i\in I}a_i$, for a family
$\famm{a_i}{i\in I}$ of elements of~$S$. Then
$\cc(a)=\bigvee_{i\in I}\cc(a_i)$.

\item $\cc(a\wedge b)=\cc(a)\wedge\cc(b)$, for all $a$, $b\in S$.
\end{enumerate}
\end{lemma}

\begin{proof}
(i) $\cc(a)\leq p$ \iff\ $p^\bot\leq\bv{a=0}$, \iff\ $p^\bot(a)=0$,
\iff\ $p(a)=a$, \iff\ $a\in pS$.

(ii) For all $q\in\BB{S}$,\index{pzzroj@$\BB{S}$} $q\leq\bv{p(a)=0}$ \iff\
$qp(a)=0$, \iff\
$qp\leq\bv{a=0}$, \iff\ $q\leq\bv{a=0}\vee p^\bot$. Hence
   \[
   \bv{p(a)=0}=\bv{a=0}\vee p^\bot.
   \]
Therefore, $\cc(p(a))=\bv{p(a)=0}^\bot=p\wedge\cc(a)$.

(iii) By (i), for any\index{pzzroj@$\BB{S}$} $p\in\BB{S}$,
$\cc(a)\leq p$ \iff\ $a\in pS$, \iff\ $a_i\in pS$ for all $i$
(by Lemma~\ref{L:ProjCont}(ii)), \iff\ $\cc(a_i)\leq p$ for all $i$, \iff\
$\bigvee_{i\in I}\cc(a_i)\leq p$. The conclusion of~(iii) follows.

(iv) It suffices to prove the inequality
   \begin{equation}\label{Eq:ccawedb=0}
   \bv{a\wedge b=0}=\bv{a=0}\vee\bv{b=0}.
   \end{equation}
Since $a\wedge b\leq a$, $b$, the inequality
$\bv{a\wedge b=0}\geq\bv{a=0}\vee\bv{b=0}$ is obvious.
Conversely, put $p=\bv{a\wedge b=0}$. By general comparability, there are
$q$, $r\in\BB{S}$\index{pzzroj@$\BB{S}$} such that $q(a)\leq q(b)$,
$r(b)\leq r(a)$, and $p=q\vee r$. It follows that
   \begin{align*}
   0&=q(a\wedge b)&&(\text{because }q\leq p)\\
   &=q(a)\wedge q(b)&&(\text{by Lemma~\ref{L:ProjCont}(i)})\\
   &=q(a)&&(\text{because }q(a)\leq q(b)),
   \end{align*}
hence $q\leq\bv{a=0}$. Similarly, $r\leq\bv{b=0}$, so
$p\leq\bv{a=0}\vee\bv{b=0}$. This completes the proof of \eqref{Eq:ccawedb=0}.
\end{proof}

\section{Least and largest difference functions}\label{S:lLDiff}

\begin{quote}
\em Standing hypothesis: $S$ is a \prm\ satisfying the following additional
properties:
\begin{itemize}
\item[(1)] $S$ is antisymmetric.

\item[(2)] $S$ has general comparability.

\item[(3)] $S$ is\index{refinement monoid!Boolean-valued partial ---}
Boolean-valued.

\item[(4)] Every element of~$S$ is the sum of a\index{directly finite}
directly finite element and a purely infinite element.
\end{itemize}
\end{quote}

\begin{lemma}\label{L:Perpdfpi}
For all $a\in S$, there exists $p\in\BB{S}$\index{pzzroj@$\BB{S}$} such that
$p(a)$ is\index{directly finite} directly finite and $p^\bot(a)$ is purely
infinite.
\end{lemma}

\begin{proof}
By assumption on $S$, there are elements $x$ and $y$ in~$S$ such that
$a=x+y$, $x$ is purely infinite, and $y$ is\index{directly finite} directly
finite. By general comparability, there exists
$p\in\BB{S}$\index{pzzroj@$\BB{S}$} such that
$p(x)\leq p(y)$ and $p^\bot(y)\leq p^\bot(x)$. Since $y$
is\index{directly finite} directly finite and $p(x)\leq p(y)\leq y$, $p(x)$
is\index{directly finite} directly finite. But $p(x)$ is purely infinite,
thus $p(x)=0$, and so
$p(a)=p(y)$ is\index{directly finite} directly finite. Since
$p^\bot(x)\geq p^\bot(y)$ with $p^\bot(x)$ purely infinite,
$p^\bot(a)=p^\bot(x)+p^\bot(y)=p^\bot(x)$ by Lemma~\ref{L:a+b=b(pi)}.
Therefore, $p^\bot(a)$ is purely infinite.
\end{proof}

\begin{corollary}\label{C:a/infty}
For any $a\in S$, the following assertions hold:
\begin{enumerate}
\item There exists a largest purely infinite element $u$ of~$S$ such that
$u\leq a$.

\item The element $u$ is also the largest $s\in S$ such that $s\ll a$.

\item There exists a unique $v\in S$ such that $a=u+v$ and $u\perp v$.

\item The element $v$ is\index{directly finite} directly finite.
\end{enumerate}
\end{corollary}

\begin{proof}
By Lemma~\ref{L:Perpdfpi}, there exists $p\in\BB{S}$\index{pzzroj@$\BB{S}$}
such that $p(a)$ is\index{directly finite} directly finite and $p^\bot(a)$
is purely infinite. Set
$u=p^\bot(a)$ and
$v=p(a)$. Observe that $u\leq a$ and $u$ is purely infinite, so $u\ll a$.

(ii) For any $s\in S$, $s\ll a$ implies that $p(s)\ll p(a)$. Since $p(a)$ is
\index{directly finite}directly finite, $p(s)=0$, and thus $s=p^\bot(s)\leq
p^\bot(a)=u$.

(i) For any purely infinite $t\in S$ such that $t\leq a$, it follows from
Lemma~\ref{L:a+b=b(pi)} that $t\ll a$, whence $t\leq u$ by part~(ii).

(iii), (iv) We already have $a=u+v$ with $u\perp v$ and $v$
\index{directly finite} directly finite. For any $w\in S$, if $a=u+w$ with
$u\perp w$, then
$u+w=u+v$, thus, by refinement (and since $u\perp v,w$), $v=w$.
\end{proof}

\begin{notation}
For any $a\in S$, we shall denote by $\di{a}$\index{adzzi@$\di{a}$|ii}
the largest purely infinite
element $u$ of~$S$ such that $u\leq a$.
\end{notation}

\begin{lemma}\label{L:a+b/infty}
Let $a$, $b\in S$ such that $a+b$ is defined. Then
   \[
   \di{a+b}=\di{a}+\di{b}.
   \]
\end{lemma}

\begin{proof}
First, $\di{a}+\di{b}$ is purely infinite and below $a+b$, thus
$\di{a}+\di{b}\leq\di{a+b}$.

Conversely, put $c=\di{a+b}$. Then $c+a+b=a+b$, thus, by
canceling the directly finite\index{directly finite} parts of $a$ and $b$
(use Lemma~\ref{L:DfimCanc}), $c+a+\di{b}=a+\di{b}$, thus, again,
$c+\di{a}+\di{b}=\di{a}+\di{b}$. In particular, $c\leq\di{a}+\di{b}$.
\end{proof}

Our next result involves the \emph{least difference} function introduced in
Definition~\ref{D:LsDiff}.

\begin{proposition}\label{P:bsminusa}
Let $a\leq b$ in $S$; then $b\sd a$ exists.
\end{proposition}

\begin{proof}
By Lemma~\ref{L:Perpdfpi}, there exists $q\in\BB{S}$\index{pzzroj@$\BB{S}$}
such that $q(a)$ is directly finite\index{directly finite} and $q^\bot(a)$
is purely infinite. Let
$c_0\in S$ such that
\begin{equation}\label{Eq:qac0b}
   q(a)+c_0=q(b).
   \end{equation}
Put $p=\bv{q^\bot(b)\leq q^\bot(a)}$. Since $qq^\bot=0$, the inequality
$q\leq p$ holds. Observe also the following equality:
   \begin{equation}\label{Eq:pqbotab}
   pq^\bot(b)=pq^\bot(a).
   \end{equation}
Since $a\leq b$ and $q^\bot(a)$ is
purely infinite, it follows from Lemma~\ref{L:BVtr} that
$p^\bot q^\bot(a)\rem p^\bot q^\bot(b)$, that is, since $p^\bot\leq q^\bot$,
   \begin{equation}\label{Eq:pbotatrb}
   p^\bot(a)\rem p^\bot(b).
   \end{equation}
In particular, we obtain the relation
   \begin{equation}\label{Eq:pbotllb}
   p^\bot(a)\ll p^\bot(b).
   \end{equation}
Since $c_0\leq q(b)\leq p(b)$, $c=c_0+p^\bot(b)$ is defined.
So we obtain that
   \begin{align*}
   b&=p^\bot(b)+pq^\bot(b)+q(b)\\
   &=p^\bot(a)+p^\bot(b)+pq^\bot(a)+q(a)+c_0
   &&\text{(by \eqref{Eq:qac0b}, \eqref{Eq:pqbotab}, and \eqref{Eq:pbotllb})}\\
   &=a+c_0+p^\bot(b)\\
   &=a+c.
   \end{align*}

Furthermore, let $x\in S$ such that $b\leq a+x$. So $q(b)\leq q(a)+q(x)$,
that is, $q(a)+c_0\leq q(a)+q(x)$. Thus, since $q(a)$
is\index{directly finite} directly finite and by Lemma~\ref{L:DfimCanc}, we
obtain
   \begin{equation}\label{Eq:c0q(x)}
   c_0\leq q(x).
   \end{equation}
Furthermore, $p^\bot(b)\leq p^\bot(a)+p^\bot(x)$, thus, by
\eqref{Eq:pbotatrb}, we obtain that
   \begin{equation}\label{Eq:pbotbbotx}
   p^\bot(b)\leq p^\bot(x).
   \end{equation}
By adding \eqref{Eq:c0q(x)} and \eqref{Eq:pbotbbotx} together, we thus
obtain that $c\leq q(x)+p^\bot(x)\leq p(x)+p^\bot(x)=x$. So we have verified
that $c=b\sd a$.
\end{proof}

A similar result holds for the existence of the ``largest difference''.

\begin{proposition}\label{P:blminusa}
Let $a\leq b$ in $S$. Then there exists a largest element $c$ of~$S$ such that
$a+c\leq b$, and then $b=a+c$.
\end{proposition}

The element $c$ of the statement above will be denoted by $b-a$, the
\emph{largest difference}\index{cabszzmsms@$c=b-a$|ii} of $b$ and $a$.

\begin{proof}
By the definition of the algebraic preordering, there exists $d\in S$ such
that $a+d=b$. So $c=\di{a}+d$ is defined (because $\di{a}\leq a$) and
$c\leq b$. {}From $\di{a}\ll a$ it follows that $a+c=b$. If
$x\in S$ is such that $a+x\leq b$, then, by pseudo-cancellation
(see Lemma~\ref{L:GCimpPC}), $x\leq d+y$ for some $y\ll a$, so
$x\leq d+\di{a}=c$.
\end{proof}

\begin{corollary}\label{C:+MJ}
Let $a$, $b\in S$, let $X$ be a nonempty subset of~$S$.
We assume that $a+x$ is defined for all $x\in X$.
\begin{enumerate}
\item If $b=\bigwedge X$, then $a+b=\bigwedge(a+X)$.

\item If $b=\bigvee X$ and $a+X$ is majorized, then $a+b=\bigvee(a+X)$.
\end{enumerate}
\end{corollary}

\begin{proof}
(i) Pick $x\in X$. Since $a+x$ is defined and $b\leq x$, $a+b$ is defined.
Furthermore, $a+b\leq a+X$. Conversely, let $c\leq a+X$. By
Lemma~\ref{L:meetjoinS}, $c'=a\vee c$ is defined, and $a\leq c'\leq a+X$.
By Proposition~\ref{P:bsminusa}, $c'\sd a\leq X$, so
$c'\sd a\leq b$. Therefore, by adding $a$ on both sides of this inequality,
we obtain that $c\leq c'=a+(c'\sd a)\leq a+b$.

(ii) Pick a majorant $c$ of $a+X$. In particular, $c\geq a$. By
Proposition~\ref{P:blminusa}, $c-a\geq X$, so
$c-a\geq b$. Since $c=a+(c-a)$, $a+b$ is defined and $a+b\leq c$.
This holds for any majorant $c$ of $a+X$. Since $a+b$ is itself a majorant
of $a+X$, it is the supremum of $a+X$.
\end{proof}

\chapter{Continuous dimension scales}

\section[Basic properties]
{Basic properties; the monoids $\ZZ_\gamma$, $\RR_\gamma$, and $\two_\gamma$}
\label{S:BasicRZ2}
\index{Zzzgamma@$\ZZ_\gamma$}
\index{Rzzgamma@$\RR_\gamma$}
\index{Tzzgamma@$\two_\gamma$}

The fundamental definition underlying this chapter is the following.

\begin{definition}\label{D:DimInt}
A \emph{continuous dimension scale}\index{continuous dimension scale|ii} is a
\pcm\index{partial commutative monoid}\ $S$ which satisfies the following
axioms.
\begin{itemize}
\item[(M1)] $S$ has refinement (see Definition~\ref{D:RefPpty}), and the
algebraic preordering on $S$ is antisymmetric.\index{mzzonesix@(M1--6)|ii}

\item[(M2)] Every nonempty subset of~$S$ admits an infimum. Equivalently,
every majorized subset of~$S$ admits a supremum.

\item[(M3)] $S$ has general comparability
(see Definition~\ref{D:GenCompMon}).

\item[(M4)] $S$ is Boolean-valued
\index{refinement monoid!Boolean-valued partial ---}
(see Definition~\ref{D:BoolVal}).

\item[(M5)] Every element $a$ of~$S$ can be written $a=x+y$, where $x$ is
directly finite\index{directly finite} (Definition~\ref{D:DirFinMon}) and
$y$ is purely infinite (Definition~\ref{D:PurInf}).

\item[(M6)] Let $a$, $b$ be purely infinite elements of~$S$.
If $a\rem b$ (see Definition~\ref{D:Remov}), then the set of all
purely infinite elements $x$ of~$S$ such that $a\rem x$ and $x^\bot=b^\bot$
(see \eqref{Eq:DefXbot}, page~\pageref{Eq:DefXbot}) has a least element.
\end{itemize}
A continuous dimension scale $S$ is \emph{bounded}, if it has a largest element.
\end{definition}

All axioms (M1)--(M5) have been considered in Chapter~\ref{Ch:PartMon}.
Axiom~(M6) is a newcomer, whose importance will appear in
Section~\ref{S:alphap}.

We shall first give an alternative axiomatization of continuous dimension scales. In
order to prepare for this, we first prove the following result, which
extends the result of Lemma~\ref{L:GenCompAx}.

\begin{proposition}\label{P:AltAx}
Let $S$ be a \prm\ satisfying the following properties:
\begin{itemize}
\item[(1)] $S$  is antisymmetric.

\item[(2)] Any two elements of $S$ have a meet.

\item[(3)] $S$ satisfies Axiom~\textup{(M5)}.
\end{itemize}
Then the following assertions are equivalent:
\begin{enumerate}
\item $S$ satisfies the following axioms:
\begin{itemize}
\item[(N1)]
$\forall a,b$, $\exists c,x,y$ such that $a=c+x$, $b=c+y$, and $x\perp y$.
\index{nzzonethree@(N1--3)|ii}

\item[(N2)] $S=a^\bot+a^{\bot\bot}$, for all $a\in S$.

\item[(N3)] $b\sd a$ exists, for all $a$, $b\in S$ such that $a\leq b$.
\end{itemize}

\item $S$ is Boolean-valued
\index{refinement monoid!Boolean-valued partial ---} and it satisfies general
comparability.
\end{enumerate}
\end{proposition}

\begin{proof}
(i)$\Rightarrow$(ii) Let $S$ satisfy (N1), (N2), and (N3).
The fact that $S$ satisfies general comparability follows from
Lemma~\ref{L:GenCompAx}. Now we prove that $S$ is
\index{refinement monoid!Boolean-valued partial ---}
Boolean-valued. So let $a$,
$b\in S$. By (N3), $c=a\sd(a\wedge b)$ exists. For any projection $p$
of~$S$,
   \begin{align}
   p(a)\leq p(b)&&\text{\iff}&&p(a)&=p(a)\wedge p(b)&&\notag\\
   &&\text{\iff}&&p(a)&=p(a\wedge b)&&
   \text{(by Lemma~\ref{L:ProjCont}(i))}\notag\\
   &&\text{\iff}&&p(a)\sd p(a\wedge b)&=0&&\notag\\
   &&\text{\iff}&&p(c)&=0&&\text{(by Lemma~\ref{L:SdProj})}.\label{Eq:p(c)=0}
   \end{align}
By (N2), $S=c^\bot+c^{\bot\bot}$. By Lemma~\ref{L:XbotId}(ii), $c^\bot$
and $c^{\bot\bot}$ are ideals of $S$, thus, since
$c^\bot\cap c^{\bot\bot}=\set{0}$, $S=c^\bot\oplus c^{\bot\bot}$, hence, by
Lemma~\ref{L:bva=0}, $p=\bv{c=0}$ exists. Therefore, by \eqref{Eq:p(c)=0},
$\bv{a\leq b}$ exists, and $\bv{a\leq b}=\bv{c=0}$.

(ii)$\Rightarrow$(i) Suppose that $S$ is Boolean-valued
\index{refinement monoid!Boolean-valued partial ---}
and satisfies general
comparability. We verify that $S$ satisfies (N1)--(N3).

(N1) By general comparability, there exists
$p\in\BB{S}$\index{pzzroj@$\BB{S}$} such that
$p(a)\leq p(b)$ and $p^\bot(b)\leq p^\bot(a)$. Let $x$, $y\in S$ such that
$p(a)+y=p(b)$ and $p^\bot(b)+x=p^\bot(a)$. Furthermore, since $p(a)\leq p(b)$
and $b=p(b)+p^\bot(b)$, the element $c=p(a)+p^\bot(b)$ is defined. {}From
$x\leq p^\bot(a)$ and $y\leq p(b)$ it follows that $x\perp y$. Finally,
   \[
   \begin{aligned}
   a&=p(a)+p^\bot(a)&&=c+x,\\
   b&=p(b)+p^\bot(b)&&=c+y.
   \end{aligned}
   \]
Hence we have obtained (N1).

(N2) follows immediately from Lemma~\ref{L:cc(a)}.

(N3) follows immediately from Proposition~\ref{P:bsminusa}.
\end{proof}

\begin{corollary}\label{C:AltAx}
Let $S$ be a \pcm.\index{partial commutative monoid} Then $S$ is a continuous dimension scale\index{continuous dimension scale} \iff\ it satisfies the axioms
\textup{(M1)}, \textup{(M2)},
\textup{(M5)}, \textup{(M6)}, \textup{(N1)}, \textup{(N2)}, and
\textup{(N3)}.
\end{corollary}

\begin{remark}\label{Rk:al1order}
It follows from Corollary~\ref{C:AltAx} that for a \pcm\index{partial
commutative monoid}\ $S$, to be a continuous dimension scale\index{continuous dimension scale} is equivalent to the conjunction of the
\emph{second-order} axiom (M2) and a finite list of
\emph{first-order} axioms.
\end{remark}

As a corollary of this alternative description of continuous dimension scales,\index{continuous dimension scale} we observe the following.

\begin{definition}\label{D:DirProd}
The \emph{direct product}\index{direct product (of \pcm s)|ii}
of a family $\famm{S_i}{i\in I}$ of
\pcm s\index{partial commutative monoid} is
obtained by endowing the ordinary cartesian product $S=\prod_{i\in I}S_i$
with the partial addition defined by
   \[
   \famm{a_i}{i\in I}+\famm{b_i}{i\in I}=\famm{a_i+b_i}{i\in I},
   \]
for all $\famm{a_i}{i\in I}$, $\famm{b_i}{i\in I}\in S$.
\end{definition}

Of course, Definition~\ref{D:DirProd} is an obvious generalization
of the definition of the product of finitely many
\pcm s\index{partial commutative monoid} introduced in
Remark~\ref{R:EmbdSSi}.

It is trivial that the direct product of any family of
\pcm s\index{partial commutative monoid} is a \pcm. Far less trivial is the
following preservation result.

\begin{lemma}\label{L:FinDPDI}
Any direct product of a family of continuous dimension scales\index{continuous dimension scale} is a continuous dimension scale.
\end{lemma}

\begin{proof}
We use the characterization of continuous dimension scales\index{continuous dimension scale}
obtained in Corollary~\ref{C:AltAx}. The proof is relatively long but very
easy, so we will not give the details of it but rather the basic idea.
A key point is to verify that the operations
$x\mapsto x^\bot$, $x\mapsto x^{\bot\bot}$, $(x,y)\mapsto x\wedge y$,
$(x,y)\mapsto y\sd x$ (for $x\leq y$), and the relations $x\leq y$,
$x\rem y$, $x\perp y$, and $y\in x^{\bot\bot}$ can be ``read
componentwise'', that is, for example, if $x=\famm{x_i}{i\in I}$ and
$y=\famm{y_i}{i\in I}$, then $x\leq y$
\iff\ $x_i\leq y_i$ for all $i$, $x^\bot=\prod_{i\in I}x_i^\bot$,
and so on. Once these simple facts are established, the verification of
the axioms (M1), (M2), (M5), (M6), (N1), (N2), and (N3) is routine.
\end{proof}

For a continuous dimension scale\index{continuous dimension scale} $S$ and elements $a$ and
$b$ in a lower subset $T$ of $S$, the orthogonality of $a$ and $b$ means the
same in $S$ and in $T$. We capture this pattern in a definition.

\begin{definition}\label{D:Abs}\hfill
\begin{enumerate}
\item A \emph{statement} $\varphi(x_1,\dots,x_n)$ in the language of
\pcm s\index{partial commutative monoid}
is \emph{absolute},\index{absolute (statement, function)|ii}
if for any continuous dimension scale\index{continuous dimension scale}
$S$, every lower subset
$T$ of~$S$, and all elements $a_1$, \dots, $a_n\in T$, $S$ satisfies
$\varphi(a_1,\dots,a_n)$ \iff~$T$ satisfies $\varphi(a_1,\dots,a_n)$.

\item A \emph{definable function}\index{definable function|ii}
$y=\varphi(x_1,\dots,x_n)$ in the
language of \pcm s\index{partial commutative monoid} is \emph{absolute}, if
for any continuous dimension scale\index{continuous dimension scale} $S$, every lower subset
$T$ of $S$, and all elements $a_1$, \dots, $a_n\in T$,
$\varphi(a_1,\dots,a_n)$ is defined in $S$ \iff\ it is defined in $T$, and
then both values are equal.
\end{enumerate}
\end{definition}

\begin{lemma}\label{L:Abs}
The following statements
\begin{enumerate}
\item $x\leq y$;
\item $x\perp y$;
\item $x^\bot\subseteq y^\bot$;
\item $x^\bot=y^\bot$;
\item $x\rem y$
\end{enumerate}
and the following function
\begin{itemize}
\item[(vi)] $z=y\sd x$
\end{itemize}
are absolute.
\end{lemma}

\begin{proof}
Most items are trivial, except perhaps (v) and (vi). Let $T$ be a
lower subset of a continuous dimension scale\index{continuous dimension scale} $S$, let
$a\leq b$ in $T$. Since $S$ is a continuous dimension scale\index{continuous dimension scale}, $c=b\sd a$ is defined in $S$. By Lemma~\ref{L:LstDiff}, $b=a+c$,
thus, since $T$ is a lower subset of $S$,
$c\in T$. It readily follows that $c=b\sd a$ in $T$, which concludes the
proof of (vi) (because $c$ always exists). As $a\rem b$ \iff\ $a\leq b$ and
$b\sd a=b$, item (v) follows immediately.
\end{proof}

As an easy consequence of Corollary~\ref{C:AltAx} and Lemma~\ref{L:Abs}, we
obtain the following.

\begin{lemma}\label{L:SegDI}
Let $S$ be a continuous dimension scale,\index{continuous dimension scale} let $T$ be a
lower subset of~$S$, viewed as a partial submonoid of~$S$ \pup{see
Definition~\textup{\ref{D:PartSubm}}}. Then $T$ is a continuous dimension scale.\index{continuous dimension scale}
\end{lemma}

We refer to Lemma~\ref{L:SegDimInt} for more information on lower subsets
of continuous dimension scales.\index{continuous dimension scale}

We observe that trying to use the original axioms (M1)--(M6) for the proof
of Lemma~\ref{L:FinDPDI} would have been much more difficult, since we
would have needed to understand the projections of the product
$\prod_{i\in I}S_i$. By using Corollary~\ref{C:AltAx}, the proof is
still somewhat tedious, but essentially trivial.

We present another way to produce continuous dimension scales\index{continuous dimension scale}.

\begin{lemma}\label{L:DirUn}
Let $I$ be an upwards directed partially ordered set, let $(S_i)_{i\in I}$
be a family of continuous dimension scales\index{continuous dimension scale} such that $S_i$
is a lower subset of
$S_j$, for all $i\leq j$ in $I$. Then $\bigcup_{\in I}S_i$ is a continuous dimension scale.\index{continuous dimension scale}
\end{lemma}

\begin{proof}
The set $S$ is, of course, endowed with the union of all the
\pcm\ operations on all the $S_i$-s. By using Lemma~\ref{L:Abs}, it is easy
to verify that $S$ satisfies (M1), (M5), (N1), (N2), and (N3).

Let $X$ be a nonempty subset of $S$; so
there exists $i\in I$ such that $X\cap S_i\neq\es$. Denote by $a_j$ the
meet of $X\cap S_j$, for all $j\geq i$ in $I$. Then $i\leq j\leq k$ implies
that $a_j\geq a_k$, in particular, all the $a_j$-s belong to $S_i$, and the
meet of all the $a_j$-s in $S_i$ is also the meet of $X$ in $S$. Hence $S$
satisfies (M2).

Let $a\rem b$ in $\DI{S}$. There exists $i\in I$ such that $a$, $b\in S_i$.
It follows from Lemma~\ref{L:Abs} that the statement $a\rem b$ holds in all
$S_j$ with $j\geq i$, thus, since $S_j$ is a continuous dimension scale,\index{continuous dimension scale} the set of all elements $x\in\DI{S_j}$
such that $a\rem x$ and $x^\bot=b^\bot$ (we use Lemma~\ref{L:Abs}) has a
least element, say, $c_j$. It follows again from Lemma~\ref{L:Abs} that
$c_j=c_i$, for all $j\geq i$; denote by $c$ this element, then $c\in\DI{S}$
and $c$ is minimum in the set of all elements $x\in\DI{S}$ such that $a\rem
x$ and $x^\bot=b^\bot$. Hence $S$ satisfies (M6). By
Corollary~\ref{C:AltAx}, $S$ is a continuous dimension scale\index{continuous dimension scale}.
\end{proof}

We now provide fundamental examples of continuous dimension scales\index{continuous dimension scale}.

For an ordinal $\gamma$, the monoids $\ZZ_\gamma$, $\RR_\gamma$, and
$\two_\gamma$ are defined in the Introduction:\label{Pg:ZR2gam}
\index{Zzzgamma@$\ZZ_\gamma$}\index{Rzzgamma@$\RR_\gamma$}%
\index{Tzzgamma@$\two_\gamma$}%
   \begin{align*}
   \ZZ_\gamma&=\ZZ^+\cup\setm{\aleph_\xi}{0\leq\xi\leq\gamma},\\
   \RR_\gamma&=\RR^+\cup\setm{\aleph_\xi}{0\leq\xi\leq\gamma},\\
   \two_\gamma&=\set{0}\cup\setm{\aleph_\xi}{0\leq\xi\leq\gamma}.
   \end{align*}
We call the elements of $\ZZ^+$, $\RR^+$, and $\set{0}$ the
\emph{finite elements} of $\ZZ_\gamma$, $\RR_\gamma$, and $\two_\gamma$,
respectively.
We endow each of the sets $\ZZ_\gamma$, $\RR_\gamma$, and $\two_\gamma$
with the addition that extends the natural addition on the finite elements
and on the alephs (so $\aleph_\alpha+\aleph_\beta=\aleph_\beta$ if
$\alpha\leq\beta$), and such that every finite element is absorbed by every
aleph. Hence $\ZZ_\gamma$, $\RR_\gamma$, and $\two_\gamma$ are monoids (and
not just partial ones). We also observe that the finite elements of
$\ZZ_\gamma$, $\RR_\gamma$, and $\two_\gamma$ are precisely the directly
finite\index{directly finite} ones.

Therefore, the algebraic ordering on each of the structures $\ZZ_\gamma$,
$\RR_\gamma$, and $\two_\gamma$ is the natural \emph{total} ordering.

\begin{proposition}\label{P:ZRTwoThDimInt}
For every ordinal $\gamma$, the monoids
$\ZZ_\gamma$,\index{Zzzgamma@$\ZZ_\gamma$}
$\RR_\gamma$,\index{Rzzgamma@$\RR_\gamma$} and
$\two_\gamma$\index{Tzzgamma@$\two_\gamma$} are totally ordered continuous dimension scales.\index{continuous dimension scale}
\index{Zzzgamma@$\ZZ_\gamma$}\index{Rzzgamma@$\RR_\gamma$}%
\index{Tzzgamma@$\two_\gamma$}%
\end{proposition}

\begin{proof}
All axioms (M1)--(M6) are trivially satisfied, except perhaps refinement.
Let $S$ be one of the structures $\ZZ_\gamma$, $\RR_\gamma$, or
$\two_\gamma$.
\index{Zzzgamma@$\ZZ_\gamma$}\index{Rzzgamma@$\RR_\gamma$}%
\index{Tzzgamma@$\two_\gamma$}%
Since every element of~$S$ is either cancellable\index{cancellable (element)}
or purely infinite, $S$ satisfies the following weak form of the
pseudo-cancellation property:
   \[
   \forall a,b,c,\quad a+c=b+c\Rightarrow\exists x\ll c,a\leq b+x.
   \]
By Proposition~1.23 of \index{Wehrung, F.}\cite{Wehr92a}, since $S$ is
totally ordered, it has refinement.
\end{proof}

Of course, one could also have verified refinement directly, but at the
expense of a few more calculations. For example, start with the observation
that $\ZZ^+$, $\RR^+$, and $\set{0}$ have refinement; note the easy fact that
adjoining a new infinity element to any refinement monoid
\index{refinement monoid} produces another
refinement monoid; use ordinal induction.

We shall also see in Proposition~\ref{P:lgrpDI} that the positive cone
of any Dedekind complete lattice ordered group is a continuous dimension scale.\index{continuous dimension scale} Furthermore,
Proposition~\ref{P:ZRTwoThDimInt} will be considerably extended in
Theorem~\ref{T:C(O,K)DimInt}.

\section{Dedekind complete lattice-ordered groups}\label{S:Dcpllgrp}

We recall that every Dedekind complete partially ordered group is abelian,
see \cite[Theorem~28]{Birk}.\index{Birkhoff, G.} Much more is true.

\begin{proposition}\label{P:lgrpDI}
Let $G$ be a Dedekind complete lattice-ordered group. Then $G^+$ is a
continuous dimension scale.\index{continuous dimension scale}
\end{proposition}

\begin{proof}
It is trivial that the algebraic preordering on $G^+$ is antisymmetric.
It is well-known that $G^+$ (or, more generally, the positive cone
of any lattice-ordered group) satisfies refinement, see, for example,
\index{Bigard, A.}\index{Keimel, K.}\index{Wolfenstein, S.}%
\cite[Th\'eor\`eme~1.2.16]{BKW}. Hence $G^+$ satisfies (M1).

Axiom~(M2) is just a reformulation of the fact that $G$ is Dedekind complete.

Axioms (M5) and (M6) are trivially satisfied, because every element of $G^+$
is\index{directly finite} directly finite.

To verify that $G^+$ satisfies Axiom~(M3), it suffices to verify that it
satisfies the assumptions of Lemma~\ref{L:GenCompAx}. For $a$, $b\in G^+$,
if we put $c=a\wedge b$, $x=a-c$, and $y=b-c$, then $a=c+x$, $b=c+y$, and
$x\wedge y=0$. This takes care of Assumption~(i).

For $x$, $y\in G$, we define $x\perp y$ to hold, if $|x|\wedge|y|=0$, and we
put
   \[
   X^\bot=\setm{g\in G}{g\perp x\text{ for all }x\in X},
   \]
the \emph{polar}\index{polar (subset)|ii} of $X$. Since $G$ is Dedekind
complete, the direct factors of $G$ are exactly the polar subsets of $G$, see
\index{Bigard, A.}\index{Keimel, K.}\index{Wolfenstein, S.}%
\cite[Th\'eor\`eme~11.2.4]{BKW}. In particular,
$G=\set{a}^\bot+\set{a}^{\bot\bot}$, for all
$a\in G$. Assumption~(ii) follows.

Finally, we verify that $G^+$ satisfies Axiom~(M4). Let $a$, $b\in G^+$. We
put $c=a-a\wedge b$. For a polar subset $H$
of $G$, $G$ is the orthogonal sum of $H$ and $H^\bot$ (see, for example,
\cite[Theorem~27]{Birk}).\index{Birkhoff, G.} In  particular, the
projection $p_H$ of $G$ onto
$H$ (relatively to $H^\bot$) is an idempotent homomorphism of
lattice-ordered groups. By Proposition~\ref{P:CharProj}, the projections of
$G^+$ are exactly the restrictions to $G^+$ of the maps of the form $p_H$,
for a polar subset $H$ of $G$. For any such $H$, $p_H(a)\leq p_H(b)$ \iff\
$p_H(c)=0$, that is, $H\subseteq\set{c}^\bot$. Hence, the restriction of
$p_{\set{c}^\bot}$ to $G^+$ is the largest projection $p$ of $G^+$ such that
$p(a)\leq p(b)$.
\end{proof}

\begin{lemma}\label{L:RefSM1}
Let $S$ be a conical \prm. If $S$ is cancellative\index{cancellative} and
satisfies Axiom~\textup{(M2)}, then $\Ref S$ \pup{see
Proposition~\textup{\ref{P:PrmRm}}} is the positive cone of a Dedekind
complete lattice-ordered group.
\end{lemma}

\begin{proof}
By Proposition~\ref{P:PrmRm}, $\Ref S$ is a
\index{refinement monoid} refinement monoid, and it is
generated by $S$ as a monoid. By Lemma~\ref{L:BasicRefS}(ii), $\Ref S$ is a
cancellative\index{cancellative} \cm, thus it is the positive cone $G^+$ of
a directed, partially preordered abelian group~$G$. Since $S$ is conical, so
is $\Ref S$ (see Lemma~\ref{L:BasicRefS}(iii)), and hence $G$ is, in fact,
partially ordered. Since $G^+=\Ref S$ has refinement,
$G$ is an interpolation group, see Proposition~2.1
\index{Goodearl, K.\,R.}in~\cite{Gpoag}.

\setcounter{claim}{0}
\begin{claim}\label{Cl:vwSG}
Let $a$, $b$, $c\in S$ such that $c$ is the infimum of $\set{a,b}$ in $S$.
Then $c$ is also the infimum of $\set{a,b}$ in $G$.
\end{claim}

\begin{cproof}
Let $x\in G$ be a minorant of $\set{a,b}$. Since $G$ has interpolation,
there exists $y\in G$ such that $0$, $x\leq y\leq a$, $b$. So $y\in S$,
hence $x\leq y\leq c$ by the assumption on $c$.
\end{cproof}

\begin{claim}\label{Cl:LatOrd}
$G$ is lattice-ordered.
\end{claim}

\begin{proof}
We put
$X=\setm{x\in G^+}{s\wedge x\text{ exists in }G,\text{ for all }s\in S}$.
It follows from Claim~\ref{Cl:vwSG} that $X$ contains $S$. Let $a$,
$b\in X$, we prove that $a+b\in X$. Let $s\in S$, put $s'=s-s\wedge a$. For
any $x\in G$,
   \begin{align*}
   x\leq s,\,a+b&\text{ \iff\ }x\leq s,\,s+b,\,a+b\\
   &\text{ \iff\ }x\leq(s\wedge a)+s',\,(s\wedge a)+b\\
\intertext{(because $t\mapsto t+b$ is an order-automorphism of $G$)}
   &\text{ \iff\ }x\leq(s\wedge a)+(s'\wedge b),
   \end{align*}
so $s\wedge(a+b)=(s\wedge a)+(s'\wedge b)$. So $X$ is closed under addition,
whence $X=G^+$. Then, replacing $S$ by $G^+$ in the definition of $X$
yields, by a similar argument, that $a\wedge b$ is defined, for all
$a$, $b\in G$.
\end{proof}

\begin{claim}\label{Cl:VWSG}
Let $X$ be a nonempty subset of~$S$ and $a\in S$. If $a$ is the supremum of
$X$ in $S$, then it is also the supremum of $X$ in $G$.
\end{claim}

\begin{cproof}
Let $b\in G$ be a majorant of $X$. Then, using Claim~\ref{Cl:LatOrd},
$a\wedge b$ is also a majorant of $X$, but $0\leq a\wedge b\leq a$, thus
$a\wedge b\in S$. Since $a$ is the supremum of $X$ in $S$,
$a=a\wedge b\leq b$.
\end{cproof}

\begin{claim}
$2S$ satisfies Axiom~\textup{(M2)}.
\end{claim}

\begin{cproof}
Let $X=\setm{c_i}{i\in I}$ be a nonempty subset of $2S$, majorized by some
element of $2S$, say, $a+b$, where $a$, $b\in S$. We prove that $X$ admits a
supremum in $G$.

By (i), $2S$ satisfies
refinement, thus, for all $i\in I$, there are $a_i\leq a$ and $b_i\leq b$
in $S$ such that $c_i=a_i+b_i$. Since $\setm{a_i}{i\in I}$ is a nonempty
subset of~$S$, majorized by $a$, it has, by assumption, a supremum in $S$,
say, $u$. Observe that $a_i\leq u$ (in $S$) for all $i$, and put
$a_i^*=u-a_i$. Then $\setm{b_i-a_i^*\wedge b_i}{i\in I}$ is a nonempty subset
of~$S$, majorized by $b$, thus it has a supremum, say, $v$, in $S$.

For all $i\in I$, $b_i\leq a_i^*\wedge b_i+v\leq a_i^*+v$, thus, by adding
$a_i$ to this inequality, we obtain that $c_i\leq u+v$.

So, to conclude the proof, it suffices to prove that $u+v$ is the least
common majorant of $X$. So let $x$ be a majorant of $X$.
It follows from Claim~\ref{Cl:VWSG} that $u$ is the supremum
of $\setm{a_i}{i\in I}$ in $G$. Then $a_i\leq c_i\leq x$, for all $i\in I$,
thus $u\leq x$. Put $y=u-x$; observe that $y\in G^+$. For $i\in I$,
   \[
   c_i=a_i+b_i\leq x=u+y=a_i+a_i^*+y,
   \]
so $b_i\leq a_i^*+y$. Since $b_i\leq b_i+y$ and $G$ is lattice-ordered,
$b_i\leq(a_i^*+y)\wedge(b_i+y)=a_i^*\wedge b_i+y$, so
$b_i-a_i^*\wedge b_i\leq y$. This holds for all $i$, thus, by
Claim~\ref{Cl:VWSG}, $v\leq y$, so $u+v\leq u+y=x$.
\end{cproof}

The conclusion of the proof is then easy: $2^nS$ is a lower subset of
$2^{n+1}S$, for each $n<\omega$, and all the sets $2^nS$ satisfy (M2), thus
their union, namely, $G^+$, also satisfies (M2). Therefore, $G$ is Dedekind
complete.
\end{proof}

\section[Continuous functions]
{Continuous functions on extremally disconnected\\
topological spaces}\label{S:cmpBooleanspc}

We shall first present a very general result about continuous functions from
extremally disconnected topological spaces to totally ordered sets with their
\index{interval topology}interval topology. Some particular cases of this
result are well-known. For example, if $\Omega$ is a complete Boolean
space,\index{Boolean space} then $\CC(\Omega,\RR)$ is a Dedekind complete
lattice ordered group, see S\"atze 1 and 3 in \cite{Naka41}
\index{Nakano, H.} and Theorem~14 in
\index{Stone, M.\,H.}\cite{Ston49}. We also refer to
\index{Goodearl, K.\,R.}\cite[Lemma~9.1]{Gpoag} for a version of
Proposition~\ref{P:ExtLsc} for $\Omega$ basically disconnected and $K$ an
arbitrary closed interval of $\RR$.

We recall a basic result of general topology, see
\cite[Theorem~10]{Halm}.\index{Halmos, P.\,R.}

\begin{proposition}\label{P:ExtrDisc}\hfill
A topological space is the ultrafilter space of a complete Boolean
algebra\index{Boolean algebra!complete ---} \iff\ it is a complete Boolean
space.\index{Boolean space}
\end{proposition}

\begin{proposition}\label{P:ExtLsc}
Let $\Omega$ be an extremally disconnected topological space, let~$K$ be a
complete totally ordered set. We endow $K$ with its\index{interval topology}
interval topology. Let
$f\colon\Omega\to K$ be a \lsc\ map. Then the map $f^*\colon\Omega\to K$
defined by the rule
   \[
   f^*(x)=\bigwedge_{V\in\Nh(x)}\bigvee f[V],\quad\text{for all }x\in\Omega,
   \]
is continuous, and it is the least continuous map $g\colon\Omega\to K$ such
that $f\leq g$ \pup{with respect to the componentwise ordering of
$\CC(\Omega,K)$}.
\end{proposition}

\begin{proof}
Obviously, $f\leq f^*$.
For any $\alpha\in K$, we define subsets of $\Omega$ by
   \[
   F_\alpha=\setm{x\in\Omega}{f(x)\leq\alpha},\
   G_\alpha=\setm{x\in\Omega}{f^*(x)\leq\alpha},\
   G^*_\alpha=\setm{x\in\Omega}{f^*(x)<\alpha}.
   \]
We record a few basic facts about the sets $F_\alpha$, $G_\alpha$,
$G^*_\alpha$.

\setcounter{claim}{0}
\begin{claim}\label{Cl:BasicFG}
For all $\alpha\in K$, the following assertions hold:

\begin{enumerate}
\item $F_\alpha$ is closed.

\item $G_\alpha\subseteq F_\alpha$.

\item $\overset{\,\circ}{F}_\alpha=\overset{\circ}{G}_\alpha$.
\end{enumerate}
\end{claim}

\begin{cproof}
(i) follows from the lower semicontinuity of $f$.

(ii) follows from the fact that $f\leq f^*$.

(iii) By (ii),
$\overset{\,\circ}{F}_\alpha$ contains $\overset{\circ}{G}_\alpha$. Put
$V=\overset{\,\circ}{F}_\alpha$. For any $x\in V$, we have $V\in\Nh(x)$ and
$\bigvee f[V]\leq\alpha$, hence $f^*(x)\leq\alpha$, that is, $x\in G_\alpha$.
Hence $V$ is contained in~$G_\alpha$, thus, since $V$ is open, in
$\overset{\circ}{G}_\alpha$.
\end{cproof}

To prove that $f^*$ is continuous, it is sufficient to
prove that $G_\alpha^*$ is open and that $G_\alpha$ is closed, for any
$\alpha\in K$.

We start with $G_\alpha^*$. For any $x\in G_\alpha^*$, there exists
$V\in\Nh(x)$ such that $\bigvee f[V]<\alpha$, that is, there exists
$\beta<\alpha$ such that $V\subseteq F_\beta$. Since $V$ is open, it
follows from Claim~\ref{Cl:BasicFG}(iii) that $V\subseteq G_\beta$, whence
$V\subseteq G_\alpha^*$.

Let now $x\in\oll{G_\alpha}$, we prove that $x\in G_\alpha$. Suppose
first that $\alpha$ is right isolated, that is,
$[\alpha,\gamma]=\set{\alpha,\gamma}$ for some $\gamma>\alpha$. Then
$G_\alpha=G^*_\gamma$, and so $G_\alpha$ is open. Since $\Omega$ is
extremally disconnected, $\oll{G_\alpha}$ is clopen. On the other hand,
it follows from Claim~\ref{Cl:BasicFG}(i, ii) that
$\oll{G_\alpha}\subseteq F_\alpha$. Therefore, by
Claim~\ref{Cl:BasicFG}(iii),
$\oll{G_\alpha}\subseteq\overset{\,\circ}{F}_\alpha=
\overset{\circ}{G}_\alpha$, so $x\in G_\alpha$.

Suppose now that $\alpha$ is not right isolated, and put
$U_\beta=\overset{\circ}{F_{\beta}}$. Let $V\in\Nh(x)$. By assumption,
$V\cap G_\alpha\ne\es$, so there exists
$y\in V$ such that $f^*(y)\leq\alpha$. By the definition of $f^*$, there
exists $W\in\Nh(y)$ such that $\bigvee f[W]\leq\beta$,
that is, $W\subseteq F_\beta$. Since $W$ is open, $W\subseteq U_\beta$ as
well, so $y\in U_\beta$ since $y\in W$. Therefore, every neighborhood of $x$
meets $U_\beta$, that is, $x\in\oll{U_\beta}$. However, since $\Omega$ is
extremally disconnected and $U_\beta$ is open, $\oll{U_\beta}$ is open as
well, thus $\oll{U_\beta}\in\Nh(x)$. Furthermore,
$U_\beta\subseteq F_\beta$, thus
$\oll{U_\beta}\subseteq F_\beta$, thus, by Claim~\ref{Cl:BasicFG}(iii),
$\oll{U_\beta}\subseteq G_\beta$. In particular,
$f^*(x)\leq\beta$. This holds for all $\beta>\alpha$, whence
$f^*(x)\leq\alpha$, that is, $x\in G_\alpha$.

Therefore, in both cases, $G_\alpha$ is closed.
So we have proved the continuity of $f^*$.

\begin{claim}\label{Cl:fcontf*}
If $f$ is continuous, then $f=f^*$.
\end{claim}

\begin{cproof}
We prove in fact that for any $x\in\Omega$, if $f$ is continuous at~$x$, then
$f(x)=f^*(x)$. We have already observed that $f(x)\leq f^*(x)$.

To prove the converse inequality, suppose first that $f(x)$ is right isolated,
and let $\gamma>f(x)$ such that $[f(x),\gamma]=\set{f(x),\gamma}$. Since
$f(x)<\gamma$ and since $f$ is continuous, there exists $V\in\Nh(x)$
such that for all $y\in V$, $f(y)<\gamma$, that is, $f(y)\leq f(x)$. Thus
$f^*(x)\leq f(x)$, hence $f^*(x)=f(x)$.

Suppose now that $f(x)$ is not right isolated. Since $f$ is continuous at
$x$, for any $\gamma>f(x)$, there exists $V\in\Nh(x)$ such that
$f(y)<\gamma$ holds for all $y\in V$. Hence $f^*(x)\leq\bigvee
f[V]\leq\gamma$. This holds for all $\gamma>f(x)$ and $f(x)$ is not right
isolated, thus, again, $f^*(x)\leq f(x)$, hence $f^*(x)=f(x)$.
\end{cproof}

If $g\geq f$ is a continuous map from $\Omega$ to~$K$, then, by
Claim~\ref{Cl:fcontf*}, $g=g^*\geq f^*$; the minimality statement follows.
\end{proof}

We introduce a few convenient notations.

\begin{notation}\label{Not:kappa+}
We denote by $\On$\index{Ozzn@$\On$|ii}
the proper class of all ordinals, and we extend the
definitions of $\ZZ_\gamma$, $\RR_\gamma$, and $\two_\gamma$
\index{Zzzgamma@$\ZZ_\gamma$}\index{Rzzgamma@$\RR_\gamma$}%
\index{Tzzgamma@$\two_\gamma$}%
(see Section~\ref{S:BasicRZ2}), defining further proper classes $\Zn$, $\Rn$,
and $\Cn$ as follows:
\index{Zzzinfty@$\ZZ_\infty$|ii}\index{Rzzinfty@$\RR_\infty$|ii}%
\index{Tzzinfty@$\two_\infty$|ii}%
   \begin{align*}
   \Zn&=\ZZ^+\cup\setm{\aleph_\alpha}{\alpha\in\On},\\
   \Rn&=\RR^+\cup\setm{\aleph_\alpha}{\alpha\in\On},\\
   \Cn&=\set{0}\cup\setm{\aleph_\alpha}{\alpha\in\On}.
   \end{align*}
For $\kappa\in\Cn$, we define $\kappa^+$ as the successor cardinal
of $\kappa$ if $\kappa$ is an infinite cardinal, and we put $0^+=\aleph_0$.
That is, $\kappa^+$ is the immediate successor of $\kappa$ in $\Cn$.
\end{notation}

\begin{notation}\label{Not:fresU}
Let $\Omega$ be a set, let $U$ be a subset of $\Omega$, let $K$ be a set
with a distinguished zero element $0$, let $f\colon\Omega\to K$. We denote
by $f\rfloor_U$ the map from $\Omega$ to~$K$ defined by the rule
   \[
   f\rfloor_U(x)=
   \begin{cases}
   f(x)&(\text{if }x\in U),\\
   0&(\text{otherwise}),
   \end{cases}
   \qquad\text{for all }x\in\Omega.
   \]
\end{notation}

\begin{notation}\label{Not:OmiKi}
Let $\Omega$ be a topological space, written as a disjoint union
$\Omega=\bigsqcup_{i<n}\Omega_i$, where $n\in\NN$ and
$\Omega_0$, \dots,
$\Omega_{n-1}$ are clopen subsets of $\Omega$. Let $K_0$, \dots, $K_{n-1}$ be
topological spaces. We define the set
   \[
   \CC(\Omega_0,K_0;\Omega_1,K_1;\dots;\Omega_{n-1},K_{n-1})
   \]
as the set of all maps $f\colon\Omega\to\bigcup_{i<n}K_i$ such that
$f|_{\Omega_i}\in\CC(\Omega_i,K_i)$, for all $i<n$.
\end{notation}

Of course, $\CC(\Omega_0,K_0;\Omega_1,K_1;\dots;\Omega_{n-1},K_{n-1})$ is
naturally isomorphic to the direct product $\prod_{i<n}\CC(\Omega_i,K_i)$.
However, we find the present notation more convenient for such statements as
Proposition~\ref{P:Gendelta}.

\begin{theorem}\label{T:C(O,K)DimInt}
Let $\Omega$ be an extremally disconnected topological space, written as a
disjoint union $\Omega=\Omega_{\I}\sqcup\Omega_{\II}\sqcup\Omega_{\III}$, for
clopen subsets $\Omega_{\I}$, $\Omega_{\II}$, $\Omega_{\III}$ of $\Omega$.
Let $\gamma$ be an ordinal. Then the space
   \[
\CC(\Omega_{\I},\ZZ_\gamma;\Omega_{\II},\RR_\gamma;\Omega_{\III},\two_\gamma)
   \]
\index{Zzzgamma@$\ZZ_\gamma$}\index{Rzzgamma@$\RR_\gamma$}%
\index{Tzzgamma@$\two_\gamma$}%
is a continuous dimension scale.\index{continuous dimension scale}
\end{theorem}

\begin{proof}
Let $K$ be one of the totally ordered sets $\ZZ_\gamma$,
$\RR_\gamma$, or $\two_\gamma$, endowed with its\index{interval topology}
interval topology and its natural monoid structure. We prove that
$S=\CC(\Omega,K)$ is a continuous dimension scale.\index{continuous dimension scale} The
proof of the general case of Theorem~\ref{T:C(O,K)DimInt} follows since the
direct product of any family of continuous dimension scales\index{continuous dimension scale} is again a continuous dimension scale, see Lemma~\ref{L:FinDPDI}.

We first observe that $S$ is a \emph{total} (as opposed to partial)
\cm. It is obvious that the algebraic preordering of
$S$ is antisymmetric and that $S$ has a largest element. Furthermore, let
$\famm{f_i}{i\in I}$ be a nonempty family of elements of~$S$. Define
$f\colon\Omega\to K$ as the componentwise join of all the $f_i$, for $i\in I$.
Then $f$ is \lsc. By Proposition~\ref{P:ExtLsc}, there exists a least
continuous map $f^*\colon\Omega\to K$ such that $f\leq f^*$. Then $f^*$ is the
supremum of $\setm{f_i}{i\in I}$ in $S$. Hence $S$ satisfies Axiom~(M2).

We now observe that $f\rfloor_U$ (see Notation~\ref{Not:fresU}) belongs to
$S$, for any $f\in S$ and any clopen subset $U$ of $\Omega$.
If $f$ is the constant function with value $\alpha$, for $\alpha\in K$, then
we shall write $\alpha\rfloor_U$ instead of $f\rfloor_U$, and we shall of
course write $\alpha$ instead of $\alpha\rfloor_\Omega$.

We shall use the following notation. For $f$, $g\in S$, we put
\index{fgbzzvvfg@$\vbv{f\leq g}$, $\vbv{f=g}$, $\vbv{f<g}$|ii}
   \begin{align*}
   \vbv{f\leq g}&=\setm{x\in\Omega}{f(x)\leq g(x)},\\
   \vbv{f=g}&=\setm{x\in\Omega}{f(x)=g(x)},\\
   \vbv{f<g}&=\setm{x\in\Omega}{f(x)<g(x)}.
   \end{align*}
We observe that $\vbv{f\leq g}$ and $\vbv{f=g}$ are closed, while
$\vbv{f<g}$ is open.

For any $f\in S$, put $U_f=\vbv{f<\aleph_0}$. Then $U_f$ is
open, thus, since $\Omega$ is extremally disconnected, $\oll{U_f}$ is
clopen. We put $\widetilde{f}=f\rfloor_{\Omega\setminus\oll{U_f}}$.
Observe that $\widetilde{f}$ has only infinite values; in particular, it is
purely infinite.

\setcounter{claim}{0}
\begin{claim}\label{Cl:CancDf}\label{Cl:CanDfT}
Let $f\in S$. Then the following are equivalent:
\begin{enumerate}
\item $f$ is cancellable\index{cancellable (element)} in $S$.

\item $f$ is directly finite\index{directly finite} in $S$.

\item $\widetilde{f}=0$.
\end{enumerate}
\end{claim}

\begin{cproof}
(i)$\Rightarrow$(ii) is trivial.

(ii)$\Rightarrow$(iii) Suppose that $\widetilde{f}>0$. Then the clopen set
$V=\Omega\setminus\oll{U_f}$ is nonempty, and $\aleph_0\rfloor_V+f=f$, so
$f$ is not\index{directly finite} directly finite.

(iii)$\Rightarrow$(i) If $\widetilde{f}=0$, then $U_f$ is\index{dense} dense
in $\Omega$. Now let $g$, $h\in S$ such that $f+g=f+h$. Then $g(x)=h(x)$ for
all $x\in U_f$ (because the values of $f$ on $U_f$ are finite), thus, since
$U_f$ is dense in $\Omega$, $g=h$. Hence $f$ is\index{cancellable (element)}
cancellable.
\end{cproof}

We thus obtain Axiom~(M5).

\begin{claim}\label{Cl:pi+df}
Every element of~$S$ can be written under the form $g+h$, where $g$ is
cancellable\index{cancellable (element)} and $h$ is purely infinite.
\end{claim}

\begin{cproof}
Let $f\in S$. By Claim~\ref{Cl:CancDf}, $g=f\rfloor_{\oll{U_f}}$
is\index{cancellable (element)} cancellable. We put $h=\widetilde{f}$. Since
$f=g+h$ and $h$ is purely infinite, the conclusion of the claim holds.
\end{cproof}

Now a weak form of pseudo-cancellation.

\begin{claim}
$S$ satisfies the following statement:
   \[
   \forall a,b,c,\ a+c\leq b+c\Rightarrow\exists x,
   \ x+c=c\text{ and }a\leq b+x.
   \]
\end{claim}

\begin{cproof}
This follows immediately from Claims \ref{Cl:CanDfT} and \ref{Cl:pi+df}, by
putting $x=\tilde{c}$.
\end{cproof}

As a consequence, the equality $\widetilde{f}=\di{f}$ holds (see
Corollary~\ref{C:a/infty}), for any $f\in S$, that is,
$\widetilde{f}$ is, indeed, the largest $g\in S$ such that $g\ll f$.

Furthermore, we observe that for $f$, $g\in S$, the componentwise meet
$f\wedge g$ of $\set{f,g}$ belongs to $S$, and it is the meet of $\set{f,g}$
in $S$. Obviously, $(f+h)\wedge(g+h)=(f\wedge g)+h$, for all $f$, $g$,
$h\in S$. By using Proposition~1.23 of \index{Wehrung, F.}\cite{Wehr92a}, it
follows that $S$ satisfies the refinement property. So, $S$ satisfies
Axiom~(M1).

Now we characterize the projections of~$S$. For any clopen set $U$ of
$\Omega$, the map $f\mapsto f\rfloor_U$ defines obviously a projection of~$S$.

\begin{claim}\label{Cl:SGenComp}
$S$ has general comparability.
\end{claim}

\begin{cproof}
Let $f$, $g\in S$. Put $V=\oll{\vbv{f<g}}$. Since $\Omega$ is extremally
disconnected, $V$ is clopen. It is obvious that $p_V(f)\leq p_V(g)$
and that $p_{\Omega\setminus V}(g)\leq p_{\Omega\setminus V}(f)$.
\end{cproof}

\begin{claim}\label{Cl:ProjPU}
The projections of~$S$ are exactly the $p_U$, where $U$ is a clopen subset
of $\Omega$.
\end{claim}

\begin{cproof}
Let $p$ be a projection of~$S$. Put $u=p(\aleph_\gamma)$ and
$v=p^\bot(\aleph_\gamma)$. Then $\aleph_\gamma=u+v$ and
$u\wedge v=0$, thus there exists a clopen subset $U$ of $\Omega$ such that
$u=\aleph_\gamma\rfloor_U$. So $p(f)\leq f\rfloor_U$, for all $f\in S$.
Conversely, $f=p(f)+p^\bot(f)\leq p(f)+f\rfloor_{\Omega\setminus U}$, thus
$f\rfloor_U\leq p(f)$. So
$p(f)=f\rfloor_U$, for all $f\in S$.
\end{cproof}

\begin{claim}\label{Cl:Compbvfleg}
$S$ satisfies Axiom~\textup{(M4)}.
\end{claim}

\begin{cproof}
Let $f$, $g\in S$, we prove that there exists a largest projection $p$ of
$S$ such that $p(f)\leq p(g)$. For $U$ a clopen subset of
$\Omega$, $p_U(f)\leq p_U(g)$ \iff\ $U\subseteq F$, where we put
$F=\vbv{f\leq g}$. Since $\Omega$ is extremally
disconnected, $\overset{\,\circ}{F}$ is clopen, hence, by
Claim~\ref{Cl:ProjPU}, $p_{\overset{\,\circ}{F}}$ is the largest projection
$p$ of~$S$ such that $p(f)\leq p(g)$.
\end{cproof}

By Claims \ref{Cl:ProjPU} and \ref{Cl:Compbvfleg}, it follows that
for any $f\in S$, $\cc(f)=p_U$, where we put $U=\oll{\vbv{0<f}}$.

It remains to verify (M6).

\begin{claim}\label{Cl:Charrem}
For all $f$, $g\in\DI{S}$, $f\rem g$ \iff\ $f\leq g$ and\linebreak
$\vbv{0<g}\cap\vbv{f=g}$ is nowhere dense.
\end{claim}

\begin{cproof}
We use the alternate form of $\rem$ given in Lemma~\ref{L:GCIVPtr}(ii),
justified by Lemma~\ref{L:GCimpPC} and Claim~\ref{Cl:SGenComp}.

We put $A=\vbv{0<g}\cap\vbv{f=g}$. Observe that
$\vbv{0<g}=\vbv{\aleph_0\leq g}$ (because $g\in\DI{S}$), hence it is
clopen. Hence $A$ is closed.

Suppose first that $f\rem g$. Of course, it follows that $f\leq g$.
Towards a contradiction, suppose that $A$ is
not nowhere dense. Since $\Omega$ is extremally disconnected,
$U=\overset{\,\circ}{A}$ is nonempty and clopen. Put
$h=g\rfloor_{\Omega\setminus U}$. Then $f+h=g$, thus, since $f\rem g$, $g=h$,
so $g\rfloor_U=0$, a contradiction since $U$ is nonempty and contained in
$\vbv{0<g}$.

Conversely, suppose that $f\leq g$ and that
$A$ is nowhere dense. Let $h\in S$ such that
$f+h=g$. Let $x\in\vbv{0<g}\setminus A$. Since $f(x)+h(x)=g(x)$ with
$f(x)<g(x)$ and $g(x)\geq\aleph_0$, $h(x)=g(x)$. So $g$ and $h$ agree on
$\vbv{0<g}\setminus A$, with $A$ nowhere dense and $\vbv{0<g}$ clopen, hence
$g$ and $h$ agree on $\vbv{0<g}$. It is obvious that both $g$ and $h$ are
zero on $\vbv{g=0}$, so, finally, $g=h$. This proves that $f\rem g$.
\end{cproof}

Towards a proof of (M6), let $f$, $g\in\DI{S}$ with $f\rem g$, and set
$A=\vbv{0\nobreak <\nobreak g}\cap\vbv{f\nobreak =\nobreak g}$.
We define a map $\ol{f}\colon\Omega\to K$ by
the rule
   \[
   \ol{f}(x)=
   \begin{cases}
   0&(\text{if }x\in\vbv{g=0})\\
   f(x)&(\text{if }x\in A)\\
   f(x)^+&(\text{if }x\in\vbv{0<g}\setminus A),
   \end{cases}
   \qquad\text{for any }x\in\Omega.
   \]
In the display above, $f(x)$ is an element of $\Cn$, and $f(x)^+$ denotes
the successor of $f(x)$ in $\Cn$, see Notation~\ref{Not:kappa+}.

At this point, we observe the obvious inequality $f\leq\ol{f}\leq g$.

\begin{claim}\label{Cl:olfusc}
The map $\ol{f}$ is \usc.
\end{claim}

\begin{cproof}
Let $x\in\Omega$, and put $\alpha=f(x)$. Since $\ol{f}|_{\vbv{g=0}}=0$ and
$\vbv{g=0}$ is clopen, $\ol{f}$ is continuous at every point of $\vbv{g=0}$.
Now suppose that $g(x)>0$. Suppose first that $x\in A$. So
$g(x)=f(x)=\alpha$, thus $V=\vbv{0<g}\cap\vbv{g=\alpha}$ is an open
neighborhood of $x$. Let $y\in V$. If $y\in A$, then
$\ol{f}(y)=f(y)=g(y)\leq\alpha$. If $y\notin A$, then $f(y)<g(y)\leq\alpha$,
thus $\ol{f}(y)\leq f(y)^+\leq\alpha$. Therefore,
$f\rfloor_V\leq\ol{f}\rfloor_V\leq\alpha$. Since $f(x)=\alpha$, it follows
that $\ol{f}$ is continuous at $x$.

Now suppose that $x\notin A$. Then $V=\vbv{f<\alpha^+}$ is an open
neighborhood of~$x$, and $f(y)\leq\alpha$ and $\ol{f}(y)\leq\alpha^+$ for
all $y\in V$. So, if $\beta\in K$ such that $\ol{f}(x)<\beta$, that is,
$\alpha^+<\beta$, then $\ol{f}(y)<\beta$ for any $y\in V$. Hence $\ol{f}$ is
\usc\ at~$x$.
\end{cproof}

By Claim~\ref{Cl:olfusc} and Proposition~\ref{P:ExtLsc} (used for the dual
partially ordered set of~$K$), for any \usc\ map $k\colon\Omega\to K$, there
exists a largest element $k_*$ of~$S$ such that $k_*\leq k$, and $k_*$ is
given by the formula
   \[
   k_*(x)=
   \bigvee_{V\in\Nh(x)}\bigwedge k[V],\quad\text{for all }x\in\Omega.
   \]
By Claim~\ref{Cl:olfusc}, we can apply this to $k=\ol{f}$, thus obtaining an
element $h=\ol{f}_*$ of~$S$.
Since the range of $\ol{f}$ is contained in $\two_\gamma$,
\index{Tzzgamma@$\two_\gamma$}
so is the range of $h$, whence $h\in\DI{S}$. Furthermore, $f\leq\ol{f}\leq g$
and, since $f$ is continuous, $f=f_*$ (see Proposition~\ref{P:ExtLsc}), thus
   \[
   f=f_*\leq\ol{f}_*=h\leq\ol{f}\leq g.
   \]
It follows from the definition of $\ol{f}$ that $\ol{f}(x)>0$ \iff\
$g(x)>0$, for any $x\in\Omega$. Since $\vbv{0<g}$ is clopen, it follows that
$\vbv{0<g}=\vbv{0<h}$, whence $\cc(g)=\cc(h)$.

\begin{claim}
The relation $f\rem h$ holds.
\end{claim}

\begin{cproof}
We have seen that $f\leq h$.
Put $B=\vbv{0<h}\cap\vbv{f=h}$. Towards a contradiction, assume that
$\overset{\!\circ}{B}\ne\es$. Since $A$ is closed and nowhere dense,
$U=\overset{\!\circ}{B}\setminus A$ is a nonempty open subset of $\vbv{0<g}$.
Pick $x\in U$ such that
$\alpha=f(x)$ is minimum. Then $V=U\cap\vbv{f<\alpha^+}$
is an open neighborhood of~$x$. For any $y\in V$, $f(y)\leq\alpha$ and
$y\in U$, thus, by minimality of $\alpha$, $f(y)=\alpha$.
This holds for any $y\in V$, and $V\subseteq\vbv{0<g}\setminus A$, thus
$\ol{f}(y)=\alpha^+$ for any $y\in V$. Since $V$ is an open neighborhood of
$x$, $h(x)=\alpha^+>\alpha=f(x)$, which contradicts $x\in B$. By
Claim~\ref{Cl:Charrem}, $f\rem h$.
\end{cproof}

To conclude the proof of Theorem~\ref{T:C(O,K)DimInt}, it suffices to
prove that if $k$ is any element of $\DI{S}$ such that $f\rem k$ and
$k^\bot=g^\bot$, then $h\leq k$. Observe first that from the assumption
$k^\bot=g^\bot$, Lemma~\ref{L:cc(a)}(i) implies that $\cc(k)=\cc(g)$, and so
$\vbv{0<g}=\vbv{0<k}$.

Since $f\rem k$,
$f\leq k$ and $B=\vbv{0<k}\cap\vbv{f=k}$ is nowhere dense.
For any $x\in\vbv{0<f}\setminus B$, the inequality $f(x)<k(x)$ holds, thus
$h(x)\leq\ol{f}(x)\leq f(x)^+\leq k(x)$. Since $B$ is nowhere dense and
$\vbv{0<f}$ is open, $h|_{\vbv{0<f}}\leq k|_{\vbv{0<f}}$. Now let
$x\in\vbv{f=0}$. If $g(x)=0$, then $h(x)=0\leq k(x)$. If $g(x)>0$, then,
since $\vbv{0<g}=\vbv{0<k}$, $k(x)>0$, so $h(x)=\aleph_0\leq k(x)$.
Therefore, we have proved that $h\leq k$.
\end{proof}

\section{Completeness of the Boolean algebra of projections}
\label{S:cBaProj}\index{Boolean algebra!complete ---}

\begin{quote}
\em Standing hypothesis: $S$ is a continuous dimension scale\index{continuous dimension scale}.
\end{quote}

We first prove the following lemma.

\begin{lemma}\label{L:OrthCpl}
Let $a$, $b\in S$ and $X\subseteq S$ such that $b=\bigvee X$. If
$a\perp X$, then $a\perp b$.
\end{lemma}

\begin{proof}
The statement $a\perp X$ means that $X\subseteq a^\bot$, that is, by
Lemma~\ref{L:cc(a)}, $X\subseteq\bv{a=0}S$. But by
Lemma~\ref{L:ProjCont}(ii), the range of any projection of~$S$ is closed under
suprema. In particular, $b\in\bv{a=0}S$.
\end{proof}

We prove here an important structural result about\index{pzzroj@$\BB{S}$}
$\BB{S}$.

\begin{proposition}\label{P:B(S)cBa}
The Boolean algebra\index{Boolean algebra!complete ---}
$\BB{S}$\index{pzzroj@$\BB{S}$} is complete. If $\famm{p_i}{i\in I}$ is any
family of projections of~$S$, then, for all $x\in S$,

\begin{enumerate}
\item If $I\ne\es$, then
$\bigl(\bigwedge_{i\in I}p_i)(x)=\bigwedge_{i\in I}p_i(x)$.

\item $\bigl(\bigvee_{i\in I}p_i)(x)=\bigvee_{i\in I}p_i(x)$.
\end{enumerate}

Furthermore, if $p=\bigwedge_{i\in I}p_i$, then $pS=\bigcap_{i\in I}p_iS$.

\end{proposition}

\begin{proof}
The cases of (i) and (ii) where $I$ is finite follow from
Proposition~\ref{P:pvwq(x)}. Therefore, by replacing in (i) (resp., in (ii))
the family $\famm{p_i}{i\in I}$ by the family of all nonempty finite meets
(resp., finite joins) of the $p_i$-s, we can assume without loss of generality
that $I$ is an upward directed partially ordered set and that $i\leq j$ in
$I$ implies $p_i\geq p_j$ (resp., $p_i\leq p_j$).

Let us suppose that
   \begin{equation}\label{Eq:vecpinc}
   i\leq j\text{ implies that }p_i\leq p_j,\quad\text{for all }i,\,j\in I.
   \end{equation}
We consider only the nontrivial case where $I\ne\es$. The
supremum $p(x)=\bigvee_{i\in I}p_i(x)$ is defined and
$p(x)\leq x$, for all $x\in S$. Since $I\ne\es$,
$q(x)=\bigwedge_{i\in I}p_i^\bot(x)$ is also defined.

\setcounter{claim}{0}
\begin{claim}
The maps $p$ and $q$ are endomorphisms of~$S$.
\end{claim}

\begin{cproof}
It is obvious that $p(0)=0$.
Let $x$, $y$, $z\in S$ such that $z=x+y$. We compute:
   \begin{align*}
   p(x)+p(y)&=\bigvee_{i\in I}p_i(x)+\bigvee_{j\in I}p_j(y)\\
   &=\bigvee_{(i,j)\in I\times J}(p_i(x)+p_j(y))
   &&(\text{by Corollary~\ref{C:+MJ}(ii)})\\
   &=\bigvee_{k\in I}(p_k(x)+p_k(y))
   &&(\text{because }I\text{ is upward directed}\\
   &&&\text{\hphantom{because is upward}} \text{and by
\eqref{Eq:vecpinc}})\\
   &=\bigvee_{k\in I}p_k(x+y)\\
   &=p(x+y).
   \end{align*}
A similar argument, based on Corollary~\ref{C:+MJ}(i), proves that $q$ is an
endomorphism of~$S$.
\end{cproof}

\begin{claim}
$pS\perp qS$.
\end{claim}

\begin{cproof}
Let $x$, $y\in S$. For all $i\in I$, $p_i(x)\perp p_i^\bot(y)$, thus, since
$p_i^\bot(y)\geq q(y)$, $p_i(x)\perp q(y)$. This holds for all $i\in I$,
thus, by Lemma~\ref{L:OrthCpl}, $p(x)\perp q(y)$.
\end{cproof}

\begin{claim}\label{Cl:pperpq}
$p$ and $q$ are projections of~$S$, and $q=p^\bot$.
\end{claim}

\begin{cproof}
By the definition of a projection (Definition~\ref{D:Proj}) and by
Claim~\ref{Cl:pperpq}, it remains to prove that $x=p(x)+q(x)$ for all
$x\in S$.

For all $i\in I$, $x=p_i(x)+p_i^\bot(x)$. Since $q\leq p_i^\bot$,
$p_i(x)+q(x)$ is defined and $p_i(x)+q(x)\leq x$.
This holds for all $i$, thus, by Corollary~\ref{C:+MJ}(ii),
$p(x)+q(x)$ is defined and $p(x)+q(x)\leq x$.

For all $i\in I$, $p_i(x)\leq p(x)\leq x$, thus
$x\sd p(x)\leq x\sd p_i(x)\leq p_i^\bot(x)$. By taking the infimum over all
$i$, we obtain that $x\sd p(x)\leq q(x)$; whence $x\leq p(x)+q(x)$.

Therefore, $x=p(x)+q(x)$.
\end{cproof}

It follows easily (use Lemma~\ref{L:pLeqOrtq}) that $p=\bigvee_{i\in I}p_i$
and that $q=\bigwedge_{i\in I}p_i^\bot$. This proves simultaneously (i) and
(ii).

Finally, suppose that $p=\bigwedge_{i\in I}p_i$. We prove that
$pS=\bigcap_{i\in I}p_iS$. For $I=\es$, this is trivial ($p=1$), so suppose
that $I\ne\es$. We have seen above that $p(x)=\bigwedge_{i\in I}p_i(x)$, for
all $x\in S$. Hence, $x\in pS$ \iff\ $p(x)=x$, \iff\ $p_i(x)=x$ for all $i$
(because $p_i(x)\leq x$ for all $i$), \iff\ $x\in\bigcap_{i\in I}p_iS$.
\end{proof}

\begin{proposition}\label{P:Xbotbot}
Let $X$ be a subset of~$S$. Then $S=X^\bot\oplus X^{\bot\bot}$.
\end{proposition}

\begin{proof}
By using Proposition~\ref{P:B(S)cBa}, we define
a projection $p$ of~$S$ by the formula
   \[
   p=\bigvee_{x\in X}\cc(x).
   \]
Then $p^\bot=\bigwedge_{x\in X}\cc(x)^\bot=\bigwedge_{x\in X}\bv{x=0}$,
thus, by Proposition~\ref{P:B(S)cBa} and Lemma~\ref{L:cc(a)}(i),
   \[
   p^\bot S=\bigcap_{x\in X}\bv{x=0}S=\bigcap_{x\in X}x^\bot=X^\bot.
   \]
It follows from this that $pS=(p^\bot S)^\bot=X^{\bot\bot}$.
\end{proof}

\section{The elements $\scal{p}{\alpha}$}\label{S:alphap}

\begin{quote}
\em Standing hypothesis: $S$ is a continuous dimension scale\index{continuous dimension scale}.
\end{quote}

We shall now define a certain doubly indexed family of elements of $\DI{S}$.
These elements represent in some sense the ``layers'' of $\DI{S}$, and a
process of ``measuring'' $\DI{S}$ against these elements will allow us to
pin down the dimension theory of $\DI{S}$.

\begin{notation}
\index{pkszzcal@$\scal{p}{\kappa}$|ii}
For any\index{pzzroj@$\BB{S}$} $p\in\BB{S}$, we define inductively a
transfinite sequence of elements $\scal{p}{\kappa}$
\index{pzzscalalp@$\scal{p}{\kappa}$|ii} of $\DI{S}$, for certain
elements $\kappa$ of $\Cn$, as follows.
\begin{enumerate}
\item $\scal{p}{0}=0$.

\item Let $\kappa\in\Cn$, and suppose that $\scal{p}{\kappa}$ is defined,
and that it is purely infinite. We put
   \[
   X=\Setm{x\in\DI{S}}{\scal{p}{\kappa}\rem x\text{ and }\cc(x)=p}.
   \]
If $X$ is nonempty, then, by Axiom~(M6), it has a least
element. We denote this element by $\scal{p}{\kappa^+}$.
If $X=\es$, then we say that $\scal{p}{\kappa^+}$ is undefined.

\item Let $\lambda$ be a limit cardinal.
Suppose that $\scal{p}{\alpha}$ has been defined for all $\alpha<\lambda$ in
$\Cn$, and that these elements form an increasing, majorized sequence of
elements of $\DI{S}$. Then we put
   \[
   \scal{p}{\lambda}=\bigvee_{\alpha<\lambda}\scal{p}{\alpha}.
   \]
Otherwise, we say that $\scal{p}{\lambda}$ is undefined.

\end{enumerate}

For any\index{pzzroj@$\BB{S}$} $p\in\BB{S}$, we define $\Lambda_p$
\index{Lzzambdap@$\Lambda_p$|ii} as the
class of all $\alpha\in\Cn$ such that $\scal{p}{\alpha}$ is defined.
\end{notation}

Observe that $\scal{0}{\kappa}=0$ for all $\kappa$. In particular, it
follows that $\Lambda_0=\Cn$. The following lemma summarizes some elementary
properties of the elements~$\scal{p}{\alpha}$.

\begin{lemma}\label{L:Basicekp}\hfill
\begin{enumerate}
\item $\Lambda_p$ is a proper initial segment
of $\Cn$, for all\index{pzzrojst@$\BBp{S}$} $p\in\BBp{S}$.

\item $\scal{p}{\alpha}<\scal{p}{\beta}$ for all
$p\in\BBp{S}$\index{pzzrojst@$\BBp{S}$} and all $\alpha<\beta$ in
$\Lambda_p$.

\item $\scal{p}{\alpha}\rem\scal{p}{\beta}$ for all
$p\in\BB{S}$\index{pzzroj@$\BB{S}$} and all $\alpha<\beta$ in $\Lambda_p$.

\item $\scal{p}{0}=0$, and $\cc(\scal{p}{\alpha})=p$ for all
$p\in\BB{S}$\index{pzzroj@$\BB{S}$} and all
$\alpha\in\Lambda_p\setminus\set{0}$.

\item Let $p$, $q\in\BB{S}$\index{pzzroj@$\BB{S}$} such that $p\leq q$. Then
$\Lambda_q\subseteq\Lambda_p$ and $\scal{p}{\alpha}=p(\scal{q}{\alpha})$,
for all $\alpha\in\Lambda_q$.

\item Let $p\in\BB{S}$\index{pzzroj@$\BB{S}$} and let $\alpha$ be an
infinite cardinal number such that $\scal{p}{\alpha}$ is defined.
Then $\scal{p}{\alpha^+}$ is defined if and only if there exists
$x\in\DI{S}$ such that $\scal{p}{\alpha}\rem x$, and then,
$\scal{p}{\alpha^+}$ is the least such $x$.
\end{enumerate}
\end{lemma}

\begin{proof}
(i) and (ii) are obvious.

(iii) is an easy consequence of Lemma~\ref{L:TrLeqTr}.

(iv) By induction on $\alpha\in\Lambda_p\setminus\set{0}$. The assertion
$\scal{p}{0}=0$ holds by definition. If $\alpha=\beta^+$ for some
$\beta\in\Cn$, then, by the definition of $\scal{p}{\alpha}$,
$\cc(\scal{p}{\alpha})=p$.

The limit step is an easy consequence of Lemma~\ref{L:Basiccc}.

(v) We prove the statement by induction on $\alpha\in\Cn$. It is trivial
for $\alpha=0$. The limit step is an easy consequence of
Lemma~\ref{L:ProjCont}(ii).

Now suppose that $\alpha=\beta^+$, for $\beta\in\Cn$.
Since $\scal{q}{\beta^+}$ is defined, $\scal{q}{\beta}$ is defined,
thus, by the induction hypothesis,
   \begin{equation}\label{Eq:timbetpbetq}
   \scal{p}{\beta}\text{ is defined, and }\scal{p}{\beta}=p(\scal{q}{\beta}).
   \end{equation}
We put $e=p(\scal{q}{\alpha})$. Since
$\scal{q}{\beta}\rem\scal{q}{\alpha}$, it follows from
Lemma~\ref{L:ProjTr}(i) that
$p(\scal{q}{\beta})\rem p(\scal{q}{\alpha})$, that is,
$\scal{p}{\beta}\rem e$. Furthermore, by (iv) above
and by Lemma~\ref{L:Basiccc}(ii), $\cc(e)=p\wedge q=p$. Hence,
   \begin{equation}\label{Eq:Defeap}
   \scal{p}{\alpha}\text{ is defined, and }\scal{p}{\alpha}\leq e.
   \end{equation}
So, $\scal{p}{\alpha}\leq p(\scal{q}{\alpha})$, thus, since
$qp^\bot(\scal{q}{\alpha})\leq p^\bot(\scal{q}{\alpha})$, the element
$e'=\scal{p}{\alpha}+qp^\bot(\scal{q}{\alpha})$ is defined, and
$e'\leq \scal{q}{\alpha}$. Note, further, that $e'$ is purely infinite.
Furthermore, by Lemma~\ref{L:ProjTr}(i), the following relations hold:
   \begin{align*}
   \scal{p}{\beta}&\rem\scal{p}{\alpha},\\
   qp^\bot(\scal{q}{\beta})&\rem qp^\bot(\scal{q}{\alpha}).
   \end{align*}
Hence, by Lemma~\ref{L:ProjTr}(ii), we obtain that
   \begin{equation}\label{Eq:betpqpbotq}
   \scal{p}{\beta}+qp^\bot(\scal{q}{\beta})\rem e'.
   \end{equation}
{}From $\scal{q}{\beta}\in qS$, it follows that
$qp^\bot(\scal{q}{\beta})=
p^\bot q(\scal{q}{\beta})=p^\bot(\scal{q}{\beta})$. By using
\eqref{Eq:timbetpbetq}, we obtain that \eqref{Eq:betpqpbotq} can be written
as
   \begin{equation}\label{Eq:timbetqe'}
   \scal{q}{\beta}\rem e'.
   \end{equation}
By Lemma~\ref{L:Basiccc} and by (iv) above, $\cc(e')=q$, hence,
\eqref{Eq:timbetqe'} implies that $\scal{q}{\alpha}\leq e'$.
Hence, by taking the image under $p$ of each side, we obtain that
$e\leq\scal{p}{\alpha}$. Therefore, by \eqref{Eq:Defeap},
$\scal{p}{\alpha}=e=p(\scal{q}{\alpha})$.

(vi) We define subsets $X$ and $Y$ of $\DI{S}$ by
   \begin{align*}
   X&=\Setm{x\in\DI{S}}{\scal{p}{\alpha}\rem x\text{ and }\cc(x)=p},\\
   Y&=\Setm{x\in\DI{S}}{\scal{p}{\alpha}\rem x}.
   \end{align*}
So, $\scal{p}{\alpha^+}$ is defined \iff\ $X$ is nonempty, which implies
that $Y$ is nonempty.

Conversely, suppose that $Y$ is nonempty. For all $x\in Y$,
$\scal{p}{\alpha}\leq x$, thus, since $\alpha>0$ and by (iv), $p\leq\cc(x)$.
Thus, by Lemma~\ref{L:Basiccc}(ii), $\cc(p(x))=p$. Furthermore,
$\scal{p}{\alpha}\in pS$, thus, by Lemma~\ref{L:ProjTr}(i),
$\scal{p}{\alpha}\rem p(x)$, that is, $p(x)\in Y$. Therefore, $pY\subseteq
X$. In particular, $X\ne\es$, so $\scal{p}{\alpha^+}$ is defined. For all
$x\in Y$, $p(x)\in X$, thus $\scal{p}{\alpha^+}\leq p(x)$; whence
$\scal{p}{\alpha^+}\leq x$. Therefore, $\scal{p}{\alpha^+}$ is also the
least element of~$Y$.
\end{proof}

It follows immediately from the definition of the
$(p,\alpha)\mapsto\scal{p}{\alpha}$ operation that
$\scal{p}{\bigvee_{i\in I}\alpha_i}=
\bigvee_{i\in I}\scal{p}{\alpha_i}$ provided the left hand side is defined.
Our next lemma shows that a similar ``linearity'' with respect to the
variable $p$ holds, thus showing a ``bilinearity'' property of the
operation $(p,\alpha)\mapsto\scal{p}{\alpha}$.

\begin{lemma}\label{L:*bilin}
Let $\alpha\in\Cn$, let $\famm{p_i}{i\in I}$ be a family of elements
of\index{pzzroj@$\BB{S}$}
$\BB{S}$. Put $p=\bigvee_{i\in I}p_i$. We make the following
assumptions:
\begin{enumerate}
\item $\scal{p_i}{\alpha}$ is defined for all $i\in I$.

\item $\setm{\scal{p_i}{\alpha}}{i\in I}$ is majorized.
\end{enumerate}
Then $\scal{p}{\alpha}$ is defined, and
$\scal{p}{\alpha}=\bigvee_{i\in I}\scal{p_i}{\alpha}$.

\end{lemma}

\begin{proof}
We argue by induction on $\alpha$.
The supremum $x=\bigvee_{i\in I}\scal{p_i}{\alpha}$ is, by assumption
(ii), defined. Furthermore, the result of Lemma~\ref{L:*bilin} is obvious
for $\alpha=0$. Now suppose that $\alpha>0$.

Suppose first that $\alpha$ is a limit cardinal. By
Lemma~\ref{L:Basicekp}(ii), $\scal{p_i}{\beta}\leq x$ for any cardinal
number $\beta<\alpha$. Therefore, by the induction hypothesis,
$\scal{p}{\beta}$ is defined and
$\scal{p}{\beta}=\bigvee_{i\in I}\scal{p_i}{\beta}\leq x$. This holds for
all $\beta<\alpha$, thus, by definition,
$\scal{p}{\alpha}$ is defined and $\scal{p}{\alpha}\leq x$, thus, since the
converse inequality is obvious, $\scal{p}{\alpha}=x$.

Now we assume that $\alpha=\beta^+$, for some $\beta\in\Cn$.
Put $y=\di{x}$. For any $i\in I$, since $\scal{p_i}{\alpha}$ is purely
infinite, it follows from Corollary~\ref{C:a/infty}(i) that
$\scal{p_i}{\alpha}\leq y$. Furthermore, it follows from
Lemma~\ref{L:Basicekp}(v) that
$p_i(\scal{p_i}{\alpha})=\scal{p_i}{\alpha}$, thus
$\scal{p_i}{\alpha}\leq p_i(y)\leq p(y)$. Hence, by
Lemma~\ref{L:Basicekp}(iv), for all $i\in I$,
$p_i=\cc(\scal{p_i}{\alpha})\leq\cc(p(y))$. Since this holds for all $i\in
I$ and since $\cc(p(y))\leq p$, we obtain the equality
   \begin{equation}\label{Eq:ccp(y)}
   \cc(p(y))=p.
   \end{equation}
For all $i\in I$, since the relations
$\scal{p_i}{\beta}\rem\scal{p_i}{\alpha}$ and $\scal{p_i}{\alpha}\leq p(y)$
hold, we deduce from Lemma~\ref{L:TrLeqTr} that the following relation
holds:
   \begin{equation}\label{Eq:pibremx}
   \scal{p_i}{\beta}\rem p(y).
   \end{equation}
In particular, $\scal{p_i}{\beta}\leq p(y)$, thus, by the induction
hypothesis, $\scal{p}{\beta}$ is defined and $\scal{p}{\beta}\leq p(y)$.
Furthermore, it follows from \eqref{Eq:pibremx} and from
Lemma~\ref{L:ProjTr}(i) that the relation
$p_i(\scal{p_i}{\beta})\rem p_i(y)$ holds for all $i\in I$. By
Lemma~\ref{L:Basicekp}(v),
$p_i(\scal{p_i}{\beta})=p_i(\scal{p}{\beta})=\scal{p_i}{\beta}$, hence
$p_i(\scal{p}{\beta})\rem p_i(y)$. This holds for all $i\in I$, thus, by
Lemma~\ref{L:ProjTr}(ii), $p(\scal{p}{\beta})\rem p(y)$. Again by
Lemma~\ref{L:Basicekp}(v), $p(\scal{p}{\beta})=\scal{p}{\beta}$, so that
$\scal{p}{\beta}\rem p(y)$. Therefore, by \eqref{Eq:ccp(y)},
$\scal{p}{\alpha}$ is defined and $\scal{p}{\alpha}\leq p(y)\leq y\leq x$.
Since the inequality $x\leq\scal{p}{\alpha}$ is obvious,
$x=\scal{p}{\alpha}$.
\end{proof}

\section{The dimension function $\mu$}\label{S:DimFctmu}

\begin{quote}
\em Standing hypothesis: $S$ is a continuous dimension scale\index{continuous dimension scale}.
\end{quote}

For every $x\in S$, we put\index{pzzrojst@$\BBp{S}$}
   \begin{equation}\label{Eq:UxpSpxap}
   U^{(x)}=\Setm{p\in\BBp{S}}{\exists\alpha\in\Lambda_p\text{ such that }
   p\pdi{x}=\scal{p}{\alpha}}.
   \end{equation}
\index{Uzz(x)@$U^{(x)}$|ii}%
The following result is the main fact about the dimension theory of
$\DI{S}$.

\begin{lemma}\label{L:Uxdense}
The set $U^{(x)}$ is a coinitial lower subset of\index{pzzrojst@$\BBp{S}$}
$\BBp{S}$, for all $x\in S$.
\end{lemma}

\begin{proof}
By replacing $x$ by $\di{x}$, we may assume without loss of generality
that $x$ is purely infinite.
The fact that $U^{(x)}$ is a lower subset of
$\BBp{S}$\index{pzzrojst@$\BBp{S}$} is an obvious consequence of
Lemma~\ref{L:Basicekp}. Now let us prove that $U^{(x)}$ is coinitial in
$\BBp{S}$.

So, let\index{pzzrojst@$\BBp{S}$} $p\in\BBp{S}$. We find $q\in(0,p]$ such
that $q\in U^{(x)}$. First, if $p\nleq\cc(x)$, then there exists $q\in(0,p]$
such that $q(x)=0$, thus, obviously, $q\in U^{(x)}$.

So we consider now the case where $p\leq\cc(x)$.
Since the sequence of all $\scal{p}{\alpha}$, for
$\alpha\in\Lambda_p$, is strictly increasing (see Lemma~\ref{L:Basicekp}(ii))
and continuous at limits, there exists a largest element $\alpha$ of $\Cn$
such that $\scal{p}{\alpha}$ is defined and $\scal{p}{\alpha}\leq x$. Since
$0\rem p(x)$ and $\cc(p(x))=p$
(use Lemma~\ref{L:Basiccc}), the element $\scal{p}{\aleph_0}$ is defined,
and $\scal{p}{\aleph_0}\leq p(x)$. Hence, $\alpha$ is an infinite cardinal
number.

Now we put $q=\bv{x\leq\scal{p}{\alpha}}$. By Lemma~\ref{L:BVtr},
   \begin{equation}\label{Eq:qbotalphp}
   q^\bot(\scal{p}{\alpha})\rem q^\bot(x).
   \end{equation}
If $p\leq q^\bot$, then we obtain that
$q^\bot(\scal{p}{\alpha})=q^\bot p(\scal{p}{\alpha})=\scal{p}{\alpha}$,
whence, by \eqref{Eq:qbotalphp}, $\scal{p}{\alpha}\rem x$. By
Lemma~\ref{L:Basicekp}(vi),
$\scal{p}{\alpha^+}$ is defined and $\scal{p}{\alpha^+}\leq x$, which
contradicts the definition of $\alpha$.

Hence, $p\nleq q^\bot$, that is, $r=p\wedge q$ is nonzero. Furthermore,
   \begin{align*}
   r(x)&\leq r(\scal{p}{\alpha})&&(\text{because }r\leq q)\\
   &=\scal{r}{\alpha}&&(\text{because }r\leq p).
   \end{align*}
Since, on the other hand, $\scal{p}{\alpha}\leq x$, we obtain that
$r(x)=\scal{r}{\alpha}$. Therefore, $r\in U^{(x)}$.
\end{proof}

\begin{quote}
\em For the remainder of Section~\textup{\ref{S:DimFctmu}}, we denote by
$\Omega$ the ultrafilter space of\index{pzzroj@$\BB{S}$} $\BB{S}$. By
Propositions
\textup{\ref{P:B(S)cBa}} and \textup{\ref{P:ExtrDisc}},
$\Omega$ is a complete Boolean space.\index{Boolean space}
The clopen sets of $\Omega$ are exactly the sets of the form
   \[
   \Omega_p=\setm{\fa\in\Omega}{p\in\fa},
   \]
for all\index{pzzroj@$\BB{S}$} $p\in\BB{S}$. Moreover, we shall fix an
ordinal $\gamma$ such that $\scal{p}{\aleph_\alpha}$ defined implies
that $\alpha\leq\gamma$, for all\index{pzzrojst@$\BBp{S}$} $p\in\BBp{S}$. The
existence of such a
$\gamma$ is ensured by Lemma~\textup{\ref{L:Basicekp}(i)}.
\end{quote}

We now define, for any $x\in S$ and any
$\fa\in\Omega$,\index{Tzzgamma@$\two_\gamma$}
   \begin{equation}\label{Eq:Defmu}
   \mu(x)(\fa)=\bigvee\setm{\alpha\in\two_\gamma}
   {\exists p\in\fa\text{ such that }\scal{p}{\alpha}\text{ is defined and }
   \scal{p}{\alpha}\leq x}.
   \end{equation}
\index{mzzu(x)@$\mu(x)(\fa)$|ii}%
The rather involved construction of the elements $\scal{p}{\alpha}$ will
give us more control over the function $\mu(x)$ just defined than one has
over (analogues of) the infinite dimension functions on nonsingular
injective modules constructed in\index{Boyle, A.\,K.}
\cite[Chapter XIII]{GoBo} and \index{Goodearl, K.\,R.}
\cite[Chapter 12]{GvnRR}.

\begin{lemma}\label{L:muRelev}
The function $\mu(x)$ is a continuous map from $\Omega$ to
$\two_\gamma$,\index{Tzzgamma@$\two_\gamma$} for any $x\in S$.
\end{lemma}

We recall here that $\two_\gamma$\index{Tzzgamma@$\two_\gamma$} is endowed
with its\index{interval topology} interval topology.

\begin{proof}
For any $\kappa\in\two_\gamma$,\index{Tzzgamma@$\two_\gamma$} we define
subsets $U_\kappa$ and $V_\kappa$ of $\Omega$ by the formulas
   \begin{align*}
   U_\kappa&=\setm{\fa\in\Omega}{\mu(x)(\fa)\geq\kappa},\\
   V_\kappa&=\setm{\fa\in\Omega}{\mu(x)(\fa)>\kappa}.
   \end{align*}

\begin{sclaim}
$V_\kappa$ is open, for any
$\kappa\in\two_\gamma$.\index{Tzzgamma@$\two_\gamma$}
\end{sclaim}

\begin{scproof}
Let $\fa\in V_\kappa$. By the definition of $\mu(x)$, there exist
$\alpha>\kappa$ in $\two_\gamma$ and $p\in\fa$ such that $\scal{p}{\alpha}$
is defined and $\scal{p}{\alpha}\leq x$. Thus, for any $\fb\in\Omega_p$,
$\mu(x)(\fb)\geq\alpha>\kappa$, that is, $\fb\in V_\kappa$.
\end{scproof}

To conclude the proof of Lemma~\ref{L:muRelev}, it suffices to prove that
$U_\kappa$ is closed, for any
$\kappa\in\two_\gamma$.\index{Tzzgamma@$\two_\gamma$} This is trivial for
$\kappa=0$. If $\kappa$ is a limit cardinal, then the equality
   \[
   U_\kappa=\bigcap_{\alpha<\kappa\text{ in }\two_\gamma}U_{\alpha^+}
   \]
holds, hence it is sufficient to prove that $U_{\alpha^+}$ is closed, for any
$\alpha\in\two_\gamma$. Towards this goal, we first observe that
$U_{\alpha^+}=V_\alpha$, thus, by the Claim above, $U_{\alpha^+}$ is open.
Since $\Omega$ is extremally disconnected, the closure
$\oll{U_{\alpha^+}}$ of $U_{\alpha^+}$ is clopen, thus it has the form
$\Omega_p$, for some\index{pzzroj@$\BB{S}$} $p\in\BB{S}$. If $p=0$ then
$\oll{U_{\alpha^+}}=\es$ and we are done, so suppose that $p>0$.
For any $q\in(0,p]$, $\Omega_q$ meets
$U_{\alpha^+}$, thus there exists $\fa\in\Omega_q$ such that
$\mu(x)(\fa)\geq\alpha^+$. Hence there exists $r\in(0,q]\cap\fa$ such that
$\scal{r}{\alpha^+}$ is defined and $\scal{r}{\alpha^+}\leq x$. Therefore,
the set of all $r\in(0,p]$ such that $\scal{r}{\alpha^+}$ is defined and
$\scal{r}{\alpha^+}\leq x$ is coinitial in $(0,p]$, which proves, by
Lemma~\ref{L:*bilin}, that $\scal{p}{\alpha^+}$ is defined and
$\scal{p}{\alpha^+}\leq x$. This means that $\Omega_p\subseteq
U_{\alpha^+}$. Therefore, $U_{\alpha^+}=\Omega_p$ is clopen.
\end{proof}

For all $x\in S$, we put
   \begin{equation}\label{Eq:Omega(x)}
   \Omega^{(x)}=\bigcup\setm{\Omega_p}{p\in U^{(x)}},
   \end{equation}
where $U^{(x)}$ has been defined in \eqref{Eq:UxpSpxap}. It follows from
Lemma~\ref{L:Uxdense} that \emph{$\Omega^{(x)}$ is a dense, open subset of
$\Omega$}.

\begin{lemma}\label{L:mucst}
Let $x\in S$. For any $\fa\in\Omega^{(x)}$, $\mu(x)(\fa)$ is the unique
element $\alpha$ of $\two_\gamma$\index{Tzzgamma@$\two_\gamma$} such that
   \[
   \exists p\in\fa\text{ with }
   \alpha\in\Lambda_p\text{ and }p\pdi{x}=\scal{p}{\alpha}.
   \]
\end{lemma}

\begin{proof}
Let $p\in U^{(x)}$ such that $\fa\in\Omega_p$. By the definition of
$U^{(x)}$, there exists $\alpha\in\Lambda_p$ such that
$p\pdi{x}=\scal{p}{\alpha}$. In particular, $\mu(x)(\fa)\geq\alpha$.

Let $q\in\fa$ and $\beta\in\Lambda_q$ such that
$\scal{q}{\beta}\leq x$. Then $r=p\wedge q$ belongs to $\fa$ and
$\beta\in\Lambda_r$, so $\scal{r}{\beta}\leq\scal{q}{\beta}\leq x$, thus
$\scal{r}{\beta}\leq\di{x}$, from which it follows that
$\scal{r}{\beta}=p(\scal{r}{\beta})\leq p\pdi{x}=\scal{p}{\alpha}$. Hence
$\scal{r}{\beta}=r(\scal{r}{\beta})\leq
r(\scal{p}{\beta})=\scal{r}{\alpha}$, so $\beta\leq\alpha$.
Hence $\mu(x)(\fa)\leq\alpha$, so, finally, $\mu(x)(\fa)=\alpha$.
\end{proof}

\begin{proposition}\label{P:PIdim}\hfill
\begin{enumerate}
\item $\mu(x+y)=\mu(x)+\mu(y)$, for all $x$, $y\in S$ such that $x+y$ is
defined.

\item $\mu(x)\leq\mu(y)$ \iff\ $x\leq y$, for all $x$, $y\in\DI{S}$.

\item The set $\mu\left[\DI{S}\right]$ is a lower subset of
$\CC(\Omega,\two_\gamma)$.\index{Tzzgamma@$\two_\gamma$}
\end{enumerate}
In particular, the restriction of $\mu$ to $\DI{S}$ is a lower
embedding\index{lower embedding} from $\DI{S}$ into
$\CC(\Omega,\two_\gamma)$.\index{Tzzgamma@$\two_\gamma$}
\end{proposition}

\begin{proof}
(i) Put $\Omega'=\Omega^{(x)}\cap\Omega^{(y)}$. Let $\fa\in\Omega'$, and put
$\alpha=\mu(x)(\fa)$ and $\beta=\mu(y)(\fa)$. By
Lemma~\ref{L:mucst}, there exists $p\in\fa$ such that
   \[
   p\pdi{x}=\scal{p}{\alpha}\quad\text{and}\quad p\pdi{y}=\scal{p}{\beta}.
   \]
Hence,
   \begin{align*}
   p\pdi{x+y}&=p\left(\di{x}+\di{y}\right)&&
   (\text{by Lemma~\ref{L:a+b/infty}})\\
   &=p\pdi{x}+p\pdi{y}\\
   &=\scal{p}{\alpha}+\scal{p}{\beta}\\
   &=\scal{p}{\alpha+\beta}&&(\text{by Lemma~\ref{L:Basicekp}(iii)}).
   \end{align*}
Hence, $p\in U^{(x+y)}$, and $\mu(x+y)(\fa)=\mu(x)(\fa)+\mu(y)(\fa)$.
Therefore, $\mu(x+y)$ and $\mu(x)+\mu(y)$ agree on an open dense subset of
$\Omega$, so, since they are continuous, they are equal.

(ii) By (i), $x\leq y$ implies that $\mu(x)\leq\mu(y)$. Conversely,
for any $p\in U^{(x)}\cap U^{(y)}$,
there exist $\alpha$, $\beta\in\two_\gamma$\index{Tzzgamma@$\two_\gamma$}
such that $\scal{p}{\alpha}$ and
$\scal{p}{\beta}$ are defined and equal to $p(x)$ and $p(y)$, respectively.
Since $p\ne0$, there exists $\fa\in\Omega$ such that $p\in\fa$. Then, by
Lemma~\ref{L:mucst}, $\mu(x)(\fa)=\alpha$ and $\mu(y)(\fa)=\beta$, so
$\mu(x)\leq\mu(y)$ implies, by Lemma~\ref{L:Basicekp}(ii), that
$\alpha\leq\beta$. Hence, $p(x)\leq p(y)$, that is, $p\leq\bv{x\leq y}$.
This holds for all $p$ in the coinitial subset $U^{(x)}\cap U^{(y)}$ of
$\BBp{S}$\index{pzzrojst@$\BBp{S}$} (see Lemma~\ref{L:Uxdense}), so
$x\leq y$.

(iii) Let $x\in\DI{S}$, and let
$f\in\CC(\Omega,\two_\gamma)$\index{Tzzgamma@$\two_\gamma$} such that
$f\leq\mu(x)$. We find $y\leq x$ in $\DI{S}$ such that $\mu(y)=f$.

We put\index{pzzrojst@$\BBp{S}$}
$U=\setm{p\in\BBp{S}}{f|_{\Omega_p}\text{ is constant}}$.

\setcounter{claim}{0}
\begin{claim}\label{Cl:Udense}
$U$ is a coinitial lower subset of\index{pzzrojst@$\BBp{S}$} $\BBp{S}$.
\end{claim}

\begin{cproof}
It is obvious that $U$ is a lower subset of\index{pzzrojst@$\BBp{S}$}
$\BBp{S}$.

Let\index{pzzrojst@$\BBp{S}$} $p\in\BBp{S}$, we find $q\in(0,p]\cap U$.
Let $\alpha$ be the minimum value of $f$ on $\Omega_p$.
Then $\Omega'=\setm{\fa\in\Omega}{f(\fa)<\alpha^+}$ is, by continuity of
$f$, an open subset of $\Omega$, and $\Omega'\cap\Omega_p\ne\es$. Let
$q\in(0,p]$ such that $\Omega_q\subseteq\Omega'\cap\Omega_p$. Then
$f|_{\Omega_q}$ is constant with value $\alpha$.
\end{cproof}

Now let $\setm{p_i}{i\in I}$ be a maximal\index{antichain} antichain of
$U\cap U^{(x)}$. By Claim~\ref{Cl:Udense} above, $\setm{p_i}{i\in I}$ is also
a maximal antichain of\index{pzzrojst@$\BBp{S}$} $\BBp{S}$. Let $\alpha_i$
denote the constant value of $f$ on $\Omega_{p_i}$, for all $i\in I$. If
$\fa\in\Omega_{p_i}$, then $\alpha_i=f(\fa)\leq\mu(x)(\fa)$. Since $p_i\in
U^{(x)}$, the equality
$p_i(x)=\scal{p_i}{\mu(x)(\fa)}$ holds. Hence, $\scal{p_i}{\alpha_i}$ is
defined and $\scal{p_i}{\alpha_i}\leq p_i(x)$.

It follows that $\setm{\scal{p_i}{\alpha_i}}{i\in I}$ is a majorized (by
$x$) subset of $\DI{S}$, thus it has a supremum, say, $y$. Note that
$y\leq x$. Furthermore, $\mu(y)(\fa)=\alpha_i=f(\fa)$, for all $i\in I$ and
all $\fa\in\Omega_{p_i}$. Hence, $\mu(y)$ and $f$ agree on a dense open
subset of $\Omega$, so, since they are continuous, $\mu(y)=f$.
\end{proof}

The following trivial property of $\mu$ will later prove very useful.

\begin{proposition}\label{P:muCommProj}
The equality $\mu(p(x))=\mu(x)\rfloor_{\Omega_p}$ holds, for all $x\in S$ and
all\index{pzzroj@$\BB{S}$} $p\in\BB{S}$.
\end{proposition}

\section{Projections on the directly finite elements}\label{S:ProjDF}
\index{directly finite}
We start with the following easy but fundamental result.

\begin{lemma}\label{L:SegDimInt}
Let $S$ be a continuous dimension scale,\index{continuous dimension scale} let $T$ be a
lower subset of~$S$, viewed as a partial submonoid of~$S$ \pup{see
Definition~\textup{\ref{D:PartSubm}}}. Then the following assertions hold:
\begin{enumerate}
\item $T$ is closed under all projections of $S$, and $p|_T\in\BB{T}$
\index{pzzroj@$\BB{S}$} for
all $p\in\BB{S}$.

\item Every projection of $T$ extends to a projection of~$S$.

\item If $T$ is dense\index{dense} in $S$, then every
projection of $T$ extends to a unique projection of~$S$.
\end{enumerate}
\end{lemma}

\begin{proof}
We recall that by Lemma~\ref{L:SegDI}, $T$ is a continuous dimension scale.\index{continuous dimension scale}

Furthermore, $p(x)\leq x$ for any projection $p$ of~$S$ and any $x\in S$.
Since $T$ is a lower subset of $S$ and $p(x)\leq x$ for any projection $p$
of~$S$ and any $x\in S$, (i) holds.

Now we prove (ii).
So let\index{pzzroj@$\BB{S}$} $p\in\BB{T}$. By the definition of a
projection, $T=pT\oplus p^\bot T$. Moreover, it follows from
Proposition~\ref{P:Xbotbot} that
$S=(pT)^\bot\oplus(pT)^{\bot\bot}$, where orthogonals are computed in $S$.
Therefore, there exists a projection $q$ of $S$ such that
$qS=(pT)^{\bot\bot}$ and $q^\bot S=(pT)^\bot$. {}From $pT\subseteq qS$ it
follows that
   \begin{equation}\label{Eq:qx=xpT}
   q(x)=x,\quad\text{for all }x\in pT.
   \end{equation}
Since $p^\bot T$ and $pT$ are orthogonal in $T$ and $T$ is a lower subset of
$S$, $p^\bot T$ and $pT$ are orthogonal in $S$. Hence
$p^\bot T\subseteq(pT)^\bot=q^\bot S$, which implies that
   \begin{equation}\label{Eq:qx=0pbT}
   q(x)=0,\quad\text{for all }x\in p^\bot T.
   \end{equation}
{}From \eqref{Eq:qx=xpT} and \eqref{Eq:qx=0pbT} it follows that $q|_T=p$.

Now we prove (iii). It suffices to prove that if $p$,
$q\in\BB{S}$\index{pzzroj@$\BB{S}$} and
$p(x)\leq q(x)$ for all $x\in T$, then $p\leq q$. Suppose otherwise. Then
$pq^\bot>0$, thus there exists a nonzero element $a$ in $pq^\bot S$. By the
assumption of (iii), we can suppose without loss of generality that $a\in
T$. So $p(a)>0$ while $q(a)=0$ with $a\in T$, a contradiction.
\end{proof}

The following series of results allows to relate the structure of the lower
subset of directly finite\index{directly finite} elements of a continuous dimension scale\index{continuous dimension scale} to the dimension function $\mu$
introduced in Section~\ref{S:DimFctmu}. We first introduce a definition.

\begin{definition}
Let $S$ be a \pcm.\index{partial commutative monoid} An element $a$
of~$S$ is \emph{multiple-free} (or, in some references, \emph{abelian}),
\index{multiple-free (element)|ii} if
$2x\leq a$ implies that $x=0$, for all $x\in S$. We denote by $S_\mf$
\index{Szzmf@$S_\mf$|ii} the subset of $S$ consisting of all multiple-free
elements.
\end{definition}

It is not hard to verify that in any \prm, multiple-free elements are
cancellable\index{cancellable (element)} (see Definition~\ref{D:DirFinMon}).
In continuous dimension scales,\index{continuous dimension scale} this is also an immediate
consequence of Axiom~(M5) and Lemma~\ref{L:DfimCanc}. Multiple-free elements
of a lattice-ordered group are called
\emph{singular} in \cite{BKW}.
\index{Bigard, A.}\index{Keimel, K.}\index{Wolfenstein, S.}%

We recall that $\ZZ_0=\ZZ^+\cup\set{\aleph_0}$
and $\RR_0=\RR^+\cup\set{\aleph_0}$, see Section~\ref{S:BasicRZ2}.

\begin{notation}\label{Not:Cfin}
Let $\Omega$ be a topological space, and let $K$ be either $\ZZ_0$ or
$\RR_0$, endowed with the\index{interval topology} interval topology. We
denote by $\CC_{\fin}(\Omega,K)$\index{Czzfin@$\CC_{\fin}(\Omega,K)$|ii}
the set of all continuous maps $f\colon\Omega\to K$
such that $f^{-1}\set{\aleph_0}$ is nowhere dense. We extend
Notation~\ref{Not:OmiKi} to $\CC_{\fin}$, thus defining, for topological
spaces $\Omega_{\I}$ and $\Omega_{\II}$,
   \[
   \CC_{\fin}(\Omega_{\I},\ZZ_0;\Omega_{\II},\RR_0)=
   \setm{f\in\CC(\Omega_{\I},\ZZ_0;\Omega_{\II},\RR_0)}
   {f^{-1}\set{\aleph_0}\text{ is nowhere dense}}.
   \]
\end{notation}

\begin{proposition}\label{P:EmbCancDImf}
Let $S$ be a stably finite\index{stably finite} continuous dimension scale\index{continuous dimension scale}
\pup{see Definition~\textup{\ref{D:DirFinMon}}}. We suppose that
$S_\mf$ is dense\index{dense} in $S$. Let $\Omega$ be the the ultrafilter
space of\index{pzzroj@$\BB{S}$} $\BB{S}$. Then there exists a map
$\delta\colon S\to\CC_{\fin}(\Omega,\ZZ_0)$ that satisfies the following
conditions:
\begin{enumerate}
\item $\delta$ is a lower embedding\index{lower embedding} \pup{see
Definition~\textup{\ref{D:LowEmb}}}.

\item $\delta(p(x))=\delta(x)\rfloor_{\Omega_p}$, for all $x\in S$ and
all\index{pzzroj@$\BB{S}$}
$p\in\BB{S}$.
\end{enumerate}
\end{proposition}

\begin{proof}[Outline of proof]
Observe that $S$ is cancellative\index{cancellative} (see
Lemma~\ref{L:DfimCanc}). By Lemma~\ref{L:RefSM1}, $\Ref S$ is the positive
cone of a Dedekind complete lattice-ordered group, say, $G$. By
Proposition~\ref{P:lgrpDI},
$G^+$ is a continuous dimension scale,\index{continuous dimension scale} and by
Lemma~\ref{L:SegDimInt}(i, iii), the restriction map\index{pzzroj@$\BB{S}$}
$\BB{G^+}\to\BB{S}$,
$p\mapsto p|_S$ is an isomorphism. Hence it suffices to prove that the
conclusion of Proposition~\ref{P:EmbCancDImf} holds in case $S=G^+$, that
is, $S$ is a \emph{total} monoid.

By Th\'eor\`eme~13.5.6 of
\index{Bigard, A.}\index{Keimel, K.}\index{Wolfenstein, S.}%
\cite{BKW}, there exist a complete
Boolean space\index{Boolean space} $\Omega'$ and a lower
embedding\index{lower embedding}
$\delta$ of $G^+$ into $\CC_{\fin}(\Omega',\ZZ_0)$. This map $\delta$ is
defined as an ``evaluation map'' on the Stone space\index{Stone space}
$\Omega'$, which implies that the condition~(ii) above is satisfied (see
p.~272 in
\index{Bigard, A.}\index{Keimel, K.}\index{Wolfenstein, S.}%
\cite{BKW} for the definition of $\delta$). Furthermore, $\Omega'=\sigma G$
is the ultrafilter space of the complete Boolean algebra
\index{Boolean algebra!complete ---} of polar subsets of
$G$, thus, $\Omega'\cong\Omega$ canonically, so we may replace $\Omega'$ by
$\Omega$.
\end{proof}

In the case where there are no nontrivial multiple-free elements, we get the
following.

\begin{proposition}\label{P:EmbCancDIdiv}
Let $S$ be a stably finite\index{stably finite} continuous dimension scale\index{continuous dimension scale} with no nontrivial multiple-free element.
Then there exists a map
$\delta\colon S\to\CC_{\fin}(\Omega,\RR_0)$ that satisfies the
following conditions:
\begin{enumerate}
\item $\delta$ is a lower embedding.\index{lower embedding}

\item $\delta(p(x))=\delta(x)\rfloor_{\Omega_p}$, for all $x\in S$ and
all\index{pzzroj@$\BB{S}$}
$p\in\BB{S}$.
\end{enumerate}
\end{proposition}

\begin{proof}[Outline of proof]
The proof of Proposition~\ref{P:EmbCancDIdiv} follows the lines of the proof
of Proposition~\ref{P:EmbCancDImf}, with the following changes. The Dedekind
complete lattice-ordered group $G$ has no nontrivial multiple-free
element, thus, by Theorem~11.2.13 of
\index{Bigard, A.}\index{Keimel, K.}\index{Wolfenstein, S.}%
\cite{BKW}, it is divisible. The rest
of the proof is the same as for Proposition~\ref{P:EmbCancDImf}, by using
Corollaire~13.4.2 of
\index{Bigard, A.}\index{Keimel, K.}\index{Wolfenstein, S.}%
\cite{BKW} instead of Th\'eor\`eme~13.5.6 of \cite{BKW}.
\end{proof}

As experience proves, it is often useful to state explicitly the definition
of the embedding $\delta$ of Propositions~\ref{P:EmbCancDImf} and
\ref{P:EmbCancDIdiv}. The definition that we present here is equivalent to
the one given by S.\,J. Bernau's embedding theorem for Archimedean
lattice-ordered groups, see \index{Bernau, S.\,J.}\cite{Bern65}, or
\index{Anderson, M.}\index{Feil, T.}%
\cite[Theorem~2.4]{AnFe88}. It is convenient to first define the concept of a
\emph{finitary unit} in a continuous dimension scale.
\index{continuous dimension scale}

\begin{definition}\label{D:FinUn}
Let $S$ be a continuous dimension scale.\index{continuous dimension scale} A \emph{finitary
unit}\index{finitary unit|ii}
of $S$ is a dense\index{dense}\index{antichain}
antichain $E$ of $S_\fin$ such that for any $e\in E$, either $e$ is 
multiple-free or there is
no nonzero multiple-free element below $e$.
\end{definition}

\begin{lemma}\label{L:Basis}
Every continuous dimension scale\index{continuous dimension scale} has a finitary unit.
\end{lemma}

\begin{proof}
Let $U$ denote the set of all elements $x\in S_\fin\setminus\set{0}$ that
are either multiple-free or without nonzero multiple-free element below. Let
$a\in S_\fin\setminus\set{0}$. If there is no nonzero
multiple-free element below $a$, then $a\in U$. If there exists a nonzero
multiple-free element $e\leq a$, observe that $e\in U$. So $U$ is
dense\index{dense} in $S_\fin$, and the finitary units of $S$ are exactly
the maximal antichains\index{antichain} of $U$.
\end{proof}

Now let $G$ be a Dedekind complete lattice-ordered group. Every polar
subset of $G$ is an orthogonal direct summand of~$G$, thus
$p\mapsto pG^++(-pG^+)$ defines an isomorphism from
$\BB{G^+}$\index{pzzroj@$\BB{S}$} onto the Boolean algebra
\index{Boolean algebra!complete ---} of polar subsets of $G$. We denote
again by $\Omega$ the ultrafilter space of
\index{pzzroj@$\BB{S}$} $\BB G^+$. Let $E$ be a finitary
unit of the continuous dimension scale\index{continuous dimension scale} $G^+$ (see
Definition~\ref{D:FinUn}). We put
   \begin{equation}\label{Eq:Defdelta}
   \delta(x)(\fa)=\bigvee\Setm{m/n}{(m,n)\in\ZZ^+\times\NN\text{ and }
   \bv{me\leq nx}\in\fa\text{ for all }e\in E},
   \end{equation}
\index{dzzelta(x)@$\delta(x)(\fa)$|ii}
for all $x\in G^+$ and $\fa\in\Omega$. Then $\delta(x)$ is a continuous map
from $\Omega$ to $\RR_0$, for all $x\in\Omega$, and $\delta$ is the desired
embedding.

Unlike the map $\mu$ given in \eqref{Eq:Defmu}, the map $\delta$ is
not intrinsic, for it depends on the choice of a finitary unit of~$G^+$.
Furthermore, it is not apparent through \eqref{Eq:Defdelta} that under the
assumptions of Proposition~\ref{P:EmbCancDImf}, the map $\delta$ takes its
values in $\CC_{\fin}(\Omega,\ZZ_0)$ (rather than just in
$\CC_{\fin}(\Omega,\RR_0)$). However, it is possible to prove that under
those assumptions, since $E$ is a finitary unit of $G^+$, the following
equality holds,
   \begin{equation}\label{Eq:Defdelta2}
   \delta(x)(\fa)=\bigvee\Setm{n\in\ZZ^+}
   {\bv{ne\leq x}\in\fa\text{ for all }e\in E},
   \end{equation}
for all $x\in G^+$ and all $\fa\in\Omega$. Hence the map $\delta(x)$ is
integer-valued. Similarly, the proof that the range of $\delta$ is, in the
context of Proposition~\ref{P:EmbCancDImf}, a lower embedding\index{lower
embedding}, uses the assumption that $E$ is a finitary unit of $G^+$.

\begin{definition}\label{D:SI,II,III}
For a general continuous dimension scale\index{continuous dimension scale} $S$, we define
ideals $S_{\I}$,
$S_{\II}$, and $S_{\III}$ of~$S$ as follows:
   \begin{align*}
   S_{\I}&=S_{\mf}^{\bot\bot};\\
   S_{\II}&=S_{\mf}^\bot\cap S_{\fin}^{\bot\bot};\\
   S_{\III}&=S_{\fin}^{\bot}.
   \end{align*}
We say that $S$ is \emph{Type $\I$} (resp., \emph{Type~$\II$},
\emph{Type~$\III$}),
\index{Type I, II, III (continuous dimension scale)|ii} if
$S_{\II}=S_{\III}=\set{0}$ (resp., $S_{\I}=S_{\III}=\set{0}$,
$S_{\I}=S_{\II}=\set{0}$).
\end{definition}

It follows from Proposition~\ref{P:Xbotbot} that the equality
   \[
   S=S_{\I}\oplus S_{\II}\oplus S_{\III}
   \]
holds (see Notation~\ref{Not:SumSubsets}). Observe, in particular, that
$S_\fin\subseteq S_{\fin}^{\bot\bot}=S_{\I}\oplus S_{\II}$.
We denote by $p_{\I}$ (resp., $p_{\II}$, $p_{\III}$) the projection
of~$S$ on $S_{\I}$ (resp., $S_{\II}$, $S_{\III}$).
So $p_{\I}\oplus p_{\II}\oplus p_{\III}=1$ in
\index{pzzroj@$\BB{S}$} $\BB{S}$, hence
$\Omega=\Omega_{\I}\sqcup\Omega_{\II}\sqcup\Omega_{\III}$,
where we put
   \[
   \Omega_{\I}=\Omega_{p_{\I}},\qquad\Omega_{\II}=\Omega_{p_{\II}},\qquad
   \Omega_{\III}=\Omega_{p_{\III}}.
   \]
Observe that $\Omega_{\I}$, $\Omega_{\II}$, and $\Omega_{\III}$ are clopen
subsets of $\Omega$.

By using Propositions \ref{P:EmbCancDImf} and \ref{P:EmbCancDIdiv}, we
obtain lower embeddings\index{lower embedding}
   \[
   \delta_{\I}\colon S_{\I}\cap S_\fin\hookrightarrow
   \CC_{\fin}(\Omega_{\I},\ZZ_0)\quad\text{and}\quad
   \delta_{\II}\colon S_{\II}\cap S_\fin\hookrightarrow
   \CC_{\fin}(\Omega_{\II},\RR_0)
   \]
such that $\delta_i(p(x))=\delta_i(x)\rfloor_{\Omega_p}$, for all
$i\in\set{\I,\II}$, all $x\in S_i$, and all $p\in\BB{S_i}$. Now we identify
$\BB(S_{\I}\oplus S_{\II})$\index{pzzroj@$\BB{S}$} with $\BB(S_\fin)$,
\emph{via} Lemma~\ref{L:SegDimInt}. Hence, by combining $\delta_{\I}$ and
$\delta_{\II}$, we obtain the following result.

\begin{proposition}\label{P:Gendelta}
Let $S$ be a continuous dimension scale.\index{continuous dimension scale} Then there exists
a map $\delta\colon S_\fin\to
\CC(\Omega_{\I},\ZZ_0;\Omega_{\II},\RR_0;\Omega_{\III},\set{0})$
\pup{see Definition~\textup{\ref{D:SI,II,III}}}
that satisfies the following conditions:
\begin{enumerate}
\item $\delta$ is a lower embedding.\index{lower embedding}

\item The values of $\delta(x)$ are finite on an open dense subset of
$\Omega$, for every $x\in S_\fin$.

\item $\delta(p(x))=\delta(x)\rfloor_{\Omega_p}$, for all $x\in S_\fin$ and
all\index{pzzroj@$\BB{S}$} $p\in\BB{S}$.
\end{enumerate}
\end{proposition}

The map $\delta$ depends on the choice of a finitary unit
$E$ of $S_\fin$, and it is then given by the formula
\eqref{Eq:Defdelta}, for all $x\in S_\fin$ and all
$\fa\in\Omega_{\I}\cup\Omega_{\II}$.

\section{Embedding arbitrary continuous dimension scales}\label{S:EmbDI}
\index{continuous dimension scale}
\begin{quote}
\em Standing hypotheses: $S$ is a continuous dimension scale\index{continuous dimension scale}, $S_{\I}$, $S_{\II}$,
$S_{\III}$, $\Omega_{\I}$, $\Omega_{\II}$, $\Omega_{\III}$,
$p_{\I}$, $p_{\II}$, $p_{\III}$ are as in Section~\textup{\ref{S:ProjDF}}.

Let $\gamma$ be the ordinal and
$\mu\colon S\to\CC(\Omega,\two_\gamma)$\index{Tzzgamma@$\two_\gamma$} the
dimension function defined in Section~\textup{\ref{S:DimFctmu}}.

We put $\ol{S}=\CC
(\Omega_{\I},\ZZ_\gamma;\Omega_{\II},\RR_\gamma;\Omega_{\III},\two_\gamma)$.
\index{Zzzgamma@$\ZZ_\gamma$}\index{Rzzgamma@$\RR_\gamma$}%
\index{Tzzgamma@$\two_\gamma$}\index{Szzbar@$\ol{S}$|ii}%
We pick a lower embedding\index{lower embedding}
$\delta\colon S_\fin\hookrightarrow\ol{S}$ as in
Proposition~\textup{\ref{P:Gendelta}}.
\end{quote}

\begin{lemma}\label{L:LocDF}
Let $\famm{a_i}{i\in I}$ be a majorized family of pairwise orthogonal
directly finite\index{directly finite} elements of~$S$. Then
$a=\oplus_{i\in I}a_i$ is directly finite.
\end{lemma}

\begin{proof}
Put $p_i=\cc(a_i)$, for all $i\in I$. Observe that $a_i=p_i(a)$, for all $i$.
Let $x\in S$ such that $a+x=a$. It follows that for any $i\in I$,
$a_i+p_i(x)=a_i$, thus, since $a_i$ is\index{directly finite} directly
finite, $p_i(x)=0$, that is,
$x\in p_i^\bot S$. Put $p=\bigvee_{i\in I}p_i$. By Lemma~\ref{L:Basiccc},
$p=\cc(a)$, thus, by using Lemma~\ref{L:cc(a)} and
Proposition~\ref{P:B(S)cBa},
   \[
   x\in\bigcap_{i\in I}p_i^\bot S=p^\bot S=a^\bot.
   \]
But $x\leq a$, therefore, $x=0$.
\end{proof}

\begin{lemma}\label{L:valmuinf}
Let $x\in\DI{S}$ and $\fa\in\Omega$. If $\cc(x)\in\fa$, then
$\mu(x)(\fa)\geq\aleph_0$.
\end{lemma}

\begin{proof}
For all $p\in(0,\cc(x)]$, there exist, by Lemma~\ref{L:Uxdense}, an
infinite cardinal number $\alpha$ and $q\in(0,p]$ such that
$q(x)=\alpha\cdot q$. So, for any $\fa\in\Omega_q$,
$\mu(x)(\fa)=\alpha\geq\aleph_0$. Therefore, the set of all
$\fa\in\Omega_{\cc(x)}$ such that $\mu(x)(\fa)\geq\aleph_0$ is dense in
$\Omega_{\cc(x)}$. Since $\mu(x)$ is continuous, the conclusion of
Lemma~\ref{L:valmuinf} follows.
\end{proof}

\begin{lemma}\label{L:dflepi}
Let $a\in S_\fin$ and $b\in\DI{S}$. If $\cc(a)\leq\cc(b)$, then $a\leq b$
and $\delta(a)\leq\mu(b)$.
\end{lemma}

\begin{proof}
Put $p=\bv{b\leq a}$ and $q=\bv{a\leq b}$. Then $p(b)\leq p(a)$, with $p(a)$
directly finite\index{directly finite} and $p(b)$ purely infinite, thus
$p(b)=0$, that is,
$p\wedge\cc(b)=0$. By assumption, $p\wedge\cc(a)=0$, that is, $p(a)=0$, thus
$p\leq q$. By general comparability, $p\vee q=1$, so, in fact, $q=1$.
Therefore, $a\leq b$.

Now put $r=\cc(b)$. For any $\fa\in\Omega_r$, it follows from
Lemma~\ref{L:valmuinf} that $\delta(a)(\fa)\leq\aleph_0\leq\mu(b)(\fa)$.
For any $\fa\in\Omega_{r^\bot}$,
   \[
   \delta(a)(\fa)=\delta(a)\rfloor_{\Omega_{r^\bot}}(\fa)=
   \delta(r^\bot(a))(\fa)=0,
   \]
because $r^\bot(a)=0$. Therefore, $\delta(a)(\fa)\leq\mu(b)(\fa)$, for any
$\fa\in\Omega$.
\end{proof}

\begin{lemma}\label{L:InSegf}
Let $f\in\ol{S}$ be\index{directly finite} directly finite, let $b\in\DI{S}$
such that $f\leq\mu(b)$. Then there exists a directly finite $a\leq b$ in $S$
such that $f=\delta(a)$.
\end{lemma}

\begin{proof}
We put $p=\cc(b)\wedge(p_{\I}\oplus p_{\II})$, and we define a subset $U$ of
$\BBp{S}$\index{pzzroj@$\BB{S}$} by
   \[
   U=\setm{q\in(0,p]}{\exists x\in S_\fin,\ \cc(x)=q\text{ and }
   f\rfloor_{\Omega_q}\leq\delta(x)}.
   \]

\setcounter{claim}{0}
\begin{claim}\label{Cl:weirdU}
$U$ is coinitial in $(0,p]$.
\end{claim}

\begin{cproof}
By Lemma~\ref{L:dflepi}, every directly finite\index{directly finite}
element of $pS$ lies below
$b$, thus, since $b$ is purely infinite and by Lemma~\ref{L:leqStillDef},
$pS_\fin$ is a total monoid.

Let $q\in(0,p]$, we prove that $U\cap(0,q]$ is nonempty. Since
$0<q\leq p_{\I}\oplus p_{\II}$ and
$(p_{\I}\oplus p_{\II})S=S_{\fin}^{\bot\bot}$, $qS$ has
a\index{directly finite} directly finite, nonzero element $y$. Observe that
$\delta(y)$ is a nonzero element of
$\CC_{\fin}(\Omega,\RR_0)$, thus, since $\delta(y)$ is continuous, there
exists
$q'\in(0,q]$ such that $\delta(y)(\fa)>0$ for all $\fa\in\Omega_{q'}$.
Without loss of generality, we may assume that $\cc(y)=q'$.

Suppose that $n\delta(y)\leq f$ for all $n\in\ZZ^+$. Then
$f(\fa)\geq\aleph_0$, for any $\fa\in\Omega_{q'}$, so
$f+\aleph_0\rfloor_{\Omega_{q'}}=f$, which contradicts the assumption that
$f$ is\index{directly finite} directly finite. Hence, there exists a largest
nonnegative integer $n$ such that $n\delta(y)\leq f$. Since $\delta(y)$
vanishes outside
$\Omega_{q'}$, there exists $r\in(0,q']$ such that
$f\rfloor_{\Omega_r}\leq(n+1)\delta(y)$. Therefore, $r$ belongs to $U$ (with
witness $x=(n+1)r(y)$).
\end{cproof}

By Claim~\ref{Cl:weirdU}, there exists a maximal\index{antichain} antichain
$W$ of $[0,p]$ such that $W\subseteq U$. For each $q\in W$, pick
a\index{directly finite} directly finite
$x_q\in S$ such that $\cc(x_q)=q$ and $f\rfloor_{\Omega_q}\leq\delta(x_q)$.
Observe that $x_q\leq b$ for all $q$, thus $x=\oplus_{q\in W}x_q$ is defined
and $x\leq b$. Furthermore, by Lemma~\ref{L:LocDF}, $x$ is
\index{directly finite} directly
finite. By Lemma~\ref{L:ProjCont}(ii), $q(x)=x_q$ for all $q\in W$. For any
$q\in W$ and any $\fa\in\Omega_q$,
   \[
   f(\fa)=f\rfloor_{\Omega_q}(\fa)\leq\delta(x_q)(\fa)=\delta(q(x))(\fa)
   =\delta(x)\rfloor_{\Omega_q}(\fa)=\delta(x)(\fa).
   \]
Since $\bigcup_{q\in W}\Omega_q$ is dense in $\Omega_p$ and both $f$ and
$\delta(x)$ are continuous and vanish outside $\Omega_p$, it follows that
$f\leq\delta(x)$. Since $\delta$ is a lower embedding\index{lower
embedding}, there exists $a\leq x$ in $S_\fin$ such that $f=\delta(a)$.
\end{proof}

We now wish to define a homomorphism of partial monoids
$\eps\colon S\to\ol{S}$\index{ezzpsilon@$\eps$ (canonical $S\to\ol{S}$)|ii}
by the rule
   \begin{equation}\label{Eq:DefEps}
   \eps(x+y)=\delta(x)+\mu(y),\qquad\text{for all }x\in S_\fin
   \text{ and all }y\in\DI{S}.
   \end{equation}
Since $\mu(y)$ is purely infinite, $\delta(x)+\mu(y)$ is the maximum
of $\delta(x)$ and $\mu(y)$. It follows that
the existence of $\eps$ is ensured by the following Lemmas
\ref{L:abcinf} and \ref{L:acinfb}.

\begin{lemma}\label{L:abcinf}
Let $a$, $b\in S_\fin$ and $c\in\DI{S}$. If $a\leq b+c$, then
$\delta(a)\leq\delta(b)+\mu(c)$.
\end{lemma}

\begin{proof}
By the refinement property, $a=b'+c'$ for some $b'\leq b$ and $c'\leq c$. In
particular, $b'$ and $c'$ are\index{directly finite} directly finite, and,
of course,
$\cc(c')\leq\cc(c)$. By Lemma~\ref{L:dflepi}, $\delta(c')\leq\mu(c)$. It
follows that $\delta(a)=\delta(b')+\delta(c')\leq\delta(b)+\mu(c)$.
\end{proof}

\begin{lemma}\label{L:acinfb}
Let $a$, $c\in\DI{S}$ and $b\in S_\fin$. If $a\leq b+c$, then
$\mu(a)\leq\delta(b)+\mu(c)$.
\end{lemma}

\begin{proof}
Since $b$ is directly finite\index{directly finite} and $c$ is purely
infinite, $\di{b}=0$ and
$\di{c}=c$. By Lemma~\ref{L:a+b/infty}, $a=\di{a}\leq\di{b}+\di{c}=c$.
Therefore, $\mu(a)\leq\mu(c)\leq\delta(b)+\mu(c)$.
\end{proof}

This shows the existence of a unique homomorphism of partial monoids
$\eps\colon S\to\ol{S}$ satisfying the condition \eqref{Eq:DefEps}.
Observe that the following additional condition is satisfied by
$\eps$ (because it is satisfied by $\mu$ and by $\delta$, see
Propositions \ref{P:muCommProj} and \ref{P:Gendelta}):
\index{pzzroj@$\BB{S}$}
   \begin{equation}\label{Eq:epsProj}
   \eps(p(x))=\eps(x)\rfloor_{\Omega_p},\qquad
   \text{for all }x\in S\text{ and all }p\in\BB{S}.
   \end{equation}
The purpose of the following Lemmas \ref{L:epsEmb} and \ref{L:epsId} is to
prove that $\eps$ is a lower embedding.\index{lower embedding}

\begin{lemma}\label{L:epsEmb}
The map $\eps$ is an order-embedding.
\end{lemma}

\begin{proof}
Let $a$, $b\in S$ such that $\eps(a)\leq\eps(b)$, we prove
that $a\leq b$. By Corollary~\ref{C:a/infty}(iii,iv), there exists a directly
finite\index{directly finite} $b'\in S$ such that $b=b'+\di{b}$.

Suppose now that $a$ is\index{directly finite} directly finite. We put
$p=\cc\pdi{b}$. Then
$\cc(p(a))\leq p=\cc\pdi{b}$, with $p(a)$\index{directly finite} directly
finite and $\di{b}$ purely infinite, thus, by Lemma~\ref{L:dflepi}, we
obtain that
   \begin{equation}\label{Eq:paledib}
   p(a)\leq\di{b}.
   \end{equation}
Furthermore, $p^\bot\pdi{b}=0$, hence, by using \eqref{Eq:epsProj} and the
definition of $\eps$,
   \[
   \delta(p^\bot(a))=\eps(p^\bot(a))\leq
   \eps(p^\bot(b'))+\eps\left(p^\bot\pdi{b}\right)
   =\delta(p^\bot(b')),
   \]
thus, since $\delta$ is an embedding, we obtain that
   \begin{equation}\label{Eq:pbotaleb'}
   p^\bot(a)\leq p^\bot(b')\leq b'.
   \end{equation}
It follows from \eqref{Eq:paledib} and \eqref{Eq:pbotaleb'} that
$a\leq b$.

In the general case, there exists, by Corollary~\ref{C:a/infty}(iii,iv), a
directly finite\index{directly finite} $a'\in S$ such that $a=a'+\di{a}$. It
follows from the previous paragraph that the inequality
   \begin{equation}\label{Eq:a'leqb}
   a'\leq b.
   \end{equation}
holds. Furthermore, dividing by $\infty$ the inequality
$\delta(a')+\mu\pdi{a}\leq\delta(b')+\mu\pdi{b}$ and using
Proposition~\ref{P:PIdim} yields, since both $\delta(a')$ and $\delta(b')$
are finite-valued on a dense subset of $\Omega$ and both $\mu\pdi{a}$ and
$\mu\pdi{b}$ are purely infinite, that
   \begin{equation}\label{Eq:mudialedib}
   \di{a}\leq\di{b}.
   \end{equation}
By \eqref{Eq:a'leqb} and \eqref{Eq:mudialedib}, $a=a'+\di{a}\leq b+\di{b}=b$.
\end{proof}

\begin{lemma}\label{L:epsId}
The range of $\eps$ is a lower subset of $\ol{S}$.
\end{lemma}

\begin{proof}
Let $b\in S$ and let $f\in\ol{S}$ such that $f\leq\eps(b)$. We find
$a\leq b$ in $S$ such that $f=\eps(a)$.
By Corollary~\ref{C:a/infty}(iii,iv), there exists a directly finite
\index{directly finite} $b'\in S$ such that $b=b'+\di{b}$.

We start with the case where $f$ is\index{directly finite} directly finite.
Since
$f\leq\delta(b')+\mu\pdi{b}$ and since $\ol{S}$ satisfies the refinement
property, there are $g$, $h\in\ol{S}$ such that $f=g+h$, $g\leq\delta(b')$, and
$h\leq\mu\pdi{b}$. Since $\delta$ is a lower embedding\index{lower
embedding}, there exists
$x\leq b'$ in $S_\fin$ such that $g=\delta(x)$. Since $f$
is\index{directly finite} directly finite and $h\leq f$, $h$ is directly
finite, thus, by Lemma~\ref{L:InSegf}, there exists a directly finite
$y\leq\di{b}$ such that
$h=\delta(y)$. Since
$x\leq b'$, $y\leq\di{b}$, and $b'+\di{b}=b$, $x+y$ is defined, $x+y\leq b$,
and $\delta(x+y)=g+h=f$.

Now the general case. Since $\ol{S}$ is a continuous dimension scale\index{continuous dimension scale}, there exists a directly finite\index{directly finite}
$f'\in\ol{S}$ such that
$f=f'+\di{f}$. By the previous paragraph, $f'=\delta(a')$ for some directly
finite\index{directly finite} $a'\leq b$. Furthermore, by dividing by
$\infty$ the inequality
$f'+\di{f}\leq\delta(b')+\mu\pdi{b}$, we obtain that $\di{f}\leq\mu\pdi{b}$,
thus, by Proposition~\ref{P:PIdim}, there exists $\oll{a}\leq\di{b}$ in
$\DI{S}$ such that
$\di{f}=\mu(\oll{a})$. Since $a'\leq b$, $\oll{a}\leq\di{b}$, and
$b+\di{b}=b$, $a'+\oll{a}$ is defined, $a'+\oll{a}\leq b$, and
$\eps(a'+\oll{a})=\delta(a')+\mu(\oll{a})=f'+\di{f}=f$.
\end{proof}

We finally arrive at the following more precise version of Theorem~C.

\begin{theorem}\label{T:EmbDimInt}
Let $S$ be a continuous dimension scale,\index{continuous dimension scale} let $\BB{S}$
\index{pzzroj@$\BB{S}$} be
the complete Boolean algebra\index{Boolean algebra!complete ---} of
projections of~$S$, and let $\Omega$ be the ultrafilter space of
$\BB{S}$, with the decomposition
$\Omega=\Omega_{\I}\sqcup\Omega_{\II}\sqcup\Omega_{\III}$ as given in
Section~\textup{\ref{S:ProjDF}}.

Then there exist an ordinal $\gamma$ and a
lower embedding\index{lower embedding}
   \[
   \eps\colon S\hookrightarrow
\CC(\Omega_{\I},\ZZ_\gamma;\Omega_{\II},\RR_\gamma;\Omega_{\III},\two_\gamma)
   \]
\index{Zzzgamma@$\ZZ_\gamma$}\index{Rzzgamma@$\RR_\gamma$}%
\index{Tzzgamma@$\two_\gamma$}%
such that $\eps(p(x))=\eps(x)\rfloor_{\Omega_p}$, for all
$x\in S$ and all\index{pzzroj@$\BB{S}$} $p\in\BB{S}$.

Conversely, for every ordinal $\gamma$, every complete Boolean
space\index{Boolean space}~$\Omega$, decomposed as
$\Omega=\Omega_{\I}\sqcup\Omega_{\II}\sqcup\Omega_{\III}$ with
$\Omega_{\I}$, $\Omega_{\II}$, and $\Omega_{\III}$ clopen, every lower subset
of the space
   \[
\CC(\Omega_{\I},\ZZ_\gamma;\Omega_{\II},\RR_\gamma;\Omega_{\III},\two_\gamma),
   \]
\index{Zzzgamma@$\ZZ_\gamma$}\index{Rzzgamma@$\RR_\gamma$}%
\index{Tzzgamma@$\two_\gamma$}%
endowed with its canonical structure of \pcm, is a continuous dimension scale.\index{continuous dimension scale}
\end{theorem}

\begin{remark}
Of course, as observed earlier, the following isomorphism
\index{Zzzgamma@$\ZZ_\gamma$}\index{Rzzgamma@$\RR_\gamma$}%
\index{Tzzgamma@$\two_\gamma$}%
   \begin{equation}\label{Eq:isoCC}
\CC(\Omega_{\I},\ZZ_\gamma;\Omega_{\II},\RR_\gamma;\Omega_{\III},\two_\gamma)
\cong\CC(\Omega_{\I},\ZZ_\gamma)\times\CC(\Omega_{\II},\RR_\gamma)
\times\CC(\Omega_{\III},\two_\gamma)
   \end{equation}
holds, so we could have formulated part of Theorem~\ref{T:EmbDimInt} by
using the right hand side of \eqref{Eq:isoCC} instead of its left hand side.
However, the formulation of the relation
$\eps(p(x))=\eps(x)\rfloor_{\Omega_p}$ would have then been
more cumbersome.
\end{remark}

\begin{proof}
The first statement (existence of $\eps$) follows from
the construction of~$\eps$ discussed in all
previous results of Section~\ref{S:EmbDI}.

Conversely, by Theorem~\ref{T:C(O,K)DimInt}, every monoid of the form
\index{Zzzgamma@$\ZZ_\gamma$}\index{Rzzgamma@$\RR_\gamma$}%
\index{Tzzgamma@$\two_\gamma$}%
   \[
\CC(\Omega_{\I},\ZZ_\gamma;\Omega_{\II},\RR_\gamma;\Omega_{\III},\two_\gamma)
   \]
is a continuous dimension scale,\index{continuous dimension scale} and, by
Lemma~\ref{L:SegDI}, every lower subset of a continuous dimension scale,\index{continuous dimension scale} viewed as a partial submonoid, is a
continuous dimension scale. Theorem~\ref{T:EmbDimInt} follows.
\end{proof}

Any continuous dimension scale\index{continuous dimension scale} is a \emph{partial} monoid,
which sometimes makes computations cumbersome. However, the following
corollary makes it possible to reduce most problems about continuous dimension scales to \emph{total} monoids.

\begin{corollary}\label{C:EmbDimInt}
Let $S$ be a continuous dimension scale.\index{continuous dimension scale} Then the
universal monoid $\Ref S$ of $S$ is a continuous dimension scale, and $p\mapsto
p|_S$ defines an isomorphism from
$\BB{\Ref S}$\index{pzzroj@$\BB{S}$} onto $\BB{S}$.
\end{corollary}

\begin{proof}
By universality, the map $\eps$ extends to a unique monoid
homomorphism $\tilde{\eps}$ from $\Ref S$ to $\ol{S}$.
Since $S$ is a lower subset of $\Ref S$ and $\eps$ is one-to-one, it follows
from\index{Wehrung, F.} \cite[Lemma~3.9]{WDim} that $\tilde{\eps}$ is
one-to-one. Since $\eps[S]$ is a lower subset of the refinement monoid
\index{refinement monoid} $\ol{S}$, the monoid
$\tilde{\eps}[\Ref S]$, which is equal to the submonoid of $\ol{S}$ generated
by $\eps[S]$, is also a lower subset of $\ol{S}$. In particular, $\Ref S$ is
isomorphic to a lower subset of~$\ol{S}$. The conclusion follows from
Theorem~\ref{T:C(O,K)DimInt} and Lemma~\ref{L:SegDimInt}.
\end{proof}

For the reader's convenience, we restate explicitly the construction of the
map $\eps\colon S\hookrightarrow\ol{S}$ of Theorem~\ref{T:EmbDimInt},
with $\ol{S}=
\CC(\Omega_{\I},\ZZ_\gamma;\Omega_{\II},\RR_\gamma;\Omega_{\III},\two_\gamma)$,
\index{Zzzgamma@$\ZZ_\gamma$}\index{Rzzgamma@$\RR_\gamma$}%
\index{Tzzgamma@$\two_\gamma$}%
for some large enough ordinal $\gamma$.
We first pick a finitary unit $E$ of $S$ (see Definition~\ref{D:FinUn}).
Then we embed $S$ into the universal monoid $\Ref S$ of $S$, see
Proposition~\ref{P:PrmRm}. We observe that $S$ is a lower subset of $\Ref S$
(Proposition~\ref{P:PrmRm}), and that $p\mapsto p|_S$ defines an
isomorphism from $\BB{\Ref S}$\index{pzzroj@$\BB{S}$} onto $\BB S$
(Lemma~\ref{L:SegDimInt}). The definition of the map $\delta\colon\Ref
S_\fin\to\ol{S}$ has the parameter
$E$, and it is given by the formula \eqref{Eq:Defdelta}, that is,
   \[
   \delta(x)(\fa)=\bigvee\Setm{m/n}{(m,n)\in\ZZ^+\times\NN\text{ and }
   \bv{me\leq nx}\in\fa\text{ for all }e\in E},
   \]
for all $x\in\Ref S_\fin$ and all $\fa\in\Omega_\I\cup\Omega_\II$.

The definition of $\mu$ is intrinsic (it does not depend on the finitary
unit $E$), but it requires the ordinal $\gamma$ to be chosen large enough.
It is given by the formula \eqref{Eq:Defmu}, that is,
\index{Tzzgamma@$\two_\gamma$}
   \[
   \mu(x)(\fa)=\bigvee\setm{\alpha\in\two_\gamma}
   {\exists p\in\fa\text{ such that }\scal{p}{\alpha}\text{ is defined and }
   \scal{p}{\alpha}\leq x}.
   \]
for all $x\in S$ and all $\fa\in\Omega$.

Finally, $\eps(x+y)=\delta(x)+\mu(y)$, for all
$(x,y)\in S_\fin\times\DI{S}$. We will call this map $\eps$ the
\emph{canonical embedding from $S$ into $\ol{S}$, relatively to the
finitary unit $E$}.

\section{Uniqueness of the canonical embedding}\label{S:UnCanEmb}
\begin{quote}
\em Standing hypotheses: $S$ is a continuous dimension scale\index{continuous dimension scale}, $S_{\I}$, $S_{\II}$,
$S_{\III}$, $\Omega_{\I}$, $\Omega_{\II}$, $\Omega_{\III}$,
$p_{\I}$, $p_{\II}$, $p_{\III}$ are as in Section~\textup{\ref{S:ProjDF}}.

Let $\gamma$ be an ordinal, large enough for the dimension
function $\mu\colon S\to\CC(\Omega,\two_\gamma)$
\index{Tzzgamma@$\two_\gamma$}introduced in
Section~\textup{\ref{S:DimFctmu}} to be defined.

We put $\ol{S}=\CC
(\Omega_{\I},\ZZ_\gamma;\Omega_{\II},\RR_\gamma;\Omega_{\III},\two_\gamma)$.
\index{Zzzgamma@$\ZZ_\gamma$}\index{Rzzgamma@$\RR_\gamma$}%
\index{Tzzgamma@$\two_\gamma$}%
We fix a finitary unit $E$ of $S$, and we let
$\delta\colon S_\fin\hookrightarrow\ol{S}$ and
$\eps\colon S\hookrightarrow\ol{S}$ be the canonical maps defined from $E$ in
Sections~\textup{\ref{S:ProjDF}} and \textup{\ref{S:EmbDI}}.
\end{quote}

Throughout the present section until Theorem~\ref{T:EmbDimInt}, we let
$\eps'\colon S\hookrightarrow\ol{S}$ be a lower embedding\index{lower
embedding} satisfying the conditions\index{pzzroj@$\BB{S}$}
   \begin{align}
   \eps'(p(x))&=\eps'(x)\rfloor_{\Omega_p},&&
   \text{for all }(x,p)\in S\times\BB{S},\label{Eq:Sheaf}\\
   \eps'(e)&=\eps(e),&&\text{for all }e\in E\cap S_\II.\label{Eq:ResE}
   \end{align}
We shall prove that $\eps'=\eps$.

We first embed $S$ into its universal monoid $\Ref S$. By
Corollary~\ref{C:EmbDimInt}, $\Ref S$ is a continuous dimension scale.\index{continuous dimension scale} Furthermore, the argument of the proof
of Corollary~\ref{C:EmbDimInt} shows that the unique extension of $\eps'$ to
a map from $\Ref S$ to $\ol{S}$ is a lower embedding.\index{lower embedding}
Since $p\mapsto p|_S$ defines an isomorphism from $\BB{\Ref
S}$\index{pzzroj@$\BB{S}$} onto $\BB S$ (Lemma~\ref{L:SegDimInt}),
$\Ref{\eps'}$ satisfies the condition \eqref{Eq:Sheaf}. Of course, it
obviously satisfies \eqref{Eq:ResE}. \emph{Therefore, we may assume, without
loss of generality, that $S$ is a \emph{total} monoid, that is, $S=\Ref S$}.

Next, it follows from Theorem~\ref{T:C(O,K)DimInt} that $\ol{S}$ is a
continuous dimension scale.\index{continuous dimension scale} Furthermore, by
Claim~\ref{Cl:ProjPU} of the proof of Theorem~\ref{T:C(O,K)DimInt}, the
projections of $\ol{S}$ are exactly the maps $f\mapsto f|_{\Omega_p}$,
for\index{pzzroj@$\BB{S}$} $p\in\BB S$, in particular,
$\BB{S}\cong\BB{\ol{S}}$. Thus we shall identify every projection $p$ of
$S$ with the associated projection of $\ol{S}$. Modulo this identification,
the central cover $\cc(f)$ of any $f\in\ol{S}$ is exactly the topological
closure of the set $\vbv{0<f}=\setm{\fa\in\Omega}{f(\fa)>0}$.

\begin{lemma}\label{L:ccinolS}
The equality $\cc(\eps'(x))=\cc(x)$ holds, for all $x\in S$.
\end{lemma}

\begin{proof}
Put $p=\cc(x)$. {}From $x=p(x)$ and \eqref{Eq:Sheaf} it follows that
$\eps'(x)=\eps'(x)\rfloor_{\Omega_p}$, whence $\cc(\eps'(x))\leq p$. Put
$q=p\wedge\neg\cc(\eps'(x))$. Then
$\eps'(q(x))=\eps'(x)\rfloor_{\Omega_q}=0$, thus, since $\eps'$ is an
embedding, $q(x)=0$. However, $q\leq p=\cc(x)$, whence $q=0$, that is,
$p=\cc(\eps'(x))$.
\end{proof}

\subsection{Uniqueness on the directly finite elements}\label{Sub:DirFinUn}
\index{directly finite}
We compute the values of $\eps$ on the elements of $E$. Let $\chi(p)$
\index{czzhi(p)@$\chi(p)$ ($p$ projection)|ii} denote
the characteristic function of $\Omega_p$, for any\index{pzzroj@$\BB{S}$}
$p\in\BB{S}$.

\begin{lemma}\label{L:epsonE}
The equality $\eps(e)=\chi(\cc(e))$ holds, for all $e\in E$.
\end{lemma}

\begin{proof}
Since $e$ is\index{directly finite} directly finite, $\eps(e)=\delta(e)$,
and it is given by
\eqref{Eq:Defdelta}. Put $p=\cc(e)$. {}From $p(e)=e$ and \eqref{Eq:Sheaf} it
follows that $\eps'(e)$ vanishes outside $\Omega_p$. Now let
$\fa\in\Omega_p$. For $e'\in E$ and $(m,n)\in\ZZ^+\times\NN$, the relation
$\bv{me'\leq ne}\in\fa$ always holds for $e'\neq e$ (because then, $p(e')=0$,
thus $\bv{me'=0}\in\fa$), while for $e'=e$, it is equivalent to the existence
of $q\in\fa$ such that $mq(e)\leq nq(e)$. However, for any $q\in\fa$,
$q\wedge\cc(e)$ is nonzero, thus $q(e)$ is nonzero, but it is
\index{directly finite} directly finite, thus $mq(e)\leq nq(e)$ \iff\
$m\leq n$. The conclusion follows immediately.
\end{proof}

\begin{lemma}\label{L:eps'onI}
The equality $\eps'(e)=\eps(e)$ holds, for all $e\in E$.
\end{lemma}

\begin{proof}
The conclusion holds by assumption for $e\in E\cap S_\II$. Now let
$e\in E\cap S_\I$, so $e$ is multiple-free and
$\Omega_{\cc(e)}\subseteq\Omega_\I$. Moreover, it follows from
Lemma~\ref{L:ccinolS} that $\cc(\eps'(e))=\cc(e)$, therefore, since
$\eps'(e)\in\ol{S}$ and $\Omega_{\cc(e)}\subseteq\Omega_\I$, we obtain the
inequality
   \begin{equation}\label{Eq:e'egchie}
   \eps'(e)\geq\chi(\cc(e)).
   \end{equation}
Let $p\in[0,\cc(e)]$ such that $2\chi(p)\leq\eps'(e)$. Since $\eps'$
is a lower embedding,\index{lower embedding} there exists $x\leq e$ such that
$\eps'(x)=\chi(p)$, thus (we recall that $S$ is a total monoid)
$\eps'(2x)=2\chi(p)\leq\eps'(e)$, thus, since $\eps'$ is an
embedding, $2x\leq e$, thus, since $e$ is multiple-free, $x=0$, whence $p=0$.
This holds for all $p\in[0,\cc(e)]$ such that $2\chi(p)\leq\eps'(e)$,
thus, since $\eps'(e)$ vanishes outside $\Omega_{\cc(e)}$, we get
$\eps'(e)\leq\chi(\cc(e))$. Therefore, by \eqref{Eq:e'egchie} and
Lemma~\ref{L:epsonE}, $\eps'(e)=\chi(\cc(e))=\eps(e)$.
\end{proof}

\begin{lemma}\label{L:eps'delta}
The equality $\eps'(x)=\delta(x)$ holds, for all $x\in S_\fin$.
\end{lemma}

\begin{proof}
If the result has been established for $S_\fin$, then it obviously holds
for~$S$. Hence we may assume that $S=S_\fin$, that is, since $S$ is total, $S$
is the positive cone of some Dedekind complete lattice-ordered group
(Lemma~\ref{L:RefSM1}). Now put
$\Omega'=\bigcup_{e\in E}\Omega_{\cc(e)}$, an open subset of $\Omega$. Since
every element of $S$ meets some element of $E$, it follows from
Proposition~\ref{P:Xbotbot} that $\Omega'$ is dense in $\Omega$.

Now let $x\in S$. Since $\delta(x)$ belongs to
$\CC_{\fin}(\Omega_{\I},\ZZ_0;\Omega_{\II},\RR_0)$ (see
Notation~\ref{Not:Cfin}), there exists an open dense subset $\Omega''$ of
$\Omega'$ such that $\delta(x)(\fa)$ is finite, for all $\fa\in\Omega''$.
Since both maps $\eps'(x)$ and $\delta(x)$ are continuous and $\Omega''$ is
dense, in order to conclude the proof, it suffices to establish the equality
$\eps'(x)(\fa)=\delta(x)(\fa)$, for all $\fa\in\Omega''$. Since
$\fa\in\Omega'$, there exists a unique $e\in E$ such that $e\in\fa$, hence
   \[
   \delta(x)(\fa)=\bigvee\setm{m/n}
   {(m,n)\in\ZZ^+\times\NN\text{ and }\bv{me\leq nx}\in\fa}.
   \]
Let $(m,n)\in\ZZ^+\times\NN$ such that $p=\bv{me\leq nx}$ belongs to $\fa$.
Applying $\eps'$ to the inequality $mp(e)\leq np(x)$ and using
\eqref{Eq:Sheaf}, we obtain the inequalities
   \[
   m\eps'(e)\rfloor_{\Omega_p}=m\eps'(p(e))\leq n\eps'(p(x))\leq n\eps'(x),
   \]
thus, by Lemmas~\ref{L:epsonE} and \ref{L:eps'onI},
$m\chi(p\wedge\cc(e))\leq n\eps'(x)$. Evaluate at $\fa$. Since
$p\wedge\cc(e)$ belongs to $\fa$, we obtain that $m/n\leq\eps'(x)(\fa)$. This
holds for all $(m,n)$ such that $\bv{me\leq nx}\in\fa$, whence
   \begin{equation}\label{Eq:delleqeps'}
   \delta(x)(\fa)\leq\eps'(x)(\fa).
   \end{equation}
Now the converse. {}From $\fa\in\Omega''$ it follows that there exists
$h\in\NN$ such that $\bv{he\leq x}\notin\nobreak\fa$. Let $n\in\NN$. There
exists a largest $m\in\ZZ^+$ such that $\bv{me\leq nx}\in\fa$, in fact $m<hn$.
Suppose that $(m+1)/n<\delta(x)(\fa)$. There exists
$(m',n')\in\ZZ^+\times\NN$ such that $(m+1)/n\leq m'/n'$ and $\bv{m'e\leq
n'x}\in\fa$. On the other hand, from $(m+1)n'e\leq m'ne$ it follows that
   \[
   \bv{(m+1)e\leq nx}=\bv{(m+1)n'e\leq nn'x}\geq\bv{m'ne\leq n'nx}
   =\bv{m'e\leq n'x}\in\fa,
   \]
whence $\bv{(m+1)e\leq nx}\in\fa$, a contradiction; so we have proved that
$\delta(x)(\fa)\leq(m+1)/n$.

{}From $\bv{(m+1)e\leq nx}\notin\fa$ and general comparability it follows that
the projection $q=\bv{nx\leq(m+1)e}$ belongs to $\fa$. Thus, applying $\eps'$
to the inequality $nq(x)\leq(m+1)q(e)$ and using \eqref{Eq:Sheaf} together
with Lemmas~\ref{L:epsonE} and \ref{L:eps'onI}, we obtain the inequalities
   \[
n\eps'(x)\rfloor_{\Omega_q}=n\eps'(q(x))\leq(m+1)\eps'(q(e))
\leq(m+1)\eps'(e)=(m+1)\chi(\cc(e)),
   \]
whence, evaluating at $\fa$, $n\eps'(x)(\fa)\leq m+1$, so
   \[
   \eps'(x)(\fa)\leq(m+1)/n\leq\delta(x)(\fa)+1/n.
   \]
This holds for all $n\in\NN$, thus $\eps'(x)(\fa)\leq\delta(x)(\fa)$.
By \eqref{Eq:delleqeps'}, the conclusion follows.
\end{proof}

\subsection{Uniqueness on the purely infinite elements}\label{Sub:PurInfUn}
We need to prove that $\eps'(x)=\mu(x)$, for all $x\in\DI{S}$. We recall that
we have identified the projections of $S$ and those of $\ol{S}$. For
$\alpha\in\Cn$ and\index{pzzroj@$\BB{S}$} $p\in\BB S$, we shall denote by
$\scal{p}{\alpha}_S$ (resp., $\scal{p}{\alpha}_{\ol{S}}$) the value of
$\scal{p}{\alpha}$ in $S$ (resp., in $\ol{S}$), if defined.
By $\alpha\cdot\chi(p)$,\index{czzhi(p)alp@$\alpha\cdot\chi(p)$|ii}
we denote the function defined on
$\Omega$ sending any element of $\Omega_p$ to $\alpha$ and any element of
$\Omega\setminus\Omega_p$ to $0$.

\begin{lemma}\label{L:ScalolS}
The value $\scal{p}{\alpha}_{\ol{S}}$ is defined and equal to
$\alpha\cdot\chi(p)$, for all $\alpha\in\two_\gamma$
\index{Tzzgamma@$\two_\gamma$}and\index{pzzroj@$\BB{S}$} $p\in\BB S$.
\end{lemma}

\begin{proof}
By induction on $\alpha$. The case $\alpha=0$ and the limit step are obvious.
Suppose that $\alpha=\beta^+$, for some $\beta\in\two_\gamma$.
\index{Tzzgamma@$\two_\gamma$}
It is easy to verify that $\beta\cdot\chi(p)\rem\alpha\cdot\chi(p)$ and
$\cc(\alpha\cdot\chi(p))=p$; thus, by the induction hypothesis,
$\scal{p}{\alpha}_{\ol{S}}$ is defined and lies above $\alpha\cdot\chi(p)$.
Now suppose that $\scal{p}{\alpha}_{\ol{S}}<\alpha\cdot\chi(p)$. There exists
$q\in(0,p]$ such that
   \begin{equation}\label{Eq:Paradchiq}
   q(\scal{p}{\alpha}_{\ol{S}})\leq\beta\cdot\chi(q).
   \end{equation}
In particular, from $\cc(\scal{p}{\alpha}_{\ol{S}})=p$ it follows that
$\beta>0$. However, by applying~$q$ to the relation
$\scal{p}{\beta}_{\ol{S}}\rem\scal{p}{\alpha}_{\ol{S}}$ and using
Lemma~\ref{L:ProjTr}(i), we obtain that
$\beta\cdot\chi(q)\rem q(\scal{p}{\alpha}_{\ol{S}})$. Hence, by
\eqref{Eq:Paradchiq} and Lemma~\ref{L:TrLeqTr},
$\beta\cdot\chi(q)\rem\beta\cdot\chi(q)$, whence, since $\beta>0$, we obtain
that $q=0$, a contradiction.
\end{proof}

Since $\eps'$ is a lower embedding,\index{lower embedding} the following
lemma is obvious.

\begin{lemma}\label{L:LowEmbRem}
$a\rem b$ \iff\ $\eps'(a)\rem\eps'(b)$, for all $a$, $b\in S$.
\end{lemma}

\begin{lemma}\label{L:eps'onpalph}
The equality $\eps'(\scal{p}{\alpha}_S)=\alpha\cdot\chi(p)$ holds, for
all\index{pzzroj@$\BB{S}$}
$p\in\BB S$ and all $\alpha\in\Cn$ such that $\scal{p}{\alpha}_S$ is defined.
\end{lemma}

\begin{proof}
By induction on $\alpha$. For $\alpha=0$ it is trivial. Suppose that
$\alpha>0$ and $\scal{p}{\alpha}_S$ is defined. It follows from the induction
hypothesis that $\beta\cdot\chi(p)<\eps'(\scal{p}{\alpha}_S)$, for all
$\beta<\alpha$ in $\Cn$; whence
$\alpha\cdot\chi(p)\leq\eps'(\scal{p}{\alpha}_S)$. Since $\eps'$ is a lower
embedding,\index{lower embedding} there exists $x\leq\scal{p}{\alpha}_S$ in
$\DI{S}$ such that
$\eps'(x)=\alpha\cdot\chi(p)$. The relation
$\eps'(\scal{p}{\beta}_S)\rem\eps'(x)$ holds, for all $\beta<\alpha$, thus, by
Lemma~\ref{L:LowEmbRem}, $\scal{p}{\beta}_S\rem x$. Furthermore, by
Lemma~\ref{L:ccinolS},
   \[
   \cc(x)=\cc(\eps'(x))=\cc(\alpha\cdot\chi(p))=p,
   \]
hence, by the definition of $\scal{p}{\alpha}$, we get that
$\scal{p}{\alpha}_S\leq x$, so, finally, $x=\scal{p}{\alpha}_S$ and
$\eps'(\scal{p}{\alpha}_S)=\alpha\cdot\chi(p)$.
\end{proof}

\begin{lemma}\label{L:eps'=muDI}
The equality $\eps'(x)=\mu(x)$ holds, for all $x\in\DI{S}$.
\end{lemma}

\begin{proof}
We let $\Omega^{(x)}$ denote the open dense subset of $\Omega$ defined in
\eqref{Eq:Omega(x)}. It suffices to prove that the equality
$\eps'(x)(\fa)=\mu(x)(\fa)$ holds, for all $\fa\in\Omega^{(x)}$.

Since $x$ is purely infinite and $\fa\in\Omega^{(x)}$, the value
$\alpha=\mu(x)(\fa)$ is the unique element of $\Cn$ such that
$\scal{p}{\alpha}_S$ is defined and equal to $p(x)$, for some $p\in\fa$.
Therefore, we can compute:
   \begin{align*}
   \eps'(x)(\fa)&=\eps'(p(x))(\fa)&&(\text{by \eqref{Eq:Sheaf}})\\
   &=\eps'(\scal{p}{\alpha}_S)(\fa)\\
   &=(\alpha\cdot\chi(p))(\fa)&&(\text{by Lemma~\ref{L:eps'onpalph}})\\
   &=\alpha\\
   &=\mu(x)(\fa).\tag*{\qed}
   \end{align*}
\renewcommand{\qed}{}
\end{proof}

\subsection{Uniqueness of $\eps$}\label{Sub:Uneps}
By putting together Lemmas~\ref{L:eps'delta} and \ref{L:eps'=muDI}, we obtain
the following.

\begin{corollary}\label{C:Uneps}
The equality $\eps'(x)=\eps(x)$ holds, for all $x\in S$.
\end{corollary}

In order to formulate concisely the corresponding uniqueness result, it is
convenient to extend the usual definition of a continuous dimension scale,\index{continuous dimension scale} as follows. We endow each of the proper
classes $\Zn$, $\Rn$, and $\Cn$
\index{Zzzinfty@$\ZZ_\infty$}\index{Rzzinfty@$\RR_\infty$}%
\index{Tzzinfty@$\two_\infty$}%
introduced in Notation~\ref{Not:kappa+} with its\index{interval topology}
interval topology. The latter consists, for example, of all open
sub\emph{sets} of the corresponding class, the essential fact being that for
a topological space
$\Omega$ (we emphasize that
$\Omega$ is a \emph{set}), the spaces of continuous functions
$\CC(\Omega,\Zn)$, $\CC(\Omega,\Rn)$, and
$\CC(\Omega,\Cn)$ are well-understood (anyway, any map from a set to $\Rn$ is
majorized by some $\aleph_\alpha$). Then we naturally extend
Notation~\ref{Not:OmiKi} to the case where the $K_i$-s may also be
$\Zn$, $\Rn$, or~$\Cn$.

By putting this together with Lemma~\ref{L:epsonE} and
Theorem~\ref{T:EmbDimInt}, we have obtained the following structure theorem
for continuous dimension scales.\index{continuous dimension scale}

\begin{theorem}\label{T:Uneps}
Let $S$ be a continuous dimension scale,\index{continuous dimension scale} let
$\BB{S}$\index{pzzroj@$\BB{S}$} be the complete Boolean
algebra\index{Boolean algebra!complete ---} of projections of~$S$, and let
$\Omega$ be the ultrafilter space of $\BB{S}$, with the decomposition
$\Omega=\Omega_{\I}\sqcup\Omega_{\II}\sqcup\Omega_{\III}$ as given in
Section~\textup{\ref{S:ProjDF}}. Let $E$ be a finitary unit of $S$ \pup{see
Definition~\textup{\ref{D:FinUn}}}.

Then there exists a unique lower embedding\index{lower embedding}
\index{Zzzinfty@$\ZZ_\infty$}\index{Rzzinfty@$\RR_\infty$}%
\index{Tzzinfty@$\two_\infty$}%
   \[
   \eps\colon S\hookrightarrow
   \CC(\Omega_{\I},\Zn;\Omega_{\II},\Rn;\Omega_{\III},\Cn)
   \]
such that $\eps(p(x))=\eps(x)\rfloor_{\Omega_p}$, for all
$x\in S$ and all\index{pzzroj@$\BB{S}$} $p\in\BB{S}$, and $\eps(e)$ takes
its values in $\set{0,1}$, for all $e\in E\cap S_{\II}$. Furthermore, this
embedding satisfies that
$\eps(e)$ takes its values in $\set{0,1}$, for all $e\in E$.
\end{theorem}

\section{Continuous dimension scales which are proper
classes}\label{S:PperClDI}
All the forthcoming section can
easily be formulated in such a standard class theory as the Bernays-G\"odel
system with choice, BGC. An alternative formulation consists of working in
classical set theory ZFC and identifying any statement (with parameters) with
one free variable, say, $\varphi(x)$, with the ``class'' that it
represents, namely, $\setm{x}{\varphi(x)}$; this way, the mention to classes
becomes a mere expendable commodity.

We shall encounter in Chapter~\ref{Ch:ClEsp} situations where it may
appear as artificial to restrict continuous dimension scales to be \emph{sets},
as opposed to \emph{proper classes}.\index{proper classes}
For example, with every right
self-injective regular ring $R$,
\index{ring!Von Neumann regular ---}\index{ring!right self-injective ---}%
we associate the category $\NSIR$\index{NzzSIR@$\NSIR$|ii}
of all nonsingular injective right $R$-modules. With the class
$\NSIR$ is associated a class that meets all attributes of a
continuous dimension scale, except that it is not a set (see Section~\ref{S:RSIReg}).
We shall call such objects \emph{Continuous Dimension Scales} (with
capitals), and we shall define them shortly. We first do this for monoids.

\begin{definition}\label{D:Monoid}
A \emph{Monoid}\index{Monoid (monoid that might be a proper class)|ii}
is a class $M$, endowed with an associative
binary operation $+$ and a zero element $0$. A \emph{\PCM}
\index{Partial Commutative Monoid (possibly a proper class)|ii}
is a class $S$,
endowed with a commutative, associative (in the sense of
Definition~\ref{D:PartCM}) partial binary operation $+$, with a zero
element $0$.
\end{definition}

Hence the definition of a Monoid (resp., \PCM) extends the one
of a monoid (resp., \pcm), by allowing proper classes.

The problem in defining Continuous Dimension Scales is not that easy to
solve. Indeed, we wish our ``Continuous Dimension Scales'' to satisfy a
version of the main embedding theorem, Theorem~\ref{T:EmbDimInt}. More
precisely, we wish every ``Continuous Dimension Scale'' to embed as a
lower subclass into a (proper class) monoid of the form
\index{Zzzinfty@$\ZZ_\infty$}\index{Rzzinfty@$\RR_\infty$}%
\index{Tzzinfty@$\two_\infty$}%
   \begin{equation}\label{Eq:InfSpace}
   \CC(\Omega_{\I},\Zn;\Omega_{\II},\Rn;\Omega_{\III},\Cn),
   \end{equation}
(see Notation~\ref{Not:kappa+}),
for pairwise disjoint complete Boolean
spaces\index{Boolean space} $\Omega_{\I}$, $\Omega_{\II}$,
and~$\Omega_{\III}$. We shall now state the new axioms defining
Continuous Dimension Scales. Of course, our definition is modeled on
Definition~\ref{D:DimInt} and Corollary~\ref{C:AltAx}.

\begin{definition}\label{D:DimIntClass}
A \emph{Continuous Dimension Scale}
\index{Continuous Dimension Scale (continuous dimension scale that might be a
proper class)|ii} is a \PCM~$S$ which satisfies the following axioms.
\begin{itemize}
\item[(M1)] $S$ has refinement, and the algebraic preordering on $S$ is
antisymmetric.

\item[(M2)] Every nonempty subset of~$S$ admits an infimum.

\item[(N1)]
$\forall a,b$, $\exists c,x,y$ such that $a=c+x$, $b=c+y$, and $x\perp y$.

\item[(N2)] $S=a^\bot+a^{\bot\bot}$, for all $a\in S$
(\emph{where $x\in a^\bot$ means, of course, that $x\wedge a=0$, while
$x\in a^{\bot\bot}$ means that $x\in u^\bot$ for any $u\in a^\bot$}).

\item[(N3)] $b\sd a$ exists, for all $a$, $b\in S$ such that $a\leq b$.

\item[(M5)] Every element $a$ of~$S$ can be written $a=x+y$, where $x$ is
directly finite\index{directly finite} and $y$ is purely infinite.

\item[(M6)] Let $a$, $b$ be purely infinite elements of~$S$.
If $a\rem b$, then the set of all purely infinite elements $x$ of~$S$ such
that $a\rem x$ and $x^\bot=b^\bot$ has a least element.

\item[\Mh] The class $(a]=\setm{x\in S}{x\leq a}$ is a set,
for all $a\in S$.\index{Mzzh@\Mh|ii}

\item[\Ml] There exists a dense\index{dense} lower \emph{subset} $U$ of $S$.
\index{Mzzl@\Ml|ii}
(\emph{We will call $U$ a \emph{generating lower subset} of $S$.})
\end{itemize}
\end{definition}

Axiom \Mh\ is there to ensure that the ``infinity'' in \eqref{Eq:InfSpace}
does not exceed the class of all ordinals. Axiom \Ml\ is there to ensure that
the ``base spaces'' $\Omega_{\I}$, $\Omega_{\II}$, $\Omega_{\III}$ in
\eqref{Eq:InfSpace} are sets (as opposed to proper classes).
We emphasize that we require no condition on \emph{subclasses} of $S$,
lest this might pave the way to undesirable set-the\-o\-ret\-i\-cal
paradoxes. In fact, since the axioms defining Continuous Dimension Scales
are requirements on either elements or subsets of $S$, we obtain the
following result.

\begin{proposition}\label{P:DimIntClass}
Let $S$ be a \PCM.\index{Partial Commutative Monoid (possibly a proper
class)} Then $S$ is a Continuous Dimension Scale
\index{Continuous Dimension Scale (continuous dimension scale that might be a
proper class)}
\iff\ every lower subset
of $S$ is a continuous dimension scale
\index{continuous dimension scale} and $S$ satisfies both \Mh\ and \Ml.
\end{proposition}

\begin{proof}
If $S$ is a Continuous Dimension Scale,
\index{Continuous Dimension Scale (continuous dimension scale that might be a
proper class)} then every lower subset of $S$ is a continuous dimension scale\index{continuous dimension scale}: the proof is \emph{mutatis mutandis} the
same as for Lemma~\ref{L:SegDI}.

Conversely, suppose that every lower subset of $S$ is a continuous dimension scale
\index{continuous dimension scale} and
$S$ satisfies both \Mh\ and \Ml. Every subset $X$ of $S$ is contained in a
generating lower subset $\ol{X}$ of $S$: namely, take
   \[
   \ol{X}=\bigcup\setm{(x]}{x\in X}.
   \]
The rest of the proof is similar to the proof of Lemma~\ref{L:DirUn}.
\end{proof}

Observe, in particular, that \emph{every generating lower subset of a
continuous dimension scale is a continuous dimension scale}.\index{continuous dimension scale} We
also obtain the following extension of Theorem~\ref{T:C(O,K)DimInt}.

\begin{corollary}\label{C:DimIntClass}
Let $\Omega$ be a complete Boolean space,\index{Boolean space} written as a
disjoint union $\Omega=\Omega_{\I}\sqcup\Omega_{\II}\sqcup\Omega_{\III}$, for
clopen subsets $\Omega_{\I}$, $\Omega_{\II}$, $\Omega_{\III}$ of $\Omega$.
Then the Monoid
\index{Monoid (monoid that might be a proper class)}
\index{Zzzinfty@$\ZZ_\infty$}\index{Rzzinfty@$\RR_\infty$}%
\index{Tzzinfty@$\two_\infty$}%
   \[
\CC(\Omega_{\I},\Zn;\Omega_{\II},\Rn;\Omega_{\III},\Cn)
   \]
is a Continuous Dimension Scale.
\end{corollary}

Everything is now ready for the proof of our general embedding theorem
for Continuous Dimension Scales.

\begin{theorem}\label{T:GenEmbDI}
Let $S$ be a Continuous Dimension Scale, let $E$ be a finitary unit
of~$S$. Let $U$ be a generating lower subset of $S$ containing $E$,
let~$\Omega$ be the ultrafilter space of\index{pzzroj@$\BB{S}$} $\BB{U}$,
with the decomposition
$\Omega=\Omega_{\I}\sqcup\Omega_{\II}\sqcup\Omega_{\III}$ as given in
Section~\textup{\ref{S:ProjDF}}. Then there exists a unique lower
embedding\index{lower embedding}
\index{Zzzinfty@$\ZZ_\infty$}\index{Rzzinfty@$\RR_\infty$}%
\index{Tzzinfty@$\two_\infty$}%
   \[
   \eps\colon S\hookrightarrow
   \CC(\Omega_{\I},\Zn;\Omega_{\II},\Rn;\Omega_{\III},\Cn)
   \]
such that $\eps(p(x))=\eps(x)\rfloor_{\Omega_p}$, for all
$x\in S$ and all\index{pzzroj@$\BB{S}$} $p\in\BB{S}$, and $\eps(e)$ takes
its values in $\set{0,1}$, for all $e\in E$ \pup{or, which is equivalent,
for all $e\in E\cap U_{\II}$}.
\end{theorem}

\begin{proof}
Let $\mathcal{C}$ be the class of all lower subsets $T$ of $S$ containing
$U$. In particular, for all $T\in\mathcal{C}$, $U$ is a generating lower
subset of $T$, thus, by Lemma~\ref{L:SegDimInt}, $p\mapsto p|_U$ defines an
isomorphism from\index{pzzroj@$\BB{S}$} $\BB{T}$ onto $\BB{U}$. Let
$p\mapsto p^T$ denote its inverse. Therefore, the ultrafilter space
$\Omega^T$ of $\BB{T}$ is homeomorphic to $\Omega$, \emph{via} the map
   \[
   \fa\mapsto\fa|_U=\setm{p|_U}{p\in\fa},\text{ for all }\fa\in\Omega^T.
   \]
Let the projections of $U$ act on $T$, by defining $p(x)=p^T(x)$, for any
$x\in T$ and $p\in\BB{U}$. Hence, by carrying the structure
of $\Omega^T$ to $\Omega$ \emph{via} this isomorphism and then applying
Theorem~\ref{T:Uneps} to the continuous dimension scale\index{continuous dimension scale} $T$
with the finitary unit $E$, we obtain that there exists a unique lower
embedding\index{lower embedding}
\index{Zzzinfty@$\ZZ_\infty$}\index{Rzzinfty@$\RR_\infty$}%
\index{Tzzinfty@$\two_\infty$}%
   \[
   \eps_T\colon T\hookrightarrow\ol{S},
   \]
where $\CC(\Omega_{\I},\Zn;\Omega_{\II},\Rn;\Omega_{\III},\Cn)$,
such that $\eps_T(p(x))=\eps_T(x)\rfloor_{\Omega_p}$, for all $x\in T$ and
all\index{pzzroj@$\BB{S}$}
$p\in\BB{U}$, and $\eps_T(e)$ takes its values in $\set{0,1}$, for any
$e\in E$.

Furthermore, for elements $T_1$, $T_2$ of $\mathcal{C}$ such that
$T_1\subseteq T_2$, the restriction of $\eps_{T_2}$ to $T_1$ satisfies the
requirements of $\eps_{T_1}$. Hence, by the uniqueness statement of
Theorem~\ref{T:Uneps}, $\eps_{T_2}$ extends $\eps_{T_1}$.

Let $\eps$ denote the union of all the maps $\eps_T$, for $T\in\mathcal{C}$.
It follows from Axiom \Ml\ that the union of all the elements of $\mathcal{C}$
is $S$, thus $\eps$ is a map from $S$ to $\ol{S}$. It obviously satisfies the
required conditions. This concludes the ``existence'' part.

If $\eps'\colon S\hookrightarrow\ol{S}$ is another lower
embedding\index{lower embedding} satisfying the conditions of the conclusion
of Theorem~\ref{T:GenEmbDI}, then, by the uniqueness statement of
Theorem~\ref{T:Uneps} applied to $T$, the restriction of $\eps'$ to $T$
equals $\eps_T$, for all $T\in\mathcal{C}$; whence $\eps'=\eps$. This
concludes the ``uniqueness'' part.
\end{proof}

\chapter{Espaliers}\label{Ch:MCLs}

\section{The axioms}\label{S:LattAxioms}

We shall now give the fundamental lattice-theoretical definition underlying
the whole paper, the definition of an \emph{espalier}. This
definition will consist of a list of simple axioms, numbered from (L1) to
(L8). Interspersed between these axioms, we shall also list some very
elementary properties of espaliers. The role of each of these comments will
also be to prepare for the formulation of the axioms that follow them.

\begin{definition}\label{D:MeasChLatt}
An \emph{espalier}\index{espalier|ii} is a structure $(L,\leq,\perp,\sim)$,
where $(L,\leq)$ is a partially ordered set, $\perp$ is a binary relation on
$L$, and $\sim$ is an equivalence relation on $L$, subject to the following
axioms:
\begin{itemize}
\item[(L1)] Every nonempty subset of $L$ has an infimum. Equivalently, every
majorized subset of $L$ has a supremum.\index{Lzzoneeight@(L1--8)|ii}

\begin{quote}
\em In particular, $L$ has a smallest element, that we shall denote by~$0$.
For $a$, $b\in L$, the meet $a\wedge b$ is always defined,
while the join $a\vee b$ is defined \iff\ the pair $\set{a,b}$ is majorized.
\end{quote}

\item[(L2)] For all $a$, $b$, $c\in L$, the following statements hold:
\begin{enumerate}
\item $a\perp 0$.

\item if $a\perp b$, then $b\perp a$.

\item if $a\leq b$ and $b\perp c$, then $a\perp c$.

\item if $\set{a,b,c}$ is majorized, $a\perp b$, and $(a\vee b)\perp c$,
then $a\perp(b\vee c)$.

\item if $a\perp a$, then $a=0$.

\end{enumerate}

\begin{quote}
\em We can then define in $L$ a partial binary operation $\oplus$, by
putting $c=a\oplus b$ \iff\ $c=a\vee b$ and $a\perp b$. So,
\textup{(i)--(iv)} above means exactly that $(L,\oplus,0)$ is a \pcm.

We say that a family $\famm{a_i}{i\in I}$ of elements of $L$
is \emph{orthogonal}, if it is majorized and
$\oplus_{i\in J}a_i$ is defined for every finite subset $J$ of $I$. We then
define $\oplus_{i\in I}a_i=\bigvee_{i\in I}a_i$.
\end{quote}

\smallskip
\item[(L3)] For all $a$, $b\in L$, if $a\leq b$, then there exists $x\in L$
such that $a\oplus x=b$.

\begin{quote}
\em Since $a\oplus x=a\vee x$, the converse of Axiom~\textup{(L3)} is, of
course, trivial.
\end{quote}

\item[(L4)] Let $a\in L$, let $\famm{b_i}{i\in I}$ be an orthogonal family of
elements of $L$. If $a\perp(\oplus_{i\in J}b_i)$, for all finite
$J\subseteq I$, then $a\perp(\oplus_{i\in I}b_i)$.

\item[(L5)] $x\sim 0$ implies that $x=0$, for all $x\in L$.

\item[(L6)] The relation $\sim$ is \emph{unrestrictedly refining},
\index{unrestrictedly refining relation|ii} that is,
for every $a\in L$ and every orthogonal family $\famm{b_i}{i\in I}$ of
elements of $L$, if $a\sim\oplus_{i\in I}b_i$, then there exists a
decomposition $a=\oplus_{i\in I}a_i$ such that $a_i\sim b_i$ for all $i\in I$.

\item[(L7)] The relation $\sim$ is \emph{unrestrictedly additive},
\index{unrestrictedly additive!--- relation|ii} that is,
for all orthogonal families $\famm{a_i}{i\in I}$ and $\famm{b_i}{i\in I}$ of
elements of $L$, if $a_i\sim b_i$ for all $i\in I$, then
$\oplus_{i\in I}a_i\sim\oplus_{i\in I}b_i$.

\item[(L8)] (the \emph{parallelogram rule})\index{parallelogram rule|ii}
For all $a$, $b$, $x$, $y\in L$ such that $a\vee b$ is defined,
   \[
   \bigl((a\wedge b)\oplus x=a\text{ and }b\oplus y=a\vee b\bigr)
   \Longrightarrow x\sim y.
   \]
\end{itemize}
An espalier\index{espalier!bounded ---|ii} is \emph{bounded}, if it has a
largest element.
\end{definition}

The following result makes it possible to create new espaliers from old ones.
We leave the straightforward proof to the reader.

\begin{proposition}\label{P:LSPEsp}\hfill
\begin{enumerate}
\item For any espalier\index{espalier} $(L,\leq,\perp,\sim)$, any lower
subset $K$ of $L$, endowed with the restrictions of $\leq$, $\perp$, and
$\sim$, is an espalier.

\item Let $(L_i,\leq_i,\perp_i,\sim_i)_{i\in I}$ be a family of espaliers.
Then the product $L=\prod_{i\in I}L_i$, endowed with the componentwise $\leq$,
$\perp$, $\sim$, is an espalier.
\end{enumerate}
\end{proposition}

In the context of Proposition~\ref{P:LSPEsp}(i), we shall say that $K$ is a
\emph{lower subespalier}\index{espalier!lower sub ---|ii} of $L$. In the
context of Proposition~\ref{P:LSPEsp}(ii), we shall say that $L$ is the
\emph{direct product} of the family $(L_i)_{i\in I}$ of espaliers. If $K$ and
$L$ are espaliers, we shall say that a map $\varphi\colon K\to L$
is a \emph{lower embedding}\index{lower embedding} of espaliers, if it is an
isomorphism from $K$ onto a lower subespalier of $L$. The verification of
the following lemma is straightforward.

\begin{lemma}\label{L:LowEmbEsp}
Let $K$ and $L$ be espaliers,\index{espalier} let $\varphi\colon K\to L$ be
a map. Then $\varphi$ is a lower embedding\index{lower embedding} \iff\ the
following conditions hold:
\begin{enumerate}
\item the range of $\varphi$ is a lower subset of $L$.

\item $x\leq_Ky$ \iff\ $\varphi(x)\leq_L\varphi(y)$, for all $x$, $y\in K$.

\item $x\perp_Ky$ \iff\ $\varphi(x)\perp_L\varphi(y)$, for all $x$, $y\in K$.

\item $x\sim_Ky$ \iff\ $\varphi(x)\sim_L\varphi(y)$, for all $x$, $y\in K$.
\end{enumerate}
\end{lemma}

\begin{quote}
\em For the remainder of Section~\textup{\ref{S:LattAxioms}}, we shall
fix an espalier\index{espalier}\linebreak $(L,\leq,\perp,\sim)$.
\end{quote}

We start up with elementary properties of orthogonal families.

\begin{lemma}\label{L:ModPair}
For all $a$, $b\in L$ such that $a\perp b$, the following holds:
\begin{enumerate}
\item $a\wedge b=0$.

\item $x=(x\oplus b)\wedge a$, for all $x\leq a$. \pup{Here, this means that
$(a,b)$ is a \emph{modular pair}, see \cite{GLT2}\index{Gr\"atzer, G.} for
the general definition of those.}

\end{enumerate}

\end{lemma}

\begin{proof}
As in \cite[Lemma~1.1]{SMae55}.\index{Maeda, S.}
\end{proof}

\begin{corollary}\label{C:RelComp}
Let $a$, $b$, $c\in L$.
\begin{enumerate}
\item If $(a,b,c)$ is orthogonal, then $(a\oplus c)\wedge(b\oplus c)=c$.

\item If $a\leq b\leq c$, then there exists $x\in L$ such that $b\wedge x=a$
and $b\vee x=c$. \pup{That is, every closed interval of $L$ is a
\emph{relatively complemented} lattice.}
\index{lattice!relatively complemented ---}
\end{enumerate}

\end{corollary}

\begin{proof}
(i) Apply Lemma~\ref{L:ModPair} to the pair $(b\oplus c,a)$ and to
$c\leq b\oplus c$.

(ii) By Axiom~(L3), there are $u$, $v$ such that $a\oplus u=b$ and
$b\oplus v=c$. So $c=a\oplus u\oplus v$. By (i), $x=a\oplus v$ satisfies the
required conditions.
\end{proof}

By using Axiom~(L4), it is easy to prove the following result (see also
Theorem~1.2 of \cite{SMae55}).\index{Maeda, S.}

\begin{lemma}\label{L:LattAssoc}
Let $I$ and $J$ be sets, let
$\pi\colon I\twoheadrightarrow J$ be a surjective map, let
$\famm{a_i}{i\in I}$ be a family of elements of $L$, and let $a\in S$.
Then the following are equivalent:

\begin{enumerate}
\item $a=\oplus_{i\in I}a_i$.

\item For all $j\in J$, the family $\famm{a_i}{i\in\pi^{-1}\set{j}}$ is
orthogonal, and, if we denote its join by~$b_j$, then $a=\oplus_{j\in J}b_j$.
\end{enumerate}
\end{lemma}

Another useful elementary orthogonality property of
espaliers\index{espalier} is the following.

\begin{lemma}\label{L:aperpJJX}
Let $a\in L$, let $X$ be a majorized subset of $L$.
If $a\perp\bigvee Y$ for all finite $Y\subseteq X$, then $a\perp\bigvee X$.
\end{lemma}

\begin{proof}
For finite $X$, this is trivial. Now suppose that $X$ is infinite. Write
$X=\setm{b_\xi}{\xi<\kappa}$, where $\kappa$ is the cardinality of $X$, and
put $b=\bigvee X$. We argue by induction on $\kappa$. Put
$\oll{b}_\xi=\bigvee_{\eta<\xi}b_\eta$, for all $\xi<\kappa$.
Observe that $\oll{b}_0=0$. For all
$\xi<\kappa$, there exists, by Axiom~(L3), $c_\xi\in L$ such that
$\oll{b}_\xi\oplus c_\xi=\oll{b}_{\xi+1}$. It follows easily that
$\oll{b}_\xi=\oplus_{\eta<\xi}c_\eta$ for all $\xi<\kappa$, while
$b=\oplus_{\eta<\kappa}c_\eta$. Furthermore, it follows from the induction
hypothesis that $a\perp\oll{b}_\xi$ for all $\xi<\kappa$, whence
$a\perp\oplus_{\eta\in J}c_\eta$ for every finite subset $J$ of $\kappa$. By
Axiom~(L4), it follows that $a\perp\oplus_{\eta<\kappa}c_\eta$, that is,
$a\perp b$.
\end{proof}

\begin{notation}
For $a$, $b\in L$, let $a\lesssim b$ hold,
\index{abzzlesssim@$a\lesssim b$ (in espaliers)|ii}
if $a\sim x$ for some $x\leq b$.
\end{notation}

\begin{lemma}\label{L:DotJoinSim}
Let $a$, $b$, $c\in L$.\index{perspectivity@perspectivity ($\sim$)}
\begin{enumerate}
\item If $a\vee c=b\vee c$ and $a\wedge c=b\wedge c$, then $a\sim b$.
\emph{(That is, if $a$ and $b$ are perspective, then $a\sim b$.)}

\item If $a\oplus c=b\oplus c$, then $a\sim b$.

\item If $a\vee c\leq b\vee c$ and $a\wedge c\leq b\wedge c$, then
$a\lesssim b$.
\end{enumerate}

\end{lemma}

\begin{proof}
(i) We put $u=a\wedge c=b\wedge c$, $v=a\vee c=b\vee c$.
Let $x$, $y$, $z\in L$ such that
   \[
   u\oplus x=a;\qquad u\oplus y=b;\qquad c\oplus z=v.
   \]
By the parallelogram rule, $x\sim z$ and $y\sim z$. Thus $x\sim y$, so
   \[
   a=u\oplus x\sim u\oplus y=b.
   \]

(ii) This is, by Lemma~\ref{L:ModPair}(i), a particular case of (i).

(iii) By using Axiom~(L3), there are $u$, $v$, $x$, $y$, $z\in L$ such that
   \begin{gather*}
   a=(a\wedge c)\oplus u,\quad a\vee c=c\oplus v,\\
   b\wedge c=(a\wedge c)\oplus x,\quad
   b=(b\wedge c)\oplus y,\quad b\vee c=(a\vee c)\oplus z.
   \end{gather*}
It follows that $b\vee c=c\oplus v\oplus z$, thus, by Axiom~(L8),
$y\sim v\oplus z$, thus, by Axiom~(L6), there are $v'\sim v$ and $z'\sim z$
such that $y=v'\oplus z'$. Moreover, it also follows from Axiom~(L8) that
$u\sim v$, thus $u\sim v'$. {}From $v'\leq y$, $a\wedge c\leq b\wedge c$,
and $y\perp(b\wedge c)$ it follows (by Axiom~(L2)) that $(a\wedge c)\perp v'$,
thus, since $a\wedge c$, $v'\leq b$, $(a\wedge c)\oplus v'$ is defined and
below $b$. Therefore, by using Axiom~(L7),
$a=(a\wedge c)\oplus u\sim(a\wedge c)\oplus v'\leq b$.
\end{proof}

We observe that Lemma~\ref{L:DotJoinSim}(iii) is stronger than Axiom~(iii) in
Definition~2.2 of~\cite{Fill65}.\index{Fillmore, P.\,A.}

\begin{quote}
\em {}From now on, we denote by $\DD(a)$, or $\DD_L(a)$
\index{dzzelta@$\DD(a)$|ii} if there is any
ambiguity on $L$, the $\sim$-equivalence class of $a$, for every $a\in L$.
Furthermore, we denote by $S$ the range of $\DD$,
and we call it the
\emph{dimension range}\index{dimension range|ii} of $L$, in notation,
$S=\Drng L$.\index{Dzzrng@$\Drng L$}
\end{quote}

We endow $S$ with the partial binary operation $+$ defined
by
   \begin{multline*}
   \gamma=\alpha+\beta,\text{ if there are }a,\,b,\,c\in L\text{ such that }\\
   \alpha=\DD(a),\ \beta=\DD(b),\ \gamma=\DD(c),\text{ and }c=a\oplus b,
   \end{multline*}
for all $\alpha$, $\beta$, $\gamma\in L$. The fact that $+$ is indeed
well-defined follows from the finite case of Axiom~(L7), namely, if
$a\oplus b$ and $a'\oplus b'$ are defined and $a\sim a'$ and $b\sim b'$, then
$a\oplus b\sim a'\oplus b'$.

We denote $\DD(0)$ by $0$.

\begin{proposition}\label{P:Spcm}
$(S,+,0)$ is a \pcm.
\end{proposition}

\begin{proof}
Only the verification of the associativity of $+$ is not completely
trivial. Given $a$, $b$, $c\in L$ such that $(\DD(a)+\DD(b))+\DD(c)$ is
defined, there exist $x$, $y\in L$ such that $\DD(a)+\DD(b)=\DD(x)$ and
$\DD(x)+\DD(c)=\DD(y)$. Then $x=a'\oplus b'$ for some $a'\sim a$ and
$b'\sim b$ in $L$, while $y=x'\oplus c'$ for some $x'\sim x$ and $c'\sim c$
in $L$. By the finite case of Axiom~(L6), $x'=a''\oplus b''$ for some
$a''\sim a'$ and $b''\sim b'$ in $L$. Axiom~(L2) then implies that
$y=a''\oplus(b''\oplus c')$, and therefore, $\DD(a)+(\DD(b)+\DD(c))$
is defined, and equal to
   \[
   \DD(a)+\DD(b''\oplus c'')=\DD(y)=(\DD(a)+\DD(b))+\DD(c).
   \]
Since $+$ is obviously commutative, this implies that $+$ is associative as
well.
\end{proof}

\begin{quote}
\em The dimension range $S$ will always be endowed
with its \emph{algebraic} preordering $\leq$, see
Definition~\textup{\ref{D:AlgPr}}. Hence
$\DD(a)\leq\DD(b)$ \iff\ $a\lesssim b$, for all $a$, $b\in L$.
\end{quote}

\begin{lemma}\label{L:DelVmeas}
The following assertions hold:
\begin{enumerate}
\item Let $\alpha$, $\beta\in S$ and let $c\in L$ such that
$\DD(c)=\alpha+\beta$. Then there are $a$, $b\in L$ such that $c=a\oplus b$
while $\DD(a)=\alpha$ and $\DD(b)=\beta$.

\item Let $a\in L$ and let $\xi\in S$. If $\xi\leq\DD(a)$, then there exists
$x\leq a$ in $L$ such that $\DD(x)=\xi$.
\end{enumerate}
\end{lemma}

\begin{proof}
(i) By the definition of $\alpha+\beta$, there are orthogonal $u$, $v\in L$
such that $\DD(c)=\DD(u\oplus v)$ while $\DD(u)=\alpha$ and $\DD(v)=\beta$.
So $c\sim u\oplus v$, thus, by the finite case of Axiom~(L6), there are $a$,
$b\in L$ such that $a\sim u$, $b\sim v$, and $c=a\oplus b$. Observe that
$\DD(a)=\DD(u)=\alpha$ and $\DD(b)=\DD(v)=\beta$.

(ii) By the definition of the algebraic preordering of $S$, there exists
$\eta\in S$ such that $\DD(a)=\xi+\eta$. The conclusion follows from (i).
\end{proof}

Proposition~\ref{P:Spcm} and Lemma~\ref{L:DelVmeas} are particular cases of
more general results in Chapter~4\index{Wehrung, F.} of~\cite{WDim}.

\begin{proposition}\label{P:SchBern}
The algebraic preordering on $S$ is antisymmetric. That is, if $a\lesssim b$
and $b\lesssim a$, then $a\sim b$, for all $a$, $b\in L$.
\end{proposition}

\begin{proof}
Similar to the proof of Theorem~41 of \cite{Kapl68},\index{Kaplansky, I.}
see also Lemma~6.1.3 of\index{Murray, F.\,J.}\index{von Neumann, J.}
\cite{MuNe36}. We give the details here for convenience.

Put $a_0=a$ and $b_0=b$.
By induction on $n<\omega$, we construct $a_{n+1}$, $b_{n+1}$, $x_n$,
$y_n\in L$ such that the following relations hold:
   \begin{gather}
   a_n=b_{n+1}\oplus x_n,\quad b_n=a_{n+1}\oplus y_n,\label{Eq:anbn+1}\\
   a_n\sim a_{n+1},\quad b_n\sim b_{n+1}.\label{Eq:ansiman+1}
   \end{gather}
Since $a_0\lesssim b_0$ and $b_0\lesssim a_0$, this is easy to satisfy for
$n=0$. Now the induction step. Since $b_{n+1}\sim b_n=a_{n+1}\oplus y_n$,
there are $a_{n+2}$ and $y_{n+1}$ in $L$ such that
   \[
   b_{n+1}=a_{n+2}\oplus y_{n+1},
   \quad a_{n+2}\sim a_{n+1},\quad y_{n+1}\sim y_n.
   \]
Similarly, we obtain elements $b_{n+2}$, $x_{n+1}$ in $L$ such that
   \[
   a_{n+1}=b_{n+2}\oplus x_{n+1},
   \quad b_{n+2}\sim b_{n+1},\quad x_{n+1}\sim x_n.
   \]
This completes the induction step for \eqref{Eq:anbn+1} and
\eqref{Eq:ansiman+1}. Observe that all the $x_n$ (resp., all the $y_n$) are
mutually $\sim$-equivalent. But then, the family consisting of all the
$x_{2n+1}$ and $y_{2n+2}$, for $n<\omega$, is orthogonal and majorized by
$a_1$. Thus there exists $c\in L$ such that
   \begin{equation}\label{Eq:a1crem}
   a_1=c\oplus(\oplus_{n<\omega}x_{2n+1})\oplus(\oplus_{n<\omega}y_{2n+2}).
   \end{equation}
Since $b=a_1\oplus y_0$, we obtain the equality
   \begin{equation}\label{Eq:bcrem}
   b=c\oplus(\oplus_{n<\omega}x_{2n+1})\oplus(\oplus_{n<\omega}y_{2n}).
   \end{equation}
Since $y_{2n}\sim y_{2n+2}$ for all $n$, it follows from
\eqref{Eq:a1crem}, \eqref{Eq:bcrem}, and Axiom~(L7) that $a_1\sim b$. But
$a\sim a_1$; whence $a\sim b$.
\end{proof}

Of course, it follows that $S$ is \emph{conical}. However, a direct proof of
the conicality of~$S$ is immediate from Axiom~(L5).

\begin{proposition}\label{P:Ref}
The \pcm\ $(S,+,0)$ satisfies the refinement property
\pup{see Definition~\textup{\ref{D:RefPpty}}}.
\end{proposition}

\begin{proof}
Let $\alpha$, $\alpha'$, $\beta$, $\beta'\in S$ such that
$\alpha+\alpha'=\beta+\beta'$. By the definition of the
addition of $S$ and the finite version of Axiom~(L6), there are $a$,
$a'$, $b$, $b'\in L$ such that $\DD(a)=\alpha$, $\DD(a')=\alpha'$,
$\DD(b)=\beta$, $\DD(b')=\beta'$, and $a\oplus a'=b\oplus b'$.
We put $e=a\oplus a'=b\oplus b'$. We observe that $a\vee b$ is defined.
There are elements $u$, $u'$, $v$, $v'$, $w$ in $L$ such that
   \begin{align*}
   (a\wedge b)\oplus u&=a; & b\oplus u'=a\vee b;\\
   (a\wedge b)\oplus v&=b; & a\oplus v'=a\vee b;\\
   &\qquad(a\vee b)\oplus w=e.
   \end{align*}
By Axiom~(L8), $u\sim u'$ and $v\sim v'$. Furthermore,
   \[
   b\oplus b'=e=(a\vee b)\oplus w=b\oplus u'\oplus w,
   \]
thus, by Lemma~\ref{L:DotJoinSim}(ii), $b'\sim u'\oplus w$. Hence
$b'=\oll{u}\oplus w_2$, for some $\oll{u}\sim u'$ (thus
$\oll{u}\sim u$) and some $w_2\sim w$. Similarly,
$a'=\oll{v}\oplus w_1$, for some
$\oll{v}\sim v$ and some $w_1\sim w$. Hence we have obtained the
following refinement matrix:
   \[
   \begin{tabular}{|c|c|c|}
   \cline{2-3}
   \multicolumn{1}{l|}{} & $\beta$ & $\beta'$\tvi\\
   \hline
   $\alpha$ & $\DD(a\wedge b)$ & $\DD(u)=\DD(\oll{u})$\tvi\\
   \hline
   $\alpha'$ & $\DD(\oll{v})=\DD(v)$ & $\DD(w_1)=\DD(w_2)$\tvi\\
   \hline
   \end{tabular}
   \]
This completes the proof.
\end{proof}

{}From Propositions \ref{P:Ref} and \ref{P:SchBern}, we deduce immediately
the following.

\begin{corollary}\label{C:ShasM0}
The \pcm\ $S$ satisfies Axiom~\textup{(M1)}.
\end{corollary}

In particular, by Proposition~\ref{P:ProjBa}, the set
$\BB{S}$\index{pzzroj@$\BB{S}$} of all projections of~$S$, endowed with the
ordering $\leq$ given by $p\leq q$
\iff\ $pS\subseteq qS$ (see Lemma~\ref{L:pLeqOrtq}), is
a\index{Boolean algebra} Boolean algebra. We shall also refer to the
projections of $S$ as the \emph{projections of $L$}. They operate also on
$L$ in a projection-like manner (up to equivalence)---see Lemmas
\ref{L:BaspcdotL} and \ref{L:ppbotonL} below.

For the remainder of Section~\ref{S:LattAxioms}, we shall analyze in some
detail the pairs $(a,b)$ of elements of $L$ such that $\DD(a)\perp\DD(b)$.
Observe that the conditions $a\perp b$ (in $L$) and $\DD(a)\perp\DD(b)$
(in $S$) are \emph{a priori} unrelated. For $a$, $b\in L$,
$\DD(a)\perp\DD(b)$ \iff\ the only element $x\in L$ such that
$x\lesssim a,b$ is $x=0$.

\begin{corollary}\label{C:Ref}
Let $a$, $b$, $c\in L$ such that $\DD(a)\perp\DD(c)$ and
$\DD(b)\perp\DD(c)$. If $\set{a,b}$ is majorized, then
$\DD(a\vee b)\perp\DD(c)$.
\end{corollary}

\begin{proof}
Let $b'\in L$ such that $a\oplus b'=a\vee b$.
In particular, $\DD(a\vee b)=\DD(a)+\DD(b')$.
It follows from Lemma~\ref{L:DotJoinSim}(iii) that $b'\lesssim b$, thus, by
assumption, $\DD(b')\perp\DD(c)$. Therefore, by Proposition~\ref{P:Ref} and
Lemma~\ref{L:XbotId}(i), $\DD(a)+\DD(b')\perp\DD(c)$, that is,
$\DD(a\vee b)\perp\DD(c)$.
\end{proof}

\begin{notation}
For $a$, $b$, $c\in L$, we define $c=a\boxplus b$
\index{cabzzabplusb@$c=a\boxplus b$|ii}
to mean that $c=a\vee b$ and $\DD(a)\perp\DD(b)$.
\end{notation}

Note, in particular, that $c=a\boxplus b$ implies that $a\wedge b=0$. Much
more is true, see Proposition~\ref{P:aboxplusb}.

\begin{lemma}\label{L:aboxbplusb}
Let $a$, $a'$, $b$, $c\in L$. If $c=a\boxplus b=a'\oplus b$ and if
$a'\leq a$, then $a'=a$, so $a\perp b$.
\end{lemma}

\begin{proof}
Let $b'\in L$ such that $c=a\oplus b'$. Since $c=a\vee b$ with $a\wedge
b=0$, it follows from Lemma~\ref{L:DotJoinSim}(i) that $b\sim b'$. Let $v$,
$v'\in L$ such that
   \[
   b=(b\wedge b')\oplus v\quad\text{and}\quad b'=(b\wedge b')\oplus v'.
   \]
Let $u\in L$ such that $a'\oplus u=a$. Then
$c=a'\oplus b=a'\oplus(b\wedge b')\oplus v$, so that
$c=a\oplus b'=a'\oplus(b\wedge b')\oplus u\oplus v'$. Therefore, by
Lemma~\ref{L:DotJoinSim}(ii), $u\oplus v'\sim v$, thus $u\lesssim v$. So
$u\leq a$ and $u\lesssim v\leq b$, hence, by assumption, $u=0$. It follows
that $a=a'$. In particular, $a\perp b$.
\end{proof}

We obtain the following important tool, Proposition~\ref{P:aboxplusb}. It is
an analogue of Axiom~(D) in \cite{Loom55}\index{Loomis, L.\,H.} and of
Axiom~$(2,\eps)$ in\index{Maeda, S.}
\cite{SMae55}. It also holds in the ``cardinal lattices''
\index{lattice!cardinal ---} considered in\index{Fillmore, P.\,A.}
\cite{Fill65}, as Lemma~2.7 of \cite{Fill65} shows. However, the proof of
Lemma~2.7 of \cite{Fill65} cannot be applied here, because there is no
``orthocomplement'' in our axiom system for espaliers.\index{espalier}

\begin{proposition}\label{P:aboxplusb}
Let $a$, $b\in L$. If $\DD(a)\perp\DD(b)$ and $\set{a,b}$ is majorized, then
$a\perp b$; so $a\boxplus b=a\oplus b$.
\end{proposition}

\begin{proof}
By Axiom~(L1), $c=a\vee b$ exists.
Let $a'\in L$ such that $c=a'\oplus b$. By Lemma~\ref{L:DotJoinSim}(i),
$a\sim a'$. Since $a$, $a'\leq c$, there exists $a_1=a\vee a'$ in $L$, but
$\DD(a)=\DD(a')\perp\DD(b)$, so, by Corollary~\ref{C:Ref},
$\DD(a_1)\perp\DD(b)$. So, $c=a_1\boxplus b=a'\oplus b$, with $a'\leq a_1$.
Hence, by Lemma~\ref{L:aboxbplusb}, $a_1\perp b$. Since $a\leq a_1$, it
follows from Axiom~(L2)(iii) that $a\perp b$.
\end{proof}

We now extend Proposition~\ref{P:aboxplusb} to arbitrary families of
elements of $L$.

\begin{definition}\label{D:StrongOrth}
A family $\famm{a_i}{i\in I}$ is \emph{strongly orthogonal},
\index{strongly orthogonal (family)|ii} if
it is majorized and $\DD(a_i)\perp\DD(a_j)$, for all $i\neq j$ in $I$.
\end{definition}

\begin{corollary}\label{C:aiStrOrth}
Every strongly orthogonal family of elements of $L$ is orthogonal.
\end{corollary}

\begin{proof}
It suffices to prove the result for $I$ finite. We argue by induction on
the cardinality of $I$. Pick $i\in I$. By the induction hypothesis,
$\famm{a_j}{j\neq i}$ is orthogonal, and, by Corollary~\ref{C:Ref},
$\DD(a_i)\perp\DD(\oplus_{j\neq i}a_j)$. Hence, by
Proposition~\ref{P:aboxplusb}, $a_i\perp\oplus_{j\neq i}a_j$ in $L$, that is,
$\famm{a_j}{j\in I}$ is orthogonal.
\end{proof}

\begin{corollary}\label{C:aboxplusb}
Let $a$, $b$, $x$, $y$, $c\in L$. If $c=a\boxplus b=x\oplus y$ with
$\DD(x)\perp\DD(b)$ and $\DD(y)\perp\DD(a)$, then $x=a$ and $y=b$.
\end{corollary}

\begin{proof}
Since $a\vee x$ is defined (and $a\vee x\leq c$) and $\DD(a)$,
$\DD(x)\perp\DD(b)$, it follows from Corollary~\ref{C:Ref} that
$\DD(a\vee x)\perp\DD(b)$. So, $c=(a\vee x)\boxplus b$.
Since $c=a\oplus b$ (by Proposition~\ref{P:aboxplusb}) and $a\leq a\vee x$,
it follows from Lemma~\ref{L:aboxbplusb} that $a=a\vee x$, that is, $x\leq a$.
Similarly, $y\leq b$. Let $x'$, $y'\in L$ such that $a=x\oplus x'$ and
$b=y\oplus y'$. Then
   \[
   c=a\oplus b=x\oplus y\oplus x'\oplus y'=c\oplus x'\oplus y',
   \]
whence $x'=y'=0$. Therefore, $x=a$ and $y=b$.
\end{proof}

\begin{corollary}\label{C:DelPerp}
Let $a\in L$, let $X$ be a majorized subset of $L$.
If $\DD(a)\perp\DD(x)$ for all $x\in X$, then $\DD(a)\perp\DD(\bigvee X)$.
\end{corollary}

\begin{proof}
Put $b=\bigvee X$. It follows from Corollary~\ref{C:Ref} that
   \begin{equation}\label{Eq:DDaperpJJY}
   \DD(a)\perp\DD\left(\bigvee Y\right),\quad\text{for all finite }Y\subseteq X.
   \end{equation}
Let $t\in(0,b]$. Then $t\not\perp b$, thus, by Lemma~\ref{L:aperpJJX},
$t\not\perp\bigvee Y$ for some finite $Y\subseteq X$. But $\set{t,\bigvee Y}$
is majorized (by $b$), thus, by Proposition~\ref{P:aboxplusb},
$\DD(t)\not\perp\DD(\bigvee Y)$. Therefore, by \eqref{Eq:DDaperpJJY},
$\DD(t)\nleq\DD(a)$, that is, $t\not\lesssim a$. So we have proved that
$\DD(a)\perp\DD(b)$.
\end{proof}

\begin{corollary}\label{C:ProjOnL}
Let $a\in L$ and let \index{pzzroj@$\BB{S}$}$p\in\BB{S}$. Then there exists
a largest element $u$ of
$[0,a]$ such that $\DD(u)\in pS$. Furthermore, $\DD(u)=p(\DD(a))$.
\end{corollary}

\begin{proof}
By Lemma~\ref{L:OplusId}, $pS=(pS)^{\bot\bot}$. Hence, by
Corollary~\ref{C:DelPerp}, the supremum $u$ of the set $X$ of all elements
$x$ of $[0,a]$ such that $\DD(x)\in pS$ belongs to~$X$.

{}From $u\leq a$ it follows that $\DD(u)\leq\DD(a)$, thus, since
$\DD(u)\in pS$,
$\DD(u)=p(\DD(u))\leq p(\DD(a))$. Conversely, $p(\DD(a))\leq\DD(a)$, thus,
by Lemma~\ref{L:DelVmeas}(ii),
there exists $v\leq a$ such that $p(\DD(a))=\DD(v)$. But $\DD(v)\in pS$, so
$v\leq u$, and thus $p(\DD(a))=\DD(v)\leq\DD(u)$. Finally, by
Proposition~\ref{P:SchBern}, $p(\DD(a))=\DD(u)$.
\end{proof}

We shall denote by $p\cdot a$\index{pzzcdota@$p\cdot a$|ii}
the element $u$ of Corollary~\ref{C:ProjOnL},
and we shall repeatedly use the properties $p\cdot a\leq a$ and
$\DD(p\cdot a)=p(\DD(a))$, for all $a\in L$ and all\index{pzzroj@$\BB{S}$}
$p\in\BB{S}$.

We shall also put $p\cdot L=\setm{p\cdot x}{x\in L}$. Note that $p\cdot L$
is a lower subset of~$L$: indeed, if $a\in L$ and $b\leq p\cdot a$, then
$\DD(b)\in pS$, so $p\cdot b=b$.

We gather up various elementary properties of the map $(p,a)\mapsto p\cdot a$
in the following Lemmas \ref{L:pLls} and \ref{L:BaspcdotL}.

\begin{lemma}\label{L:pLls}
Let\index{pzzroj@$\BB{S}$} $p\in\BB{S}$. Then the following assertions hold:
\begin{enumerate}
\item $p\cdot L=\DD^{-1}[pS]$.

\item $p\cdot L$ is a lower subset of $(L,\lesssim)$.

\item $p\cdot L$ is a lower subset of $(L,\leq)$.

\item $p\cdot L$ is closed under majorized suprema.

\end{enumerate}
\end{lemma}

\begin{proof}
(i) is an immediate consequence of Corollary~\ref{C:ProjOnL}.
The assertions (ii), (iii) follow immediately.

(iv) follows immediately from (i), Corollary~\ref{C:DelPerp}, and the fact
that $pS=(p^\bot S)^\bot$.
\end{proof}

\begin{lemma}\label{L:BaspcdotL}
Let $a$, $b\in L$, let $p$,\index{pzzroj@$\BB{S}$} $q\in\BB{S}$. Then the
following assertions hold:
\begin{enumerate}
\item $a\leq b$ implies that $p\cdot a\leq p\cdot b$.

\item $a\lesssim b$ implies that $p\cdot a\lesssim p\cdot b$.

\item $a\sim b$ implies that $p\cdot a\sim p\cdot b$.

\item $p\leq q$ implies that $p\cdot a\leq q\cdot a$.
\end{enumerate}
\end{lemma}

\begin{proof}
(i) $p\cdot a\leq a\leq b$ and $\DD(p\cdot a)\in pS$, thus
$p\cdot a\leq p\cdot b$ by the definition of~$p\cdot b$.

(ii) It follows from $a\lesssim b$ that $\DD(a)\leq\DD(b)$, thus
$\DD(p\cdot a)=p(\DD(a))\leq p(\DD(b))=\DD(p\cdot b)$, that is,
$p\cdot a\lesssim p\cdot b$.

(iii) follows immediately from (ii) and from Proposition~\ref{P:SchBern}.

(iv) Since $p\leq q$, $\DD(p\cdot a)\in pS\subseteq qS$, and so, since
$p\cdot a\leq a$, it follows that $p\cdot a\leq q\cdot a$.
\end{proof}

\begin{proposition}\label{P:piaOrth}
Let $\famm{p_i}{i\in I}$ be an orthogonal family of elements
of\index{pzzroj@$\BB{S}$}
$\BB{S}$ and let $a\in L$. Then the family $\famm{p_i\cdot a}{i\in I}$ is
orthogonal in $L$.
\end{proposition}

\begin{proof}
For $i\neq j$ in $I$, $\DD(p_i\cdot a)\in p_iS$ and
$\DD(p_j\cdot a)\in p_jS$, thus, since $p_ip_j=0$,
$\DD(p_i\cdot a)\perp\DD(p_j\cdot a)$. The result follows then from
Corollary~\ref{C:aiStrOrth}.
\end{proof}

\begin{lemma}\label{L:ppbotonL}
Let $a\in L$ and let \index{pzzroj@$\BB{S}$}$p\in\BB{S}$. Then the following
assertions hold:
\begin{enumerate}
\item $a=p\cdot a\boxplus p^\bot\cdot a$
\pup{thus, by Proposition~\textup{\ref{P:aboxplusb}},
$a=p\cdot a\oplus p^\bot\cdot a$}.

\item Let $x\in p\cdot L$ and $y\in p^\bot\cdot L$ such that $a=x\vee y$. Then
$x=p\cdot a$ and $y=p^\bot\cdot a$.
\end{enumerate}

\end{lemma}

\begin{proof}
(i) Since the join $a'=p\cdot a\vee p^\bot\cdot a$ is defined
(and $a'\leq a$) and since $\DD(p\cdot a)\perp\DD(p^\bot\cdot a)$, it
follows from Proposition~\ref{P:aboxplusb} that
$a'=p\cdot a\oplus p^\bot\cdot a$. Since
   \[
   \DD(a)=p(\DD(a))+p^\bot(\DD(a))=\DD(p\cdot a)+\DD(p^\bot\cdot a)
   =\DD(p\cdot a\oplus p^\bot\cdot a)
   \]
and by Axiom~(L6), there are $u\sim p\cdot a$ and
$v\sim p^\bot\cdot a$ such that $a=u\oplus v$. Since
$\DD(u)=\DD(p\cdot a)\in pS$ and $u\leq a$, we have $u\leq p\cdot a$.
Likewise, $v\leq p^\bot\cdot a$, thus $a=u\oplus v\leq a'$, whence
$a=a'=p\cdot a\boxplus p^\bot\cdot a$.

(ii) By assumption, $a=x\boxplus y$. By Proposition~\ref{P:aboxplusb},
$a=x\oplus y$. By Corollary~\ref{C:aboxplusb} and by (i), $x=p\cdot a$ and
$y=p^\bot\cdot a$.
\end{proof}

\begin{proposition}\label{P:ShasGC}
$S$ has general comparability.
\end{proposition}

\begin{proof}
We prove the two following claims.

\setcounter{claim}{0}

\begin{claim}\label{Cl:abot+abotbot=S}
$S=\la^\bot+\la^{\bot\bot}$, for all $\la\in S$.
\end{claim}

\begin{cproof}
Let $\la$, $\lx\in S$.
Pick $a$, $x\in L$ such that $\DD(a)=\la$ and $\DD(x)=\lx$. By
Corollary~\ref{C:DelPerp}, there exists a largest element $u\leq x$ such
that $\DD(u)\perp\la$. Let $v\in L$ such that $u\oplus v=x$. So
$\lx=\DD(u)+\DD(v)$, with $\DD(u)\in\la^\bot$. For any $t\leq v$ such that
$\DD(t)\in\la^\bot$, the inequality $t\leq u$ holds by the definition of~$u$,
but $t\leq v$, so $t=0$ since $u\perp v$. Hence $\DD(v)\in\la^{\bot\bot}$.
\end{cproof}

\begin{claim}\label{Cl:abcxyperp}
For all $a$, $b\in L$, there are $u$, $v$, $x$, $y\in L$ such that
$a=u\oplus x$ and $b=v\oplus y$ while $u\sim v$ and $\DD(x)\perp\DD(y)$.
\end{claim}

\begin{cproof}
An easy application of Zorn's Lemma yields a subset
$X=\setm{(a_i,b_i)}{i\in I}$ of $(0,a]\times(0,b]$
which is maximal with respect to the following properties:
\begin{enumerate}
\item Both families $\famm{a_i}{i\in I}$ and $\famm{b_i}{i\in I}$ are
orthogonal.

\item $a_i\sim b_i$ for all $i\in I$.
\end{enumerate}

Put $u=\bigvee_{i\in I}a_i$ and $v=\bigvee_{i\in I}b_i$. By Axiom~(L7),
$u\sim v$. Pick $x$ and $y$ such that $a=u\oplus x$ and $b=v\oplus y$. If
$\DD(x)\not\perp\DD(y)$, then there are nonzero $x'\leq x$ and $y'\leq y$
such that $x'\sim y'$. But then, $X\cup\set{(x',y')}$
still satisfies (i) and (ii) above, which contradicts the maximality of $X$.
Hence $\DD(x)\perp\DD(y)$.
\end{cproof}

By Lemma~\ref{L:GenCompAx}, general comparability follows from Claims
\ref{Cl:abot+abotbot=S} and \ref{Cl:abcxyperp}.
\end{proof}

We recall that general comparability is also Axiom~(M3). Note that general
comparability in $S$ can be stated as the following property of $L$, which
we also refer to as \emph{general comparability}: for any $a$, $b\in L$,
there exists $p\in\BB{S}$\index{pzzroj@$\BB{S}$} such that $p\cdot a\lesssim
p\cdot b$ and $p^\bot\cdot b\lesssim p^\bot\cdot a$.

\begin{lemma}\label{L:p(aopplusb)}
Let $a\in L$, let $\famm{a_i}{i\in I}$ be a family of elements of $L$, let
$p\in\BB{S}$.
\begin{enumerate}
\item If $I\ne\es$ and $a=\bigwedge_{i\in I}a_i$, then
$p\cdot a=\bigwedge_{i\in I}(p\cdot a_i)$.

\item If $a=\bigvee_{i\in I}a_i$, then
$p\cdot a=\bigvee_{i\in I}(p\cdot a_i)$.
\end{enumerate}
\end{lemma}

\begin{proof}
(i) It is clear that $p\cdot a\leq p\cdot a_i$ for all $i$. Let $b\in L$ such
that $b\leq p\cdot a_i$ for all $i\in I$. In particular, $b\leq a_i$ for all
$i$, thus $b\leq a$. In addition, since $I\ne\es$, it follows from
Lemma~\ref{L:pLls}(iii) that $b\in p\cdot L$, so $b=p\cdot b\leq p\cdot a$.

(ii) The equality $a_i=(p\cdot a_i)\vee(p^\bot\cdot a_i)$ holds for all
$i\in I$, so
   \[
   a=\bigvee_{i\in I}(p\cdot a_i)\vee
   \bigvee_{i\in I}(p^\bot\cdot a_i).
   \]
By Lemma~\ref{L:pLls}(iv), $\bigvee_{i\in I}(p\cdot a_i)\in p\cdot L$
and $\bigvee_{i\in I}(p^\bot\cdot a_i)\in p^\bot\cdot L$.
Therefore, by Lemma~\ref{L:ppbotonL},
$p\cdot a=\bigvee_{i\in I}(p\cdot a_i)$ and
$p^\bot\cdot a=\bigvee_{i\in I}(p^\bot\cdot a_i)$.
\end{proof}

\begin{proposition}\label{P:B(L)cBa}
The Boolean algebra\index{Boolean algebra!complete ---}
$\BB{S}$\index{pzzroj@$\BB{S}$} is complete.
\end{proposition}

\begin{proof}
It suffices to prove that every orthogonal family $\famm{p_i}{i\in I}$ of
elements of $\BB{S}$ admits a supremum. By Proposition~\ref{P:piaOrth},
the family $\famm{p_i\cdot a}{i\in I}$ is orthogonal, for all $a\in L$.
Furthermore, if $b\in L$ such that $a\sim b$, then, by
Lemma~\ref{L:BaspcdotL}(iii), $p_i\cdot a\sim p_i\cdot b$ for all $i\in I$,
thus, by Axiom~(L7),
   \[
   \oplus_{i\in I}(p_i\cdot a)\sim\oplus_{i\in I}(p_i\cdot b).
   \]
Hence, we can define a map $p\colon S\to S$ by the rule
   \[
   p(\DD(a))=\DD(\oplus_{i\in I}(p_i\cdot a)),\quad\text{for all }a\in L.
   \]
It is obvious that $p(0)=0$.
Let $\la$, $\lb\in S$ such that $\la+\lb$ is defined. So
$\la+\lb=\DD(a\oplus b)$, for some $a\in\la$ and $b\in\lb$ such that
$a\perp b$. So $p(\la)=\DD(a')$ and $p(\lb)=\DD(b')$, where $a'$ and $b'$
are defined by
   \[
   a'=\oplus_{i\in I}(p_i\cdot a)\quad\text{and}\quad
   b'=\oplus_{i\in I}(p_i\cdot b).
   \]
In particular, $a'\leq a$ and $b'\leq b$, so $a'\perp b'$, and, by using
Lemma~\ref{L:p(aopplusb)}(ii),
   \[
   a'\oplus b'=\oplus_{i\in I}(p_i\cdot a\oplus p_i\cdot b)
   =\oplus_{i\in I}(p_i\cdot(a\oplus b)),
   \]
whence
   \[
   p(\la)+p(\lb)=\DD(a'\oplus b')=\DD(\oplus_{i\in I}(p_i\cdot(a\oplus b)))
   =p(\la+\lb).
   \]
So $p$ is an endomorphism of $(S,+,0)$.

Now let $\la\in S$. Pick $a\in\la$, and pick $u\in L$ such that
$a=\oplus_{i\in I}(p_i\cdot a)\oplus u$. For all $v\leq u$, if
$\DD(v)\in p_iS$, then $v\leq p_i\cdot a$, so $v\leq u\wedge(p_i\cdot a)=0$.
Hence $\DD(u)\in(p_iS)^\bot$. We infer that $\DD(u)\in(pS)^\bot$.
Indeed, let $\lx\in S$. Pick $x$ such that $\DD(x)=\lx$.
So $p(\lx)=\DD(\oplus_{i\in I}(p_i\cdot x))$ by the definition of $p$.
But $\DD(p_i\cdot x)\perp\DD(u)$ for
all $i\in I$, thus, by Corollary~\ref{C:DelPerp}, $p(\lx)\perp\DD(u)$.
Hence $\DD(u)\in(pS)^\perp$, with $\la=p(\la)+\DD(u)$. It follows that $p$
is a projection of~$S$.

It is clear that $p_i\leq p$ for all $i\in I$. Let
$q\in\BB{S}$\index{pzzroj@$\BB{S}$} such that
$p_i\leq q$ for all $i\in I$. Let $a\in L$. Then $p_i\cdot a\leq q\cdot a$,
for all $i\in I$, thus $\oplus_{i\in I}(p_i\cdot a)\leq q\cdot a$. Taking
the image under $\DD$ of both sides yields that $p(\DD(a))\leq q(\DD(a))$.
This holds for all $a\in L$, whence $p\leq q$. So we have verified that
$p=\bigvee_{i\in I}p_i$.
\end{proof}

As a consequence, we obtain that $S$ satisfies Axiom~(M4) (observe that
$p\cdot a\lesssim b$ \iff\ $p(\DD(a))\leq p(\DD(b))$).

\begin{proposition}\label{P:BoolValL}
For all $a$, $b\in L$, there exists a largest
$p\in\BB{S}$\index{pzzroj@$\BB{S}$} such that
$p\cdot a\lesssim b$.
\end{proposition}

\begin{proof}
Put\index{pzzroj@$\BB{S}$} $X=\setm{q\in\BB{S}}{q\cdot a\lesssim b}$, and
put $p=\bigvee X$. Let $\famm{p_i}{i\in I}$ be a maximal orthogonal family
of elements of $X$. For all $i\in I$, let $b_i\leq b$ such that $p_i\cdot
a\sim b_i$. In particular, $\DD(b_i)=\DD(p_i\cdot a)\in p_iS$, so $b_i\leq
p_i\cdot b$. By Corollary~\ref{C:aiStrOrth}, the family $\famm{b_i}{i\in I}$
is orthogonal, and by Proposition~\ref{P:piaOrth}, the family
$\famm{p_i\cdot a}{i\in I}$ is orthogonal. Hence,
   \begin{align*}
   \DD(p\cdot a)&=p(\DD(a))
   &&(\text{by the definition of }p\cdot a)\\
   &=\DD(\oplus_{i\in I}(p_i\cdot a))
   &&(\text{by the proof of Proposition~\ref{P:B(L)cBa}})\\
   &=\DD(\oplus_{i\in I}b_i)
   &&(\text{by Axiom~(L7)})\\
   &\leq\DD(b).\tag*{\qed}
   \end{align*}
\renewcommand{\qed}{}
\end{proof}

\begin{notation}\label{Not:bvalsmb}
For $a$, $b\in L$, we put $\bv{a\lesssim b}=\bv{\DD(a)\leq\DD(b)}$. That is,
in accordance to Definition~\ref{D:BoolVal}, $\bv{a\lesssim b}$ is the
largest projection $p$ of~$S$ such that $p\cdot a\lesssim b$.

In a similar spirit as for Notation~\ref{Not:bva=b}, we put
$\bv{a\sim b}=\bv{a\lesssim b}\wedge\bv{b\lesssim a}$, which by
Proposition~\ref{P:SchBern} is the largest
$p\in\BB{S}$\index{pzzroj@$\BB{S}$} such that
$p\cdot a\sim p\cdot b$.
\end{notation}

\section{Purely infinite elements; trim sequences}\label{S:PItrim}

\begin{quote}
\em Standing hypotheses: $L$ is an espalier,\index{espalier} $S$ is the
dimension range of $L$, and
$\DD\colon L\twoheadrightarrow S$ is the canonical map.
\end{quote}

Following Definition~\ref{D:PurInf}, we say that an element $a$ of $L$ is
\emph{purely infinite},\index{purely infinite (=idem-multiple)|ii}
if $\DD(a)$ is purely infinite in $S$. This occurs
if and only if $a=a'\oplus a''$ for some $a'$, $a''\sim a$ in $L$.

Purely infinite elements are also, in some references, called
\emph{idempotent}, or, as in \index{Tarski, A.}\cite{Tars},
\emph{idem-multiple}.

Similarly, following Definition~\ref{D:DirFinMon}, we say that an element
$a$ of $L$ is\index{directly finite} \emph{directly finite}, if $\DD(a)$ is
directly finite\index{directly finite} in $S$. Observe that $a$ is directly
finite\index{directly finite} if and only if
$a\sim b\leq a$ implies that $b=a$, for any $b\in L$.

\begin{definition}\label{D:HomogSeq}
A family $\famm{a_i}{i\in I}$ of elements of $L$
is \emph{homogeneous},\index{homogeneous|ii} if it is orthogonal and
$a_i\sim a_j$, for all $i$, $j\in I$.

A homogeneous family $\famm{a_i}{i\in I}$ is \emph{trivial}, if $a_i=0$, for
all $i$; equivalently, $a_i=0$ for \emph{some} $i\in I$.
\end{definition}
\begin{lemma}\label{L:CharPurInf}
Let $a\in L$. Then the following are equivalent:
\begin{enumerate}
\item $a$ is purely infinite.

\item There exists a homogeneous sequence $\famm{x_n}{n<\omega}$ such that
$a=\oplus_{n<\omega}x_n$.

\end{enumerate}
\end{lemma}

\begin{proof}
(i)$\Rightarrow$(ii) There are $a'$ and $a''$ such that $a=a'\oplus a''$
and $a\sim a'\sim a''$.
By using Axiom~(L6), it is then easy to construct inductively sequences
$\famm{a_n}{n<\omega}$ and $\famm{a'_n}{n<\omega}$ such that
$a_0=a$, $a_n=a_{n+1}\oplus a'_n$, and $a_n\sim a_{n+1}\sim a'_n$,
for all $n$. So $\famm{a'_n}{n<\omega}$ is a homogeneous sequence whose
join, $a'$, belongs to $[0,a]$. Since $a\sim a'_0\leq a'\leq a$, we obtain
that $a\sim a'$ by Proposition~\ref{P:SchBern}.
Hence, the conclusion follows from Axiom~(L6).

(ii)$\Rightarrow$(i) Put $a'=\oplus_nx_{2n}$ and $a''=\oplus_nx_{2n+1}$. By
Axiom~(L7), $a\sim a'\sim a''$. Since $a=a'\oplus a''$, $a$ is purely
infinite.
\end{proof}

\begin{lemma}\label{L:CharDf}
Let $a\in L$. The following are equivalent:
\begin{enumerate}
\item $a$ is\index{directly finite} directly finite.

\item There is no nontrivial purely infinite element below $a$.

\item The interval $[0,a]$ has no infinite nontrivial homogeneous sequence.
\end{enumerate}
\end{lemma}

\begin{proof}
(i)$\Rightarrow$(ii)
Let $b\in[0,a]$ be a purely infinite element, and let $x\in L$ such that
$a=b\oplus x$. Then $\DD(a)=\DD(b)+\DD(x)=2\DD(b)+\DD(x)=\DD(b)+\DD(a)$.
Since $a$ is\index{directly finite} directly finite, $\DD(b)=0$, whence
$b=0$.

(ii)$\Leftrightarrow$(iii) follows immediately from Lemma~\ref{L:CharPurInf}.

(iii)$\Rightarrow$(i)
Let $x\in L$ such that $\DD(a)+\DD(x)=\DD(a)$. So $a=a'\oplus x'$, for some
$a'\sim a$ and $x'\sim x$. It is then easy to construct, by induction (and
Axiom~(L6)), sequences $\famm{a_n}{n<\omega}$ and $\famm{x_n}{n<\omega}$
of elements of $L$ such that $a_0=0$, $a_n\sim a$, $x_n\sim x$,
and $a_n=a_{n+1}\oplus x_n$ for all $n$. In particular, the sequence
$\famm{x_n}{n<\omega}$ is homogeneous, thus, by assumption, $x_n=0$ for all
$n$. Therefore, $x=0$. So $a$ is\index{directly finite} directly finite.
\end{proof}

We deduce from this that $S$ satisfies Axiom~(M5).

\begin{proposition}\label{P:DfPiDec}
For all $a\in L$, there are $b$, $c\in L$ such that $b$ is
purely infinite, $c$ is\index{directly finite} directly finite, and
$a=b\oplus c$.
\end{proposition}

\begin{proof}
Let $\famm{x_i}{i\in I}$ be a maximal orthogonal family of nonzero purely
infinite elements of $[0,a]$. We put $b=\oplus_{i\in I}x_i$. For all
$i\in I$, there exists a decomposition $x_i=x'_i\oplus x''_i$ where
$x'_i\sim x''_i\sim x_i$. Put $b'=\oplus_{i\in I}x'_i$ and
$b''=\oplus_{i\in I}x''_i$. By Axiom~(L7), $b'\sim b''\sim b$. Since
$b=b'\oplus b''$, $b$ is purely infinite.

Let $c\in L$ such that $a=b\oplus c$. Suppose that $c$ is
not\index{directly finite} directly finite. Then, by Lemma~\ref{L:CharDf},
there exists a purely infinite element $x$ such that $0<x\leq c$. But then,
enlarging the family
$\famm{x_i}{i\in I}$ by $x$ yields an orthogonal family of nonzero purely
infinite elements of
$[0,a]$, which contradicts the maximality of $\famm{x_i}{i\in I}$.
So, $c$ is\index{directly finite} directly finite.
\end{proof}

We can then reformulate Proposition~\ref{P:bsminusa} in the language of
lattices.

\begin{proposition}\label{P:bsminusaL}
Let $a$, $b\in L$ such that $a\lesssim b$. Then there exists $c\in L$ such
that $c\leq b$ and $\DD(c)=\DD(b)\sd\DD(a)$.
\end{proposition}

\begin{proof}
Propositions \ref{P:SchBern}, \ref{P:Ref}, \ref{P:ShasGC},
\ref{P:BoolValL}, and \ref{P:DfPiDec} establish the hypotheses of
Proposition~\ref{P:bsminusa}. Thus, $\DD(b)\sd\DD(c)$ exists in $S$. Since
this element lies below $\DD(b)$, there exists $c\in[0,b]$ such that
$\DD(c)=\DD(b)\sd\DD(a)$.
\end{proof}

The following important definition involves both the lattice structure and
the dimension function. It is the key to proving the existence of majorized
suprema in~$S$.

\begin{definition}\label{D:trim}\hfill
\begin{enumerate}
\item Let $a$, $b\in L$. We write\index{abbtzzrim@$a\trim b$|ii}
$a\trim b$, if there exists $c\in L$ such
that $a\oplus c=b$ and $\DD(c)=\DD(b)\sd\DD(a)$.

\item Let $\kappa$ be an ordinal. A $\kappa$-sequence
$\famm{a_\xi}{\xi<\kappa}$ of elements of $L$ is \emph{trim},
\index{trim (sequence)|ii} if the following conditions hold:
\begin{enumerate}
\item $a_\xi\trim a_{\xi+1}$ for all $\xi$ such that $\xi+1<\kappa$.

\item For any limit ordinal $\lambda<\kappa$,
the sequence $\setm{a_\xi}{\xi<\lambda}$ is majorized, and
$\bigvee_{\xi<\lambda}a_\xi\trim a_\lambda$.
\end{enumerate}
\end{enumerate}
\end{definition}

\begin{lemma}\label{L:TrimIns}
Let $a$, $b\in L$ such that $a\leq b$, let $\lx\in S$ such that
$\DD(a)\leq\lx\leq\DD(b)$. Then there exists $x\in L$ such that
$a\trim x\leq b$ and $\DD(x)=\lx$.
\end{lemma}

\begin{proof}
Let $c\in L$ such that $a\oplus c=b$. So
$\DD(a)\leq\lx\leq\DD(b)=\DD(a)+\DD(c)$, thus $\lx\sd\DD(a)\leq\DD(c)$.
Hence there exists $y\leq c$ such that
$\DD(y)=\lx\sd\DD(a)$. Now we put $x=a\oplus y$. Then
$\DD(x)=\DD(a)+(\lx\sd\DD(a))=\lx$, $a\leq x\leq b$, and
$\DD(x)\sd\DD(a)=\lx\sd\DD(a)=\DD(y)$. So, $a\trim x$.
\end{proof}

In the statement of the following Lemma~\ref{L:TrimSeq}, a \emph{lifting} of
a family $\famm{\la_i}{i\in I}$ of elements of~$S$ is a family
$\famm{a_i}{i\in I}$ of elements of $L$ such that $\DD(a_i)=\la_i$ for all
$i\in I$.

\begin{lemma}\label{L:TrimSeq}
Let $\kappa$ be an ordinal.
\begin{enumerate}
\item For all $b\in L$, every increasing $\kappa$-sequence of elements of
$[0,\DD(b)]$ has a trim lifting in $[0,b]$.

\item For any majorized trim sequences $\famm{x_\xi}{\xi<\kappa}$ and
$\famm{y_\xi}{\xi<\kappa}$ of elements of $L$,
   \[
   x_\xi\sim y_\xi\text{ for all }\xi<\kappa\quad\text{implies that}\quad
   \bigvee_{\xi<\kappa}x_\xi\sim\bigvee_{\xi<\kappa}y_\xi.
   \]

\item For every majorized trim lifting $\famm{a_\xi}{\xi<\kappa}$ of
a $\kappa$-sequence $\famm{\la_\xi}{\xi<\kappa}$ of elements of~$S$,
   \[
   \DD\biggl(\bigvee_{\xi<\kappa}a_\xi\biggr)=\bigvee_{\xi<\kappa}\la_\xi.
   \]
\end{enumerate}

\end{lemma}

\begin{proof}
We argue by transfinite induction on $\kappa$. The result is vacuous  for
$\kappa=0$. Suppose that we have proved the lemma for all ordinals
$\kappa'<\kappa$, with $\kappa>0$.

(i) Let $\famm{\la_\xi}{\xi<\kappa}$ be an increasing $\kappa$-sequence of
elements of $[0,\DD(b)]$.
We construct inductively elements $a_\xi$ of $[0,b]$, for $\xi<\kappa$.

For $\xi=0$, pick any element $a_0$ of $[0,b]$ such that $\DD(a_0)=\la_0$.

Suppose we have constructed $a_\xi\leq b$ such that $\DD(a_\xi)=\la_\xi$,
with $\xi+1<\kappa$. By Lemma~\ref{L:TrimIns}, there exists
$a_{\xi+1}\leq b$ such that $a_\xi\trim a_{\xi+1}$ and
$\DD(a_{\xi+1})=\la_{\xi+1}$.

Suppose finally that $\lambda<\kappa$ is a limit ordinal and that
$\famm{a_\xi}{\xi<\lambda}$ is a trim lifting of
$\famm{\la_\xi}{\xi<\lambda}$ in $[0,b]$. We put
   \[
   \oll{a}_\lambda=\bigvee_{\xi<\lambda}a_\xi.
   \]
Now, we observe that $\DD(a_\xi)=\la_\xi\leq\la_\lambda$ for all
$\xi<\lambda$. Since $\famm{a_\xi}{\xi<\lambda}$ is trim and majorized, it
follows from (iii) of the induction hypothesis that
$\DD(\oll{a}_\lambda)\leq\la_\lambda$. By applying once again
Lemma~\ref{L:TrimIns}, we obtain $a_\lambda\leq b$ such that
$\oll{a}_\lambda\trim a_\lambda$ and
$\DD(a_\lambda)=\la_\lambda$.

By the definition of a trim sequence, $\famm{a_\xi}{\xi<\kappa}$ is a trim
lifting of $\famm{\la_\xi}{\xi<\kappa}$ in $[0,b]$.\smallskip

(ii) We construct inductively elements $x'_\xi\in L$, for $\xi<\kappa$, of
$L$, as follows. We put $x'_0=x_0$. If $\xi+1<\kappa$, then
$x_\xi\trim x_{\xi+1}$, so there exists $x'_{\xi+1}$ such that
$x_{\xi+1}=x_\xi\oplus x'_{\xi+1}$ and
$\DD(x'_{\xi+1})=\DD(x_{\xi+1})\sd\DD(x_\xi)$.
If $\lambda<\kappa$ is a limit ordinal, we put
$\oll{x}_\lambda=\bigvee_{\xi<\lambda}x_\xi$. Since
$\oll{x}_\lambda\trim x_\lambda$, there exists $x'_\lambda$ such that
$\oll{x}_\lambda\oplus x'_\lambda=x_\lambda$ and
$\DD(x'_\lambda)=\DD(x_\lambda)\sd\DD(\oll{x}_\lambda)$.
It follows that $x_\xi=\oplus_{\eta\leq\xi}x'_\eta$ for all $\xi<\kappa$.
In particular, $\bigvee_{\xi<\kappa}x_\xi=\oplus_{\xi<\kappa}x'_\xi$.

Let $\famm{y'_\xi}{\xi<\kappa}$ be constructed from
$\famm{y_\xi}{\xi<\kappa}$ the same way $\famm{x'_\xi}{\xi<\kappa}$ is
constructed from $\famm{x_\xi}{\xi<\kappa}$. So $x'_0=x_0\sim y_0=y'_0$.
Let $\xi$ such that $\xi+1<\kappa$. Since $x_\xi\sim y_\xi$ and
$x_{\xi+1}\sim y_{\xi+1}$,
   \[
   \DD(x'_{\xi+1})=\DD(x_{\xi+1})\sd\DD(x_\xi)=\DD(y_{\xi+1})\sd\DD(y_\xi)
   =\DD(y'_{\xi+1}).
   \]
Let $\lambda<\kappa$ be a limit ordinal. By (iii) of the induction
hypothesis, $\oll{x}_\lambda\sim\oll{y}_\lambda$. Since
$x_\lambda\sim y_\lambda$, we obtain that
   \[
   \DD(x'_\lambda)=\DD(x_\lambda)\sd\DD(\oll{x}_\lambda)=
   \DD(y_\lambda)\sd\DD(\oll{y}_\lambda)=\DD(y'_\lambda).
   \]
Hence we have proved that $x'_\xi\sim y'_\xi$ for all $\xi<\kappa$. Hence,
by Axiom~(L7), $\oplus_{\xi<\kappa}x'_\xi\sim\oplus_{\xi<\kappa}y'_\xi$,
that is, $\bigvee_{\xi<\kappa}x_\xi\sim\bigvee_{\xi<\kappa}y_\xi$.
\smallskip

(iii) Let $\famm{a_\xi}{\xi<\kappa}$ be a majorized trim lifting of
$\famm{\la_\xi}{\xi<\kappa}$. We put $a=\bigvee_{\xi<\kappa}a_\xi$.
So, $\la_\xi=\DD(a_\xi)\leq\DD(a)$ for all $\xi<\kappa$. Now let $\lb\in S$
such that $\la_\xi\leq\lb$ for all $\xi<\kappa$, and let $b\in L$ such that
$\DD(b)=\lb$. By (i), $\famm{\la_\xi}{\xi<\kappa}$ has a trim lifting
$\famm{a'_\xi}{\xi<\kappa}$ in $[0,b]$. In particular,
$\DD\Bigl(\bigvee_{\xi<\kappa}a'_\xi\Bigr)\leq\DD(b)$. However, by (ii),
$\DD(a)=\DD\Bigl(\bigvee_{\xi<\kappa}a'_\xi\Bigr)$; whence $\DD(a)\leq\lb$.
So $\DD(a)=\bigvee_{\xi<\kappa}\la_\xi$.
\end{proof}

\begin{corollary}\label{C:ShasM1}
$S$ satisfies Axiom~\textup{(M2)}.
\end{corollary}

\begin{proof}
We prove that every majorized subset $X$ of~$S$ has a supremum.
By Proposition~\ref{P:SchBern}, Proposition~\ref{P:ShasGC}, and
Lemma~\ref{L:meetjoinS}, every majorized finite subset of
$S$ has a supremum.

So it remains to conclude in case $X$ is infinite. We argue by induction on
the cardinality of $X$. Write $X=\setm{\la_\xi}{\xi<\kappa}$, where $\kappa$
is the cardinality of $X$. By the finite case and the induction hypothesis,
for all $\xi<\kappa$, the set $\setm{\la_\eta}{\eta\leq\xi}$ has a supremum,
say, $\lb_\xi$. Since $X$ is majorized, so is $\setm{\lb_\xi}{\xi<\kappa}$,
that is, there exists $b\in L$ such that $\lb_\xi\leq\DD(b)$ for all
$\xi<\kappa$. By Lemma~\ref{L:TrimSeq}(i), the family
$\famm{\lb_\xi}{\xi<\kappa}$ has a trim lifting in $[0,b]$, say,
$\famm{b_\xi}{\xi<\kappa}$. Put $c=\bigvee_{\xi<\kappa}b_\xi$.
By Lemma~\ref{L:TrimSeq}(iii), $\DD(b)$ is the supremum of
$\setm{\lb_\xi}{\xi<\kappa}$, that is, the supremum of~$X$.
\end{proof}

\section{Axiom (M6)}\label{S:Shas(M6)}

\begin{quote}
\em Standing hypotheses: $L$ is an espalier,\index{espalier} $S$ is the
dimension range\index{dimension range} of $L$, and $\DD\colon
L\twoheadrightarrow S$ is the canonical map.
\end{quote}

At this point, what remains to do in order to conclude the proof of
Theorem~A is to establish that $S$ satisfies Axiom~(M6). We shall devote
Section~\ref{S:Shas(M6)} to this.

In accordance with Definition~\ref{D:Remov}, we state the following
definition.

\begin{definition}\label{D:RemovL}
Let $a$, $b\in L$. We say that $a$ is \emph{removable} from $b$,
\index{removable (in espaliers)|ii} in notation
$a\srem b$\index{abrszzem@$a\srem b$|ii}, if $\DD(a)\rem\DD(b)$ in $S$.
Equivalently, $a\srem b$, if
$a\lesssim b$, and $b\lesssim a\oplus x$ implies that $b\lesssim x$, for all
$x\in S$.
\end{definition}

\begin{notation}\label{Not:CardElt}
Let $\kappa$ be a cardinal number, let $a$, $b\in L$.
\begin{enumerate}
\item Let $\kappa\cdot a\sim b$ be the statement that there exists a
homogeneous $\kappa$-sequence $\famm{a_\xi}{\xi<\kappa}$ such that
   \[
   \oplus_{\xi<\kappa}a_\xi=b\quad\text{and}\quad a\sim a_0.
   \]
\item Let $\kappa\cdot a\lesssim b$ be the statement that there exists a
homogeneous $\kappa$-sequence $\famm{a_\xi}{\xi<\kappa}$ such that
   \[
   \oplus_{\xi<\kappa}a_\xi\leq b\quad\text{and}\quad a\sim a_0.
   \]
\end{enumerate}
\end{notation}

For example, $1\cdot a\sim b$ (resp., $1\cdot a\lesssim b$) means that
$a\sim b$ (resp., $a\lesssim b$). Another example is that $2\cdot a\sim a$
\iff\ $a$ is purely infinite.

\begin{lemma}\label{L:LocInfQuot}
Let $a$, $b\in L\setminus\set{0}$, let $\beta$ be an infinite cardinal number.
If $\beta\cdot a\lesssim b$, then there exist an infinite
cardinal number $\gamma\geq\beta$ and a projection $p$ of~$S$ such
that $p\cdot a>0$ and $\gamma\cdot(p\cdot a)\sim p\cdot b$.
\end{lemma}

\begin{proof}
We start with a homogeneous family of $\beta$ elements of $[0,b]$ all
equivalent to $a$ (modulo $\sim$), and enlarge it to a maximal such family,
say, $\vec a=\famm{a_\xi}{\xi<\gamma}$, where $\gamma\geq\beta$ is an
infinite cardinal number. Let $b'\in L$ such that
   \[
   b=b'\oplus(\oplus_{\xi<\gamma}a_\xi).
   \]
By general comparability, there exists $p\in\BB{S}$\index{pzzroj@$\BB{S}$}
such that
$p\cdot b'\lesssim p\cdot a$ and $p^\bot\cdot a\lesssim p^\bot\cdot b'$. By
the maximality of $\vec a$, $a\not\lesssim b'$, hence $p\cdot a>0$.
Now we put
   \[
   b^*=b'\oplus(\oplus_{0<\xi<\gamma}a_\xi).
   \]
Since $\gamma$ is an infinite cardinal and by Axiom~(L7), $b\sim b^*$.
Moreover, $p\cdot b'\lesssim p\cdot a_0$, so, by Lemmas
\ref{L:BaspcdotL} and \ref{L:p(aopplusb)},
   \[
   p\cdot b\sim p\cdot b^*\lesssim p\cdot(\oplus_{\xi<\gamma}a_\xi)
   \leq p\cdot b.
   \]
Hence, by using Proposition~\ref{P:SchBern} and Lemma~\ref{L:p(aopplusb)},
   \begin{equation}
   p\cdot b\sim p\cdot(\oplus_{\xi<\gamma}a_\xi)=
   \oplus_{\xi<\gamma}(p\cdot a_\xi).\tag*{\qed}
   \end{equation}
\renewcommand{\qed}{}
\end{proof}

\begin{lemma}\label{L:PurInfAl0}
$\aleph_0\cdot a\sim a$, for all purely infinite $a\in L$.
\end{lemma}

\begin{proof}
By Lemma~\ref{L:CharPurInf}, we have $a=\oplus_{n<\omega}x_n$ for some
homogeneous sequence $\famm{x_n}{n<\omega}$. Let
$\omega=\bigsqcup_{n<\omega}I_n$ be an infinite partition of $\omega$, with
all the $I_n$ infinite. Put
$a_n=\oplus_{k\in I_n}x_k$, for all $n<\omega$. By Axiom~(L7), $a_n\sim a$
for all $n$. The proof is concluded by the observation that
$a=\oplus_{n<\omega}a_n$.
\end{proof}

By replacing a bijection from $\omega\times\omega$ onto $\omega$ by a
bijection from $\kappa\times\kappa$ onto $\kappa$, for any infinite cardinal
$\kappa$, in the proof above, we easily obtain the following result.

\begin{lemma}\label{L:PurInfKappa}
Let $a$, $b\in L$, let $\kappa$ be an infinite cardinal number. If
$b\sim\kappa\cdot a$, then $b\sim\kappa\cdot b$.
\end{lemma}

\begin{notation}
For $a\in L$, we put $\cc(a)=\cc(\DD(a))$ (see Definition~\ref{D:cc(a)}).
\end{notation}

So, $\cc(a)=\bv{a\sim 0}^\bot$, for all $a\in L$. In view of
Lemma~\ref{L:ppbotonL}, $\cc(a)$ is the smallest projection $p$ of $S$ such
that $p\cdot a=a$.

\begin{lemma}\label{L:LocM5}
Let $a$, $b\in L$ purely infinite such that $a\srem b$ and $b\neq 0$. Then
there exists a purely infinite $e\in L$ such that
\begin{enumerate}
\item $e\leq b$ and $e\not\lesssim a$.

\item $e\lesssim c$, for all purely infinite $c\leq b$ such that $a\srem c$
and $\cc(e)\leq\cc(c)$.

\end{enumerate}
\end{lemma}

\begin{proof}
We put $F=\setm{x\leq b}{x\text{ is purely infinite and }x\not\lesssim a}$.
If $b\notin F$, then
$b\lesssim a$, thus, since $a\rem b$, $b=0$, a contradiction. So, $b\in F$.
For all $x\in F$, we denote by $\nu(x)$ the least infinite cardinal number
$\alpha$ such that $\alpha\cdot y\not\lesssim x$ for all $y\in F$. By
Lemma~\ref{L:PurInfAl0}, $\nu(x)\geq\aleph_1$ for all $x\in F$. We pick
$e\in F$ such that $\nu(e)=\alpha$ is minimum, and we prove that this $e$
satisfies the required conditions. Of course, (i) holds since $e\in F$.

Let $c\in L$ be purely infinite such that $a\srem c\leq b$ and
$\cc(e)\leq\cc(c)$. Note that $\cc(e)>0$ (otherwise, $e=0$, a contradiction).

\setcounter{claim}{0}
\begin{claim}\label{Cl:qcbetaqc}
For all $p\in(0,\cc(c)]$ and for every infinite cardinal number
$\beta<\alpha$, there exists $q\in(0,p]$ such that
$q\cdot c\sim\beta\cdot(q\cdot c)$.
\end{claim}

\begin{cproof}
If $p\cdot c\lesssim p\cdot a$, then, since $p\cdot a\srem p\cdot c$
(Lemma~\ref{L:ProjTr}(i)), $p\cdot c=0$, which is impossible since
$0<p\leq\cc(c)$. So, $p\cdot c\not\lesssim p\cdot a$, so $p\cdot c\in F$.
In particular, $\nu(p\cdot c)\geq\alpha>\beta$, so, by the definition of
$\nu(p\cdot c)$, there exists $d\in F$ such that
$\beta\cdot d\lesssim p\cdot c$. Note that $p\cdot d=d$.
Hence, by Lemma~\ref{L:LocInfQuot}, there are
$q\in\BBp{S}$\index{pzzrojst@$\BBp{S}$} and an infinite cardinal number
$\gamma\geq\beta$ such that
$\gamma\cdot(q\cdot d)\sim qp\cdot c$ and $q\cdot d>0$.
In particular, $qp\cdot d=q\cdot d>0$, so we may replace $q$ by $qp$,
and then $\gamma\cdot(q\cdot d)\sim q\cdot c$. Therefore, by
Lemma~\ref{L:PurInfKappa}, $\gamma\cdot(q\cdot c)\sim q\cdot c$, with
$\gamma\geq\beta$, thus, since $\gamma\geq\beta$ and by
Proposition~\ref{P:SchBern}, $\beta\cdot(q\cdot c)\sim q\cdot c$.
\end{cproof}

\begin{claim}\label{Cl:qelsmqc}
For all $p\in(0,\cc(c)]$, there exists $q\in(0,p]$ such that
$q\cdot e\lesssim q\cdot c$.
\end{claim}

\begin{cproof}
By general comparability, there exists a decomposition $p=p'\oplus p''$
(in\index{pzzroj@$\BB{S}$}
$\BB{S}$) such that $p'\cdot c\lesssim p'\cdot e$ and
$p''\cdot e\lesssim p''\cdot c$. If $p''>0$, then we may take $q=p''$. So
suppose that $p''=0$, so $p\cdot c\lesssim p\cdot e$. Since $p\cdot c$ is
purely infinite, there exist, by Lemmas \ref{L:LocInfQuot} and
\ref{L:PurInfAl0}, $q\in\BB{S}$\index{pzzroj@$\BB{S}$} and an infinite
cardinal $\beta$ such that
$\beta\cdot(qp\cdot c)\sim qp\cdot e$ and $qp\cdot c>0$. After replacing
$q$ by $qp$, we have $q\in(0,p]$, with
   \begin{equation}\label{Eq:betqcqe}
   \beta\cdot(q\cdot c)\sim q\cdot e
   \end{equation}
and $q\cdot c>0$. Since $q\cdot a\srem q\cdot c$, we must have
$q\cdot c\not\lesssim q\cdot a$, whence $q\cdot c\not\lesssim a$, and so
$q\cdot c\in F$. Hence $\beta<\nu(e)=\alpha$. Hence, by
Claim~\ref{Cl:qcbetaqc}, there exists $r\in(0,q]$ such that
$\beta\cdot(r\cdot c)\sim r\cdot c$. Therefore, by \eqref{Eq:betqcqe},
$r\cdot e\sim\beta\cdot(r\cdot c)\sim r\cdot c$.
\end{cproof}

By Claim~\ref{Cl:qelsmqc} and by Proposition~\ref{P:BoolValL},
$\cc(c)\cdot e\lesssim c$. However, by assumption, $\cc(e)\leq\cc(c)$, thus
$\cc(c)\cdot e=e$ (see Lemma~\ref{L:cc(a)}), so $e\lesssim c$.
\end{proof}

And now, Axiom~(M6) (recall that $\DD(c)^\bot=\DD(b)^\bot$ if and only if
$\cc(c)=\cc(b)$; see Definition~\ref{D:cc(a)}).

\begin{proposition}\label{P:Shas(M6)}
Let $a$, $b\in L$ be purely infinite such that $a\srem b$. Then there
exists a purely infinite $e\leq b$ such that
\begin{enumerate}
\item $a\srem e$ and $\cc(e)=\cc(b)$.

\item $e\lesssim c$, for all purely infinite $c\in L$ such that $a\srem c$
and $\cc(c)=\cc(b)$.
\end{enumerate}

\end{proposition}

\begin{proof}
We first claim that it will suffice to find a purely infinite element
$e\leq b$ satisfying (i) and the the statement
\begin{itemize}
\item[(i$'$)] $e\lesssim c$, for all purely infinite $c\leq b$ such that
$a\srem c$ and $\cc(c)=\cc(b)$.
\end{itemize}

Indeed, suppose that (i) and (ii$'$) are satisfied. Let $c\in L$ such that
$a\srem c$ and $\cc(c)=\cc(b)$. Then $\DD(a)\srem\DD(b),\DD(c)$, and so
$\DD(a)\srem\DD(b)\wedge\DD(c)$ by Corollary~\ref{C:2-5.?}. There exists
$d\leq b$ such that $\DD(d)=\DD(b)\wedge\DD(c)$, whence
$a\srem d\lesssim c$. Moreover, $d$ is purely infinite by
Lemma~\ref{L:2-5.?}, and $\cc(d)=\cc(b)\wedge\cc(c)=\cc(b)$ by
Lemma~\ref{L:Basiccc}(ii). Since any element $e\lesssim d$ would then
satisfy $e\lesssim c$, the claim is proved.

Let $P$ be the set of all pairs $(p,x)\in(\BB{S})\times
L$\index{pzzroj@$\BB{S}$} such that
$x$ is purely infinite and the following conditions hold:
\begin{itemize}
\item[(a)] $p\leq\cc(b)$ and $x\leq p\cdot b$.

\item[(b)] $p\cdot a\srem x$ and $\cc(x)=p$.

\item[(c)] For all purely infinite $y\leq p\cdot b$, the conditions
$p\cdot a\srem y$ and $\cc(y)=p$ imply that $x\lesssim y$.
\end{itemize}

Let $\setm{(p_i,x_i)}{i\in I}$ be a subset of $P$, maximal with the
property that the $p_i$ are nonzero and pairwise orthogonal. We observe
that since $\setm{x_i}{i\in I}$ is majorized (by $b$) and since the $p_i$
are pairwise orthogonal, it follows from Corollary~\ref{C:aiStrOrth} that
$\famm{x_i}{i\in I}$ is an orthogonal family of $L$. We put
   \[
   p=\bigvee_{i\in I}p_i\quad\text{and}\quad x=\oplus_{i\in I}x_i.
   \]

\setcounter{claim}{0}
\begin{claim}
The pair $(p,x)$ belongs to $P$.
\end{claim}

\begin{cproof}
Observe that $\cc(x)=p\leq\cc(b)$ and $x\leq p\cdot b$.
Since all the $x_i$ are purely infinite, $x$ is purely infinite.
Furthermore, $p_i\cdot a\srem x_i\leq x$ for all $i$, thus,
by Lemma~\ref{L:TrLeqTr},
$p_i\cdot a\srem x$. This holds for all $i$, thus, by
Lemma~\ref{L:ProjTr}(ii) and Proposition~\ref{P:B(S)cBa}, $p\cdot a\srem x$.

Let $y\leq p\cdot b$ be a purely infinite element of $L$ such that
$p\cdot a\srem y$ and $\cc(y)=p$. For all $i\in I$,
$p_i\cdot a\srem p_i\cdot y$ and, by Lemma~\ref{L:Basiccc}(ii),
$\cc(p_i\cdot y)=p_i$, so $p_i\cdot x=x_i\lesssim p_i\cdot y$. This holds
for all $i$, thus $x=p\cdot x\lesssim y$.
\end{cproof}

So, it suffices to prove that $p=\cc(b)$. Until the end of the proof, we
suppose otherwise. Put $q=\cc(b)p^\bot>0$. Since $0<q\leq\cc(b)$ and
$a\srem b$, the relation $q\cdot b\not\lesssim q\cdot a$ holds. Since
$q\cdot a\srem q\cdot b$, there exists, by Lemma~\ref{L:LocM5}, a purely
infinite $x^*\leq q\cdot b$ such that
\begin{itemize}
\item[($\alpha$)] $x^*\not\lesssim q\cdot a$;

\item[($\beta$)] $x^*\lesssim y^*$, for all purely infinite
$y^*\leq q\cdot b$ such that
$q\cdot a\srem y^*$ and $\cc(x^*)\leq\cc(y^*)$.
\end{itemize}

Now put $r=\bv{x^*\lesssim a}$. So, by definition,
$r\cdot x^*\lesssim r\cdot a$. If $r\cdot x^*=x^*$, then
$x^*\lesssim r\cdot a$, but $x^*\lesssim q\cdot b$, so
$x^*=q\cdot x^*\lesssim qr\cdot a\lesssim q\cdot a$, which contradicts
($\alpha$) above. So, $r\cdot x^*\neq x^*$, thus, since $x^*\in\cc(b)S$,
the projection $p'=\cc(b)r^\bot$ is nonzero.

Now, from $pq=0$ it follows that $p\cdot x^*=0$, so $p\leq r$, that is,
$p\perp p'$. Hence, $p'=qr^\bot$. Moreover,
$\cc(b)\cc(x^*)^\bot\leq\cc(x^*)^\bot\leq r$, whence, taking complements in
$\cc(b)$, $p'\leq\cc(x^*)$. In particular,
$\cc(p'\cdot x^*)=p'$.

By general comparability, there exists $g\leq\cc(b)$ in
$\BB{S}$\index{pzzroj@$\BB{S}$} such that
$g\cdot a\lesssim g\cdot x^*$ and $g^\bot\cdot x^*\lesssim g^\bot\cdot a$.
Then $\cc(b)g^\bot\leq r$, whence $p'\leq g$, and so
$p'\cdot a\lesssim p'\cdot x^*$.
By the definition of $r$ and of $p'$, $s\cdot x^*\not\lesssim s\cdot a$, for
all $s\in(0,p']$, thus, by Lemma~\ref{L:BVtr}, $p'\cdot a\srem p'\cdot x^*$.

Consider a purely infinite $c'\leq p'\cdot b$ such that $p'\cdot a\srem c'$
and $\cc(c')=p'$. Since $c'\leq p'\cdot b$ and $qr\perp p'$, the element
$c^*=c'\oplus qr\cdot b$ is defined and
$c^*\leq(p'\vee qr)\cdot b=q\cdot b$. Furthermore, since $q\leq\cc(b)$,
$\cc(c^*)=p'\vee qr=q$.
Since $p'\cdot a\srem c'$ and $a\srem b$, it follows from
Lemma~\ref{L:ProjTr} that $q\cdot a\srem c^*$.

Therefore, by part ($\beta$) of the definition of $x^*$, $x^*\lesssim c^*$,
thus $p'\cdot x^*\lesssim p'\cdot c^*=c'$. So we have proved that
$(p',p'\cdot x^*)\in P$, with $p'$ nonzero and orthogonal to all the
$p_i$ for $i\in I$, which contradicts the maximality of
$\setm{(p_i,x_i)}{i\in I}$. So, $p=\cc(b)$.
\end{proof}

Proposition~\ref{P:Shas(M6)} concludes the proof of Theorem~A. A more
complete form of Theorem~A is the following.

\begin{theorem}\label{T:DimEsp}
Let $(L,\leq,\perp,\sim)$ be an espalier.\index{espalier} Then the quotient
$\Drng L=L/{\sim}$\index{Dzzrng@$\Drng L$} can be endowed with a partial
addition $+$, defined by the rule
   \[
   \DD(c)=\DD(a)+\DD(b),\quad\text{for all }a,\,b,\,c\in L\text{ such that }
   c=a\oplus b,
   \]
that makes it a continuous dimension scale,\index{continuous dimension scale} with zero
element $\DD(0)$.
\end{theorem}

\section{D-universal classes of espaliers}\label{S:DUniv}
\index{espalier}

One of the questions that we shall regularly encounter throughout the study
of various classes of espaliers,\index{espalier} in Chapter~\ref{Ch:ClEsp},
will be what are the possible dimension ranges\index{dimension range} of
members of a given class of espaliers.

\begin{definition}\label{D:DUniv}
A class $\mathcal{E}$ of espaliers\index{espalier} is
\emph{D-universal},\index{D-universal|ii} if every continuous dimension scale\index{continuous dimension scale} admits a lower embedding\index{lower
embedding} into the dimension range\index{dimension range} of some member of
$\mathcal{E}$.
\end{definition}

We recall that the class of espaliers\index{espalier} is closed under
so-called \emph{lower subespaliers}, and also under \emph{direct products} of
espaliers, see Proposition~\ref{P:LSPEsp}.
The following lemma records some elementary facts about these notions. We
leave its easy proof to the reader.

\begin{lemma}\label{L:BasicDrng}\hfill
\begin{enumerate}
\item Let $K$ and $L$ be espaliers,\index{espalier} let $\varphi\colon K\to
L$ be a lower embedding\index{lower embedding} \pup{see
Lemma~\textup{\ref{L:LowEmbEsp}}}. Then the rule
$\DD_K(x)\mapsto\DD_L(\varphi(x))$ defines a lower embedding\index{lower
embedding} from
$\Drng K$\index{Dzzrng@$\Drng L$} into $\Drng L$.

\item Let $(L,\leq,\perp,\sim)$ be an espalier,\index{espalier} let $S$ be a
lower subset of\index{Dzzrng@$\Drng L$} $\Drng L$, put
   \[
   K=\setm{x\in L}{\DD_L(x)\in S}.
   \]
Then $K$ is a lower subespalier\index{espalier!lower sub ---} of $L$, and the
rule $\DD_K(x)\mapsto\DD_L(x)$ defines an isomorphism
from\index{Dzzrng@$\Drng L$} $\Drng K$ onto
$S$.

\item Let $(L_i)_{i\in I}$ be a family of espaliers,\index{espalier} let
$L=\prod_{i\in I}L_i$ be its direct product. Then the rule
$(\DD_{L_i}(x_i))_{i\in I}\mapsto\DD_L((x_i)_{i\in I})$ defines an isomorphism
from $\prod_{i\in I}\Drng L_i$\index{Dzzrng@$\Drng L$} onto $\Drng L$.
\end{enumerate}
\end{lemma}

In the context of Lemma~\ref{L:BasicDrng}(i), we shall of course write
\index{Dzzrnggf@$\Drng\varphi$} $\Drng\varphi\colon\Drng K\to\Drng L$ to
denote the map that sends $\DD_K(x)$ to $\DD_L(\varphi(x))$, for every
$x\in K$.

As a consequence of Lemma~\ref{L:BasicDrng}, the dimension
ranges\index{dimension range} of members of D-universal\index{D-universal}
classes of espaliers\index{espalier} can be nearly anything reasonable.

\begin{proposition}\label{P:DUniv}
Let $\mathcal{E}$ be a D-universal\index{D-universal} class of
espaliers.\index{espalier}
\begin{enumerate}
\item If every bounded lower subespalier\index{espalier!lower sub ---} of
every member of $\mathcal{E}$ belongs to $\mathcal{E}$, then every bounded
continuous dimension scale\index{continuous dimension scale} is isomorphic to the dimension
range\index{dimension range} of some bounded member of $\mathcal{E}$.

\item If every lower subespalier\index{espalier!lower sub ---} of every
member of $\mathcal{E}$ belongs to $\mathcal{E}$, then every continuous dimension scale\index{continuous dimension scale} is isomorphic to the dimension
range\index{dimension range} of some member of $\mathcal{E}$.
\end{enumerate}
\end{proposition}

\begin{proof}
(i) Let $S$ be a bounded continuous dimension scale,\index{continuous dimension scale}
denote by $\la$ the largest element of $S$. Since $\mathcal{E}$ is
D-universal,\index{D-universal} there exists
$L\in\mathcal{E}$ such that $S$ is (isomorphic to) a lower subset
of\index{Dzzrng@$\Drng L$}
$\Drng L$. Let $a\in L$ such that $\DD_L(a)=\la$, put $K=(a]$, a lower
subespalier\index{espalier!lower sub ---} of $L$. It follows from the
assumption and Lemma~\ref{L:BasicDrng}(i) that $K$ belongs to
$\mathcal{E}$ and $\Drng K$\index{Dzzrng@$\Drng L$} is isomorphic to $S$.
Observe that $K$ is bounded.

(ii) Let $S$ be a continuous dimension scale.\index{continuous dimension scale} Since
$\mathcal{E}$ is D-universal,\index{D-universal} there exists
$L\in\mathcal{E}$ such that $S$ is (isomorphic to) a lower subset
of\index{Dzzrng@$\Drng L$} $\Drng L$. Put $K=\setm{x\in L}{\DD_L(x)\in S}$.
It follows from the assumption and Lemma~\ref{L:BasicDrng}(ii) that $K$
belongs to $\mathcal{E}$ and\index{Dzzrng@$\Drng L$} $\Drng K$ is isomorphic
to $S$.
\end{proof}

The following result gives us a sufficient condition for
D-universality.\index{D-universal}

\begin{lemma}\label{L:DUniv}
Let $\mathcal{E}$ be a class of espaliers\index{espalier} satisfying the
following conditions:
\begin{enumerate}
\item $\mathcal{E}$ is closed under finite direct products.

\item For every ordinal $\gamma$ and every
complete Boolean space\index{Boolean space}
$\Omega$, there are $L_{\I}$, $L_{\II}$, $L_{\III}\in\mathcal{E}$ such
that $\CC(\Omega,\ZZ_\gamma)$ has a lower embedding\index{lower embedding}
into $\Drng L_{\I}$,\index{Dzzrng@$\Drng L$}
$\CC(\Omega,\RR_\gamma)$ has a lower embedding\index{lower embedding} into
$\Drng L_{\II}$,\index{Dzzrng@$\Drng L$} and
$\CC(\Omega,\two_\gamma)$ has a lower embedding\index{lower embedding} into
$\Drng L_{\III}$.\index{Dzzrng@$\Drng L$}
\index{Zzzgamma@$\ZZ_\gamma$}\index{Rzzgamma@$\RR_\gamma$}%
\index{Tzzgamma@$\two_\gamma$}%
\end{enumerate}
Then $\mathcal{E}$ is D-universal.\index{D-universal}
\end{lemma}

\begin{proof}
It follows from Lemma~\ref{L:BasicDrng} and assumptions (i), (ii) above
that for every ordinal $\gamma$ and any complete Boolean
spaces\index{Boolean space} $\Omega_{\I}$,
$\Omega_{\II}$, and
$\Omega_{\III}$, there exists $L\in\mathcal{E}$ such that the continuous dimension scale
\index{continuous dimension scale}%
\index{Zzzgamma@$\ZZ_\gamma$}\index{Rzzgamma@$\RR_\gamma$}%
\index{Tzzgamma@$\two_\gamma$}%
   \[
   \CC(\Omega_{\I},\ZZ_\gamma)\times
   \CC(\Omega_{\II},\RR_\gamma)\times\CC(\Omega_{\III},\two_\gamma)
   \]
embeds into\index{Dzzrng@$\Drng L$} $\Drng L$. The conclusion follows from
Theorem~\ref{T:EmbDimInt}.
\end{proof}

The results of Section~\ref{S:DUniv} will make it possible to prove further
results of D-universality.\index{D-universal}

\section{Existence of large constants}\label{S:LrgConst}

Taking account of the various examples of espaliers\index{espalier}
discussed in the Introduction, one is led to the conjecture that the
appearance of large cardinal values in the functional representation of the
dimension range\index{dimension range} of an espalier\index{espalier} should
be closely related to the existence of certain large orthogonal sums within
the espalier. (See also the proof of Lemma~\ref{L:LocM5}.)
Moreover, in the construction of espaliers of different types, we
will need to know what ingredients will ensure that the dimension
range\index{dimension range} of an example will be as large as desired. In
the present section, we provide some answers to the above questions.

\begin{quote}
\em Standing hypotheses: $L$ is an espalier,\index{espalier} $S$ is the
dimension range\index{dimension range} of $L$, and $\DD\colon
L\twoheadrightarrow S$ is the canonical map. Moreover, $\Omega$,
$\Omega_{\I}$, $\Omega_{\II}$,
$\Omega_{\III}$ are as in Section~\textup{\ref{S:ProjDF}}.

Let $\gamma$ be the ordinal and $\mu\colon S\to\CC(\Omega,\two_\gamma)$
\index{Tzzgamma@$\two_\gamma$}the
dimension function defined in Section~\textup{\ref{S:DimFctmu}}.

We put $\ol{S}=\CC(\Omega_{\I},\ZZ_\gamma;
\Omega_{\II},\RR_\gamma;\Omega_{\III},\two_\gamma)$.
\index{Zzzgamma@$\ZZ_\gamma$}\index{Rzzgamma@$\RR_\gamma$}%
\index{Tzzgamma@$\two_\gamma$}%
We pick a lower embedding\index{lower embedding} $\delta\colon
S_\fin\hookrightarrow\ol{S}$ as in Proposition~\textup{\ref{P:Gendelta}}. Let
$\varepsilon\colon S\hookrightarrow \ol{S}$ be the corresponding lower
embedding\index{lower embedding} defined in Section~\textup{\ref{S:EmbDI}}.
\end{quote}

\begin{lemma}\label{L:homogrem}
Let $a$, $b\in L$ be purely infinite elements with $a\lesssim b$.
Suppose that there is an infinite cardinal $\kappa$ such that
$\kappa\cdot b\lesssim b$ but $\kappa\cdot c\not\lesssim a$ for
all nonzero $c\in L$. Then $a\srem b$.
\end{lemma}

\begin{proof} Suppose $x\in L$ and $b\lesssim a\oplus x$. Because
of general comparability in~$L$ (Proposition~\ref{P:ShasGC} and comment
following), there is some $p\in\BB S$\index{pzzroj@$\BB{S}$} such that
$p\cdot a\lesssim p\cdot x$ and $p^\perp\cdot x\lesssim p^\perp\cdot a$. Now
   \[
   p^\perp\cdot b\lesssim (p^\perp\cdot a)\oplus (p^\perp\cdot x)
   \lesssim p^\perp\cdot(2\cdot a)\sim p^\perp\cdot a
   \]
by Lemmas~\ref{L:BaspcdotL} and \ref{L:p(aopplusb)}. Moreover, these lemmas
imply $\kappa\cdot(p^\perp\cdot b)\lesssim p^\perp\cdot b$, and so we
have $\kappa\cdot(p^\perp\cdot b) \lesssim a$. Our assumptions on
$a$ now imply that $p^\perp\cdot b=0$, and hence $b=p\cdot b$ by
Lemma~\ref{L:ppbotonL}.

Since $p\cdot a$ is purely infinite and $p\cdot a\lesssim p\cdot x$, we have
$(p\cdot a)\oplus (p\cdot x) \sim p\cdot x$, and so
   \[
   b= p\cdot b\lesssim (p\cdot a)\oplus (p\cdot x) \sim p\cdot x \leq x.
   \]
Therefore $a\srem b$.
\end{proof}

In many examples of espaliers,\index{espalier} the orthogonality relation
coincides with disjointness in the (partial) lattice:
\index{lattice!partial ---}
$a\perp b\Longleftrightarrow a\wedge b=0$. Let us abbreviate this condition
by the symbol $(\perp\,=\wedge0)$.

\begin{lemma}\label{L:orthocount}
Assume $(\perp\,=\wedge0)$. Let $x\in L$ be purely infinite, and
let $\eta$ be an infinite cardinal such that $x$ is not equal to
any orthogonal sum of more than $\eta$ nonzero elements. Let $y\in L$
and let $\beta\geq\eta$ be a cardinal number such that $\beta\cdot x\sim y$.

Then $y$ does not majorize any orthogonal sum of more than
$\beta$ nonzero elements. In particular, $\alpha\cdot
u\not\lesssim y$ for all $\alpha>\beta$ and all nonzero $u\in L$.
\end{lemma}

\begin{proof} By assumption, $y= \oplus_{i\in I} y_i$ with
$|I|=\beta$ and each $y_i\sim x$. Suppose that $y\geq\oplus_{j\in
J} z_j$ where $|J|>\beta$ and all $z_j\ne 0$. Since $z_j\wedge
y\ne 0$, we have $z_j\not\perp y$. Axiom~(L4) then yields a finite subset
$I_j\subset I$ such that $z_j\not\perp \oplus_{i\in I_j} y_i$. Since the set
of finite subsets of $I$ has cardinality $\beta$, the fibres of the map
$j\mapsto I_j$ cannot all have cardinality at most $\beta$. Hence,
there exist a subset $J'\subseteq J$ with $|J'|>\beta$ and a
finite subset $I'\subset I$ such that $I_j=I'$ for all $j\in J'$.
Thus, the element $y^*= \oplus_{i\in I'} y_i$ satisfies
$z_j\not\perp y^*$, and so
$z_j\wedge y^* \ne 0$, for all $j\in J'$,  because of
$(\perp\,=\wedge0)$. Consequently, $y^*$ majorizes an orthogonal
sum of more than $\beta$ nonzero elements, and after adjoining
an additional element if necessary, we may assume that
$y^*$ equals such an orthogonal sum. However, $y^*\sim x$ because
$x$ is purely infinite, and so Axiom~(L6) implies that $x$ is an
orthogonal sum of more than $\beta$ nonzero elements. This
contradicts our hypotheses. \end{proof}

\begin{lemma}\label{L:largedelta}
Assume $(\perp\,=\wedge0)$. Let $x\in L$ be purely infinite,
put $p=\cc(x)$, and let $\sigma$ be an ordinal such that $x$ is not equal to
any orthogonal sum of more than~$\aleph_\sigma$ nonzero
elements. Let $y\in L$ and let $\tau$ be an ordinal with
$\aleph_{\sigma+\tau}\cdot x\sim y$. Then $\scal{p}{\aleph_\tau}$
is defined, and $\scal{p}{\aleph_\tau}\leq\Delta(y)$.
\end{lemma}

\begin{proof} We proceed by induction on $\tau$.

Assume first that $\tau=0$. The set $X=\setm{a\in S|_\infty}{\cc(a)=p}$ is
nonempty, as it contains $\Delta(x)$. Since the element $\scal{p}{0}=0$ is
removable from any element of $L$, the element $\scal{p}{\aleph_0}$ is
defined as the least element of
$X$. Thus, $\scal{p}{\aleph_0}\leq\Delta(x)\leq\Delta(y)$.

Next, suppose that $\tau=\rho+1$ for some ordinal $\rho$. There is
some $z\leq y$ such that $\aleph_{\sigma+\rho}\cdot x\sim z$. By
induction, $\scal{p}{\aleph_\rho}$ is defined, and
$\scal{p}{\aleph_\rho}\leq\Delta(z)$. Now $\scal{p}{\aleph_\rho}=
\Delta(a)$ for some purely infinite $a\leq z$, and Lemma~\ref{L:orthocount}
shows that $\aleph_{\sigma+\tau}\cdot u\not\lesssim z$ for all nonzero $u\in
L$. On the other hand, since $\aleph_{\sigma+\tau}\cdot x\sim y$, we have
$\aleph_{\sigma+\tau}\cdot y\sim y$. Hence, $z\srem y$ by
Lemma~\ref{L:homogrem}, and so $a\srem y$. Therefore
$\scal{p}{\aleph_\rho} \rem \Delta(y)$. Since $\Delta(y)$ is a
purely infinite element with central cover $p$, it follows that
$\scal{p}{\aleph_\tau}$ is defined and majorized by $\Delta(y)$.

Finally, suppose that $\tau$ is a limit ordinal. For each ordinal
$\rho<\tau$, there exists $y_\rho\leq y$ such that
$\aleph_{\sigma+\rho}\cdot x\sim y_\rho$. By induction,
$\scal{p}{\aleph_\rho}$ is defined and $\scal{p}{\aleph_\rho}\leq
\Delta(y_\rho)\leq\Delta(y)$. Therefore $\scal{p}{\aleph_\tau}$ is defined,
and $\scal{p}{\aleph_\tau}= \bigvee_{\rho<\tau}
\scal{p}{\aleph_\rho}\leq\Delta(y)$.
\end{proof}

\begin{proposition}\label{P:existconstant}
Assume $(\perp\,=\wedge0)$. Let $x$, $y\in L$ be purely infinite elements
such that $\cc(x)=\cc(y)=1$, and let
$\sigma$, $\tau$ be ordinals, such that $x$ is not equal to
any orthogonal sum of more than $\aleph_\sigma$ nonzero elements, and
$\aleph_{\sigma+\tau}\cdot x\sim y$. Then the following statements hold:
\begin{enumerate}
\item $\mu(\Delta(y))(\fa)\geq\aleph_\tau$ for all $\fa\in\Omega$.
\item There exists a purely infinite element $u_\tau\in L$ such
that $\mu(\Delta(u_\tau))$ equals the constant function with value
$\aleph_\tau$.
\item Set $L_\tau= [0,u_\tau]\subseteq L$, and restrict $\le$,
$\perp$, $\sim$ from $L$ to $L_\tau$. Then $L_\tau$ is an
espalier,\index{espalier} and\index{Dzzrng@$\Drng L$} $\Drng L_\tau\cong
\CC(\Omega_{\I},\ZZ_\tau;
\Omega_{\II},\RR_\tau; \Omega_{\III},\two_\tau)$.
\end{enumerate}
\end{proposition}

\begin{proof} (i) In view of Lemma~\ref{L:largedelta},
$\scal{1}{\aleph_\tau}$ is defined and majorized by $\Delta(y)$.
Since $1\in\fa$ for all $\fa\in\Omega$, we get
$\mu(\Delta(y))(\fa)\geq\aleph_\tau$ for all $\fa$.

(ii) Because of (i), the lower embedding\index{lower embedding}
$\varepsilon: S
\hookrightarrow \ol{S}$ sends $\Delta(y)$ to a function
$f\in\ol{S}$ with $f(\fa)\geq\aleph_\tau$ for all $\fa\in\Omega$.
In particular, $\tau\leq\gamma$, and $\ol{S}$ contains the
constant function $t_\tau$ with $t_\tau(\fa)= \aleph_\tau$ for all
$\fa\in\Omega$. Since $\varepsilon$ is a lower embedding\index{lower
embedding}, there is some $w_\tau\in S$ such that $\varepsilon(w_\tau)=
t_\tau$. Note that $w_\tau$ is purely infinite, because $t_\tau$ is. Hence,
$\mu(w_\tau)= \varepsilon(w_\tau)= t_\tau$. It just
remains to note that $w_\tau= \Delta(u_\tau)$ for some purely
infinite element $u_\tau\in L$.

(iii) That $L_\tau$ is an espalier\index{espalier} follows from
Proposition~\ref{P:LSPEsp}(i). It is clear that\index{Dzzrng@$\Drng L$}
$\Drng L_\tau$ is isomorphic to the submonoid $S_\tau= [0,\Delta(u_\tau)]
\subseteq S$. Since
$\varepsilon$ is a lower embedding,\index{lower embedding} it maps $S_\tau$
isomorphically onto $\CC(\Omega_{\I},\ZZ_\tau;
\Omega_{\II},\RR_\tau; \Omega_{\III},\two_\tau)$.
\end{proof}

In case $(\perp\,=\wedge0)$ does not hold, it is not clear whether large
orthogonal sums are sufficient to imply large constants. For use in that
situation, we record the following more elementary approach.

\begin{lemma}\label{L:buildpkappa}
Let $\kappa$ be an infinite cardinal, and
$\famm{b_\xi}{\aleph_0\leq\xi\leq\kappa}$ a family of purely infinite
elements of $S$ \pup{indexed by infinite cardinals}.
Set $p=\cc(b_{\aleph_0})$. For all infinite cardinals
$\xi<\eta\leq\kappa$, assume that $b_\xi\leq b_\eta$ but $q(b_\eta)
\nleq q(b_\xi)$ for all nonzero projections $q\leq p$. Then
$\scal{p}{\kappa}$ is defined, and $\scal{p}{\kappa}\leq b_\kappa$.
\end{lemma}

\begin{proof} We show, by induction on $\xi$, that $\scal{p}{\xi}$
is defined and majorized by~$b_\xi$, for all infinite cardinals
$\xi\leq\kappa$. Since $b_{\aleph_0}$ is a purely infinite element
with central cover $p$, it is clear from the definition that
$\scal{p}{\aleph_0}$ is defined and
$\scal{p}{\aleph_0}\leq\nobreak b_{\aleph_0}$.

Next, suppose that $\xi$ is an infinite cardinal less than $\kappa$,
such that the element $a=\scal{p}{\xi}$ is defined and $a\leq b_\xi$.
Note that the elements $a$ and $p(b_{\xi^+})$ both have central cover
$p$. By assumption, $q(b_{\xi^+}) \nleq q(a)$ for all nonzero
projections $q\leq p$, whence Corollary~\ref{C:remequiv} implies that
$\scal{p}{\xi}= a \rem p(b_{\xi^+})$. Thus, $\scal{p}{\xi^+}$
is defined and
$\scal{p}{\xi^+}\leq p(b_{\xi^+})\leq b_{\xi^+}$.

Finally, if $\xi$ is a limit cardinal less than or equal to $\kappa$,
such that $\scal{p}{\eta}$ is defined and majorized by $b_\eta$ for all
infinite cardinals $\eta<\xi$, then $\scal{p}{\eta} \leq b_\xi$ for
all~$\eta$, whence $\scal{p}{\xi}= \bigvee_{\aleph_0\leq\eta<\xi}
\scal{p}{\eta}$ is defined and $\scal{p}{\xi} \leq b_\xi$.
\end{proof}

\begin{proposition}\label{P:secondexistconstant}
Let $\tau$ be an ordinal and $\famm{x_\alpha}{\alpha\leq\tau}$ a
family of purely infinite elements of $L$ with central cover $1$. For all
ordinals $\alpha<\beta\leq\tau$, assume that $x_\alpha\lesssim x_\beta$
but $q\cdot x_\beta \not\lesssim q\cdot x_\alpha$ for all nonzero
projections\index{pzzroj@$\BB{S}$} $q\in \BB S$. Then the following
statements hold:
\begin{enumerate}
\item $\mu(\Delta(x_\tau))(\fa)\geq\aleph_\tau$ for all $\fa\in\Omega$.
\item There exists a purely infinite element $u_\tau\in L$ such
that $\mu(\Delta(u_\tau))$ equals the constant function with value
$\aleph_\tau$.
\item Set $L_\tau= [0,u_\tau]\subseteq L$, and restrict $\le$,
$\perp$, $\sim$ from $L$ to $L_\tau$. Then $L_\tau$ is an
espalier,\index{espalier} and\index{Dzzrng@$\Drng L$} $\Drng L_\tau\cong
\CC(\Omega_{\I},\ZZ_\tau;
\Omega_{\II},\RR_\tau; \Omega_{\III},\two_\tau)$.
\end{enumerate}
\end{proposition}

\begin{proof}
(i) Set $b_{\aleph_\alpha}= \Delta(x_\alpha)$ for all ordinals
$\alpha\leq\tau$. Then
$\famm{b_\xi}{\aleph_0\leq \xi\leq \aleph_\tau}$ is a family of
purely infinite elements of $S$ with central cover $1$, such that for all
infinite cardinals $\xi<\eta\leq\aleph_\tau$, we have $b_\xi\leq b_\eta$
but $q(b_\eta) \nleq q(b_\xi)$ for all nonzero\index{pzzroj@$\BB{S}$} $q\in
\BB S$. Thus, by Lemma~\ref{L:buildpkappa}, $\scal{1}{\aleph_\tau}$ is
defined and
$\scal{1}{\aleph_\tau} \leq \Delta(x_\tau)$. Since $1\in\fa$ for all
$\fa\in\Omega$, we get
$\mu(\Delta(x_\tau))(\fa)\geq\aleph_\tau$ for all $\fa$.

(ii) and (iii) follow from (i) just as in Proposition \ref{P:existconstant}.
\end{proof}

\chapter{Classes of espaliers}\label{Ch:ClEsp}
\index{espalier}

\section{Abstract measure theory; Boolean espaliers}\label{S:AMeasTh}

\begin{definition}\label{D:BoolEsp}
An espalier $(L,\leq,\perp,\sim)$ is
\emph{Boolean},\index{espalier!Boolean ---|ii} if $L$ is a Boolean lattice
\index{lattice!Boolean ---}
and $x\perp y$ \iff\ $x\wedge y=0$, for all $x$, $y\in L$.
\end{definition}

Of course, the underlying Boolean algebra\index{Boolean algebra!complete ---}
of a Boolean espalier is complete. For a Boolean
algebra\index{Boolean algebra} $B$, we will denote by $\perp_B$ the canonical
orthogonality relation of $B$, that is,
$x\perp_By$ \iff\ $x\wedge y=0$, for all $x$, $y\in B$. We say that a family
$(a_i)_{i\in I}$ of elements of $B$ is \emph{disjoint}, if $a_i\wedge a_j=0$
for all $i\neq j$ in $I$, and then we let $\oplus_{i\in I}a_i$ denote its
join.

Many of the axioms defining the class of espaliers do not need checking in the
Boolean case.

\begin{proposition}\label{P:BoolEsp}
Let $B$ be a complete Boolean algebra\index{Boolean algebra},
let $\sim$ be a binary relation on
$B$. Then $(B,\leq_B,\perp_B,\sim)$ is a Boolean espalier \iff\
the following conditions hold:
\begin{itemize}
\item[(B0)] $x\sim 0$ implies that $x=0$, for all $x\in B$.
\index{bzzzerotwo@(B0--2)|ii}

\item[(B1)] The binary relation $\sim$ is unrestrictedly refining,
\index{unrestrictedly refining relation} that is,
for every $a\in B$ and every disjoint family $(b_i)_{i\in I}$ of elements of
$B$, if $a\sim\oplus_{i\in I}b_i$, then there exists a decomposition
$a=\oplus_{i\in I}a_i$ such that $a_i\sim b_i$ for all $i\in I$.

\item[(B2)] The binary relation $\sim$ is unrestrictedly additive,
\index{unrestrictedly additive!--- relation}
that is, for all disjoint families $(a_i)_{i\in I}$ and $(b_i)_{i\in I}$ of
elements of $B$, if $a_i\sim b_i$ for all $i\in I$, then
$\oplus_{i\in I}a_i\sim\oplus_{i\in I}b_i$.

\end{itemize}
\end{proposition}

We leave to the reader the straightforward proof of
Proposition~\ref{P:BoolEsp}.

Boolean espaliers can often be constructed from the following objects.

\begin{definition}\label{D:BpreEsp}
A \emph{Boolean pre-espalier}\index{espalier!Boolean pre ---|ii} is a pair
$(B,\sim)$, where $B$ is a Boolean algebra\index{Boolean algebra} and $\sim$
is an equivalence relation on $B$ satisfying Axioms~(B0) and (B1).
\end{definition}

We observe that the underlying Boolean algebra\index{Boolean algebra} of a
Boolean pre-espalier need not be complete. We recall (see, for example,
T. Jech \cite{Jech78})\index{Jech, T.} that for every Boolean algebra~$B$,
there exists a unique (up to isomorphism) complete Boolean algebra,
\index{Boolean algebra!complete ---} that we shall denote by
$\ol{B}$ and call the
\emph{completion} of $B$, such that $B$ is\index{dense} dense in~$\ol{B}$.
The following result makes it possible to extend to $\ol{B}$ any
Boolean pre-espalier structure on $B$.

\begin{lemma}\label{L:ExtPreEsp}
Let $(B,\sim)$ be a Boolean pre-espalier.
Define a binary relation $\sim^*$ on
$\ol{B}$ by the rule
   \begin{multline}\label{Eq:Defsim*}
   x\sim^*y\Longleftrightarrow\text{there are are decompositions }
   x=\oplus_{i\in I}x_i,\ y=\oplus_{i\in I}y_i\\
   \text{ such that }
   x_i,\,y_i\in B\text{ and }x_i\sim y_i,\text{ for all }i\in I,
   \end{multline}
\index{xzzsim*y@$x\sim^*y$|ii}%
for all $x$, $y\in\ol{B}$.
Then $(\ol{B},\sim^*)$ is a Boolean espalier.\index{espalier!Boolean ---}
Furthermore, $\sim^*$ is the smallest equivalence relation $\sim'$ on
$\ol{B}$ containing $\sim$ such that
$(\ol{B},\sim')$ is an espalier.\index{espalier}
\end{lemma}

We shall call $\sim^*$ the \emph{espalier closure}
\index{espalier!--- closure} of $\sim$.

\begin{proof}
It is clear that every equivalence relation $\sim'$ on $\ol{B}$ containing
$\sim$, such that $(\ol{B},\sim')$ is an espalier,\index{espalier} also
contains $\sim^*$, hence it suffices to prove that
$(\ol{B},\sim^*)$ is an espalier.\index{espalier}

Since every element of $\ol{B}$ can be written $\oplus_{i\in I}x_i$, where
all the $x_i$-s belong to $B$, the binary relation $\sim^*$ is reflexive. It
is obviously symmetric. Now let $a$, $b$, $c\in\ol{B}$ such that $a\sim^*b$
and $b\sim^*c$. There are decompositions of the form
   \[
   a=\oplus_{i\in I}a_i,\ b=\oplus_{i\in I}b'_i=\oplus_{j\in J}b''_j,\
   c=\oplus_{j\in J}c_j,
   \]
with $a_i\sim b'_i$ in $B$, for all $i\in I$, and $b''_j\sim c_j$, for all
$j\in J$. For any $i\in I$, $a_i\sim b'_i=\oplus_{j\in J}(b'_i\wedge b''_j)$,
thus, since $\sim$ satisfies (B1), there exists a decomposition
$a_i=\oplus_{j\in J}a_{i,j}$ with $a_{i,j}\sim b'_i\wedge b''_j$, for all
$j\in J$. For $j\in J$, since
$c_j\sim b''_j=\oplus_{i\in I}(b'_i\wedge b''_j)$ and by (B1), there exists a
decomposition $c_j=\oplus_{i\in I}c_{i,j}$ such that
$b'_i\wedge b''_j\sim c_{i,j}$, for all $i\in I$. Therefore,
$a_{i,j}\sim c_{i,j}$, for all $(i,j)\in I\times J$, and
$a=\oplus_{(i,j)\in I\times J}a_{i,j}$ and
$c=\oplus_{(i,j)\in I\times J}c_{i,j}$; whence $a\sim^*c$. Therefore, $\sim^*$
is an equivalence relation on $\ol{B}$. It is obvious that $\sim^*$ satisfies
(B0).

Now let $a\sim^*\oplus_{i\in I}b_i$ in $\ol{B}$. By definition, there are
decompositions $a=\oplus_{j\in J}a'_j$ and
$\oplus_{i\in I}b_i=\oplus_{j\in J}b'_j$ such that $a'_j\sim b'_j$, for all
$j\in J$. For $j\in J$, since $a'_j\sim b'_j=\oplus_{i\in I}(b_i\wedge b'_j)$,
there exists a decomposition $a'_j=\oplus_{i\in I}a_{i,j}$ such that
$a_{i,j}\sim b_i\wedge b'_j$, for all $i\in I$. Observe that
$a=\oplus_{(i,j)\in I\times J}a_{i,j}$; put $a_i=\oplus_{j\in J}a_{i,j}$, for
all $i\in I$. Thus $a=\oplus_{i\in I}a_i$, and, by the definition of
$\sim^*$, $a_i\sim^*\oplus_{j\in J}(b_i\wedge b'_j)=b_i$, for all $i\in I$.
Therefore, $\sim^*$ satisfies (B1).

Finally let $a=\oplus_{i\in I}a_i$ and $b=\oplus_{i\in I}b_i$ with
$a_i\sim^*b_i$, for all $i\in I$. By definition, for all $i\in I$, there are
decompositions $a_i=\oplus_{j\in J_i}a_{i,j}$ and
$b_i=\oplus_{j\in J_i}b_{i,j}$ such that $a_{i,j}\sim b_{i,j}$, for all
$i\in I$ and all $j\in J_i$. Put $J=\bigcup_{i\in I}(\set{i}\times J_i)$,
then $a=\oplus_{(i,j)\in J}a_{i,j}$ and $b=\oplus_{(i,j)\in J}b_{i,j}$,
whence $a\sim^*b$. Therefore, $\sim^*$ satisfies (B2).
\end{proof}

For a Boolean espalier\index{espalier!Boolean ---} $(B,\sim)$ and a set
$I$, we let the permutation group $\SS_I$\index{SzSI@$\SS_I$|ii}
of~$I$ act on the Boolean algebra
\index{Boolean algebra}
$B^I$ by translation: namely,
   \[
   (\sigma x)(i)=x(\sigma^{-1}(i)),\text{ for all }x\in B^I\text{ and all }
   \sigma\in\SS_I.
   \]
Next, let $\sim_I$ be the equivalence relation on $B^I$ associated
with this action and $\sim$, that is,
   \[
   x\sim_Iy\Leftrightarrow\exists\sigma\in\SS_I\text{ such that }
   y(i)\sim x(\sigma(i)),\text{ for all }x,\,y\in B^I.
   \]
Since $\SS_I$ acts on $B^I$ by automorphisms (of the Boolean algebra
\index{Boolean algebra!complete ---} $B^I$)
and $(B,\sim)$ is an espalier,\index{espalier} it is easy to
see that $(B^I,\sim_I)$ is a
Boolean pre-espalier.\index{espalier!Boolean pre ---} Since $B^I$ is already
complete, the espalier closure\index{espalier!--- closure} of
$(B^I,\sim_I)$ is an equivalence relation on $B^I$, that we shall denote by
$\sim^I$. For $i\in I$, we denote by $\varphi_i\colon B\hookrightarrow B^I$
the canonical map, that is, $\varphi_i(x)(j)$ is equal to $x$ if $i=j$, to
$0$ otherwise, for all $x\in B$ and all $j\in J$.

\begin{lemma}\label{L:Basicphii}
The following statements hold, for any $i\in I$:
\begin{enumerate}
\item The map $\varphi_i$ is a lower embedding\index{lower embedding} of
espaliers.\index{espalier}

\item The map\index{Dzzrnggf@$\Drng\varphi$} $\Drng\varphi_i$ is a lower
embedding\index{lower embedding} with dense\index{dense} image from
$B/{\sim}$ into $B^I/{\sim^I}$.
\end{enumerate}
\end{lemma}

\begin{proof}
(i) It is obvious that $\varphi_i$ is a $(\leq,\perp)$-isomorphism from $B$
onto a lower subset of $B^I$, and that $x\sim y$ implies that
$\varphi_i(x)\sim^I\varphi_i(y)$, for all $x$, $y\in B$. Now suppose that
$\varphi_i(x)\sim^I\varphi_i(y)$. There are decompositions of the form
$x=\oplus_{j\in J}x_j$ and $y=\oplus_{j\in J}y_j$ in $B$ such that
$\varphi_i(x_j)\sim_I\varphi_i(y_j)$, for all $j\in J$. Hence $x_j\sim y_j$,
for all $j\in J$, whence, since $(B,\sim)$ satisfies (B2), $x\sim y$.
This completes the proof of~(i).

Let $\la\in B^I/{\sim^I}$ be nonzero; so $\la=\DD(a)$, for some
$a\in B^I\setminus\set{0}$. Since $a$ has a nonzero component, there are
$b\in B\setminus\set{0}$ and $\sigma\in\SS_I$ such that
$\sigma\varphi_i(b)\leq a$. Hence,\index{Dzzrnggf@$\Drng\varphi$}
   \[
   0<(\Drng\varphi_i)(\DD(b))=\DD(\varphi_i(b))=\DD(\sigma\varphi_i(b))
   \leq\DD(a)=\la,
   \]
which completes the proof of (ii).
\end{proof}

For a Boolean algebra\index{Boolean algebra} $B$, we define a
cardinal number $\wdt(B)$\index{wzzdt@$\wdt(B)$|ii} by
\index{antichain}
   \[
   \wdt(B)=\sup\setm{|X|}{X\text{ is an antichain of }B}.
   \]
Furthermore, for a set $I$, we put\index{szzupp@$\supp(x)$|ii}
   \[
   \supp(x)=\setm{i\in I}{x(i)\neq 0},\text{for all }x\in B^I.
   \]
The coming set of lemmas, from \ref{L:PresCard} to \ref{L:RemGExt}, is aimed
at constructing Boolean espaliers\index{espalier!Boolean ---} whose
dimension ranges\index{dimension range} have large constants. Instead of
accomodating the results of Section~\ref{S:LrgConst} to the present context,
we propose direct proofs, probably of more interest to the Boolean
algebra-oriented reader.

The following lemma expresses, essentially, the well-known fact that
large enough cardinals are preserved from the ground universe $V$ to the
Boolean-valued\index{VzzVB@$V^B$} universe $V^B$.

\begin{lemma}\label{L:PresCard}
Let $(B,\sim)$ be a Boolean espalier,\index{espalier!Boolean ---} let
$I$ be a set. If $x\sim^Iy$, then $|\supp(y)|\leq|\supp(x)|\cdot\wdt(B)$,
for all $x$,
$y\in B^I$.
\end{lemma}

\begin{proof}
Put $\kappa=|\supp(x)|\cdot\wdt(B)$. Since $x\sim^Iy$, there are
decompositions of the form $x=\oplus_{j\in J}x_j$ and $y=\oplus_{j\in J}y_j$
such that $x_j\sim_Iy_j$ for all $j\in J$. By decomposing further the
$x_j$-s and the $y_j$-s as disjoint sums of elements of
$D=\bigcup_{i\in I}\varphi_i[B\setminus\set{0}]$, we may assume, without loss
of generality, that both $x_j$ and $y_j$ belong to $D$, for all $j\in J$.
Hence, for all $j\in J$, there are $a(j)$, $b(j)\in I$ and $\ol{x}_j$,
$\ol{y_j}\in B\setminus\set{0}$ such that $x_j=\varphi_{a(j)}(\ol{x}_j)$ and
$y_j=\varphi_{b(j)}(\ol{y}_j)$. We observe that
   \[
   \supp(x)=\setm{a(j)}{j\in J}\text{ and }\supp(y)=\setm{b(j)}{j\in J}.
   \]
Put $J_i=\setm{j\in J}{a(j)=i}$, for all $i\in\supp(x)$. Since the family
$(x_j)_{j\in J_i}$ is the image under $\varphi_i$ of the
antichain\index{antichain}
$(\ol{x}_j)_{j\in J_i}$ of $B$, the inequality
$|J_i|\leq\wdt(B)$ holds. Therefore,
   \begin{equation}
   |\supp(y)|\leq|J|=\sum_{i\in\supp(x)}|J_i|\leq\kappa.\tag*{\qed}
   \end{equation}
\renewcommand{\qed}{}
\end{proof}

For a Boolean algebra\index{Boolean algebra} $B$, a set $I$,
and a subset $X$ of $I$, we put
$X\cdot 1=(x_i)_{i\in I}$, where $x_i=1_B$ if $i\in X$ and $x_i=0_B$ if
$i\notin X$; hence $X\cdot 1\in B^I$.

\begin{lemma}\label{L:X1simY1}
Let $(B,\sim)$ be a Boolean espalier,\index{espalier!Boolean ---} let $I$ be
a set. Then $|X|=|Y|$ implies that $X\cdot 1\sim^IY\cdot 1$, for all $X$,
$Y\subseteq I$.
\end{lemma}

\begin{proof}
Let $\sigma\colon X\twoheadrightarrow Y$ be a bijection. If $X$ is finite,
then $\sigma$ can be extended to a permutation $\tau$ of $I$, thus
$X\cdot 1\sim^I\tau(X\cdot 1)=Y\cdot 1$.
Suppose now that $X$ is infinite.
There exists a partition $X=X_0\sqcup X_1$ of $X$ such that
$|X_0|=|X_1|=|X|$. Put $Y_i=\sigma[X_i]$, for $i<2$. Then the restriction of
$\sigma$ from $X_i$ onto $Y_i$ can be extended to a permutation $\tau_i$ of
$I$, for all $i<2$. Therefore,
   \begin{equation}
   X\cdot 1=(X_0\cdot 1)\oplus(X_1\cdot 1)\sim^I
   \tau_0(X_0\cdot 1)\oplus\tau_1(X_1\cdot 1)=(Y_0\cdot 1)\oplus(Y_1\cdot 1)
   =Y\cdot 1.\tag*{\qed}
   \end{equation}
\renewcommand{\qed}{}
\end{proof}

The following lemma makes it possible to find Boolean
espaliers\index{espalier!Boolean ---} with long
$\srem$-chains.

\begin{lemma}\label{L:RemGExt}
Let $(B,\sim)$ be a Boolean espalier,\index{espalier!Boolean ---} let
$\alpha$, $\beta$, $\fm$ be infinite cardinals such that
$\wdt(B)\leq\alpha<\beta\leq\fm$. Then the relation
$\alpha\cdot 1\srem\beta\cdot 1$ holds in the espalier\index{espalier}
$(B^\fm,\sim^\fm)$.
\end{lemma}

\begin{proof}
Let $x$, $y\in B^\fm$ such that $\beta\cdot 1=x\oplus y$ and
$x\sim^\fm\alpha\cdot 1$, we prove that $y\sim^\fm\beta\cdot 1$. It follows
from Lemma~\ref{L:PresCard} that $|\supp(x)]\leq\alpha$, thus, putting
$Y=\beta\setminus\supp(x)$, we obtain that $|Y|=\beta$ and $Y\cdot 1\leq y$.
{}From Lemma~\ref{L:X1simY1} it follows that $Y\cdot 1\sim\beta\cdot 1$,
whence $\beta\cdot 1\lesssim y$. Since $y\leq\beta\cdot 1$ and by
Lemma~\ref{P:SchBern}, it follows that $y\sim^\fm\beta\cdot 1$.
\end{proof}

\begin{lemma}\label{L:LotsofI}
Let $\Omega$ be a complete Boolean
space,\index{Boolean space} let $\gamma$ be an ordinal. Then there exists a
Boolean espalier
\index{espalier!Boolean ---}
$(B,\sim)$ such that $B/{\sim}\cong\CC(\Omega,\ZZ_\gamma)$.
\index{Zzzgamma@$\ZZ_\gamma$}
\end{lemma}

\begin{proof}
Denote by $D$ the complete Boolean algebra\index{Boolean algebra} of clopen
subsets of $\Omega$. Let~$\theta$ be an ordinal such that
$\wdt(D)\leq\aleph_\theta$, put
$\fm=\omega_{\theta+\gamma}$, endow the direct power $D^\fm$ with the
previously introduced $\sim^\fm$ defined from the espalier $(D,=)$. We
consider the map $\varphi_0\colon D\hookrightarrow D^\fm$
introduced earlier. It follows from Lemma~\ref{L:Basicphii} that $\varphi_0$
is a lower embedding\index{lower embedding} of espaliers\index{espalier} and
$\Drng\varphi_0$\index{Dzzrnggf@$\Drng\varphi$} is a lower
embedding\index{lower embedding} of continuous dimension scales\index{continuous dimension scale} with dense\index{dense} image. In
particular,\index{pzzroj@$\BB{S}$}
$\BB{D^\fm}\cong\BB{D}\cong D$. Since~$\Omega$ is isomorphic to the
ultrafilter space of $D$, it follows from Theorems~\ref{T:EmbDimInt} and
\ref{T:DimEsp} that $D^\fm/{\sim^\fm}$ has a lower embedding\index{lower
embedding} into an espalier
\index{espalier} of the form
   \[
   \ol{S}=
\CC(\Omega_{\I},\ZZ_\alpha;\Omega_{\II},\RR_\alpha;\Omega_{\III},\two_\alpha),
   \]
for some ordinal $\alpha$ and a partition
$\Omega=\Omega_{\I}\sqcup\Omega_{\II}\sqcup\Omega_{\III}$ of $\Omega$ into
clopen sets. However, the continuous dimension scale\index{continuous dimension scale}
$D\cong\CC(\Omega,\set{0,1})$ has a lower embedding\index{lower embedding}
into $D^\fm/{\sim^\fm}$, thus
$\Omega_{\II}=\Omega_{\III}=\es$ and $\Omega_{\I}=\Omega$. Finally, it
follows from Lemma~\ref{L:RemGExt} that the $(\gamma+1)$-sequence
$(\DD(\omega_{\theta+\xi}\cdot 1))_{\xi\leq\gamma}$ is $\rem$-increasing in
$D^\fm/{\sim^\fm}$, but all the members of this sequence have central
cover $1$, thus the image of $D^\fm/{\sim^\fm}$ in $\ol{S}$ contains a
function whose values are all above $\aleph_\gamma$, in particular,
$\CC(\Omega,\ZZ_\gamma)$
\index{Zzzgamma@$\ZZ_\gamma$}
has a lower embedding\index{lower embedding} into $D^\fm/{\sim^\fm}$. The
argument of Proposition~\ref{P:DUniv} shows then that there exists a bounded
lower subespalier\index{espalier!lower sub ---}
$(B,\sim)$ of $(D^\fm,\sim^\fm)$ such that $B/{\sim}$ is
isomorphic to $\CC(\Omega,\ZZ_\gamma)$.
\end{proof}

\begin{lemma}\label{L:LotsofII}
Let $\Omega$ be a complete Boolean
space,\index{Boolean space} let $\gamma$ be an ordinal. Then there exists a
Boolean espalier
\index{espalier!Boolean ---}
$(B,\sim)$ such that $B/{\sim}\cong\CC(\Omega,\RR_\gamma)$.
\index{Rzzgamma@$\RR_\gamma$}
\end{lemma}

\begin{proof}
The proof of Lemma~\ref{L:LotsofII} requires some familiarity with
forcing and
\index{Boolean algebra!complete ---} complete Boolean algebras, in
particular, the random real extension and two-step iterated forcing, see
\index{Jech, T.}\cite{Jech78}. We denote by $\Bo$\index{Bzzo@$\Bo$|ii}
the Boolean algebra of all Borel subsets
of the real unit interval $[0,1]$ modulo null sets, the \emph{random
algebra},\index{random algebra|ii}
and by $\mm\colon\Bo\to[0,1]$ the Lebesgue measure on $\Bo$.
Furthermore, let $x\sim y$ hold, if $\mm(x)=\mm(y)$, for all $x$, $y\in\Bo$.
Let $\alpha=\sum_{i\in I}\alpha_i$, for a family $(\alpha_i)_{i\in I}$ of
elements of $[0,1]$, mean that $\alpha$ is the supremum over all finite
subsets $J$ of $I$ of $\sum_{i\in J}\alpha_i$.
We need a couple of standard facts on the measure $\mm$, summed up in the
following claims.
\goodbreak

\setcounter{claim}{0}
\begin{claim}\label{Cl:1}\hfill
\begin{enumerate}
\item $\mm(\oplus_{i\in I}x_i)=\sum_{i\in I}\mm(x_i)$, for any disjoint
family $(x_i)_{i\in I}$ of elements of~$\Bo$.

\item Let $x\in\Bo$ and let $(\alpha_i)_{i\in I}$ be a family of elements of
$[0,1]$. If $\mm(x)=\sum_{i\in I}\alpha_i$, then there exists a decomposition
$x=\oplus_{i\in I}x_i$ in $\Bo$ such that $\mm(x_i)=\alpha_i$, for all
$i\in I$.
\end{enumerate}
\end{claim}

\begin{cproof}
(i) Observe that the assumptions imply that the set $\setm{i\in I}{x_i>0}$ is
countable; the conclusion follows from countable additivity of Lebesgue
measure.

(ii) Again, $\setm{i\in I}{\alpha_i>0}$ is countable, so we may assume
without loss of generality that $I=\omega$. It is then easy to construct
inductively a nondecreasing sequence $(a_i)_{i<\omega}$ of elements of
$[0,1]$ satisfying the conditions $\mm(x\cap[0,a_0])=\alpha_0$ and
$\mm(x\cap[a_n,a_{n+1}])=\alpha_{n+1}$, for all $n<\omega$. Put
$a=\sup_{n<\omega}a_n$, then $x_0=x\setminus[a_0,a]$ and
$x_n=x\cap[a_n,a_{n+1}]$, for all $n<\omega$, satisfy the desired conclusion.
\end{cproof}

As an immediate corollary, we obtain the following claim.

\begin{claim}\label{Cl:2}
The structure $(\Bo,\sim)$ is a Boolean
espalier,\index{espalier!Boolean ---} with dimension
range\index{dimension range} isomorphic to $[0,1]$.
\end{claim}

Now we work under the assumptions of Lemma~\ref{L:LotsofII}. Denote by $C$
the complete Boolean algebra
\index{Boolean algebra!complete ---}
of clopen subsets of $\Omega$. We consider the
quotient of the Scott-Solovay $C$-valued universe\index{VzzVB@$V^B$} $V^C$
of set theory under the equivalence relation that identifies names $\dx$ and
$\dy$ \iff\
$\bv{\dx=\dy}=1$. We still denote by $V^C$ the quotient, endowed with its
natural Boolean value, so, now, $\bv{\lx=\ly}=1$ \iff\ $\lx=\ly$, for all
$\lx$, $\ly\in V^C$. Furthermore, let $\lBo$\index{bzzlo@$\lBo$|ii}
denote the (equivalence
class of the) $C$-valued name for $\Bo$, and put $D=C*\lBo$, the two-step
iterated forcing of $C$ by the random algebra of $V^C$. Hence $D$ is the
complete Boolean algebra
\index{Boolean algebra!complete ---}
of all $\lx\in V^C$\index{VzzVB@$V^B$} such that $\bv{\lx\in\lBo}=1$,
the partial ordering $\leq$ being defined by $\lx\leq\ly$ \iff\
$\bv{\lx\leq\ly}=1$. Hence, orthogonality in $D$ is defined by $x\perp y$
\iff\ $\bv{\lx\wedge\ly=0}=1$. Let $\sim$ be the binary relation defined on
$D$ by $\lx\sim\ly$ \iff\ $\bv{\lx\sim\ly}=1$, for all $\lx$, $\ly\in D$.

\begin{claim}\label{Cl:3}
The structure $(D,\sim)$ is a Boolean espalier.\index{espalier!Boolean ---}
\end{claim}

\begin{cproof}
It is obvious that $\sim$ satisfies (B0). Now let
$\la\sim\lb$ in~$D$, with $\lb$ decomposed as $\lb=\oplus_{i\in I}\lb_i$.
If $p\mapsto p^*$ denotes the canonical embedding from $C$ into $C*\lBo$, the
relation $p\leq\bv{\lx\leq\ly}$ is equivalent to $p^*\wedge\lx\leq\ly$, for
all $\lx$, $\ly\in D$. It is an easy exercise to deduce from this the relation
$\bv{\lb=\oplus_{i\in\check I}\lb_i}=1$, where the symbol $\check I$ denotes
the canonical name in $V^C$\index{VzzVB@$V^B$} for $I$. Moreover,
   \[
   \bv{\la\sim\lb}=\bv{(\lBo,\sim)\text{ is a Boolean espalier}}=1
   \]
and $V^C$\index{VzzVB@$V^B$} is a Boolean-valued\index{VzzVB@$V^B$} model of
set theory, in particular, $V^C$\index{VzzVB@$V^B$} satisfies
Claim~\ref{Cl:2}. Therefore,
   \[
   \bv{\exists(\lx_i)_{i\in\check I}\text{ such that }
   \la=\oplus_{i\in\check I}\lx_i\text{ and }\lx_i\sim\lb_i,\text{ for all }
   i\in\check I}=1.
   \]
Since $V^C$\index{VzzVB@$V^B$} is full and the notion of function is
absolute, there exists a family $(\la_i)_{i\in I}$ of elements of $D$ such
that
$\bv{\la_i\sim\lb_i}=1$, for all $i\in I$, and
   \[
   \bv{\la=\oplus_{i\in\check I}\la_i}=1.
   \]
Hence $\la=\oplus_{i\in I}\la_i$ and $\la_i\sim\lb_i$, for all $i\in I$.
\end{cproof}

We shall now identify the dimension range\index{dimension range} of
$(D,\sim)$. For every $\la\in D$, there exists a unique
$\lga\in V^C$\index{VzzVB@$V^B$} such that $\bv{\lga=\mm(\la)}=1$. Observe,
in particular, that
$\bv{0\leq\lga\leq 1}=1$. We put $\lga=\eps(\la)$.

We need the following standard fact.

\begin{claim}\label{Cl:4}
There exists an isomorphism of \pcm s\index{partial commutative monoid} from
$\CC(\Omega,[0,1])$ onto the set of $\lga\in V^C$\index{VzzVB@$V^B$} such
that
$\bv{0\leq\lga\leq 1}=1$, endowed with the addition defined by
$\lgc=\lga+\lgb$ \iff\ $\bv{\lgc=\lga+\lgb}=1$.
\end{claim}

\begin{cproof}
This is a particular case of a much more general statement, see, for example,
\cite[Theorem~3.10]{Wehr93}.\index{Wehrung, F.}
\end{cproof}

{}From now on we identify the elements of $\CC(\Omega,[0,1])$ with the
$C$-valued names of elements of $[0,1]$, \emph{via} the isomorphism of
Claim~\ref{Cl:4}. We obtain the following result.

\begin{claim}\label{Cl:5}
The map $\eps$ can be factored through $\sim$, to an isomorphism from
$D/{\sim}$ onto $\CC(\Omega,[0,1])$.
\end{claim}

\begin{cproof}
For $\la$, $\lb\in D$, $\eps(\la)=\eps(\lb)$ \iff\
$\bv{\mm(\la)=\mm(\lb)}=\nobreak1$, \iff\ $\bv{\la\sim\lb}=1$, \iff\
$\la\sim\lb$. If
$\lc\in D$ and $\lc=\la\oplus\lb$, then $\bv{\lc=\la\oplus\lb}=1$, thus
$\bv{\mm(\lc)=\mm(\la)+\mm(\lb)}=1$, thus $\eps(\lc)=\eps(\la)+\eps(\lb)$.
Finally, let $\lga\in\CC(\Omega,[0,1])$. There exists a unique
$\la\in V^C$\index{VzzVB@$V^B$} such that $\bv{\la=[0,\lga]}=1$, thus
$\la\in D$ and
$\eps(\la)=\lga$, so $\eps$ is surjective.
\end{cproof}

The rest of the proof proceeds like in the proof of Lemma~\ref{L:LotsofI}, of
which we shall keep the notation. It follows from Claim~\ref{Cl:5} that the
Boolean algebra
\index{Boolean algebra}
of projections of $(D,\sim)$ is isomorphic to $C$, see
Claim~\ref{Cl:ProjPU} in the proof of Theorem~\ref{T:C(O,K)DimInt}. If
$\theta$ is an ordinal such that $\wdt(D)\leq\aleph_\theta$ and we put
$\fm=\omega_{\theta+\gamma}$, then $D^\fm/{\sim^\fm}$ has a $\rem$-increasing
chain of length $\gamma+1$. Furthermore, observe that this time, since
$\CC(\Omega,[0,1])$ has a lower embedding\index{lower embedding} into
$D^\fm/{\sim^\fm}$,
$\Omega_{\I}=\Omega_{\III}=\es$ while $\Omega_{\II}=\Omega$, and the image of
$D^\fm/{\sim^\fm}$ in $\ol{S}$ contains a function with values
above~$\aleph_\gamma$. Hence, there exists a bounded lower subespalier
\index{espalier!lower sub ---}
$(B,\sim)$ of $(D^\fm,\sim^\fm)$ such that
$B/{\sim}\cong\CC(\Omega,\RR_\gamma)$.
\index{Rzzgamma@$\RR_\gamma$}
\end{proof}

\begin{lemma}\label{L:LotsofIII}
Let $\Omega$ be a complete Boolean
space,\index{Boolean space} let $\gamma$ be an ordinal. Then there exists a
Boolean espalier
\index{espalier!Boolean ---}
$(B,\sim)$ such that $B/{\sim}\cong\CC(\Omega,\two_\gamma)$.
\index{Tzzgamma@$\two_\gamma$}
\end{lemma}

\begin{proof}
In order to give a proof of Lemma~\ref{L:LotsofIII}, it is also convenient to
be familiar with forcing and
\index{Boolean algebra!complete ---}
complete Boolean algebras. We denote by $\Co$\index{Czzo@$\Co$|ii}
the Boolean algebra \index{Boolean algebra}
of all Borel subsets of the Cantor space
$\set{0,1}^\omega$ modulo meager sets, the \index{Cohen algebra|ii}
\emph{Cohen algebra}. Furthermore,
for $x\in\Co$, we define $\nn(x)=\aleph_0$ if $x>0$ while $\nn(0)=0$.
Let $x\sim y$ hold, if $\nn(x)=\nn(y)$, for all $x$, $y\in\Co$.

\setcounter{claim}{0}
\begin{claim}\label{Cl:1b}
The structure $(\Co,\sim)$ is a Boolean
espalier,\index{espalier!Boolean ---} with dimension range
\index{dimension range} isomorphic to
$\set{0,\aleph_0}$.
\end{claim}

\begin{cproof}
Every nonzero element of $\Co$ can be decomposed as a disjoint union of
two (resp., $\omega$) nonzero elements of $\Co$. Furthermore, if
$u=\oplus_{i\in I}u_i$ in $\Co$, then  $\setm{i\in I}{u_i>0}$ is countable.
It follows easily that $\sim$ satisfies (B1). It obviously satisfies (B0) and
(B2). Since every nonzero element of $\Co$ can be decomposed as a disjoint
union of two nonzero elements of $\Co$, every element of $\Co$ is purely
infinite. It follows that $\Co/{\sim}\cong\set{0,\aleph_0}$.
\end{cproof}

Now we work under the assumptions of Lemma~\ref{L:LotsofIII}. Denote by $C$
the complete Boolean algebra
\index{Boolean algebra!complete ---}
of clopen subsets of $\Omega$. We define the
(quotiented) Scott-Solovay universe $V^C$\index{VzzVB@$V^B$} of set theory
as in the proof of Lemma~\ref{L:LotsofII}. Furthermore, let
$\lCo$\index{czzlo@$\lCo$|ii} denote the (equivalence class of the)
$C$-valued name for $\Co$, and put
$D=C*\lCo$, the two-step iterated forcing of $C$ by the Cohen algebra
of\index{VzzVB@$V^B$}
$V^C$. Hence~$D$ is the complete Boolean algebra
\index{Boolean algebra!complete ---}
of all $\lx\in V^C$\index{VzzVB@$V^B$} such that $\bv{\lx\in\lCo}=1$,
the partial ordering $\leq$ being defined by $\lx\leq\ly$ \iff\
$\bv{\lx\leq\ly}=1$. Hence, orthogonality in~$D$ is defined by $x\perp y$
\iff\ $\bv{\lx\wedge\ly=0}=1$. Let $\sim$ be the binary relation defined on
$D$ by $\lx\sim\ly$ \iff\ $\bv{\lx\sim\ly}=1$, for all $\lx$, $\ly\in D$.

The proof of the following Claim~\ref{Cl:2b} is, \emph{mutatis
mutandis}, the same as the one for Claim~\ref{Cl:3} in the
proof of Lemma~\ref{L:LotsofII}.

\begin{claim}\label{Cl:2b}
The structure $(D,\sim)$ is a Boolean espalier.\index{espalier!Boolean ---}
\end{claim}

We shall now identify the dimension range\index{dimension range} of
$(D,\sim)$. For every
$\la\in D$, there exists a unique $\lga\in V^C$\index{VzzVB@$V^B$} such that
$\bv{\lga=\nn(\la)}=1$. Observe, in particular, that
$\bv{\lga\in\set{0,\aleph_0}}=1$. We put $\lga=\eps(\la)$.

The analogue of Claim~\ref{Cl:4} of Lemma~\ref{L:LotsofII} takes the
following simple form, with a much more direct proof.

\begin{claim}\label{Cl:3b}
There exists an isomorphism of \pcm s\index{partial commutative monoid} from
$\CC(\Omega,\set{0,\aleph_0})$
\pup{isomorphic to $C$} onto the set of $\lga\in V^C$\index{VzzVB@$V^B$}
such that
$\bv{\lga\in\set{0,\aleph_0}}=\nobreak1$, endowed with the addition defined by
$\lgc=\lga+\lgb$ \iff\ $\bv{\lgc=\lga+\lgb}=1$.
\end{claim}

{}From now on we identify the elements of $\CC(\Omega,\set{0,\aleph_0})$ with
the $C$-valued names of elements of $\set{0,\aleph_0}$, \emph{via} the
isomorphism of Claim~\ref{Cl:3b}. We obtain the following result.

\begin{claim}\label{Cl:4b}
The map $\eps$ can be factored through $\sim$, to an isomorphism from
$D/{\sim}$ onto $\CC(\Omega,\set{0,\aleph_0})$.
\end{claim}

\begin{cproof}
The proof that $\eps$ is an embedding for $\leq$ and for $+$ is the same as
in the proof of Claim~\ref{Cl:4} of Lemma~\ref{L:LotsofII}.
Let $\lga\in\CC(\Omega,\set{0,\aleph_0})$. There exists
$\la\in V^C$\index{VzzVB@$V^B$} such that $\bv{\la=0}=\bv{\lga=0}$, thus
$\la\in D$ and $\eps(\la)=\lga$, so
$\eps$ is surjective.
\end{cproof}

The rest of the proof proceeds like in the proof of Lemma~\ref{L:LotsofII}.
It follows from Claim~\ref{Cl:4b} that the Boolean algebra
\index{Boolean algebra} of projections of
$(D,\sim)$ is isomorphic to $C$.
If $\theta$ is an ordinal such that $\wdt(D)\leq\aleph_\theta$ and we put
$\fm=\omega_{\theta+\gamma}$, then $D^\fm/{\sim^\fm}$ has a $\rem$-increasing
chain of length $\gamma+1$. Since $\CC(\Omega,\set{0,\aleph_0})$ has a lower
embedding\index{lower embedding} into $D^\fm/{\sim^\fm}$,
$\Omega_{\I}=\Omega_{\II}=\es$ while
$\Omega_{\III}=\Omega$, and the image of $D^\fm/{\sim^\fm}$ in $\ol{S}$
contains a function with values above~$\aleph_\gamma$. Hence, there exists a
bounded lower subespalier\index{espalier!lower sub ---}
$(B,\sim)$ of $(D^\fm,\sim^\fm)$ such that
$B/{\sim}\cong\CC(\Omega,\two_\gamma)$.
\index{Tzzgamma@$\two_\gamma$}
\end{proof}

\begin{remark}\label{Rk:IIneedsforc}
Since the Boolean algebra\index{Boolean algebra} $\Co$ has an
absolute (in set-the\-o\-ret\-i\-cal sense) dense\index{dense} subalgebra,
namely, the Boolean algebra\index{Boolean algebra} $F_\omega$ of clopen
subsets of
$\set{0,1}^\omega$, the forcing could, in principle, have been eliminated
from the proof of Lemma~\ref{L:LotsofIII}: for example, one could have taken
for $D$ the completion of $C\otimes F_\omega$ (the tensor product for Boolean
algebras\index{Boolean algebra} is just the coproduct). Such an
argument would not have worked for Lemma~\ref{L:LotsofII}, because $\Bo$ of
the ground universe may not be dense\index{dense} in the $\Bo$ of a generic
extension.
\end{remark}

By Lemma~\ref{L:DUniv} and Proposition~\ref{P:DUniv}(i), we thus obtain the
following result.

\begin{theorem}\label{T:MeasDUniv}
The class of Boolean espaliers\index{espalier!Boolean ---} is
D-universal.\index{D-universal} Moreover, every bounded continuous dimension scale\index{continuous dimension scale} is isomorphic to the dimension
range\index{dimension range} of some Boolean
espalier.\index{espalier!Boolean ---}
\end{theorem}

\section[Complete meet-continuous lattices]
{Conditionally complete, meet-continuous, relatively complemented,
modular lattices}\label{S:CMSMLatt}
\index{lattice!sectionally complemented modular ---}

We first recall some basic lattice-theoretical definitions, see
\index{Gr\"atzer, G.}\cite{GLT2}. A
lattice $(L,\vee,\wedge)$ is \emph{modular},
\index{lattice!modular ---|ii}
if $x\geq z$ implies that
$x\wedge(y\vee z)=(x\wedge y)\vee z$, for all $x$, $y$, $z\in L$. We say that
$L$ is
\begin{itemize}
\item[---] \emph{complemented},
\index{lattice!complemented ---}
if it has a least element $0$, a largest
element $1$, and every $x\in L$ has a complement, that is, $y\in L$ such that
$x\wedge y=0$ and $x\vee y=1$.

\item[---] \emph{sectionally complemented},
\index{lattice!sectionally complemented ---|ii}
if it has a least element $0$ and
every sublattice of the form $[0,a]$, for $a\in L$, is complemented;

\item[---] \emph{relatively complemented},
\index{lattice!relatively complemented ---|ii}
if every sublattice of the form
$[a,b]$, for $a\leq b$ in~$L$, is complemented.
\end{itemize}

In general, these notions are unrelated. However, in the \emph{modular} case,
the following implications hold:
   \[
   \text{complemented}\Rightarrow\text{sectionally complemented}
   \Rightarrow\text{relatively complemented}.
   \]
We say that the lattice $L$ is \emph{complete}, if every subset of $L$ has
an infimum---equivalently, every subset of $L$ has a supremum. We say that
$L$ is \emph{conditionally complete}, if every nonempty bounded subset of $L$
has an infimum---equivalently, the interval $[a,b]$ is complete, for all
$a\leq b$ in $L$. We say that
\index{lattice!meet-continuous ---|ii}
$L$ is \emph{meet-continuous}, if for every $a\in L$ and every upward
directed subset $X$ of $L$ admitting a supremum, the equality
$a\wedge\bigvee X=\bigvee(a\wedge X)$ holds, where we put
$a\wedge X=\setm{a\wedge x}{x\in X}$. If the dual condition holds, $L$ is
called \index{lattice!join-continuous ---|ii}
\emph{join continuous}, and if both conditions hold,~$L$ is
\index{lattice!continuous ---|ii}
\emph{continuous}. (This definition of continuity is not equivalent
to the one presented in G. Gierz \emph{et al.}~\cite{Comp},
\index{Gierz, G.}\index{Hofmann, K.\,H.}\index{Keimel, K.}%
\index{Lawson, J.\,D.}\index{Mislove, M.}\index{Scott, D.\,S.}%
which is \index{lattice!Scott-continuous ---}
nowadays often called ``Scott continuity''.)
A \emph{continuous geometry}\index{continous geometry|ii} is any complete,
\index{lattice!complemented modular ---}
complemented, modular, continuous lattice. (This is the current
terminology; von Neumann's original definition included hypotheses of
irreducibility and lack of chain conditions. What we have called a
continuous geometry was called a ``reducible continuous geometry'' or a
``continuous geometry in the general sense'' in some of the older
literature.)

The \emph{dimension monoid}\index{dimension monoid|ii}
of a lattice $L$, see \cite{WDim},\index{Wehrung, F.} is the
commutative monoid $\Dim L$\index{dzzimL@$\Dim L$|ii} defined by generators
\index{dzzeltab@$\DD(a,b)$|ii}$\DD(a,b)$, for $a\leq b$ in
$L$, and the following relations:
\begin{itemize}
\item[(D0)] $\DD(a,a)=0$, for all $a\in L$.

\item[(D1)] $\DD(a,c)=\DD(a,b)+\DD(b,c)$, for all $a\leq b\leq c$ in $L$.

\item[(D2)] $\DD(a,a\vee b)=\DD(a\wedge b,b)$, for all $a$, $b\in L$.
\end{itemize}

It is an open problem whether the dimension monoid of an arbitrary lattice is
always a\index{refinement monoid} refinement monoid, however, this is solved
in a few important particular cases: the case of \emph{finite lattices}, of
which the dimension monoids are so-called \emph{primitive monoids}, and the
case of \emph{modular lattices}, for which an alternative presentation of
the dimension monoid is given that implies refinement. We shall concentrate
here on the latter.

In a modular lattice $L$ with zero, let $x\perp y$ hold, if $x\wedge y=0$,
for any $x$, $y\in L$, and then we define $x\oplus y=x\vee y$. The following
result is folklore, see, for example,\index{Wehrung, F.}
\cite[Proposition~8.1]{WDim}; it says, essentially, that $\oplus$ is
associative in modular lattices. We state it in a way that relates it to
the axioms defining espaliers.\index{espalier}

\begin{proposition}\label{P:oplL2}
Let $L$ be a modular lattice with zero.
\index{lattice!modular ---} Then the relation $\perp$ on $L$
satisfies \textup{(L2)}.
\end{proposition}

In a modular lattice\index{lattice!modular ---} $L$, we define the binary
relations $\sim$
\index{perspectivity@perspectivity ($\sim$)} (\emph{perspectivity}),
$\sim_2$\index{bi-perspectivity@bi-perspectivity ($\sim_2$)|ii}
(\emph{bi-perspectivity}),
\index{projectivity@projectivity ($\approx$)}
$\approx$ (\emph{projectivity}), and $\approxeq$
\index{projectivity by decomposition@projectivity by decomposition
($\approxeq$)|ii}
(\emph{projectivity by decomposition}) as follows:
   \begin{align*}
   x\sim y\ &\Longleftrightarrow\exists z\text{ such that }
   x\oplus z=y\oplus z;\\
   x\sim_2y\ &\Longleftrightarrow\exists z\text{ such that }x\sim z\sim y;\\
   x\approx y\ &\Longleftrightarrow\exists n\in\NN,\,\exists z_0,\dots,z_n
   \text{ such that }x=z_0\sim z_1\sim\cdots\sim z_n=y;\\
   x\approxeq y\ &\Longleftrightarrow\exists n\in\NN,\,
   \exists x_0\approx y_0,\dots,x_{n-1}\approx y_{n-1}\\
     &\qquad\qquad\text{ such that } x=\oplus_{i<n}x_i
   \text{ and }y=\oplus_{i<n}y_i.
   \end{align*}
In case $L$ is a relatively complemented lattice
\index{lattice!relatively complemented ---}
with zero, the dimension\index{dimension monoid}
monoid $\Dim L$\index{dzzimL@$\Dim L$} of $L$ is generated by the elements
$\DD(x)=\DD(0,x)$, for
$x\in L$ (see\index{Wehrung, F.} \cite[Proposition~9.1]{WDim}). We define the
\emph{dimension range}\index{dimension range|ii} of $L$
as\index{Dzzrng@$\Drng L$}
$\Drng L=\setm{\DD(x)}{x\in L}$. In case $L$ is also
modular,\index{Dzzrng@$\Drng L$} $\Drng L$ can be endowed with the partial
addition defined by
   \[
   \DD(a)+\DD(b)=\DD(a\oplus b),\text{ for all }a,\,b\in L
   \text{ such that }a\perp b.
   \]
Furthermore, the \pcm\ $\Drng L$\index{Dzzrng@$\Drng L$} determines
the\index{dzzimL@$\Dim L$} \cm\
$\Dim L$, and any equality of the form $\DD(a)=\DD(b)$ can be tested in a
very simple way, see Corollaries~9.4 and 9.5 in\index{Wehrung, F.}
\cite{WDim} and Proposition~\ref{P:PrmRm} of the present paper.

\begin{proposition}\label{P:D(a)=D(b)}
Let $L$ be a sectionally complemented modular
\index{lattice!sectionally complemented modular ---} lattice.
Put\index{Dzzrng@$\Drng L$} $S=\Drng L$. Then the following statements hold:
\begin{enumerate}
\item $\DD(a)=\DD(b)$ \iff\ $a\approxeq b$, for all $a$, $b\in L$.

\item $\Dim L=\Ref{S}$,\index{dzzimL@$\Dim L$} the universal monoid of $S$.
\end{enumerate}
\end{proposition}

This is the point where the theory of espaliers\index{espalier}
and continuous dimension scales\index{continuous dimension scale} comes
in. Our plan is to associate, with a sectionally complemented, modular
lattice
\index{lattice!sectionally complemented modular ---}
$L$, an espalier\index{espalier} $L^*$ such that the dimension
range of $L$, as defined above, is the dimension
range of $L^*$. The structure
$L^*$ is simply defined as $(L,\leq,\perp,\approxeq)$, for those $\perp$ and
$\approxeq$ defined above, so all we need to do is find sufficient
conditions for it to be an espalier.\index{espalier} The following lemma
sums up some of the hardest (in particular item (i)) and most useful results
of \cite{WDim}.\index{Wehrung, F.}

\begin{lemma}\label{L:LattInfRef}
Let $L$ be a conditionally complete, meet-continuous, sectionally
complemented, modular
\index{lattice!sectionally complemented modular ---} lattice. Then the
following statements hold:
\begin{enumerate}
\item $x\approxeq y$ \iff\ there are decompositions
$x=x_0\oplus x_1$ and $y=y_0\oplus y_1$ in $L$ such that $x_0\sim y_0$ and
$x_1\sim y_1$, for all $x$, $y\in L$.

\item Let $a$, $b\in L$ and let $\famm{b_i}{i\in I}$ be a family of elements
of $L$. If $a\sim\oplus_{i\in I}b_i$, then there exists a decomposition
$a=\oplus_{i\in I}a_i$ such that $a_i\sim b_i$, for all $i\in I$.

\item Let $\famm{a_i}{i\in I}$ and $\famm{b_j}{j\in J}$ be families of
elements of $L$ such that $\oplus_{i\in I}a_i=\oplus_{j\in J}b_j$. Then there
are families $\famm{x_{i,j}}{(i,j)\in I\times J}$ and
$\famm{y_{i,j}}{(i,j)\in I\times J}$ of elements of $L$ such that
   \begin{align*}
   a_i&=\oplus_{j\in J}x_{i,j},&&\text{for all }i\in I,\\
   b_j&=\oplus_{i\in I}y_{i,j},&&\text{for all }j\in J,
   \end{align*}
and $x_{i,j}\sim_2y_{i,j}$, for all $(i,j)\in I\times J$.
\end{enumerate}
\end{lemma}

\begin{proof}
(i) follows immediately from Lemma~10.4 and Theorem~13.2
in\index{Wehrung, F.} \cite{WDim}.

(ii) There exists $c\in L$ such that $a\oplus c=b\oplus c$.
Let $\tau\colon[0,b]\twoheadrightarrow[0,a]$ be the perspective mapping
\index{perspective mapping}
with axis $c$, that is, $\tau(x)=(x\oplus c)\wedge a$, for all $x\in[0,b]$.
Put $a_i=\tau(b_i)$, for all $i\in I$. Since $\tau$ is an isomorphism from
$[0,b]$ onto $[0,a]$, the equality $a=\oplus_{i\in I}a_i$ holds. Furthermore,
$a_i$ is perspective\index{perspectivity@perspectivity ($\sim$)} to $b_i$
with axis $c$, for all $i\in I$.

The proof of (iii) is virtually the same as the one of\index{Wehrung, F.}
\cite[Lemma~12.17]{WDim}, except that we replace countable families by
transfinite ones. We give the proof here for the convenience of the reader.
So, let $\lambda$ and $\kappa$ be ordinals, let
$\famm{a_\alpha}{\alpha<\lambda}$ and
$\famm{b_\beta}{\beta<\kappa}$ be orthogonal families of elements of $L$
such that $\oplus_{\alpha<\lambda}a_\alpha=\oplus_{\beta<\kappa}b_\beta$. We
put
   \begin{align*}
   \oll{a}_\alpha&=
   \oplus_{\xi<\alpha}a_\xi,&&\text{for all }\alpha\leq\lambda,\\
   \oll{b}_\beta&=\oplus_{\eta<\beta}b_\eta,&&\text{for all }\beta\leq\kappa.
   \end{align*}
Furthermore, we put $c_{\alpha,\beta}=
(\oll{a}_{\alpha+1}\wedge\oll{b}_{\beta})\vee
(\oll{a}_{\alpha}\wedge\oll{b}_{\beta+1})$ and
$d_{\alpha,\beta}=\oll{a}_{\alpha+1}\wedge\oll{b}_{\beta+1}$, for all
$\alpha<\lambda$ and all $\beta<\kappa$. Since
$c_{\alpha,\beta}\leq d_{\alpha,\beta}$ and $L$ is sectionally complemented,
there exists $z_{\alpha,\beta}\in L$ such that
$c_{\alpha,\beta}\oplus z_{\alpha,\beta}=d_{\alpha,\beta}$.

\setcounter{claim}{0}
\begin{claim}\label{Cl:IneqMC1}
The statement
$(\oll{a}_\alpha\wedge\oll{b}_\beta)\oplus(\oplus_{\eta<\beta}z_{\alpha,\eta})
=\oll{a}_{\alpha+1}\wedge\oll{b}_\beta$ holds, for all $\alpha<\lambda$
and all $\beta\leq\kappa$.
\end{claim}

\begin{cproof}
We prove the conclusion by induction on $\beta$. It is trivial for $\beta=0$.
For $\beta$ a limit ordinal, it follows easily from the meet-continuity
of~$L$ and the induction hypothesis. Now suppose having proved the statement
at step $\beta$, we prove it at step $\beta+1$. It follows from the induction
hypothesis that
$\oplus_{\eta<\beta}z_{\alpha,\eta}\leq\oll{a}_{\alpha+1}\wedge\oll{b}
_{\beta}$,
thus
   \begin{align*}
   \oll{a}_{\alpha}\wedge\oll{b}_{\beta+1}\wedge
   (\oplus_{\eta<\beta}z_{\alpha,\eta})
   &=(\oll{a}_{\alpha}\wedge\oll{b}_{\beta+1})\wedge
   (\oll{a}_{\alpha+1}\wedge\oll{b}_{\beta})\wedge
   (\oplus_{\eta<\beta}z_{\alpha,\eta})\\
   &=(\oll{a}_{\alpha}\wedge\oll{b}_{\beta})\wedge
   \oplus_{\eta<\beta}z_{\alpha,\eta}\\
   &=0,
   \end{align*}
that is,
$\oll{a}_{\alpha}\wedge\oll{b}_{\beta+1}\perp\oplus_{\eta<\beta}
z_{\alpha,\eta}$. Furthermore,
   \begin{align*}
   (\oll{a}_{\alpha}\wedge\oll{b}_{\beta+1})\oplus
   (\oplus_{\eta<\beta}z_{\alpha,\eta})
   &=(\oll{a}_{\alpha}\wedge\oll{b}_{\beta+1})\vee
   (\oll{a}_{\alpha}\wedge\oll{b}_{\beta})\vee
   (\oplus_{\eta<\beta}z_{\alpha,\eta})\\
   &=(\oll{a}_{\alpha}\wedge\oll{b}_{\beta+1})\vee
   (\oll{a}_{\alpha+1}\wedge\oll{b}_{\beta})\\
   &\mbox{\hphantom{induction hypothesis}}\text{ (by the induction
hypothesis)}\\
   &=c_{\alpha,\beta},
   \end{align*}
whence
   \begin{align*}
   d_{\alpha,\beta}&=c_{\alpha,\beta}\oplus z_{\alpha,\beta}\\
   &=\bigl((\oll{a}_{\alpha}\wedge\oll{b}_{\beta+1})\oplus
   (\oplus_{\eta<\beta}z_{\alpha,\eta})\bigr)\oplus z_{\alpha,\beta}\\
   &=(\oll{a}_{\alpha}\wedge\oll{b}_{\beta+1})\oplus
   (\oplus_{\eta<\beta+1}z_{\alpha,\eta}).\tag*{\qedc}
   \end{align*}
\renewcommand{\qedc}{}
\end{cproof}

Of course, by symmetry, the following claim is also valid.

\begin{claim}\label{Cl:IneqMC2}
The statement
$(\oll{a}_\alpha\wedge\oll{b}_\beta)\oplus(\oplus_{\xi<\alpha}z_{\xi,\beta})
=\oll{a}_{\alpha}\wedge\oll{b}_{\beta+1}$ holds, for all $\alpha\leq\lambda$
and all $\beta<\kappa$.
\end{claim}

In particular, using Claim~\ref{Cl:IneqMC1} for $\beta=\kappa$ yields that
$\oll{a}_\alpha\oplus(\oplus_{\eta<\kappa}z_{\alpha,\eta})=
\oll{a}_{\alpha}\oplus a_{\alpha}$, thus
$\oplus_{\eta<\kappa}z_{\alpha,\eta}\sim a_{\alpha}$. Thus, by item (ii)
above, there exists a decomposition of the form
$a_{\alpha}=\oplus_{\eta<\kappa}x_{\alpha,\eta}$ such that
$x_{\alpha,\eta}\sim z_{\alpha,\eta}$, for all $\eta<\kappa$. Similarly, for
all $\beta<\kappa$, there exists a decomposition of the form
$b_{\beta}=\oplus_{\xi<\lambda}y_{\xi,\beta}$ such that
$z_{\xi,\beta}\sim y_{\xi,\beta}$, for all $\xi<\lambda$. It follows that
$x_{\xi,\eta}\sim_2y_{\xi,\eta}$,
\index{bi-perspectivity@bi-perspectivity ($\sim_2$)}
for all $\xi<\alpha$ and $\eta<\beta$, and
the $x_{\xi,\eta}$-s and $y_{\xi,\eta}$-s are as desired.
\end{proof}

\begin{proposition}\label{P:L2L*esp}
Let $L$ be a conditionally complete, meet-continuous, sectionally
complemented, modular
\index{lattice!sectionally complemented modular ---} lattice. Then
$L^*=(L,\leq,\perp,\approxeq)$ is an espalier,\index{espalier}
and\index{Dzzrng@$\Drng L$}
$\Drng L=\Drng L^*$.
\end{proposition}

\begin{proof}
The verifications in $L^*$ of Axioms (L1) to (L5) and (L8) are either
trivial or immediate consequences of the assumptions and
Proposition~\ref{P:oplL2}.

Let $a$, $b$, $b_i$ (for $i\in I$) be elements of~$L$ such that $a\approxeq b$
and $b=\oplus_{i\in I}b_i$. By Lemma~\ref{L:LattInfRef}(i), there are $a'$,
$a''$, $b'$, $b''\in L$ such that $a=a'\oplus a''$, $b=b'\oplus b''$,
$a'\sim b'$, and $a''\sim b''$. Applying Lemma~\ref{L:LattInfRef}(iii) to
the equality $b'\oplus b''=\oplus_{i\in I}b_i$, we obtain decompositions
$b'=\oplus_{i\in I}x'_i$, $b''=\oplus_{i\in I}x''_i$, $b_i=y'_i\oplus y''_i$,
for all $i\in I$, such that $x'_i\sim_2y'_i$ and $x''_i\sim_2y''_i$,
\index{bi-perspectivity@bi-perspectivity ($\sim_2$)}
for all
$i\in I$. Since $a'\sim b'=\oplus_{i\in I}x'_i$ and
$a''\sim b''=\oplus_{i\in I}x''_i$, there are, by
Lemma~\ref{L:LattInfRef}(ii), decompositions $a'=\oplus_{i\in I}a'_i$ and
$a''=\oplus_{i\in I}a''_i$ such that $a'_i\sim x'_i$ and $a''_i\sim x''_i$,
for all $i\in I$. Observe that $a=a'\oplus a''=\oplus_{i\in I}a_i$, where we
put $a_i=a'_i\oplus a''_i$, for all $i\in I$. Furthermore, $a_i\approxeq b_i$,
for all $i\in I$. Hence $\approxeq$ is unrestrictedly refining
\index{unrestrictedly refining relation} (Axiom (L6)).

The proof that $\approxeq$ is unrestrictedly additive
\index{unrestrictedly additive!--- relation} (Axiom~(L7)) is
virtually the same as\index{Wehrung, F.} \cite[Proposition~13.9]{WDim},
except that countable families are replaced by arbitrary families.

The statement that\index{Dzzrng@$\Drng L$} $\Drng L=\Drng L^*$ follows
immediately from Proposition~\ref{P:D(a)=D(b)}(i).
\end{proof}

\begin{corollary}\label{C:DimMonMCSLatt}
The dimension monoid\index{dimension monoid} of any conditionally complete,
meet-continuous, sectionally complemented, modular lattice
\index{lattice!sectionally complemented modular ---}
is a \pup{total} continuous dimension scale.\index{continuous dimension scale}
\end{corollary}

\begin{proof}
By Proposition~\ref{P:L2L*esp}, $L^*$ is an espalier\index{espalier}
and\index{Dzzrng@$\Drng L$}
$\Drng L=\Drng L^*$, thus, by Theorem~\ref{T:DimEsp}, $S=\Drng
L$\index{Dzzrng@$\Drng L$} is a continuous dimension scale.\index{continuous dimension scale} By Corollary~\ref{C:EmbDimInt}, $\Ref S$ is also a continuous dimension scale,\index{continuous dimension scale} but by Proposition~\ref{P:D(a)=D(b)},
$\Dim L$\index{dzzimL@$\Dim L$} is isomorphic to $\Ref S$, thus it is a
continuous dimension scale.\index{continuous dimension scale}
\end{proof}

Hence we can complete the program of determining the dimension theory of
conditionally complete, meet-continuous, relatively complemented, modular
lattices
\index{lattice!sectionally complemented modular ---}
initiated by I. Halperin\index{Halperin, I.}\index{von Neumann, J.}
and J. von~Neumann in \cite{HaNe40}.

\begin{theorem}\label{T:DimMonRCMLatt}
The dimension monoid\index{dimension monoid} of any conditionally complete,
meet-continuous, relatively complemented modular lattice
\index{lattice!sectionally complemented modular ---}
is a \pup{total} continuous dimension scale.\index{continuous dimension scale}
\end{theorem}

\begin{proof}
Let $L$ be a conditionally complete, meet-continuous, relatively complemented
\index{lattice!sectionally complemented modular ---}
modular lattice. The interval $[a,b]$ is a conditionally complete,
meet-continuous, complemented modular lattice, for all $a\leq b$ in $L$,
thus, by Corollary~\ref{C:DimMonMCSLatt}, its dimension
monoid\index{dimension monoid} $\Dim[a,b]$\index{dzzimL@$\Dim L$} is a
continuous dimension scale.\index{continuous dimension scale} Furthermore, if $a'\leq a\leq
b\leq b'$ in $L$, then the natural map from
$\Dim[a,b]$\index{dzzimL@$\Dim L$} to
$\Dim[a',b']$\index{dzzimL@$\Dim L$} is, by\index{Wehrung, F.}
\cite[Corollary~13.5]{WDim}, a lower embedding\index{lower embedding} of
\cm s. Express the lattice
$L$ as the direct limit of the direct system $\mathcal{I}$ of its closed
intervals. Since the $\Dim$\index{dzzimL@$\Dim L$} functor preserves direct
limits, it follows from Lemma~\ref{L:DirUn} that
$\Dim L=\varinjlim_{[a,b]\in\mathcal{I}}\Dim[a,b]$\index{dzzimL@$\Dim L$} is
a continuous dimension scale.\index{continuous dimension scale}
\end{proof}

We say that a (von Neumann) regular ring $R$
\index{ring!Von Neumann regular ---} is
\emph{right continuous},\index{ring!right continuous regular ---|ii}
if the lattice\index{Lzzat@$\Lat(R)$|ii}
$\Lat(R)$ of all principal right ideals of $R$
(cf.~\cite[Theorem~2.3]{GvnRR})
is complete and meet-continuous. In
particular, every right self-injective regular ring
\index{ring!Von Neumann regular ---}\index{ring!right self-injective ---}%
is right continuous, see \index{Goodearl, K.\,R.}\cite[Corollary~13.5]{GvnRR}.
The connection between the present section and the upcoming
Section~\ref{S:RSIReg} is made possible by the following immediate
consequence of \cite[Corollary~13.4]{WDim}.\index{Wehrung, F.}

\begin{proposition}\label{P:Latt2Ring}
Let $R$ be a right continuous regular ring.
\index{ring!Von Neumann regular ---}\index{ring!right self-injective ---}%
Then the following statements hold:
\begin{enumerate}
\item The monoid $V(R)$
\index{VzzofR@$V(R)$}
of isomorphism classes of finitely generated
projective right $R$-modules is isomorphic to the dimension
monoid\index{dimension monoid}\index{dzzimL@$\Dim L$}
$\Dim\Lat(R)$ of~$\Lat(R)$.

\item Two principal right ideals $I$ and $J$ of $R$ are isomorphic \iff\
there are decompositions $I=I_0\oplus I_1$ and $J=J_0\oplus J_1$ such that
$I_0\sim J_0$ and $I_1\sim J_1$. In particular, $I\cong J$ \iff\
$I\approxeq J$ in the lattice $\Lat(R)$.
\end{enumerate}
\end{proposition}

The results of Section~\ref{S:RSIReg} below imply that the
espaliers\index{espalier} of the form $L^*$ appearing in
Proposition~\ref{P:L2L*esp} form a D-universal\index{D-universal} class.
Here is a slightly sharper statement.

\begin{theorem}\label{T:meetcontinDuniv}
The class of espaliers\index{espalier} of the form
$(L,{\le},{\perp},{\approxeq})$, where
$L$ is a complete, meet-continuous, complemented, modular
\index{lattice!complemented modular ---} lattice, is
D-universal.\index{D-universal}
\end{theorem}

\begin{proof}
It follows from Theorem~\ref{T:L(R)Duniv} below that every continuous dimension scale\index{continuous dimension scale} admits a lower embedding\index{lower
embedding} into the dimension range\index{dimension range} of an
espalier\index{espalier} of the form
$(\Lat(R),{\subseteq},{\perp},{\cong})$, for some regular,
right self-injective ring $R$.
\index{ring!Von Neumann regular ---}\index{ring!right self-injective ---}%
It follows from Proposition~\ref{P:Latt2Ring}(ii) that the relations 
$\cong$ and
$\approxeq$ on $\Lat(R)$ coincide. Since $\Lat(R)$ is a
complete, meet-continuous, complemented, modular
\index{lattice!complemented modular ---}
lattice, the conclusion follows.
\end{proof}

In particular, it
follows from Propositions~\ref{P:L2L*esp} and \ref{P:Latt2Ring} and
Theorem~\ref{T:L(R)Duniv} that every continuous dimension scale\index{continuous dimension scale} admits a lower embedding\index{lower embedding} into
\index{dzzimL@$\Dim L$} $\Dim\Lat(R)$, for some regular, right self-injective
ring~$R$.
\index{ring!Von Neumann regular ---}\index{ring!right self-injective ---}%
Observe again that $\Lat(R)$ is a complete, meet-continuous,
complemented, modular
\index{lattice!complemented modular ---} lattice.

\begin{corollary}\label{C:DimSCML}
The dimension monoids\index{dimension monoid} of conditionally complete,
meet-continuous, sectionally \pup{resp., relatively} complemented, modular
lattices
\index{lattice!sectionally complemented modular ---} are
exactly the total continuous dimension scales.\index{continuous dimension scale}
\end{corollary}

Hence, validating the possibility suggested in F.
Wehrung \index{Wehrung, F.}\cite{CGApp}, the dimension theory of
conditionally complete, meet-continuous, relatively complemented, modular
lattices is completely elucidated.

\section[Self-injective regular rings]{Self-injective regular rings and
nonsingular injective modules}
\label{S:RSIReg}\index{ring!Von Neumann regular ---}
\index{ring!right self-injective ---}

For notation, terminology, and standard results on the topics of
this section, we refer to \cite{GoBo, Gnonsing, GvnRR}.
\index{Goodearl, K.\,R.}\index{Boyle, A.\,K.}%
Throughout the section, let
$R$ denote a (von Neumann) regular (unital) ring;\index{ring!Von Neumann
regular ---} after some preliminary results, we shall assume that $R$ is
also right self-injective, that is, $R$ is injective as a right module over
itself.

Let $\Lat(R)$ denote the collection of principal right
ideals of $R$. Regularity implies that $\Lat(R)$ is a complemented
modular
\index{lattice!complemented modular ---}
lattice, in which finite suprema and infima are given by
sums and intersections, respectively (e.g., \cite[Theorem~2.3]{GvnRR}).
\index{Goodearl, K.\,R.}
Define orthogonality in $\Lat(R)$ to mean lattice disjointness: $A\perp B$
if and only if $A\cap B=0$. For an equivalence relation on $\Lat(R)$, we
shall use $\cong$, that is, isomorphism of right $R$-modules.

A small amount of category-theoretical notation will be helpful in
dealing with $R$-modules. We write
$\ModR$\index{MzzodR@$\ModR$|ii} for the
category of all right $R$-modules, and $\addA$\index{azzddA@$\addA$|ii}
for the full subcategory of
$\ModR$ whose objects are all direct summands of finite direct sums of
copies of a given module $A$. In particular, the objects of
$\addR$ are precisely the finitely generated projective right
$R$-modules. An expression such as ``$A\in\ModR$'' will abbreviate
the assertion that $A$ is an object in the category $\ModR$. We
write $\rE(A)$\index{ErzzE@$\rE(M)$}
to stand for an arbitrary injective hull\index{injective hull} of a module
$A$, and if $\kappa$ is a cardinal, $\kappa\cdot A$
\index{kzzcdotA@$\kappa\cdot A$|ii} stands for a direct sum
of $\kappa$ copies of $A$. For $A,B\in \ModR$, write $A\lesssim B$
\index{ABzzlesssimB@$A\lesssim B$ (for modules)|ii} to
mean that $A$ is isomorphic to a submodule of $B$. If $A,B\in
\addR$ and $A\lesssim B$, then regularity of $R$ implies that $A$
is isomorphic to a direct summand of $B$\index{Goodearl, K.\,R.}
\cite[Theorem~1.11]{GvnRR}.

We use the notation $V(R)$\index{VzzofR@$V(R)$} for the monoid of
isomorphism classes of objects from $\addR$ (in which the
addition operation is induced from the direct sum operation on
modules). To match our notation for dimension\index{dimension range} ranges,
we shall denote elements of $V(R)$\index{VzzofR@$V(R)$} in the form
$\Delta(A)$, rather than using a more common notation like $[A]$. This
involves a slight but unproblematic abuse of notation in case $A\in\Lat(R)$,
since the element $\Delta(A)\in V(R)$\index{VzzofR@$V(R)$} stands for the
isomorphism class of
$A$ within the class of all right $R$-modules, whereas
once we have made $\Lat(R)$ into an espalier,\index{espalier} the notation
$\Delta(A)$ will also be used for the image of $A$ in
\index{Dzzrng@$\Drng L$} $\Drng\Lat(R)$, and in the latter case $\Delta(A)$
stands for the isomorphism class of $A$ within
$\Lat(R)$.

\begin{lemma}\label{L:V(R)}
For any regular ring $R$,
\index{ring!Von Neumann regular ---} the
monoid $V(R)$\index{VzzofR@$V(R)$} is a refinement mon\-oid, the interval
$[0,\Delta(R)]\subseteq V(R)$ is a\index{refinement monoid!partial ---}
partial refinement monoid, and $V(R)$\index{VzzofR@$V(R)$} is the universal
monoid of $[0,\Delta(R)]$.
\end{lemma}

\begin{proof} Refinement in $V(R)$\index{VzzofR@$V(R)$} is given by
\index{Goodearl, K.\,R.}\cite[Theorem~2.8]{GvnRR}, and it is clear that
$S=[0,\Delta(R)]$ is a partial submonoid of\index{VzzofR@$V(R)$} $V(R)$,
hence a\index{refinement monoid!partial ---} partial refinement monoid in
its own right. Since every object in
$\addR$ is isomorphic to a finite direct sum of principal right ideals of
$R$ \index{Goodearl, K.\,R.}\cite[Proposition~2.6]{GvnRR}, every element
of\index{VzzofR@$V(R)$} $V(R)$ is a sum of elements of $S$.

Let $U$ denote the universal monoid of $S$, with canonical map
$\phi\colon S\rightarrow U$. There exists a unique homomorphism
$\psi\colon U\rightarrow V(R)$\index{VzzofR@$V(R)$} such that
$\psi\phi\colon S\rightarrow V(R)$\index{VzzofR@$V(R)$} is the inclusion map.
Given $x\in V(R)$, write
$x= \sum_{i<n} x_i$ for some elements $x_i\in S$, and set $\theta(x)=
\sum_{i<n}
\phi(x_i)\in U$; this is well defined by refinement. Hence, we
obtain a homomorphism $\theta\colon V(R)\rightarrow U$. Obviously
$\psi\theta$ is the identity map on\index{VzzofR@$V(R)$} $V(R)$, and
$\theta\psi\phi=\phi$, whence $\theta\psi$ is the identity map on $U$.
Therefore $\psi$ is an isomorphism.
\end{proof}

We next determine the projections on $V(R)$\index{VzzofR@$V(R)$} and on
$[0,\Delta(R)]$. This requires working with orthogonality in
$V(R)$\index{VzzofR@$V(R)$} (as defined in
Section~\ref{S:DirectDecompRefMon}), which is determined as follows
\cite[Proposition~2.21]{GvnRR}:\index{Goodearl, K.\,R.} For any
$A,B\in \addR$, we have
   \[
   \Delta(A)\perp\Delta(B) \Longleftrightarrow \Hom_R(A,B)=0
   \Longleftrightarrow \Hom_R(B,A)=0.
   \]
Let $\rB(R)$\index{BrzzB@$\rB(R)$|ii}
denote the set of all central idempotents in $R$; this
is a Boolean algebra \index{Boolean algebra}
whose operations are given by the rules
   \[
   e\wedge f = ef \qquad\qquad e\vee f = e+f-ef \qquad\qquad e' = 1-e
   \]
\cite[p.~83]{GvnRR}.\index{Goodearl, K.\,R.} If $R$ is right self-injective,
$\rB(R)$ is complete \cite[Proposition~9.9]{GvnRR}.\index{Goodearl, K.\,R.}

\begin{lemma}\label{L:B(R)}
Let $R$ be a regular ring.\index{ring!Von Neumann regular ---}
\enumerate
\item For any $e\in \rB(R)$, there is a projection $p_e$ on
$V(R)$\index{VzzofR@$V(R)$} such that $p_e\Delta(A)= \Delta(Ae)$ for all
$A\in \addR$. Moreover,
$p_e^\perp= p_{1-e}$.
\item The rule $e\mapsto p_e$ defines an isomorphism\index{pzzroj@$\BB{S}$}
\index{VzzofR@$V(R)$}
$\rB(R)\overset{\cong}\longrightarrow \BB V(R)$.
\item The rule $e\mapsto p_e|_{[0,\Delta(R)]}$ defines an
\index{pzzroj@$\BB{S}$}
isomorphism $\rB(R) \overset{\cong}\longrightarrow \BB[0,\Delta(R)]$.
\end{lemma}

\begin{proof} (i) It is clear that there is an endomorphism $p_e$
of $V(R)$\index{VzzofR@$V(R)$} such that $p_e\Delta(A)= \Delta(Ae)$ for all
$A\in
\addR$. Moreover, $\Delta(A)= \Delta(Ae)+ \Delta(A(1-e))$, and
$\Delta(A(1-e))\perp \Delta(Be)$ for all $B\in\addR$, so that
$\Delta(A(1-e)) \in (p_eV(R))^\perp$. Therefore\index{pzzroj@$\BB{S}$}
$p_e\in \BB V(R)$. Similarly, $p_{1-e}\in \BB V(R)$,\index{VzzofR@$V(R)$}
and we observe that\index{VzzofR@$V(R)$}
$V(R)= p_eV(R)\oplus p_{1-e}V(R)$. Therefore $p_{1-e}= p_e^\perp$.

(ii) Let $e$, $f\in \rB(R)$. If $e\leq f$, then $e=fe$, whence $p_ep_f=
p_e$ and so $p_e\leq p_f$. Conversely, if $p_e\leq p_f$, then
$p_e(x)\leq p_f(x)$ for all $x\in V(R)$\index{VzzofR@$V(R)$}
(Lemma~\ref{L:pLeqOrtq}(i)), whence
$Re\lesssim Rf$ (taking $x= \Delta(R)$). Consequently,
$Re(1-f)=0$, and so $e\leq f$. This shows that the map $e\mapsto
p_e$ is an order-embedding of $\rB(R)$ into
\index{pzzroj@$\BB{S}$}\index{VzzofR@$V(R)$}%
$\BB V(R)$.

Given\index{pzzroj@$\BB{S}$} $p\in \BB V(R)$, we have $\Delta(R)= p\Delta(R)+
p^\perp\Delta(R)$, and so $R= I\oplus J$ for some right ideals
$I$, $J$ such that $\Delta(I)= p\Delta(R)$ and $\Delta(J)=
p^\perp\Delta(R)$. There is an idempotent $e\in R$ such that $eR=
I$ and $(1-e)R= J$. Moreover,
$\Delta(I)\perp\Delta(J)$, and thus $\Hom_R(I,J) =0$.
This homomorphism group being isomorphic to $(1-e)Re$, we conclude
that $e\in \rB(R)$ \cite[Lemma~3.1]{GvnRR}.\index{Goodearl, K.\,R.} In
particular, we can now write $I= Re$ and $J= R(1-e)$. Any $A\in\addR$ is
isomorphic to a direct summand of $n\cdot R$ for some positive integer $n$,
whence
$p\Delta(A)\leq np\Delta(R)= n\Delta(Re)= p_e(n\Delta(R))$. On the
other hand,
$Ae$ is isomorphic to a direct summand of $n\cdot(Re)$, and so
$p_e\Delta(A)\leq n\Delta(Re)= p(n\Delta(R))$. Since
$pV(R)$\index{VzzofR@$V(R)$} and
$p_eV(R)$ are ideals of $V(R)$, it follows that they are equal,
and therefore $p=p_e$.

(iii) This is proved in the same manner as (ii). \end{proof}

Recall that we have defined orthogonality in $\Lat(R)$ by the rule
$A\perp B \Longleftrightarrow A\cap B =0$. When this occurs, the
right ideal $A+B$ is both the orthogonal sum of $A$ and $B$ within
$\Lat(R)$ and the module-theoretic direct sum of $A$ and $B$, so
that the two uses of the expression $A\oplus B$ coincide. However,
infinite orthogonal sums in $\Lat(R)$ (when they exist) cannot
be module-theoretic direct sums, since the direct sum of an
infinite family of nonzero modules is not finitely generated. To
distinguish these cases, let us write $\oplus^\perp_{i\in I} A_i$
for the orthogonal sum of a family
$\famm{A_i}{i\in I}$ of elements of $\Lat(R)$ and $\bigoplus_{i\in I} A_i$
for the module-theoretic direct sum. (For either to exist, the family
$\famm{A_i}{i\in I}$ must be independent.)

\begin{proposition}\label{P:L(R)}
Let $\Lat(R)$ be the lattice of principal right ideals of a
regular, right self-injective ring
\index{ring!Von Neumann regular ---}\index{ring!right self-injective ---}%
$R$. Then
$(\Lat(R),{\subseteq},{\perp},{\cong})$ is an espalier,\index{espalier} and
its dimension range\index{dimension range} is isomorphic to the
interval\index{VzzofR@$V(R)$}
$[0,\Delta(R)] \subseteq V(R)$. Consequently, both $[0,\Delta(R)]$ and
$V(R)$ are continuous dimension scales.\index{continuous dimension scale} In case $R$ is
purely infinite,\index{VzzofR@$V(R)$}\index{Dzzrng@$\Drng L$}
$\Drng\Lat(R)\cong V(R)$.
\end{proposition}

\begin{proof}
By \cite[Corollary 13.5]{GvnRR},\index{Goodearl, K.\,R.} $\Lat(R)$ is
complete and upper continuous (= meet-continuous). In particular,
Axiom (L1) holds. As shown in the proof of
\cite[Proposition~13.3]{GvnRR},\index{Goodearl, K.\,R.} arbitrary infima in
$\Lat(R)$ are given by intersections, while the supremum of a family
$\famm{A_i}{i\in I}$ of elements of $\Lat(R)$ is the unique principal right
ideal of $R$ which contains $\sum_{i\in I} A_i$ as an essential submodule.
Since $R$ is right self-injective,\index{ErzzE@$\rE(M)$}
$\bigvee_{i\in I}A_i=\rE(\sum_{i\in I} A_i)$. Hence, if $\famm{A_i}{i\in I}$
is an orthogonal family,
$\oplus^\perp_{i\in I}A_i=\rE(\bigoplus_{i\in I}A_i)$.

Axioms (L2), (L3), (L5), and (L8) are clear, (L2)(iv) and (L8)
being standard properties of submodules of arbitrary modules.
Axioms (L4), (L6), and (L7) are basic properties of\index{injective hull}
injective hulls. Therefore $\Lat(R)$ is an espalier.\index{espalier} It is
clear that\index{Dzzrng@$\Drng L$}
$\Drng\Lat(R)\cong[0,\Delta(R)]$. If $R$ is purely infinite, then
$n\cdot R\cong R$ for all positive integers $n$
\cite[Theorem~10.16]{GvnRR},\index{Goodearl, K.\,R.} in which
case\index{VzzofR@$V(R)$} $[0,\Delta(R)]= V(R)$.

Finally, $[0,\Delta(R)]$ is a continuous dimension scale\index{continuous dimension scale}
by Theorem~\ref{T:DimEsp}, and it follows from Lemma~\ref{L:V(R)} and
Corollary \ref{C:EmbDimInt} that
$V(R)$\index{VzzofR@$V(R)$} is a continuous dimension scale.\index{continuous dimension scale}
\end{proof}

For the remainder of the section, assume that $R$ is right
self-injective. Before applying Theorem~\ref{T:EmbDimInt}, we show
that the type decomposition of $R$ (see \cite[Chapter VII]{GoBo} or
\cite[Chapter 10]{GvnRR})\index{Goodearl, K.\,R.}\index{Boyle, A.\,K.}
matches the type decomposition of
$V(R)$\index{VzzofR@$V(R)$} (Definition \ref{D:SI,II,III}). Here it is
natural to work with type decompositions of modules from $\addR$, as in
\cite[Theorem~7.2]{GoBo} and\index{Goodearl, K.\,R.}\index{Boyle, A.\,K.}
\cite[Theorem~10.31]{GvnRR}.

\begin{lemma}\label{L:typeV(R)}
Let $A\in\addR$.\index{VzzofR@$V(R)$} The following statements hold:
\enumerate
\item $A$ is an abelian, directly\index{directly finite} finite, or purely
infinite module, respectively, if and only if $\Delta(A)$ is a
multiple-free, directly\index{directly finite} finite, or purely infinite
element, respectively, of $V(R)$.\index{VzzofR@$V(R)$}
\item $A$ is of Type $\I$, $\II$, $\III$, respectively, if and only if
$\Delta(A)$ lies in $V(R)_{\I}$, $V(R)_{\II}$,
$V(R)_{\III}$,\index{VzzofR@$V(R)$} respectively.
\endenumerate
\end{lemma}

\begin{proof} (i) The equivalence for directly\index{directly finite} finite
modules is clear from the definitions, and the other two equivalences follow
from \cite[Theorems 2.1, 6.2]{GoBo}.
\index{Goodearl, K.\,R.}\index{Boyle, A.\,K.}%

(ii) First, $\Delta(A)\in V(R)_{\III}=V(R)_{\fin}^\perp$
\index{VzzofR@$V(R)$} if and only if
$\Hom_R(B,A)=0$ for all directly\index{directly finite} finite $B\in
\addR$, if and only if $A$ has no nonzero directly\index{directly finite}
finite direct summands, if and only if $A$ is of Type III
\cite[p.~37]{GoBo}.\index{Goodearl, K.\,R.}\index{Boyle, A.\,K.}
Similarly, $\Delta(A)\in V(R)_{\mf}^\perp$ if and only if $\Hom_R(B,A)=0$ for
all abelian $B\in \addR$, if and only if $A$ has no nonzero
abelian direct summands. Consequently, $\Delta(A)\in V(R)_{\I}$ if
and only if every nonzero direct summand of $A$ contains a nonzero
abelian direct summand, while $\Delta(A)\in V(R)_{\III}$ if and
only if $A$ has no nonzero abelian direct summands, but every
nonzero direct summand of $A$ contains a nonzero directly\index{directly
finite} finite direct summand. Thus, the remainder of part (ii) follows from
Theorems 5.1 and 5.5 of~\cite{GoBo}.\index{Goodearl, K.\,R.}
\index{Boyle, A.\,K.}
\end{proof}

\begin{corollary}\label{C:typeV(R)}
If $R$ is of Type $\I$, $\II$, $\III$, respectively, then
$V(R)$\index{VzzofR@$V(R)$} equals
$V(R)_{\I}$, $V(R)_{\II}$, $V(R)_{\III}$,
respectively.\index{VzzofR@$V(R)$}
\end{corollary}

\begin{proof}
\cite[Theorem~5.11]{GoBo}.
\index{Goodearl, K.\,R.}\index{Boyle, A.\,K.}%
\end{proof}

\begin{theorem}\label{T:isomV(R)}
Let $R$ be a regular, right self-injective ring,
\index{ring!Von Neumann regular ---}\index{ring!right self-injective ---}%
and write
$R=R_{\I}\times R_{\II}\times R_{\III}$ where $R_{\J}$ is of Type $\J$.
Let $\Omega_{\J}$ be the ultrafilter space of $\rB(R_{\J})$.
Then there
exists an ordinal $\gamma$ such that $V(R)$\index{VzzofR@$V(R)$} is
isomorphic to a lower subset of $\CC(\Omega_{\I},\ZZ_\gamma)\times
\CC(\Omega_{\II},\RR_\gamma)\times \CC(\Omega_{\III},\two_\gamma)$.
\index{Zzzgamma@$\ZZ_\gamma$}\index{Rzzgamma@$\RR_\gamma$}%
\index{Tzzgamma@$\two_\gamma$}%
\end{theorem}

\begin{proof} Since\index{VzzofR@$V(R)$}
$V(R)\cong V(R_{\I})\times V(R_{\II})\times V(R_{\III})$, it follows from
Corollary \ref{C:typeV(R)} that each\index{VzzofR@$V(R)$}
$V(R)_{\J}\cong V(R_{\J})$. By Proposition~\ref{P:L(R)},
$V(R)$\index{VzzofR@$V(R)$} is a continuous dimension scale,\index{continuous dimension scale} and using Lemma~\ref{L:B(R)}(ii) we see that each $\Omega_{\J}$ is
homeomorphic to the ultrafilter space of\index{pzzroj@$\BB{S}$} $\BB
V(R)_{\J}$. Therefore the theorem follows from Theorem~\ref{T:EmbDimInt}.
\end{proof}

We now turn our attention to nonsingular injective modules, which
allows us to extend the above results to proper Continuous Dimension
Scales, and which will allow us to show that the
espaliers\index{espalier} of the form $\Lat(R)$ form a
D-universal\index{D-universal} class.

Let $\NSIR$\index{NzzSIR@$\NSIR$} denote the full subcategory of
$\ModR$ whose objects are all nonsingular injective right $R$-modules.
Note that if
$A\in \NSIR$, then $\addA \subseteq\NSIR$. Let $V(A)$\index{VzzofR@$V(R)$}
be the monoid of isomorphism classes of objects in $\addA$, where---as
above---we use $\Delta(B)$ to denote the isomorphism class of an
object $B$. Note that for $B$, $C\in \addA$, we have $\Delta(B)\leq
\Delta(C)$ if and only if $B\lesssim C$.

For $A\in\NSIR$, let $\Lat(A)$ denote the collection of those
submodules of $A$ which are direct summands. Then $\Lat(A)$ is a
complete, complemented, modular
\index{lattice!complemented modular ---}
lattice, with infima and suprema\index{Goodearl, K.\,R.}\index{Boyle, A.\,K.}
given just as in $\Lat(R)$ \cite[Propositions 1.3, 1.6]{GoBo}. We
define $\perp$ in $\Lat(A)$ as in $\Lat(R)$.

\begin{lemma}\label{L:endinject}
Let $A$ be a nonsingular injective right $R$-module, and set $T=
\End_R(A)$. Then $T$ is a regular, right self-injective ring,
\index{ring!Von Neumann regular ---}\index{ring!right self-injective ---}%
and\index{VzzofR@$V(R)$}
$V(T)\cong V(A)$. Consequently, $V(A)$ is a dimension\index{VzzofR@$V(R)$}
interval.\index{continuous dimension scale}

Moreover, $(\Lat(A),{\subseteq},{\perp},{\cong})$ is an
espalier,\index{espalier} isomorphic to
$(\Lat(T),{\subseteq},{\perp},{\cong})$. Consequently,
\index{Dzzrng@$\Drng L$}$\Drng \Lat(A)\cong
[0,\Delta(A)]\subseteq V(A)$.\index{VzzofR@$V(R)$}
\end{lemma}

\begin{proof} For the first statement, see, for example, \cite[Corollary
1.23]{GvnRR}.\index{Goodearl, K.\,R.} It is well known that $\addA$ is
equivalent to the category of finitely generated projective right
$T$-modules (e.g.,\index{Lam, T.\,Y.}
\cite[Theorem~18.59]{Lam}), and thus $V(A)\cong V(T)$.\index{VzzofR@$V(R)$}
Therefore, Proposition~\ref{P:L(R)} implies that $V(A)$\index{VzzofR@$V(R)$}
is a continuous dimension scale.\index{continuous dimension scale}

According to \cite[Proposition~1.8]{GoBo},
\index{Goodearl, K.\,R.}\index{Boyle, A.\,K.}%
there is a lattice isomorphism $\phi\colon\Lat(T)\rightarrow \Lat(A)$
given by the rule
$\phi(J)= JA$. Any pair of right ideals of $T$ is given by $eT$,
$fT$ for some pair $e$, $f$ of idempotents in $T$, and it is well
known that $eT\cong fT$ if and only if $eA\cong fA$ (cf.~the proof
of \cite[Proposition~2.4]{GvnRR}).\index{Goodearl, K.\,R.} Hence, $\phi$ is
an isomorphism of espaliers.\index{espalier} It is clear that
\index{Dzzrng@$\Drng L$}$\Drng \Lat(A)\cong[0,\Delta(A)]$.
\end{proof}

A major advantage of working with nonsingular injective modules is
that any set of such modules can be combined to form a new one, by
taking the injective hull\index{injective hull} of the direct sum.
Consequently, we can pass from the category $\NSIR$ to a proper
Continuous Dimension Scale which contains ``arbitrarily large'' elements.
Thus, let $V(\NSIR)$\index{VzzofNSIR@$V(\NSIR)$} denote the (proper)
Monoid
\index{Monoid (monoid that might be a proper class)} consisting of all
isomorphism classes
$\Delta(A)$ of objects $A\in \NSIR$, with addition induced by
direct sum. (To help keep set-theoretic difficulties at bay, one
might wish to pass from $\NSIR$ to an equivalent skeletal category---a
category in which isomorphic objects are equal---before forming this
Monoid.)\index{Monoid (monoid that might be a proper class)} For any
$A\in\NSIR$, the ideal of $V(\NSIR)$\index{VzzofNSIR@$V(\NSIR)$} generated by
$\Delta(A)$ equals the monoid\index{VzzofR@$V(R)$} $V(A)$; in particular,
this ideal is a set.

\begin{lemma}\label{L:VNSIR}
The Monoid\index{Monoid (monoid that might be a proper class)} $V(\NSIR)$ is
a Continuous Dimension Scale, and $V(R)$ is a generating lower subset
of\index{VzzofNSIR@$V(\NSIR)$} $V(\NSIR)$.
\end{lemma}

\begin{proof} If $S$ is a lower subset of
\index{VzzofNSIR@$V(\NSIR)$} $V(\NSIR)$, then
$S=\setm{\Delta(B_i)}{i\in I}$ for some set
$\setm{B_i}{i\in I}$ of objects from $\NSIR$. Form
$B=\rE(\bigoplus_{i\in I} B_i)$,\index{ErzzE@$\rE(M)$} and observe
that\index{VzzofR@$V(R)$}
$S\subseteq V(B)$. By Lemmas
\ref{L:endinject} and \ref{L:SegDI}, $V(B)$\index{VzzofR@$V(R)$} and $S$ are
continuous dimension scales.\index{continuous dimension scale} For any element
\index{VzzofNSIR@$V(\NSIR)$}
$a=\Delta(A)\in V(\NSIR)$, the class
$(a]$ is contained in the set $V(A)$\index{VzzofR@$V(R)$} and so it is a
set. Thus, Axiom \Mh\ is satisfied in\index{VzzofNSIR@$V(\NSIR)$} $V(\NSIR)$.

Since every object in $\addR$ is injective (being a direct summand
of some injective module $n\cdot R$), we have
$\addR\subseteq\NSIR$ and\index{VzzofR@$V(R)$}\index{VzzofNSIR@$V(\NSIR)$}
$V(R)\subseteq V(\NSIR)$. It is then clear that
$V(R)$\index{VzzofR@$V(R)$} is a lower subset of\index{VzzofNSIR@$V(\NSIR)$}
$V(\NSIR)$. Given any nonzero object $A\in \NSIR$, choose a
nonzero element
$x\in A$. By
\cite[Theorem~9.2]{GvnRR},\index{Goodearl, K.\,R.} the cyclic module $xR$ is
both projective and injective. On the one hand, this means that $xR\in\addR$
and\index{VzzofR@$V(R)$} $\Delta(xR)\in V(R)$, while on the other,
$\Delta(xR)\leq\Delta(A)$. Thus,\index{VzzofR@$V(R)$} $V(R)$ is\index{dense}
dense in \index{VzzofNSIR@$V(\NSIR)$}$V(\NSIR)$, and therefore
$V(\NSIR)$ satisfies Axiom \Ml. \end{proof}

\begin{theorem}\label{T:VNSIR}
Let $R$ be a regular, right self-injective ring,
\index{ring!Von Neumann regular ---}\index{ring!right self-injective ---}%
and write $R= R_{\I}\times R_{\II}\times R_{\III}$ where $R_{\J}$ is of Type J.
Let $\Omega_{\J}$ be the ultrafilter space of $\rB(R_{\J})$. Then
\index{VzzofNSIR@$V(\NSIR)$}
   \[
   V(\NSIR)\cong \CC(\Omega_{\I},\ZZ_\infty)\times
   \CC(\Omega_{\II},\RR_\infty)\times \CC(\Omega_{\III},\two_\infty).
   \]
\end{theorem}

\begin{proof} As observed in the proof of Theorem~\ref{T:isomV(R)},
$\Omega_{\J}$ is homeomorphic to the ultrafilter space
of\index{pzzroj@$\BB{S}$} $\BB V(R)_{\J}$\index{VzzofR@$V(R)$} for $\J=\I$,
$\II$, $\III$. Let
$E$ be a finitary unit of\index{VzzofR@$V(R)$} $V(R)_{\fin}$, and observe
that since $V(R)$\index{VzzofR@$V(R)$} is dense\index{dense}
in\index{VzzofNSIR@$V(\NSIR)$}
$V(\NSIR)$, the set $E$ is dense\index{dense} in
$V(\NSIR)_{\fin}$. Thus, $E$ is a finitary unit of $V(\NSIR)$. By
Theorem~\ref{T:GenEmbDI} and its proof, there is a lower
embedding\index{lower embedding}
   \[
   \varepsilon\colon V(\NSIR)\hookrightarrow \CC(\Omega_{\I},\ZZ_\infty)\times
   \CC(\Omega_{\II},\RR_\infty)\times \CC(\Omega_{\III},\two_\infty)
   \]
(unique with respect to our choice of $E$) such that whenever
$A\in\NSIR$ and $R\in\addA$, the restriction of $\varepsilon$
to\index{VzzofR@$V(R)$}
$V(A)$ matches the embedding given in Theorem~\ref{T:EmbDimInt}.

To see that every function in $\CC(\Omega_{\I},\ZZ_\infty)\times
\CC(\Omega_{\II},\RR_\infty)\times \CC(\Omega_{\III},\two_\infty)$ lies in
the image of $\varepsilon$, it suffices to show that for any
infinite cardinal $\kappa=\aleph_\tau$, the constant function
$t_\kappa$ with
$t_\kappa(x)= \kappa$ for all $x\in \Omega_{\I}\sqcup \Omega_{\II}\sqcup
\Omega_{\III}$ lies in the image of $\varepsilon$.

Set $B= \rE(\aleph_0\cdot R)$,\index{ErzzE@$\rE(M)$} let $\aleph_\sigma$ be
the cardinality of $B$, and set $A= \rE(\aleph_{\sigma+\tau}\cdot B)$.
In particular, $B$ contains no direct sums of more than
$\aleph_\sigma$ nonzero submodules, and\index{VzzofR@$V(R)$} $V(R)\subseteq
V(B)\subseteq V(A)$. We may assume that $B$ is an actual submodule
of $A$. By Lemma~\ref{L:SegDimInt}, restriction from
$V(A)$\index{VzzofR@$V(R)$} to
$V(R)$ provides an\index{pzzroj@$\BB{S}$} isomorphism\index{VzzofR@$V(R)$}
$\BB V(A)\overset{\cong}\longrightarrow \BB V(R)$.

According to Lemma~\ref{L:endinject}, $\Lat(A)$ is an
espalier\index{espalier} whose dimension range\index{dimension range} is
isomorphic to\index{VzzofR@$V(R)$} $V(A)$ (we have $[0,\Delta(A)]= V(A)$
because $A$ is purely infinite). Now $A$ and
$B$ are purely infinite elements of $\Lat(A)$ with central cover
$1$, and $B$ is not equal to any orthogonal sum of more than
$\aleph_\sigma$ nonzero elements. The module-theoretic statement
$\rE(\aleph_{\sigma+\tau}\cdot B)\cong A$,\index{ErzzE@$\rE(M)$} when written
in the symbolism of espaliers,\index{espalier} says that
$\aleph_{\sigma+\tau}\cdot B\sim A$. Thus, Proposition~\ref{P:existconstant}
implies that there exists a purely infinite element $C\in\Lat(A)$ such that
$\mu(\Delta(C))$ equals the constant function with value
$\aleph_\tau$. Therefore we have $\Delta(C)\in V(\NSIR)$
\index{VzzofNSIR@$V(\NSIR)$} with
$\varepsilon(\Delta(C))= t_\kappa$, which completes the proof of
the theorem. \end{proof}

The following corollary is an immediate consequence of
Theorem~\ref{T:VNSIR}, in view of the fact that\index{VzzofR@$V(R)$} $V(A)$
is a lower subset of
$V(\NSIR)$ for any $A\in\NSIR$. If the reader wishes to avoid proper
Continuous Dimension Scales, this result can be proved directly, using
the same methods employed in the theorem.

\begin{corollary}\label{C:VNSIR}
Let $R$ be a regular, right self-injective ring,
\index{ring!Von Neumann regular ---}\index{ring!right self-injective ---}%
and write $R= R_{\I}\times R_{\II}\times R_{\III}$, where $R_{\J}$ is of Type
$\J$. Let $\Omega_{\J}$ be the ultrafilter space of $\rB(R_{\J})$. Given any
ordinal $\gamma$, there exists a nonsingular injective right
$R$-module $A$ such that
\index{Zzzgamma@$\ZZ_\gamma$}\index{Rzzgamma@$\RR_\gamma$}%
\index{Tzzgamma@$\two_\gamma$}\index{VzzofR@$V(R)$}%
   \begin{equation}
   V(A)\cong \CC(\Omega_{\I},\ZZ_\gamma)\times
   \CC(\Omega_{\II},\RR_\gamma)\times
   \CC(\Omega_{\III},\two_\gamma).\tag*{\qed}
   \end{equation}
\renewcommand{\qed}{}
\end{corollary}

To show that every continuous dimension scale\index{continuous dimension scale} appears as a
lower subset of some\index{VzzofR@$V(R)$} $V(A)$, it only remains to
construct regular, right self-injective rings
\index{ring!Von Neumann regular ---}\index{ring!right self-injective ---}%
of Types I, II, III having arbitrary complete
Boolean algebras\index{Boolean algebra!complete ---}
as their Boolean algebras of central idempotents.
We shall make use of the concept of a \emph{maximal quotient ring}
\index{ring!maximal quotient ---}
(see, for example, \cite[Chapter~2]{Gnonsing},
\index{Lam, T.\,Y.}\cite[\S13]{Lam})
\index{Goodearl, K.\,R.} in part of the process.

\begin{proposition}\label{P:existB(R)I,III}
Given any\index{Boolean algebra!complete ---}
complete Boolean algebra $B$, there exist regular, right
self-injective rings
\index{ring!Von Neumann regular ---}\index{ring!right self-injective ---}%
$R_{\I}$ and $R_{\III}$ of Types I and III,
respectively, such that $\rB(R_{\I})\cong \rB(R_{\III})\cong B$.
\end{proposition}

\begin{proof}
The quickest way to obtain a Type I example is to
take $R_{\I}$ to be $B$ itself, made into a ring in the canonical
way. Then $R_{\I}$ is a commutative, regular, self-injective ring
\index{ring!Von Neumann regular ---}\index{ring!right self-injective ---}%
in which all elements are idempotent, and $\rB(R_{\I})\cong B$. For
later use, we note that since $R_{\I}$ is commutative, $B\cong
\rB(R_{\I})\cong \Lat(R_{\I})$. The self-injectivity of $R_{\I}$
implies that $R_{\I}$ is a continuous regular ring,
\index{ring!continuous regular ---}
thus yielding von Neumann's well-known result that $B$ is continuous (see
\cite[Lemma~II.4.10]{GLT2}).\index{Gr\"atzer, G.}
Since $R_{\I}$ is commutative, it is abelian, and hence
is of Type I. As a ring, $B$ has characteristic~$2$, while the
reader may prefer examples having characteristic 0. We can construct
examples which are algebras over any field $F$, as follows.

Let $X$ be the ultrafilter space of $B$, so that $B$ is isomorphic
to the Boolean algebra \index{Boolean algebra}
of clopen subsets of $X$. Let $S$ be the
ring of all locally constant functions from $X$ to $F$ (that is,
functions $f:X\rightarrow F$ such that each point of $X$ has a
neighborhood on which $f$ is constant). Observe that $S$ is a
commutative regular ring, with $B\cong \rB(S)\cong \Lat(S)$. As
noted above, $B$ is a continuous lattice;
\index{lattice!continuous ---} thus~$S$ is a
continuous regular ring.
\index{ring!continuous regular ---}
Finally, let $R_{\I}$ be the maximal\index{ring!maximal right quotient ---}
(right) quotient ring of $S$. Since $S$ is regular, it is a
nonsingular ring, and so $R_{\I}$ is regular and right
self-injective \index{Goodearl, K.\,R.}(\cite[Corollary~2.31]{Gnonsing},
\index{Lam, T.\,Y.}\cite[Theorem~13.36]{Lam}).
\index{ring!Von Neumann regular ---}\index{ring!right self-injective ---}%
Moreover, since $S$ is commutative, so is\index{Lam, T.\,Y.}
$R_{\I}$ (see \cite[Lemma~14.15]{Lam}). Therefore $R_{\I}$ is of Type I.
By \cite[Theorem~13.13]{GvnRR},\index{Goodearl, K.\,R.} all the idempotents
of $R_{\I}$ lie in $S$ (this is not hard to prove directly in the present
case). Therefore
$\rB(R_{\I})= \rB(S)\cong B$.

Similar methods, worked out by Busqu\'e \cite{Busq93},\index{Busqu{\'e}, C.}
can be applied in the Type III case. First choose a commutative, regular,
self-injective ring $R_{\I}$
\index{ring!Von Neumann regular ---}\index{ring!right self-injective ---}%
with $\rB(R_{\I})\cong B$. By
\cite[Theorem~2.5]{Busq93},\index{Busqu{\'e}, C.} there exists a regular,
right self-injective ring $R_{\III}$
\index{ring!Von Neumann regular ---}\index{ring!right self-injective ---}%
of Type III whose center is
isomorphic to $R_{\I}$. Therefore $\rB(R_{\III})\cong B$. \end{proof}

It appears that the constructions used in Proposition~\ref{P:existB(R)I,III}
do not always produce rings of Type II. We approach the Type II existence
problem lattice-the\-o\-ret\-i\-cal\-ly,
\emph{via} von Neumann's \emph{Coordinatization Theorem}
(e.g., \cite[Theorem~14.1]{JvNbook},\index{von Neumann, J.}
\cite[Chapter XI, Satz 3.2]{FMae55}).\index{Maeda, F.}

\begin{proposition}\label{P:existB(R)II}
Given any\index{Boolean algebra!complete ---} complete Boolean algebra $B$,
there exists a regular, right and left self-injective ring $R$ of Type \IIf\
with $\rB(R)\cong B$.
\index{ring!Von Neumann regular ---}
\index{ring!right and left self-injective ---}
\end{proposition}

\begin{proof}
Let $L$ be an irreducible (i.e., indecomposable)
continuous geometry\index{continous geometry} such that the (unique)
dimension function $D$ on $L$ is positive on all nonzero elements of $L$ and
the range of $D$ is the unit interval $[0,1]$. Such continuous
geometries\index{continous geometry} were constructed by von Neumann
\index{von Neumann, J.}
\cite{JvNexamples}. Alternatively, one could choose a simple, regular, right
self-injective ring
\index{ring!Von Neumann regular ---}\index{ring!right self-injective ---}%
$S$ of Type \IIf\ (see \cite[Corollary
11.10]{GoBo} and \cite[Example 10.7, Theorem~10.27]{GvnRR}
\index{Goodearl, K.\,R.}\index{Boyle, A.\,K.}%
for existence) and take
$L=\Lat(S)$. Indecomposability of $L$ then follows from indecomposability of
$S$, and the properties of $D$ follow from those of the unique rank function
$N$ on $S$ (see\index{Goodearl, K.\,R.}
\cite[Corollary 16.15]{GvnRR}), since $D$ is given by the formula
$D(xR)= N(x)$ for $x\in R$.

Next, let $L(B)$ be the (reducible) continuous
geometry\index{continous geometry} constructed from $L$ and~$B$ by Halperin
\index{Halperin, I.}
in \cite[Theorem~1]{Halp}. The center of $L(B)$ (i.e., the sublattice of
neutral elements) is isomorphic to $B$ by \cite[Theorem~2]{Halp},
\index{Halperin, I.} and $L(B)$
contains a sublattice (with the same largest element) isomorphic to $L$
\cite[Remark 2, p.~351]{Halp}.\index{Halperin, I.} For any positive integer
$n$, the largest element
$1\in L$ can be written as the supremum of $n$ independent pairwise
perspective\index{perspectivity@perspectivity ($\sim$)} elements (e.g.,
because there exist elements $x\in L$ with $D(x)= 1/n$), and so the same
occurs in $L(B)$. Consequently,
$L(B)$ has order $n$ (in von Neumann's sense) for all $n$.

In particular, since $L(B)$ has order 4, von Neumann's
Coordinatization Theorem implies that there exists a regular ring
\index{ring!Von Neumann regular ---}
$R$ such that $\Lat(R)\cong L(B)$. Since $\Lat(R)$ is thus a
continuous lattice,\index{lattice!continuous ---} $R$ is a continuous
\index{ring!continuous regular ---} regular
ring. Now $R$ is unit-regular
\cite[Corollary 13.23]{GvnRR},\index{Goodearl, K.\,R.} and hence
perspectivity\index{perspectivity@perspectivity ($\sim$)} in
$\Lat(R)$ is given by module isomorphism \cite[Corollary 4.23]{GvnRR}.
\index{Goodearl, K.\,R.}
Consequently, for each positive integer $n$, the module~$R$ is a direct
sum of $n$ pairwise isomorphic right ideals. In particular, there are no
nonzero central idempotents $e\in R$ such that the ring $eR$ is
abelian, and therefore~$R$ is right and left self-injective
\index{ring!right and left self-injective ---}\index{Goodearl, K.\,R.}%
\cite[Corollary 13.18]{GvnRR}.

Since $R$ is unit-regular, it is directly\index{directly finite} finite
\cite[Proposition~5.2]{GvnRR}.\index{Goodearl, K.\,R.} Hence,
\cite[Theorems 10.13, 10.24]{GvnRR} show that $R\cong \prod_{m=1}^\infty
R_m\times R_{\II}$ where each $R_m$ is of Type I$_m$ and $R_{\II}$ is of
Type \IIf. Since the dimension theory of $\Lat(R_m)$ takes values in
$\set{0,1/m,2/m,\dots,1}$, the module $R_m$ cannot be a direct sum of $m+1$
nonzero pairwise isomorphic right ideals (cf.~\cite[Theorem~10.10]{GoBo}
\index{Goodearl, K.\,R.}\index{Boyle, A.\,K.}%
or \cite[Corollary 11.18]{GvnRR}).\index{Goodearl, K.\,R.} Thus, all $R_m=0$,
and $R\cong R_{\II}$ is of Type \IIf.

Finally, since the center of $\Lat(R)$ is isomorphic to $\rB(R)$
\cite[Chapter VI, S\"atze 1.9, 3.5]{FMae55},\index{Maeda, F.} we conclude
that $\rB(R)\cong B$.
\end{proof}

\begin{theorem}\label{T:arbVNSIR}
Let $\Omega_{\I}$, $\Omega_{\II}$, $\Omega_{\III}$ be arbitrary complete
Boolean spaces\index{Boolean space} \pup{possibly empty}.
Then there exists a regular, right
self-injective ring $R$
\index{ring!Von Neumann regular ---}\index{ring!right self-injective ---}%
such that\index{VzzofNSIR@$V(\NSIR)$}
   \[
   V(\NSIR)\cong \CC(\Omega_{\I},\ZZ_\infty)\times
   \CC(\Omega_{\II},\RR_\infty)\times \CC(\Omega_{\III},\two_\infty).
   \]
\end{theorem}

\begin{proof} Propositions \ref{P:existB(R)I,III} and
\ref{P:existB(R)II}, together with Theorem~\ref{T:VNSIR}.
\end{proof}

Although the following result is a corollary of Theorem~\ref{T:arbVNSIR}, we
give an independent proof avoiding the use of proper Continuous Dimension
Scales.

\begin{theorem}\label{T:L(R)Duniv}
The class of espaliers\index{espalier} of the form
$(\Lat(R),{\subseteq},{\perp},{\cong})$, for regular, right
self-injective rings $R$,
\index{ring!Von Neumann regular ---}\index{ring!right self-injective ---}%
is D-universal.\index{D-universal}
\end{theorem}

\begin{proof}
Let $S$ be an arbitrary continuous dimension scale.\index{continuous dimension scale} By
Propositions \ref{P:existB(R)I,III} and
\ref{P:existB(R)II}, together with Theorems \ref{T:EmbDimInt} and
\ref{T:isomV(R)}, there exists a regular, right self-injective
ring $R'$
\index{ring!Von Neumann regular ---}\index{ring!right self-injective ---}%
such that $S$ is isomorphic to a lower subset\index{VzzofR@$V(R)$}
$S'$ of $V(R')$. Set $A=\rE(\aleph_0\cdot R')$\index{ErzzE@$\rE(M)$} and
$R=\End_{R'}(A)$. By Lemma~\ref{L:endinject}, $R$ is a regular, right
self-injective ring and\index{Dzzrng@$\Drng L$} $\Drng\Lat(R)\cong
\Drng \Lat(A)\cong [0,\Delta(A)]$. Since all finitely generated
projective right $R'$-modules are isomorphic to direct summands of
$A$, we see that $V(R')$\index{VzzofR@$V(R)$} is a lower subset of
$[0,\Delta(A)]$. Thus,
$S'$ is isomorphic to a lower subset of the dimension range\index{dimension
range} of the espalier\index{espalier} $\Lat(R)$.
\end{proof}

The results above also allow us to determine the monoids
$V(R)$\index{VzzofR@$V(R)$} in the present context, as follows.

\begin{corollary}\label{C:determineV(R)}
Let $M$ be a commutative monoid.\index{VzzofR@$V(R)$}
Then $M\cong V(R)$ for some regular,
right self-injective ring $R$
\index{ring!Von Neumann regular ---}\index{ring!right self-injective ---}%
if and only if $M$ is a continuous dimension scale\index{continuous dimension scale} containing
an order-unit.
\end{corollary}

\begin{proof}
If $R$ is a regular, right self-injective ring,
\index{ring!Von Neumann regular ---}\index{ring!right self-injective ---}%
then
$\Delta(R)$ is an order-unit in\index{VzzofR@$V(R)$} $V(R)$, and $V(R)$ is a
continuous dimension scale\index{continuous dimension scale} by Proposition~\ref{P:L(R)}.
Conversely, let~$M$ be a continuous dimension scale\index{continuous dimension scale} which
contains an order-unit $u$. By Theorem~\ref{T:L(R)Duniv} and
Proposition~\ref{P:L(R)}, there exists a regular, right self-injective ring
$R'$ such that~$M$ is isomorphic to a lower subset $M'$
of\index{VzzofR@$V(R)$} $V(R')$. Let $u'$ denote the image of
$u$ under this isomorphism; then $M'$ equals the ideal
of\index{VzzofR@$V(R)$} $V(R')$ generated by $u'$. Now $u'= \Delta(A)$ for
some
$A\in\addR'$, and it is clear that\index{VzzofR@$V(R)$} $M'=V(A)$. By
Lemma~\ref{L:endinject},
$R= \End_{R'}(A)$ is a regular, right self-injective ring
\index{ring!Von Neumann regular ---}\index{ring!right self-injective ---}%
and\index{VzzofR@$V(R)$} $V(R)\cong V(A)\cong M$.
\end{proof}

\section{Projection lattices of W*- and AW*-algebras}\label{S:ProjAW*}

Our main references for W*-algebras will be the texts by
J. Dixmier \cite{Dixbook},\index{Dixmier, J.} R.\,V. Kadison and J.\,R.
Ringrose \index{Kadison, R.\,V.}\index{Ringrose, J.\,R.}%
\cite{KadRing}, \index{Li, B.-R.} B.-R. Li \cite{BingRen}, and S. Sakai
\index{Sakai, S.}
\cite{Sakai}; for AW*-al\-ge\-bras, we rely on the text by\index{Berberian,
S.\,K.} S.\,K. Berberian \cite{Berb} and the monograph by I. Kaplansky
\cite{Kapl68}.\index{Kaplansky, I.} A
\emph{W*-algebra}\index{W*-algebra|ii} (also called a
\emph{von Neumann algebra}) can be defined as any C*-algebra which is
isomorphic (\emph{qua} C*-algebra) to a *-subalgebra of $\calB(H)$ (the
algebra of all bounded linear operators) on some (complex) Hilbert space $H$
which is closed in the strong operator topology (the topology of pointwise
convergence). Kaplansky introduced the concept of an
\emph{AW*-algebra} (abbreviating ``abstract W*-algebra'')
in order to obtain a more general class of C*-algebras defined (and
analyzed) by purely algebraic properties. Before giving the definition, we
recall a few basic concepts.
\medskip
\begin{quote}
All W*- and AW*-algebras that we consider here will be
assumed to be \emph{unital}.
\end{quote}
\medskip

A \emph{projection} in a C*-algebra $A$ is any self-adjoint
idempotent, that is, any element $p\in A$ with $p=p^2=p^*$. The \emph{right
annihilator} of a subset $S$ of $A$ is the right ideal
\index{azznr@$\annr(S)$|ii}
   \[
   \annr(S)=\setm{x\in A}{sx=0\text{ for all }s\in S}.
   \]
Finally, $A$ is said to be an \emph{AW*-algebra}\index{AW*-algebra|ii} if the
right annihilator of any subset $S$ of $A$ is a principal right ideal
generated by a projection, that is, $\annr(S)= pS$ for some (necessarily
unique) projection $p\in A$. Every W*-algebra is AW*
\index{Berberian, S.\,K.} \cite[\S4, Proposition~9]{Berb}, but not
conversely. For example, the (unital) commutative AW*-algebras are precisely
(up to isomorphism) the algebras
$\CC(X,\CCx)$ of continuous complex-valued functions on complete Boolean
spaces\index{Boolean space} $X$ \index{Berberian, S.\,K.}
\cite[\S7, Theorem~1]{Berb}; such an algebra is W* if and only if $X$ is
\emph{hyperstonian} \index{Dixmier, J.}
(\cite[Th\'eor\`eme~2]{Dix51}; cf.~\cite[Theorems 5.3.3,
5.3.4]{BingRen}).\index{Li, B.-R.} By definition, $X$ is
hyperstonian\index{hyperstonian|ii}
(cf.~\cite[D\'efinition~3]{Dix51}),\index{Dixmier, J.} if for any nonempty
open subset $U$ of $X$, there exists a Radon measure $\mu$ on $X$, vanishing
on all nowhere dense subsets of $X$, such that $\mu(U)>0$.

The set $L$ of projections of an AW*-algebra $A$ is
equipped with the partial ordering~$\leq$ defined by $p\leq q$ \iff\ $pq=p$
(equivalently, $qp=p$), for $p$, $q\in L$. 
The poset $(L,\leq)$ is a complete lattice
\cite[\S4, Proposition~1]{Berb}. Furthermore, two projections
$p$, $q\in L$ are \emph{orthogonal}, in symbols $p\perp q$, if $pq=0$
(equivalently, $qp=0$). Then the sum $p+q$ is also a projection, and it is
the join of $\set{p,q}$ in $L$: hence $p\oplus q=p+q$. Finally, two
projections $p$ and $q$ of $A$ are \emph{Murray-von~Neumann equivalent}, in
symbols $p\sim q$, if there exists $x\in A$ such that $p=xx^*$ and $q=x^*x$.
Equivalently, $pA$ and $qA$ are isomorphic as right $A$-modules, that
is, there are $x$, $y\in A$ such that $p=xy$ and $q=yx$---this equivalence is
nontrivial and contained in \cite[Theorem~27]{Kapl68}.\index{Kaplansky, I.}

A projection $p\in A$ is said to be \emph{$\sigma$-finite} (or
\emph{countably decomposable}, or \emph{orthoseparable}) if $p$ does not
majorize any uncountable orthogonal family of nonzero projections; if the
projection $1\in A$ has this property, then the algebra $A$ itself is called
$\sigma$-finite. This same terminology is also used for Boolean algebras
\index{Boolean algebra} and their elements. Let us say that a Boolean algebra
$B$ is \emph{locally $\sigma$-finite} provided every element of $B$
is a supremum of $\sigma$-finite elements. Furthermore, a Boolean
space\index{Boolean space}
$X$ is locally $\sigma$-finite, if its Boolean
algebra\index{Boolean algebra} of clopen subsets is locally
$\sigma$-finite.

\begin{proposition}\label{P:Hyp2sF}
Every hyperstonian Boolean space\index{Boolean space} is locally
$\sigma$-finite.
\end{proposition}

\begin{proof}
Let $B$ denote the ultrafilter space of a hyperstonian Boolean
space\index{Boolean space}~$X$. For a Radon measure $\mu$ on $X$, we say
that a Borel subset $A$ of $X$ is \emph{$\mu$-self-supporting},
\index{self-supporting (subset)|ii} if $\mu(A\cap U)>0$ whenever
$U\subseteq X$ is open and $A\cap U\neq\es$. Then every Borel subset $A$
of $X$ of positive measure contains a $\mu$-self-supporting compact subset
$K$ of positive measure, see \cite[\S1.9]{Frem}.\index{Fremlin, D.\,H.} If
$\mu$ vanishes on all nowhere dense subsets, then
$\mu(\overset{\,\circ}{K})=\mu(K)>0$, hence
$\overset{\,\circ}{K}$ is a $\mu$-self-supporting open subset of $K$ with
positive measure. As $\mu$ is a Radon measure, $K$ contains a
$\mu$-self-supporting \emph{clopen} subset with positive measure.

Let $D$ denote the set of elements of $B$ whose associated clopen set is
$\mu$-self-supporting with respect to some finite Radon measure $\mu$ on
$X$ vanishing on all nowhere dense subsets. It follows from the
assumption that $X$ is hyperstonian\index{hyperstonian} and the paragraph
above that every element of $B$ is a supremum of elements of $D$. But every
element of $D$ is clearly $\sigma$-finite.
\end{proof}

The converse of Proposition~\ref{P:Hyp2sF} does not hold as a rule.

\begin{example}\label{Ex:SFnonHS}
As in Section~\ref{S:AMeasTh}, we denote by $\Co$ the Boolean algebra
\index{Boolean algebra!complete ---}
of all Borel subsets of the Cantor space $\set{0,1}^\omega$ modulo meager
sets. Let $X$ denote the ultrafilter space of
$\Co$. Then $X$ is clearly $\sigma$-finite. However, there is no
nontrivial Radon measure on $X$, as shown, for example, by the argument on
pages~82--83 in \cite[Chapter~21]{Oxto}.\index{Oxtoby, J.\,C.} In
particular, $X$ is not hyperstonian.\index{hyperstonian}
\end{example}

On the other hand, the ``measure'' analogue of $C_{\omega}$, that is, the
random algebra~$B_{\omega}$ (see Section~\ref{S:AMeasTh}) has, of course,
hyperstonian\index{hyperstonian} ultrafilter space.

\begin{proposition}\label{P:AW*esp}
Let $L$ be the lattice\index{lattice} of projections of
an\index{AW*-algebra} AW*-algebra $A$, endowed with the relations $\leq$,
$\perp$, and $\sim$ defined above. Then
$L$ is an espalier.\index{espalier}
\end{proposition}

\begin{proof}
Axiom (L1) follows from \cite[Theorem~19]{Kapl68}\index{Kaplansky, I.} or
\index{Berberian, S.\,K.}\cite[\S4, Proposition~1]{Berb}. Axioms~(L2) and
(L3) are easy exercises, Axiom~(L5) is trivial.

Axiom~(L4): for $p$, $r\in L$, the element $1-p\in A$ is also a
projection, and $p\perp r$ is equivalent to $r\leq 1-p$. Now let
$\famm{q_i}{i\in I}$ be an orthogonal family of elements of~$L$ such that
$p\perp(\oplus_{i\in J}q_i)$, for all finite $J\subseteq I$. This means that
$\oplus_{i\in J}q_i\leq 1-p$, for all finite $J\subseteq I$, thus
$\oplus_{i\in I}q_i\leq 1-p$, that is, $p\perp(\oplus_{i\in I}q_i)$.

Axiom~(L6) is Axiom~(C) in\index{Kaplansky, I.} \cite[Chapter~4]{Kapl68}; see
\cite[Theorem~24]{Kapl68} or \index{Berberian, S.\,K.}
\cite[\S1, Proposition~9]{Berb}.

Axioms~(L7) and (L8) are difficult results, proved in
\cite[Theorem~52, 62]{Kapl68}\index{Kaplansky, I.} and \cite[\S20,
Theorem~1, \S13, Theorem~1]{Berb}.\index{Berberian, S.\,K.}
\end{proof}

The ``projections'' of the espalier\index{espalier} $L$ are not the
projections of
$A$ (which are the elements of $L$), but they correspond to the
\emph{central} projections of $A$. In fact, all central idempotents of $A$
are projections \index{Berberian, S.\,K.}\cite[\S3, Exercise 1]{Berb}, and so
we may use without ambiguity the notation $\rB(A)$ of Section~\ref{S:RSIReg}
to stand for the Boolean algebra\index{Boolean algebra!complete ---} of
central projections in $A$.

\begin{lemma}\label{L:projoflattproj}
Let $L$ be the lattice\index{lattice} of projections of an
\index{AW*-algebra}AW*-algebra $A$. For each $e\in\rB(A)$, there is a
projection\index{pzzroj@$\BB{S}$} $\pi_e\in \BB \Drng L$
\index{Dzzrng@$\Drng L$}such that
$\pi_e(\Delta(p))= \Delta(ep)$ for all $p\in L$. The rule $e\mapsto \pi_e$
defines an isomorphism of $\rB(A)$ onto\index{pzzroj@$\BB{S}$} $\BB\Drng L$.
\end{lemma}

\begin{proof} Set $S= \Drng L= L/{\sim}$. It is clear that for each $e\in
\rB(A)$, there is a projection $\pi_e\in \BB S$\index{pzzroj@$\BB{S}$} as
described, and
$\pi_{1-e}= \pi_e^\perp$. It is also clear that $e\leq f$ implies
$\pi_e\leq \pi_f$, for $e,f\in \rB(A)$. On the other hand, if $e\nleq f$,
the projection $g= e(1-f)$ is nonzero. Note that $\pi_g(\Delta(g))=
\Delta(g)\ne 0$, whence $\pi_g\ne 0$. Since $\pi_g\leq \pi_e\wedge
\pi_f^\perp$, it follows that $\pi_e \nleq \pi_f$. Thus, the map
$\rB(A)\rightarrow \BB S$ given by $e\mapsto \pi_e$ is an order-embedding
that respects complements. It only remains to show that this map is
surjective.

Given\index{pzzroj@$\BB{S}$} $q\in\BB S$, we have
$\Delta(1)=q(\Delta(1))+q^\perp(\Delta(1))$, and so there exist orthogonal
projections $e,f\in L$ such that $1= e\oplus f$ while $\Delta(e)=
q(\Delta(1))$ and $\Delta(f)= q^\perp(\Delta(1))$. Further, $\Delta(e)\perp
\Delta(f)$, and so there is no nonzero projection
$p\in L$ such that $p\lesssim e,f$. We next show that $e$ and $f$ are central
projections. Since $fAe= (eAf)^*$, it is enough to show that $eAf=0$. Let
$x$ be an arbitrary element of $eAf$, and let
$p_r$ and $p_l$ be the \emph{right} and \emph{left projections} of $x$,
respectively \index{Berberian, S.\,K.}(\cite[\S3, Definition 4]{Berb},
\cite[p.~28]{Kapl68}),\index{Kaplansky, I.} that is, the unique projections
such that $p_r^\perp$ and $p_l^\perp$ generate, respectively, the right and
left annihilators of
$x$. Since $xe=fx=0$, we find that $e\le p_r^\perp$ and
$f\le p_l^\perp$, that is, $p_r\le e^\perp= f$ and $p_l\le e$. However,
$p_r\sim p_l$ \index{Berberian, S.\,K.}(see
\index{Berberian, S.\,K.}\cite[\S20, Theorem~3]{Berb} or
\cite[Theorem~63]{Kapl68}),\index{Kaplansky, I.} whence
$p_r\lesssim e,f$ and so $p_r=0$. Thus $x=xp_r=0$, proving that $eAf=0$,
as desired. Consequently, $e$ and $f$ are central, as claimed. Now
for any $r\in L$, we have $\Delta(r)= \Delta(er)+ \Delta(fr)$ with
$\Delta(er)\leq q(\Delta(1))$ and $\Delta(fr)\leq q^\perp(\Delta(1))$,
whence $\Delta(er)\in q(S)$ and $\Delta(fr)\in q^\perp(S)$.
Hence, we conclude that $\pi_e(\Delta(r))= \Delta(er)= q(\Delta(r))$.
Therefore $\pi_e= q$, completing the proof.
\end{proof}

In the context of Proposition~\ref{P:AW*esp}, let us denote by $[p]$ the
$\sim$-equivalence class of a projection $p$ of $A$.
The \emph{addition} of these equivalence classes is defined by
   \[
   [p]+[q]=[p\oplus q]=[p+q],
   \text{ for any orthogonal projections }p\text{ and }q.
   \]
The dimension range\index{dimension range} of $L$ is, of
course,\index{Dzzrng@$\Drng L$} $\Drng L=L/{\sim}$, equipped with the above
partial addition.

Just as in Section~\ref{S:RSIReg}, we can define the
monoid\index{VzzofR@$V(R)$} $V(A)$ of isomorphism classes of finitely
generated projective right $A$-modules. As noted above, projections $p$,
$q\in A$ satisfy $p\sim q$ if and only if
$pA\cong qA$, and so we obtain an embedding of partial
monoids,\index{Dzzrng@$\Drng L$} $\Drng L
\hookrightarrow V(A)$, where $[p]\mapsto \Delta(pA)$. Under this embedding,
$[1]\mapsto \Delta(A)$. Any direct summand of the right module $A$ has the
form $eA$ for an idempotent $e$, and since $e$ is equivalent to a projection
$p\in A$ \cite[Theorem~26]{Kapl68},\index{Kaplansky, I.} we have $pA\cong eA$
and so $[p]\mapsto
\Delta(eA)$. Similarly, any pair of orthogonal idempotents in $A$ is
equivalent to a pair of orthogonal projections, so that any pair of elements
$u$, $v\in V(A)$\index{VzzofR@$V(R)$} such that $u+v\leq\Delta(A)$ must be
the image of a pair of elements of\index{Dzzrng@$\Drng L$}
$\Drng L$ whose sum is defined. Thus, the embedding above
maps\index{Dzzrng@$\Drng L$}
$\Drng L$ isomorphically onto the interval
$[0,\Delta(A)]\subseteq V(A)$, which we record in the theorem below.

\begin{theorem}\label{T:DimAW*}
Let $L$ be the lattice\index{lattice} of projections of
an\index{AW*-algebra} AW*-algebra $A$. Then the dimension
range\index{dimension range}
$\Drng L=L/{\sim}$\index{Dzzrng@$\Drng L$} is a
\pup{bounded} continuous dimension scale,\index{continuous dimension scale}
and\index{VzzofR@$V(R)$}\index{Dzzrng@$\Drng L$}
$\Drng L\cong[0,\Delta(A)]\subseteq V(A)$. If $A$ is a\index{W*-algebra}
W*-algebra, then the ultrafilter space
of\index{pzzroj@$\BB{S}$}\index{Dzzrng@$\Drng L$}
$\BB\Drng L$ is hyperstonian.\index{hyperstonian}
\end{theorem}

\begin{proof}
That $\Drng L$ is a continuous dimension scale\index{continuous dimension scale} follows from
Theorem~\ref{T:DimEsp}. We have just seen above that\index{Dzzrng@$\Drng L$}
$\Drng L\cong [0,\Delta(A)]$. Observe that\index{pzzroj@$\BB{S}$}
$\BB\Drng L\cong\rB(A)\cong \rB(Z)$, where $Z$ is the center of $A$. If $A$
is a\index{W*-algebra} W*-algebra, then so is $Z$, whence
$Z\cong \CC(X,\CCx)$ for some hyperstonian\index{hyperstonian} complete
Boolean space\index{Boolean space} $X$. In particular, the ultrafilter space
of\index{pzzroj@$\BB{S}$} $\BB\Drng L$\index{Dzzrng@$\Drng L$} is
homeomorphic to $X$ and thus it is hyperstonian.\index{hyperstonian}
\end{proof}

In the context of Theorem \ref{T:DimAW*}, observe that by
Proposition~\ref{P:Hyp2sF}, $\BB\Drng L$
\index{pzzroj@$\BB{S}$}\index{Dzzrng@$\Drng L$}%
is locally $\sigma$-finite in case $A$ is a\index{W*-algebra} W*-algebra.

In particular, when $L$ is the lattice\index{lattice} of projections of an
\index{AW*-algebra}AW*-algebra $A$, it follows from
Theorem~\ref{T:EmbDimInt} that the
\pcm\
$L/{\sim}$ embeds as a lower subset into a \cm\ of the form
\index{Zzzgamma@$\ZZ_\gamma$}\index{Rzzgamma@$\RR_\gamma$}%
\index{Tzzgamma@$\two_\gamma$}%
   \[
   \CC(\Omega_{\I},\ZZ_\gamma)\times\CC(\Omega_{\II},\RR_\gamma)
   \times\CC(\Omega_{\III},\two_\gamma),
   \]
for complete Boolean spaces\index{Boolean space}
$\Omega_{\I}$, $\Omega_{\II}$, $\Omega_{\III}$. Theorem \ref{T:DimAW*}
implies that these spaces must be hyperstonian\index{hyperstonian} in case
$A$ is a\index{W*-algebra} W*-algebra.

There exists a Type I, II, III decomposition for\index{AW*-algebra}
AW*-algebras (see \cite{Kapl68})\index{Kaplansky, I.} which parallels that
for regular, right self-injective rings;
\index{ring!Von Neumann regular ---}\index{ring!right self-injective ---}%
in fact, Kaplansky developed much of the Type I, II, III theory for
\emph{Baer rings}
\index{ring!Baer ---}
(rings in which the right or left annihilator of any
element is generated by an idempotent), a class of rings which includes both
AW*-algebras\index{AW*-algebra} and regular, right self-injective rings.
\index{ring!Von Neumann regular ---}\index{ring!right self-injective ---}%
We shall use some of
the terminology and results of this theory without explicit references. We
point out that an AW*-algebra\index{AW*-algebra} $A$ is called a
\emph{factor} provided the center of $A$ equals the complex field $\CCx$;
equivalently, $A$ is a factor if and only if $A$ is nonzero and $\rB(A)=
\set{0,1}$.

\begin{lemma}\label{L:bigW*factors}
Let $\gamma$ be an ordinal and $\J\in \{\I,\,\II,\,\III\}$.
There exists a W*-fac\-tor $A_{\J}$ of Type $\J$ which contains a
family $\famm{p^{\J}_\alpha}{\alpha\leq\gamma}$ of nonzero purely infinite
projections such that $p^\J_\alpha \lesssim p^\J_\beta$ but $p^\J_\beta
\not\lesssim p^\J_\alpha$ for all ordinals $\alpha<\beta\leq \gamma$.
\end{lemma}

\begin{proof}
Choose a Hilbert space $H_\gamma$ with an
orthonormal basis of cardinality~$\aleph_\gamma$, and set $A_{\I}=
\calB(H_\gamma)$. For each ordinal $\alpha\leq\gamma$, choose a projection
$p^\I_\alpha \in A_\I$ such that the closed subspace $p^\I_\alpha H_\gamma$
of $H_\gamma$ has an orthonormal basis of cardinality~$\aleph_\alpha$.
The desired properties of
$A_{\I}$ and the
$p^{\I}_\alpha$ are clear.

Next, choose W*-factors $B_{\II}$ and $B_{\III}$ of Types II and
III (e.g., \cite[Part I, \S9.4]{Dixbook},\index{Dixmier, J.}
\index{Kadison, R.\,V.}\index{Ringrose, J.\,R.}%
\cite[Chapters 6, 8]{KadRing},\index{Sakai, S.} \cite[Chapter 4]{Sakai}).
These factors can be chosen as subalgebras of
$\calB(H_0)$ for a separable Hilbert space $H_0$ (e.g.,
\cite[Remark, p.~155]{Dixbook},\index{Dixmier, J.}
\cite[Theorem~7.3.16]{BingRen}),\index{Li, B.-R.} so that they are
$\sigma$-finite\cite[Proposition~1.14.3]{BingRen}. Now let
$A_{\II}$ and $A_{\III}$ be the W*-tensor products
\index{abzzarotimes@$A\barotimes B$|ii}
$B_{\II}\barotimes A_{\I}$ and $B_{\III} \barotimes A_{\I}$. These algebras
are of Types II and III, respectively (e.g.,
\cite[Propositions~11.2.21, 11.2.26]{KadRing},
\index{Kadison, R.\,V.}\index{Ringrose, J.\,R.}\index{Sakai, S.}%
\cite[Proposition~2.6.3, Theorem~2.6.4]{Sakai}), and they are factors
\cite[Proposition~2.6.7]{Sakai}.\index{Sakai, S.}

Now let $\J=\II$ or $\III$, and set $p^\J_\alpha= 1\otimes p^\I_\alpha \in
A_\J$ for all ordinals $\alpha\leq\gamma$. It is clear that these
$p^\J_\alpha$ are purely infinite projections, and that $p^\J_\alpha
\lesssim p^\J_\beta$ for all ordinals $\alpha<\beta\leq\gamma$. Observe
that the\index{W*-algebra} W*-algebra $p^\J_\alpha A_\J p^\J_\alpha$ is
isomorphic to
$B_\J\barotimes p^\I_\alpha A_\I p^\I_\alpha$, which is in turn isomorphic
to a W*-subalgebra of $\calB(H_0\otimes p^\I_\alpha H_\gamma)$. Since
$H_0\otimes p^\I_\alpha H_\gamma$ has an orthonormal basis of cardinality
$\aleph_\alpha$, we thus see that $p^\J_\alpha$ does not majorize any
orthogonal family of more than $\aleph_\alpha$ nonzero projections. On the
other hand, $\aleph_\beta\cdot p^\I_\beta \sim p^\I_\beta$, whence
$p^\J_\beta$ majorizes an orthogonal family of $\aleph_\beta$ nonzero
projections (equivalent to itself). Therefore $p^\J_\beta \not\lesssim
p^\J_\alpha$.
  \end{proof}

We can now show that the class of projection lattices\index{lattice} of
\index{W*-algebra}W*-algebras, while not D-universal,\index{D-universal} is
at least D-universal\index{D-universal} relative to continuous dimension scales\index{continuous dimension scale} for which the ultrafilter space of the
Boolean algebra\index{Boolean algebra!complete ---} of projections is
hyperstonian.\index{hyperstonian}

\begin{theorem}\label{T:bigW*Drng}
Let $\Omega_{\I}$, $\Omega_{\II}$, $\Omega_{\III}$ be arbitrary
hyperstonian\index{hyperstonian} spaces \pup{possibly empty}, and let
$\gamma$ be an arbitrary ordinal. Then there exists a\index{W*-algebra}
W*-algebra $A$ such that
\index{Zzzgamma@$\ZZ_\gamma$}\index{Rzzgamma@$\RR_\gamma$}%
\index{Tzzgamma@$\two_\gamma$}\index{VzzofR@$V(R)$}%
   \[
   V(A)\cong \CC(\Omega_{\I},\ZZ_\gamma)\times
   \CC(\Omega_{\II},\RR_\gamma)\times \CC(\Omega_{\III},\two_\gamma).
   \]
\end{theorem}

\begin{proof} For $\J= \I,\,\II,\,\III$, choose W*-factors $A_\J$ and
families $\famm{p^\J_\alpha}{\alpha\leq\gamma}$ of purely infinite
projections as in Lemma~\ref{L:bigW*factors}. Let $C_{\J}=
\CC(\Omega_{\J},\CCx)$, which is a\index{W*-algebra} W*-algebra because
$\Omega_{\J}$ is hyperstonian,\index{hyperstonian} and note that
$\rB(C_{\J})$ is isomorphic to the Boolean algebra
\index{Boolean algebra!complete ---} of clopen subsets of
$\Omega_{\J}$. Let $D_{\J}$ be the W*-tensor product $A_{\J} \barotimes
C_{\J}$, which has Type J by the results referenced in
Lemma~\ref{L:bigW*factors}. Since $A_{\J}$ is a factor, the centers of
$C_{\J}$ and $D_{\J}$ are isomorphic\index{Sakai, S.}
\cite[Proposition~2.6.7]{Sakai}, and thus $\rB(C_{\J})\cong \rB(D_{\J})$,
\emph{via} the map $e\mapsto 1\otimes e$. Consequently, if
$L_{\J}$ is the lattice\index{lattice} of projections of~$D_{\J}$, the
ultrafilter space of\index{pzzroj@$\BB{S}$}
$\BB\Drng L_{\J}$\index{Dzzrng@$\Drng L$} is homeomorphic to $\Omega_{\J}$.

Since $D_{\J}$ is of Type J, it follows that $L_{\J}/{\sim}$ is of
Type J, that is, $L_{\J}/{\sim}= (L_{\J}/{\sim})_{\J}$ in the notation of
Definition \ref{D:SI,II,III}. Set
$q^\J_\alpha= p^\J_\alpha \otimes 1\in D_\J$ for all ordinals
$\alpha\leq\gamma$, and observe that the $q^\J_\alpha$ are purely infinite
projections with central cover $1$, such that $q^\J_\alpha \lesssim
q^\J_\beta$ for all ordinals $\alpha<\beta\leq\gamma$.

\begin{sclaim} For any ordinals $\alpha<\beta\leq\gamma$, we have
$rq^\J_\beta \not\lesssim rq^\J_\alpha$ for all nonzero central projections
$r$ in $D_\J$.
\end{sclaim}

\begin{scproof} For each point
$x\in\Omega_{\J}$, let $\pi_x\colon D_{\J}\rightarrow A_{\J} \barotimes \CCx
\cong A_{\J}$ be the W*-al\-ge\-bra homomorphism obtained by
tensoring the identity map on $A_{\J}$ with the evaluation map
$f\mapsto f(x)$ from $C_{\J}$ to $\CCx$. Observe that $\pi_x(q^\J_\alpha)=
p^\J_\alpha$ and $\pi_x(q^\J_\beta)=
p^\J_\beta$. Moreover, $r=1\otimes e$ for some projection $e\in
\CC(\Omega_\J,\CCx)$, and $\pi_x(r)= e(x)\in\set{0,1}$. If $rq^\J_\beta
\lesssim rq^\J_\alpha$, then $e(x)p^\J_\beta \lesssim e(x)p^\J_\alpha$ for
all $x\in \Omega_\J$. Since $p^\J_\beta \not\lesssim p^\J_\alpha$, we must
have $e(x)=0$ for all $x\in \Omega_\J$, and thus $r=0$. This contradiction
establishes the claim.
\end{scproof}

We now apply Proposition~\ref{P:secondexistconstant},
and conclude that
there exist projections $r_{\J}\in L_{\J}$ for each J such that
the dimension ranges\index{dimension range} of the intervals $[0,r_{\J}]$
have the following form:
\index{Zzzgamma@$\ZZ_\gamma$}\index{Rzzgamma@$\RR_\gamma$}%
\index{Tzzgamma@$\two_\gamma$}%
\begin{align*}
[0,r_{\I}]/{\sim}\ &\cong\ \CC(\Omega_{\I},\ZZ_\gamma) \\
[0,r_{\II}]/{\sim}\ &\cong\ \CC(\Omega_{\II},\RR_\gamma) \\
[0,r_{\III}]/{\sim}\ &\cong\ \CC(\Omega_{\III},\two_\gamma).
\end{align*}
Therefore the dimension range\index{dimension range} of the
lattice\index{lattice} of projections of the\index{W*-algebra} W*-algebra
   \[
   A= r_{\I}D_{\I}r_{\I}\times r_{\II}D_{\II}r_{\II}\times
   r_{\III}D_{\III}r_{\III}
   \]
has the desired form. Note that each of the projections $r_{\J}$ is purely
infinite, whence the projection $1\in A$ is purely infinite, and
consequently\index{VzzofR@$V(R)$} $[0,\Delta(A)]= V(A)$.  Therefore, in view
of Theorem~\ref{T:DimAW*}, the present theorem is proved.
\end{proof}

\begin{corollary}\label{C:W*hyperDuniv}
Let $S$ be a continuous dimension scale.
\index{continuous dimension scale} Then
$S$ admits a lower embedding\index{lower embedding} into the dimension
range\index{dimension range} of the lattice\index{lattice} of projections of
some\index{W*-algebra} W*-algebra if and only if the ultrafilter space of
$\BB S$ \index{pzzroj@$\BB{S}$} is hyperstonian.\index{hyperstonian}
\end{corollary}

\begin{proof}
The sufficiency follows from Theorems
\ref{T:EmbDimInt} and \ref{T:bigW*Drng}, and the necessity from Theorem
\ref{T:DimAW*}.
\end{proof}

In order to see that the projection lattices\index{lattice}
of\index{AW*-algebra} AW*-algebras form a D-universal\index{D-universal}
class of espaliers,\index{espalier} we need an analogue of Theorem
\ref{T:bigW*Drng} in which $\Omega_{\I}$, $\Omega_{\II}$, $\Omega_{\III}$
are arbitrary complete Boolean spaces\index{Boolean space} and $A$ is an
\index{AW*-algebra}AW*-algebra. However, there is no general theory of
AW*-tensor products available to replace the W*-tensor products $A_\J
\barotimes C_\J$ used in our proof. P.~Ara has suggested that one might be
able to use the \emph{monotone complete tensor products}
\index{monotone complete tensor products}
introduced by M.~Hamana
\cite{Hamana82I, Hamana82II}\index{Hamana, M.} instead. (We thank him for
making us aware of Hamana's work.) Rather than developing the necessary
auxiliary results about monotone complete tensor products here, we complete
the picture by taking a different route. Namely, we borrow the methods and
results of G. Takeuti\index{Takeuti, G.}
\cite{Take83} and some of the subsequent results obtained in M. Ozawa
\index{Ozawa, M.}
\cite{Ozaw86}. These methods involve forcing, more specifically, the
Scott-Solovay model $V^B$\index{VzzVB@$V^B$}%
of $B$-valued set theory (also
used in Section~\ref{S:AMeasTh}), for any
\index{Boolean algebra!complete ---} complete Boolean algebra $B$.

We give a short summary of what we shall use from\index{Ozawa, M.}
\index{Takeuti, G.} \cite{Take83,Ozaw86}. If $\lA$ is an\index{AW*-algebra}
AW*-algebra in\index{VzzVB@$V^B$} $V^B$, the
\emph{bounded global section algebra} $\tlA$ of $\lA$
\index{bounded global section algebra\ii}\index{Abzzglobs@$\tlA$|ii}%
is the set of all $\lx\in V^B$\index{VzzVB@$V^B$} such that
$\bv{\lx\in\lA}=1$ and
$\bv{\bv{\lx}_{\lA}\leq\check{n}1_{\lA}}=1$
for some constant $n\in\NN$, endowed with its canonical structure
of\index{AW*-algebra} AW*-algebra. For example,
$\lz=\lx+\ly$ \iff\ $\bv{\lz=\lx+\ly}=1$. The center of $\tlA$ contains a
copy of the bounded global section algebra $\tCCx$\index{CCtzzCCx@$\tCCx$|ii}
of the complex numbers. Observe that $\tCCx$ is isomorphic to the algebra of
continuous maps from the ultrafilter space of $B$ to $\CCx$. In case $\lA$
is an AW*-factor in\index{VzzVB@$V^B$} $V^B$, the center of $\tlA$ is
exactly (the canonical copy of) $\tCCx$, see \cite[Theorem~5]{Ozaw86}.
\index{Ozawa, M.} In particular, for
$u\in B$, the central idempotent of $\tlA$ corresponding to $u$ is the unique
element $\lu\in\tlA$ such that $u=\bv{\lu=1}$ while $\neg u=\bv{\lu=0}$.
By Lemma~\ref{L:projoflattproj}, the Boolean algebra
\index{Boolean algebra!complete ---} of projections of
$\tlA$ is also isomorphic to $B$. Thus, letting $L$ be the
espalier\index{espalier} of projections of $\tlA$ and identifying the
elements of $B$ with the projections of\index{Dzzrng@$\Drng L$} $\Drng L$,
we obtain that
$\bv{\la\leq\lb}$ has the same meaning in the present paper and in
\index{Takeuti, G.}\cite{Take83,Ozaw86}.\index{Ozawa, M.}

We apply this to the D-universality\index{D-universal} problem as follows.

\begin{theorem}\label{T:bigAW*Drng}
Let $\Omega_{\I}$, $\Omega_{\II}$, $\Omega_{\III}$ be arbitrary
complete Boolean spaces\index{Boolean space} \pup{possibly empty}, and let
$\gamma$ be an arbitrary ordinal. Then there exists an\index{AW*-algebra}
AW*-algebra~$A$ such that
\index{Zzzgamma@$\ZZ_\gamma$}\index{Rzzgamma@$\RR_\gamma$}%
\index{Tzzgamma@$\two_\gamma$}\index{VzzofR@$V(R)$}%
   \[
   V(A)\cong \CC(\Omega_{\I},\ZZ_\gamma)\times
   \CC(\Omega_{\II},\RR_\gamma)\times \CC(\Omega_{\III},\two_\gamma).
   \]
\end{theorem}

\begin{proof}
By arguing as in the proof of Theorem~\ref{T:bigW*Drng}, it suffices to
prove that for every
\index{Boolean algebra!complete ---} complete Boolean algebra $B$, every
ordinal $\gamma$, and every $\J\in\set{\I,\II,\III}$, there exist an
\index{AW*-algebra}AW*-algebra $A$ of type J with algebra of central
idempotents isomorphic to~$B$ and a
$(\gamma+1)$-sequence $\famm{p_{\alpha}}{\alpha\leq\gamma}$ of projections
of central cover $1$ such that $p_{\alpha}\lesssim p_{\beta}$ but
$\bv{p_{\beta}\lesssim p_{\alpha}}=0$, for $\alpha<\beta\leq\gamma$.

By applying Lemma~\ref{L:bigW*factors} within\index{VzzVB@$V^B$}
$V^B$, we obtain a factor $\lA$ of type $\J$ in $V^B$\index{VzzVB@$V^B$} and
a $B$-valued name $\lp$ such that the following statements hold in $V^B$
(that is, they have Boolean value $1$):
  \begin{gather}
  \lp\text{ is a map from }\check{\gamma}
  \text{ to the purely infinite elements of }\lL,\notag\\
  \lp(\lga)\lesssim\lp(\lgb),\text{ for all }
  \lga<\lgb\leq\check{\gamma},\label{Eq:palesssimpb}\\
  \lp(\lgb)\not\lesssim\lp(\lga),\text{ for all }
  \lga<\lgb\leq\check{\gamma},\label{Eq:pbnlesssimpa}
  \end{gather}
where $\lL$ denotes the espalier\index{espalier} of projections of $\lA$
within\index{VzzVB@$V^B$} $V^B$.

Now let $A=\tlA$ be the bounded global section algebra of $\lA$. It
follows from \index{Ozawa, M.}\cite[Theorem~7]{Ozaw86}
\index{Takeuti, G.}(see also
\cite[\S 2]{Take83}) that $A$ is Type J, furthermore, its algebra of central
idempotents is isomorphic to $B$. For all $\alpha\leq\gamma$, let
$\lp_{\alpha}$ be the unique $B$-valued name such that
$\bv{\lp_{\alpha}=\lp(\check{\alpha})}=1$. For $\alpha<\beta\leq\gamma$,
it follows from \eqref{Eq:palesssimpb}, \eqref{Eq:pbnlesssimpa} that
$\lp_{\alpha}\lesssim\lp_{\beta}$ and
$\bv{\lp_{\beta}\lesssim\lp_{\alpha}}=0$.
\end{proof}

Therefore, we have obtained the following result, which, together with the
other main results of the present section, completely elucidates the
dimension theory of projections of W*- and\index{AW*-algebra} AW*-algebras.

\begin{theorem}\label{T:L(A)Duniv}
The class of espaliers\index{espalier} obtained from projection
lattices\index{lattice} of\index{AW*-algebra} AW*-algebras is
D-universal\index{D-universal}.
\end{theorem}

\section{Concluding remarks}\label{S:Concluding}

The questions arising naturally from this work can be divided in two parts:
namely, those where the theory reflects about itself, and those where it
reflects about other topics.

In the first group, we shall mention the following. For given, ``practical''
examples, where we need to verify that a given structure is an espalier, the
axiom~(L7) is often a source of problems. Thus we may ask to what extent it
is possible to remove Axiom~(L7) from the definition of an espalier, thus
defining ``pre-espaliers'' (see also Definition~\ref{D:BpreEsp}). But then,
in order to extend a pre-espalier $(L,\leq,\perp,\sim)$ to an espalier, we
need to define a new binary relation $\sim^*$ on $L$ by letting $x\sim^*y$
hold, if there are decompositions $x=\oplus_{i\in I}x_i$ and
$y=\oplus_{i\in I}y_i$ such that $x_i\sim y_i$, for all $i\in I$. However,
proving the transitivity of the new relation $\sim^*$ leads to the
verification of a common refinement property, see Lemma~\ref{L:ExtPreEsp}
for the Boolean case. This problem can be formulated as follows.

\begin{problem}\label{Pb:pre-espRef}
Let $(L,\leq,\perp,\sim)$ be a structure satisfying all axioms from (L0) to
(L8) with the possible exception of (L7), and let $\famm{x_i}{i\in I}$ and
$\famm{y_j}{j\in J}$ be families of elements of $L$ such that
$\oplus_{i\in I}x_i=\oplus_{j\in J}y_j$. Are there families
$\famm{u_{i,j}}{(i,j)\in I\times J}$ and $\famm{v_{i,j}}{(i,j)\in I\times J}$
of elements of $L$ such that $x_i=\oplus_{j\in J}u_{i,j}$ (for all $i\in I$),
$y_j=\oplus_{i\in I}v_{i,j}$ (for all $j\in J$), and $u_{i,j}\sim v_{i,j}$
(for all $(i,j)\in I\times J$)?
\end{problem}

The second group of questions asks for constructing further classes of
espaliers, within other areas of mathematics. Of course, isomorphism types
of various structures are privileged, see, for example, B. J\'onsson and
A. Tarski's\index{J\'onsson, B.} appendix in\index{Tarski, A.}
\cite{Tars}. In another direction, one might ask about extensions of various
results of cancellation or unique decomposition, known for finite
structures (see \cite[Chapter~5]{MMTa})
\index{McKenzie, R.\,N.}\index{McNulty G.\,F.}\index{Taylor, W.\,F.}%
to infinite structures subjected to completeness conditions. This would in
turn yield, for example, nontrivial cancellation results for further infinite
structures, of which the main result of \cite{Jenc}\index{Jen\v{c}a, G.}
about \emph{$\sigma$-complete effect algebras} would be a prototype.

Expecting infinite generalizations of finite results \emph{via} espaliers is
reasonable as long as there are enough refinement theorems around, see,
again, \cite[Chapter~5]{MMTa}. Hence the Lov\'asz cancellation theorems, see
\cite{Lova67}\index{Lov\'asz, L.} or \cite[Section~5.7]{MMTa}, do not enter
this category, as they are established by counting arguments, in contexts
where refinement does not always hold.
We do not know of any framework that could extend Lov\'asz's results to
infinite structures subjected to completeness conditions.

\printindex

\end{document}